\documentclass{article}
\usepackage{Arxiv}

\usepackage{lipsum}
\usepackage{amsfonts}
\usepackage{graphicx}
\usepackage{epstopdf}
\usepackage{algorithmic}
\usepackage[utf8]{inputenc}
\usepackage{url}
\usepackage{amssymb}
\usepackage{amsmath}
\usepackage{amsthm}
\usepackage{verbatim}
\usepackage{indentfirst}
\usepackage{float}%use [h] to fix postion of a figure
\usepackage{xcolor}
\usepackage{comment}
\usepackage{graphicx}
\usepackage{subfigure}
\usepackage{algorithm}
\usepackage{breqn}
\usepackage{todonotes}
\usepackage{etoolbox}
\makeatletter
\patchcmd{\@addmarginpar}{\ifodd\c@page}{\ifodd\c@page\@tempcnta\m@ne}{}{}
\makeatother
\reversemarginpar
\usepackage{mathtools}
\allowdisplaybreaks
\usepackage{bbm}
\usepackage{bbding}
\usepackage{stackengine}
\newcommand\barbelow[1]{\stackunder[1.2pt]{$#1$}{\rule{.8ex}{.075ex}}}
\ifpdf
  \DeclareGraphicsExtensions{-eps-converted-to.pdf,-eps-converted-to.pdf,.png,.jpg}
\else
  \DeclareGraphicsExtensions{-eps-converted-to.pdf}
\fi
\DeclareMathOperator*{\argmax}{\textbf{argmax}}
\DeclareMathOperator*{\argmin}{\textbf{argmin}}

\def\fprod#1{\left\langle#1\right\rangle}
\def\prox#1{\mathbf{prox}_{#1}}

\def\grad{\nabla}

\def\bx{\mathbf{x}}  %{\mbox{\boldmath $\lambda$}}

\def\cC{\mathcal{C}}
\def\cD{\mathcal{D}}

\def\cF{\mathcal{F}}
\def\cG{\mathcal{G}}

\def\cJ{\mathcal{J}}

\def\cL{\mathcal{L}}

\def\cO{\mathcal{O}}
\def\cP{\mathcal{P}}

\def\cR{\mathcal{R}}
\def\cS{\mathcal{S}}

\def\cX{\mathcal{X}}
\def\cY{\mathcal{Y}}
\def\cZ{\mathcal{Z}}

\def\smskip{\smallskip}

\def\texitem#1{\par\smskip\noindent\hangindent 25pt
               \hbox to 25pt {\hss #1 ~}\ignorespaces}

\def\norm#1{\|#1\|}

\newcommand{\BEAS}{\begin{eqnarray*}}
\newcommand{\EEAS}{\end{eqnarray*}}
\newcommand{\BEA}{\begin{eqnarray}}
\newcommand{\EEA}{\end{eqnarray}}
\newcommand{\BEQ}{\begin{eqnarray}}
\newcommand{\EEQ}{\end{eqnarray}}
\newcommand{\BIT}{\begin{itemize}}
\newcommand{\EIT}{\end{itemize}}
\newcommand{\BNUM}{\begin{enumerate}}
\newcommand{\ENUM}{\end{enumerate}}

\newcommand{\BA}{\begin{array}}
\newcommand{\EA}{\end{array}}

\newcommand{\reals}{\mathbb{R}}
\newcommand{\integers}{\mathbb{Z}}

\newcommand{\Tr}{\mathop{\bf Tr}}

\newcommand{\dom}{\mathop{\bf dom}}

\newif\ifpagenumbering
\pagenumberingtrue

\pagenumberingfalse
\newsavebox{\theorembox}
\newsavebox{\lemmabox}
\newsavebox{\defnbox}
\newsavebox{\assbox}
\savebox{\theorembox}{\noindent\bf Theorem}
\savebox{\lemmabox}{\noindent\bf Lemma}
\savebox{\defnbox}{\noindent Definition}

\usepackage{soul}
\usepackage[normalem]{ulem}
\DeclareUnicodeCharacter{2212}{~}
\newtheorem{assumption}{Assumption}
\newtheorem{theorem}{Theorem}
\newtheorem{lemma}{Lemma}
\newtheorem{definition}{Definition}
\newtheorem{remark}{Remark}
\newtheorem{proposition}{Proposition}
\newtheorem{corollary}{Corollary}

\usepackage{hyperref}
\usepackage{cleveref}
\crefname{assumption}{Assumption}{Assumptions}
\crefname{theorem}{Theorem}{Theorems}
\crefname{lemma}{Lemma}{Lemmas}

\title{Robust \mg{Accelerated} Primal-Dual Methods for Computing Saddle Points
}

\author{
Xuan Zhang\\
Department of Industrial and Manufacturing Engineering\\
Pennsylvania State University
\\
University Park, PA,USA.\\
\texttt{xxz358@psu.edu}
\And
Necdet Serhat Aybat\\
Department of Industrial and Manufacturing Engineering\\
Pennsylvania State University
\\
University Park, PA,USA.\\
\texttt{nsa10@psu.edu}
\And
Mert Gürbüzbalaban\\
Department of Management Science and Information Systems\\
Rutgers University\\
Piscataway, NJ, USA\\
\texttt{mg1366@rutgers.edu}
}

\usepackage{amsopn}

\DeclareMathOperator{\diag}{diag}

\ifpdf
\hypersetup{
  pdftitle={{Robust Accelerated Proximal Primal-Dual Methods for Computing Saddle Points}},
  pdfauthor={Zhang Xuan，Necdet Serhat Aybat, Mert Gürbüzbalaban}
}
\fi

\def\sa#1{\textcolor{red}{#1}}
\def\mg#1{\textcolor{blue}{#1}}
\def\xzh#1{\textcolor{black}{#1}}

\def\xuan#1{\todo[size=footnotesize]{XUan:~#1}}

\def\newcomment#1{\textcolor{brown}{#1}}

\renewcommand{\sa}[1]{\textcolor{black}{#1}}
\renewcommand{\mg}[1]{\textcolor{black}{#1}}
\newcommand{\rev}[1]{\textcolor{black}{#1}}

\def\xuan#1{\textcolor{black}{#1}}
\newcommand{\fin}[1]{\textcolor{black}{#1}}

\makeatletter
\def\hlinewd#1{%
\noalign{\ifnum0=`}\fi\hrule \@height #1 \futurelet
\reserved@a\@xhline}
\makeatother

\begin{document}
\setlength{\abovedisplayskip}{2pt}
\setlength{\belowdisplayskip}{2pt}

\maketitle

\begin{abstract}
    %\nsa{Abstract needs to be rewritten.}
    \mg{We consider strongly-convex-strongly-concave %(SCSC)
    saddle point %(SP)
    problems
    %in the setting when
    \sa{assuming} we have access to unbiased stochastic estimates of the gradients. We propose a stochastic accelerated primal-dual~\sa{(SAPD)} algorithm and show that \sa{SAPD %iterate
    sequence, generated using constant primal-dual step sizes, linearly} converges to a neighborhood of the \sa{unique} saddle point.
    % %at a linear rate, %$\rho$
    % where the size of the neighborhood is determined by the asymptotic variance of the iterates. 
    Interpreting %the asymptotic variance 
    \fin{the size of the neighborhood} as a measure of robustness to gradient noise, we obtain explicit characterizations of robustness in terms of SAPD parameters and problem constants. Based on these characterizations, we develop computationally tractable techniques for optimizing \sa{the SAPD parameters, i.e., %of the SAPD algorithm
    %such as
    the primal and dual step sizes,
    %\xtodo{Shall we use stepsize or stepsize? I notice that both show up in paper. MG: Both of used in the literature, I made all of them "step size"}
    %dual step size
    and the momentum parameter,} to achieve a desired trade-off between the convergence rate and robustness \sa{on the \mg{Pareto curve}. This allows SAPD to enjoy fast convergence properties while being robust to noise as an accelerated method. %We also show that 
    SAPD admits convergence guarantees for the %gap
    \rev{distance} metric %where the
    \sa{with a} variance term optimal up to \sa{a logarithmic factor} --which can %be made optimal
    {be} \sa{removed} by employing a restarting strategy.
    %\nsa{gap result without compact domain is removed!}
    %Furthermore, to our knowledge, our work is the first one showing %an iteration
    %a complexity result for the gap function on smooth SCSC problems without the bounded domain assumption.
    \rev{We also discuss how convergence and robustness results extend to the convex-concave setting.}
    Finally, we illustrate %the efficiency of our approach
    \rev{our framework}
    %on bilinear SP problems and
    on distributionally robust logistic regression \fin{problem}.}}\vspace*{-2mm} %problems.}} %on a regularized bilinear sP problem with synthetic data and on a distributionally robust logistic regression problem based on real datasets.
    %Our method is optimal for bilinear SP problems.%We analyze the trade-offs between the bias and variance robustness to gradient noise as a function of the parameters (primal stepsize, dual stepsize,

     %the robust saddle point problem
  %   \newcomment{analyze the trade-offs between the convergence rate and robustness with SAPD and propose a framework to systematically .}
    
\end{abstract}
\vspace*{-3mm}
\section{Introduction}
We consider the \sa{following} saddle point~(SP) problem:
\begin{align}
\label{eq:main-problem}
    \min_{x\in \mathcal{X}} \max_{y \in \mathcal{Y}} \mathcal{L}(x,y)  \triangleq f(x) +  \Phi(x,y) - g(y),
\end{align}
where $\mathcal{X}$ and $\mathcal{Y}$ are, $n$ and $m$ dimensional inner product spaces endowed with \sa{inner product} norms \sa{ $\|x\|_{\mathcal{X}}=\sqrt{\langle x,x\rangle}_{\mathcal{X}}$ and $\|y\|_{\mathcal{Y}}=\sqrt{\langle y,y\rangle}_{\mathcal{Y}}$}, respectively; $\Phi:\mathcal{X} \times \mathcal{Y} \rightarrow \reals$ is
convex in $x$, concave in $y$ \mg{with} %and has a
\sa{a Lipschitz gradient}; $f:\mathcal{X}\rightarrow \reals\sa{\cup\{\rev{+}\infty\}}$ and $g:\mathcal{Y}\rightarrow \reals\sa{\cup\{\rev{+}\infty\}}$ are \sa{closed, convex functions}. \fin{We assume that $\cL(x,y)$ is (strongly) convex in $x$ and (strongly) concave in $y$ with %modules
\sa{moduli $\mu_x,\mu_y\geq 0$}, %for all $x\in\cX$ and $y\in\cY$
respectively.}
The SP problem \sa{in~\eqref{eq:main-problem}} has a wide range of applications; in fact, {many convex optimization problems arising in machine learning~\sa{(ML)} can be recast as
%the SP problem described above by
\sa{\eqref{eq:main-problem} through Lagrangian duality.} Prominent applications %where
with %saddle-point
\sa{SP} formulations %arise
include empirical risk minimization~\sa{(ERM)} \cite{zhang2017stochastic,wang2017exploiting}, supervised learning with non-separable losses, or regularizers \cite{xu2005maximum,palaniappan2016stochastic}, distributionally robust %formulations of empirical risk minimization problems
\sa{ERM}~\cite{namkoong2016stochastic}, and robust optimization \cite{ben2009robust}.} %\nsa{To save space, I removed constrained convex optimization example, as it also does not have a strongly concave dual part.}
% \sa{A canonical example is the constrained convex optimization problems:}
% $$
% \min_{x\in\mathcal{X}} f(x)+\sa{\phi(x)}\; \text{s.t.}\; G(x)\in -\mathcal{K},
% $$
% \mg{where}
% $\mathcal{K}$ is a closed convex cone \mg{lying} in the dual space $\mathcal{Y}^*$, \sa{$f$ is a closed strongly convex function,
% \sa{$\phi$} is convex with a Lipschitz gradient}, $G:\mathcal{X}\rightarrow\mathcal{Y}^*$ is \sa{a Lipschitz, $\mathcal{K}$-convex function} \mg{with} a Lipschitz Jacobian. Using Lagrangian duality, one can equivalently reformulate the above problem as
% $$
% \min_{x\in\mathbb{R}^n}\max_{y\in\mathcal{K}^*} f(x)+g(x)+\langle G(x),y\rangle,
% $$
% which is a special case of \eqref{eq:main-problem}. \mg{This problem arises frequently in applications ranging from \cn to \cn}.
\sa{In many of these applications,} \rev{one does} not have access to exact values of the gradients
%$\nabla_x\mathcal{L}$ and $\nabla_y\mathcal{L}$
\sa{$\grad_x\Phi$ and $\grad_y\Phi$;} but, %we have 
\rev{rather has} access to their unbiased stochastic estimates \sa{$\tilde\grad_x\Phi$ and $\tilde\grad_y\Phi$}. This would typically be the case when the gradients are estimated from a subset of data points in the \sa{big-data} regime (as in stochastic gradient and stochastic approximation methods) or if noise is injected to the gradients on purpose to protect the privacy of the user data~\cite{xie2018differentially}.

%In this paper, 
\fin{We propose a first-order method, the Stochastic Accelerated Primal-Dual~\sa{(SAPD)} algorithm, to solve \eqref{eq:main-problem} under the assumption that we have access to unbiased stochastic oracles $\tilde\grad_x\Phi$ and $\tilde\grad_y\Phi$
% are unbiased \mg{estimators} of $\nabla_x\Phi(x,y),\nabla_y\Phi(x,y)$ are known, i.e.,
% $$
% \mathbb{E}\left[\tilde{\nabla}\Phi_x(x,y;\omega^x)|x,y \right] =  \nabla\Phi_x(x,y),~
%  \mathbb{E}\left[\tilde{\nabla}\Phi_y(x,y;\omega^y)|x,y \right] = \nabla\Phi_y(x,y).
% $$
%In addition, we assume
with
%the stochastic estimates
%have
a bounded %second moment
variance, see~\cref{ASPT: unbiased noise assumption} for the details. {This setting is commonly considered in the literature and is relevant %in
\sa{to} a number of applications,}
%including the problem of
e.g., training GANs~{\cite{zhong2020improving} and robust learning~\cite{pmlr-v32-wen14}}. First, assuming that $\cL$ is strongly convex strongly concave (SCSC), we show that {SAPD %iterate
    sequence, generated using constant primal-dual step sizes, linearly} converges to a neighborhood of the {unique} saddle point.
    Interpreting 
    \fin{the size of the neighborhood} as a measure of \textit{robustness}\footnote{\fin{This definition of robustness for an algorithm is inspired by the robust control literature, and that it should not be mixed with robustness in~\cite{nemirovski2009robust}.}} to gradient noise, we propose computationally tractable techniques for optimizing the SAPD parameters to achieve a desired trade-off between the convergence rate and robustness. We also discuss how convergence and robustness results extend to the convex-concave setting with $\mu_x=\mu_y=0$.} \vspace*{-1mm}
    % This allows SAPD to enjoy fast convergence properties while being robust to noise as an accelerated method.
    % SAPD admits convergence guarantees for the 
    % \rev{distance} metric 
    % \sa{with a} variance term optimal up to \sa{a logarithmic factor} --which can 
    % {be} \sa{removed} by employing a restarting strategy.

%\looseness=-1
%, where the generator and the discriminator approximate the gradient by taking a batch of data $\mathcal{D}$ and computing the average of the gradient from every single point in data $\mathcal{D}$ as the iteration oracle.
%In addition, this setting also appears in other scenarios such as privacy-related applications where the noise is added intentionally to prevent the model from remembering possibly sensitive data and preserve privacy or when the presence of noise is inevitable due to imperfections in communication and sensing.
% the study on stochasticfirst-order methods for stochastic SPP is still quite limited.  In the stochastic setting, we assume that thereexists astochastic oracle(SO) that can provide unbiased estimators to the gradient operators∇G(x) and(−Kx,KTy).
\subsection{Related Work}
\label{sec:related}
\mg{When the coupling term $\Phi(x,y)$ is bilinear, \sa{i.e.,} $\Phi(x,y) = \langle Kx,y \rangle$ for some linear operator $K:\mathcal{X}\rightarrow
\sa{\mathcal{Y}^*}$,
%the problem
\eqref{eq:main-problem} is well-studied for both strongly-convex-strongly-concave (SCSC) problems (%where
$\mu_x, \mu_y>0$) as well as for merely-convex-merely-concave \sa{(MCMC) problems ($\mu_x=\mu_y=0$).} %where $f$ and $g$ are \sa{merely} convex. %(but not necessarily strongly convex).
The convergence results cover both the stochastic case (when only stochastic estimates of the gradients are available) %as well as
\sa{and} the deterministic case (when the gradient information is exact).}
%the convergence results of not only deterministic cases but also stochastic cases have been well studied, and those include from merely convex-merely concave (MCMC) to strongly convex- strongly concave (SCSC).
\sa{In our work,
%we focus on the general SP problem for which
we do not assume bilinear $\Phi$.}
When the coupling term $\Phi$ is non-bilinear, there exist some convergence results in the deterministic case; however, the stochastic setting remains relatively understudied. %\looseness=-1
% although some deterministic results exist, the research on the stochastic case is few. Especially for the research on SCSC objective function with stochastic oracle, the results are very limited, and ours owns the most general objective function among them.
Two standard metrics to measure the quality of a random  $(\bar{x},\bar{y})\in\cX\times\cY$ returned by a stochastic algorithm are the \emph{gap function} $\cG:\cX\times\cY\to \mathbb{R}_+$ for the %general
\fin{MCMC case} and \emph{distance metric} $\cD:\cX\times\cY\to \mathbb{R}_+$ for the SCSC case, for which there is a unique saddle point $(x^*,y^*)$, i.e.,\looseness=-1 %They are defined as
\vspace*{-1mm}
{\small
\begin{equation}
\label{eq:gap}
\cG(\bar{x},\bar{y})\triangleq\mathbb{E}\big[ \sup_{(x,y)\in \cX\times \cY} \big\{ \mathcal{L}(\bar{x}, y)  - \mathcal{L}(x, \bar{y}) \big\} \big],\quad
\rev{\cD(\bar{x},\bar{y})\triangleq\mathbb{E}[\mu_x\norm{\bar{x}-x^*}^2+\mu_y\norm{\bar{y}-y^*}^2]},\vspace*{-2mm}
\end{equation}}%
\mg{where %\rev{$\mu_x$ and $\mu_y$ are the strong convexity and concavity moduli}, and 
the expectation is taken with respect to the randomness encountered in the %choice of
\sa{generation of} the point $(\bar{x},\bar{y})\in\cX\times\cY$}. \mg{In the following discussion, we summarize existing results \sa{closely related to our setting, and discuss} our contributions.}\vspace*{-2mm}
\subsubsection{The Deterministic Case}
%\sa{A significant amount of previous work analyzes the convergence properties of algorithms designed for solving \eqref{eq:main-problem} with a bilinear primal-dual coupling term,} i.e., $\Phi(x,y) = \langle Kx,y \rangle$ for some linear operator $K:\mathcal{X}\rightarrow \mathcal{Y}$. \sa{In our work, we focus on a general SP problem for which we do not assume bilinear $\Phi$ and conclude corollaries for the particular case of bilinear $\Phi$.}
%Next, we will give a brief review of the recent work for the SP problem \sa{ \eqref{eq:main-problem}}.

\sloppy The bilinear %case
\sa{structure} has been thoroughly studied;
some well-known algorithms include %Nesterov smoothing (a.k.a.,
excessive gap technique~\cite{nesterov2005smooth,nesterov2005excessive}, primal-dual hybrid gradient (PDHG) %Chamoblle and Pock
\cite{chambolle2011first,chambolle2016ergodic} \rev{--also see~\cite{thekumparampil2022lifted} achieving the best bound.}
\sa{On the other extreme}, when $\mathcal{L}$ owns a general form and the smoothness cannot be guaranteed, \sa{primal-dual subgradient algorithms} %derived from primal-dual subgradient
have been proposed in several works, \sa{e.g., %including Nedic and Ozdaglar
\cite{nedic2009subgradient,nesterov2009primal,juditsky2011first}}.
%However, it is within expectations that non-smoothness is not able to assure a low iteration complexity. As a result, those algorithms become weaker when
\sa{The iteration complexity of these subgradient-based methods can be significantly improved when $\mathcal{L}$ has further structure.
%but more general than the bilinear case.
Indeed,
%With $\mathcal{L}$ or even only the coupling term
there are methods exploiting the structure when $\Phi$ is smooth, and $f,g$ have efficient prox maps, which include} %many algorithms have been developed to achieve a better iteration complexity. Those methods include
%\xtodo{The order of this paragraph is not good.}
Mirror-Prox(MP)~\cite{nemirovski2004prox}, Optimistic Gradient Descent Ascent~(OGDA) and Extra-gradient~(EG)~\cite{mokhtari2020unified} \sa{methods}.
%When using Euclidean geneate function, MP has the same structure with EG.
Additionally, the \sa{effect of} Lipschitz constants
%on the
\sa{along} different blocks of variables also has been explored recently; \sa{some new works account for the individual effects of $L_{xx}$, $L_{yx}$ and $L_{yy}$, i.e., the Lipschitz constants of $\nabla_x\Phi(\cdot, y)$, $\nabla_y\Phi(\cdot, y)$, and $\nabla_y\Phi(x, \cdot)$, 
respectively, instead of using the worst-case
%Lipschitz constant  42
{parameters} $L\triangleq\max\{L_{xx}, L_{xy}, L_{yx}, \sa{L_{yy}}\}$,}
 {$\mu\triangleq\min\{\mu_x,\mu_y\}$.}
%For example,
For bilinear SP problems, a lower complexity bound of ${{\Omega}\Big( \sqrt{1 + \frac{L_{yx}^2}{\mu_x\mu_y}}\cdot \ln(1/\epsilon)\Big)}$ is shown in~\cite{zhang2019lower} for a class of \sa{first-order primal-dual algorithms employing proximal-gradient steps; on the other hand, the lower bound for %first-order
gradient-based methods
is $\Omega\Big( \sqrt{\frac{L_{xx}}{\mu_x} + \frac{L_{yx}^2}{\mu_x\mu_y}+\frac{L_{yy}}{\mu_y}}\cdot \ln(1/\epsilon)\Big)$ when $f(\cdot)=g(\cdot)=0$ and $\Phi$ is SCSC \cite{zhang2019lower}.}
%Our algorithm is based on  \cite{hamedani2018primal}, extending it to a SCSC scenario, \newcomment{meanwhile providing a wider ranger of step sizes for both MCMC and SCSC cases.}
% Moreover, suppose $f(x)$ is $L_x$-smooth, $\mu_x$-strongly convex and $g(y)$ is $L_y$-smooth, $\mu_y$-strongly convex, and $\Phi(x,y)$ is bilinear, lifted primal-dual (LPD)~\cite{thekumparampil2022lifted}
% method achieves the iteration complexity of $\cO\Big(\sqrt{\frac{L_x}{\mu_x}} + \frac{L_{yx}}{\sqrt{\mu_x\mu_y}} + \sqrt{\frac{L_y}{\mu_y}} \Big)$; in addition to this work, without assuming $\Phi(x,y)$ is bilinear, 
\rev{In the rest of this discussion, we focus on {the} deterministic SCSC setup, 
%It is important to note that for deterministic SCSC setting, 
for which the results in $\cD$ metric can be converted into %duality 
gap metric $\cG$ by only increasing the logarithmic term by problem parameters, see~\cite[Appendix C]{cohen2020relative}.}
%For the SCMC setting, Thekumparampil et al. \cite{thekumparampil2019efficient} proposed an algorithm that can achieve the optimal $O(1/k^2)$ rate. Note that MCMC SP problems can be formulated as variational inequalities \cn.
%\nsa{Due to limited space, we should focus on SCSC case to cut down the intro -- there are tons of work for SCMC, MCMC cases.}

{Mokhtari \mg{\emph{et al.}} \cite{mokhtari2020unified} show that both OGDA and EG have an iteration complexity of $\mathcal{O}\Big(\frac{L}{\mu}\ln(1/\varepsilon)\Big)$
%\mg{when distance to the saddle-point is taken as a performance metric}
\sa{for $\cD$ metric defined in~\eqref{eq:gap}}. %\nsa{The metric is not distance, it is distance square.}
Gidel \mg{\emph{et al.}} \cite{gidel2018variational} also \sa{show} %s
the same rate for OGDA \mg{
%based on
\sa{from} a variational inequality} (VI) \mg{perspective}.}
\sa{In the analysis of these algorithms, %e.g., EG and OGDA, treats
primal and dual step sizes are set equal, which may
%be not powerful for many large-scale problems
lead to conservative steps whenever $L_{xx}\gg L_{yy}$, or vice versa.} \mg{For instance, in the primal-dual formulation of empirical risk minimization problems in machine learning, choosing primal and dual step sizes to be different can lead to an improved convergence rate \cite{zhang2017stochastic}.}
%\nsa{What is the metric, what does $\epsilon$ refer to?}
{There are also some multi-loop algorithms. \mg{In particular,}
Lin \emph{et al.}~\cite{lin-near-optimal} propose an \sa{inexact proximal point} algorithm, which consists of 3-nested loops; %for the gap \mg{function \sa{$\cG$} as a metric}.
\sa{each proximal step computation requires calling \fin{Nesterov's accelerated gradient descent} (AGD) iteratively to solve strongly convex smooth~(SCS) optimization subproblems with a high precision}
%The inner loop requires solving an inexact problem up to a high-precision
that \mg{can be} impractical.%\nsa{I think there is an error in the complexity statement for the subproblems -- I corrected it, please check!}
%\xtodo{1. It is right to me if the $\kappa$ in AGD is equal to $\kappa_x$. 2. The footnote has a SCS function, what does it mean?}\mtodo{ok, looks good now.}
\footnote{\sa{In each AGD call, an SCS function $h$ with condition number $\kappa_x=L/\mu_x$ is minimized 
%to achieve a target accuracy of
to compute $\bar{x}\approx \argmin_{x\in X}h(x)$ such that $\norm{\bar{x}-\Pi_{X}(\bar{x}-\grad h(\bar{x})/L)}^2\leq \frac{\epsilon}{2(10 \kappa_y)^{11}\kappa_x^{13}}$}.} The computational complexity to compute $(\bar{x},\bar{y})$ such that $\cG(\bar{x},\bar{y})\leq \epsilon$ is $\mathcal{O}\left(\frac{L}{\sqrt{\mu_x\mu_y}}\cdot\sa{\ln^3}(1/\epsilon)\right)$.
%\mtodo{in which metric?}, sub-problems are required to be solved up to an accuracy of
%$\varepsilon=\frac{2(\mu_x\mu_y)^{11}(L-\mu)\epsilon}{10^{11}L^{20}\mu^2}$)
%which is restrictive for ill-conditioned problems when $\kappa = L/\mu$ is large.}
%\nsa{This quantity in the footnote should be simplified, fors instance $L-\mu$ can be deleted to get $L^{19}$ in the denominator.}
\sa{Although %they 
\cite{lin-near-optimal} claims %that their algorithm 
to achieve the lower complexity bound provided in~\cite{zhang2019lower}, this is not the case for problems with $L_{yx}\ll L$.} %\nsa{What are the disadvantages of this algorithm?}
The algorithm in~\cite{wang-li} consists of 4-nested loops and has similar shortcomings in practice. The computational complexity to compute $(\bar{x},\bar{y})$ such that $\cD(\bar{x},\bar{y})\leq \epsilon$ is $\mathcal{O}(\sqrt{\frac{L_{xx}}{\mu_x} + \frac{L L_{xy}}{\mu_x \mu_y} + \frac{L_{yy}}{\mu_y}}\cdot\ln^3(1/\epsilon))$. \fin{More recently, after our paper has appeared, Jin~\emph{el al.}~\cite{jin2022sharper} independently obtain the iteration complexity of ${\cO}\Big(\big(\frac{L_{xx}}{\mu_x} +
\frac{L_{yx}}{\sqrt{\mu_x\mu_y}}
%+ \sqrt{\frac{L_y}{\mu_y}} 
+\frac{L_{yy}}{\mu_y}\big)\cdot\ln(1/\epsilon)\Big)$ for satisfying $\cG(\bar x,\bar{y})\leq \epsilon$.} 
%-- note that this is a weaker result in terms of the metric used, i.e., $\cG(\bar{x},\bar{y})\leq\epsilon$ implies that $\cD(\bar{x},\bar{y})=\cO(\epsilon)$; but, the reverse is not true.
%\xtodo{The disadvantage for \cite{wang-li} still need to be writen more. MG: Agreed.}
%\looseness=-1
\subsubsection{The Stochastic Case}
%\xtodo{The following part is not about SCSC.}
While the deterministic SP problem has attracted %huge
much attention, the study on \sa{the first-order stochastic methods for %stochastic non-bilinear SP problems
\eqref{eq:main-problem} is still relatively} limited. \mg{For MCMC SP problems, proximal methods} %based on mirror descent
have been developed, \sa{e.g., Stochastic Mirror-Descent(SMD)~\cite{nemirovski2009robust},
%in the literature. They include
the Stochastic Mirror-Prox (SMP)~\cite{juditsky2011solving} and %the stochastic accelerated Mirror-Prox
its accelerated version (SAMP)~\cite{chen2017accelerated}}. \sa{In~\cite{zhao2019optimal}, %Zhao considers
MCMC and {strongly-convex-merely-concave}~(SCMC) scenarios are considered under \emph{additive} unbiased noise with \fin{a} bounded variance.} %setting.
%and assuming $\mathcal{X}$ and $\mathcal{Y}$ are bounded.
%the additive unbiased noise model \mg{when $g$ is merely convex, i.e. when $\mu_y = 0$}.
%When $f$ is merely convex \mg{(i.e. when $\mu_x = 0$)}, with the assumptions that the sets $\mathcal{X}$ and $\mathcal{Y}$ are bounded and the variance of the noise is bounded, they show the convergence of the \mg{supremum} %supreme
%of the dual\mg{ity} gap in expectation.
When $\mu_x>0$, %\mg{the authors} %they develop
a multi-stage \sa{scheme achieving the best known complexity for the stochastic SP problems is proposed in~\cite{zhao2019optimal}; however, this is a two-loop method and each outer iteration requires} solving a \mg{non-trivial} %complicated
sub-problem \sa{with an increasing accuracy, which is a function of some problem parameters that may not be known in practice, e.g., Bregman diameters of $\cX,\cY$ and noise variance}.  {%In addition,
There are also some \mg{VI-based} methods~\cite{cui2016analysis,gidel2018variational} %have been developed
for \sa{the MCMC scenario}.}\looseness=-1
% \mg{In our results}, if we let $\mu_x = \mu_y = 0$, for both deterministic and stochastic \mg{cases}, we can get a similar result with theirs immediately.\todo{MG: Similar in which metric?} Moreover, if we only let $\mu_y = 0$, a similar result of the deterministic case can be found in \cite{hamedani2018primal} on which our analysis is based.

Our focus in this paper will be on the stochastic SCSC case. {Let $\delta^2$ denote the variance bound on the stochastic first-order oracle.}
%\mg{In this setting,} most of the existing results are for $f(\cdot) = g(\cdot) =0$ and  $\Phi(x,y)$ is SCSC, e.g., see~\cite{fallah2020optimal}. %strong convexity-strong concavity.
Yan \emph{et al.} \cite{yan2020optimal} consider $\min_{x\in X}\max_{y\in Y}\Phi(x,y)$ for possibly non-smooth, SCSC $\Phi$, and propose Epoch-GDA with an \rev{oracle complexity of $\mathcal{O}\left(\frac{1}{\epsilon}\ln(1/p)\right)$ for computing $(\bar x,\bar y)$ such that $\cG(\bar x, \bar y)\leq \epsilon$ with probability $1-p$}. When $\Phi$ is smooth,
%Hsieh \emph{et al.}~\cite{hsieh2019convergence} show that a
stochastic
%\mg{extra-gradient} (EG)
EG method~\cite{hsieh2019convergence} for SCSC SP problems and Stochastic Operator Extrapolation method~\cite{kotsalis2020simple} for strongly monotone VIs, both using constant step sizes, can guarantee $\cD(x_k,y_k)\leq \epsilon$ within 
%in the sense of the distance to the saddle point\mtodo{MG: Do we mean expected distance to the saddle point?} with
$\mathcal{O}\left(\frac{1}{\epsilon}\right)$ and $\cO(\kappa\ln(1/\epsilon)+\tfrac{\delta^2}{\mu\epsilon}\ln(1/\epsilon))$ iterations, respectively, where $\kappa=L/\mu$. %Note that the S-OGDA updates can also be thought as an variant of S-EG.
%In addition,
%\xtodo{The font of ''et al'' are not consistent. For example, line 81 and and line 111. MG: I fixed it, usually we need et al. being italic.}
Fallah \emph{et al.}~\cite{fallah2020optimal} propose multi-stage variants (employing restarts) of Stochastic Gradient Descent Ascent~(S-GDA) and Stochastic OGDA~(S-OGDA) that can guarantee $\cD(x_k,y_k)\leq \epsilon$ within 
$\cO(\kappa^2\ln(1/\epsilon)+\tfrac{\delta^2}{\mu\epsilon})$ and $\cO(\kappa\ln(1/\epsilon)+\tfrac{\delta^2}{\mu\epsilon})$ iterations, respectively. Unlike our paper, in both \cite{hsieh2019convergence,fallah2020optimal}, Lipschitz constant of $\grad\Phi$, i.e., $L$, is used to determine the step size, rather than exploiting the block Lipschitz structure.
\vspace*{-2mm}
\subsection{Comparison}
\renewcommand{\arraystretch}{1.5}
{
	\begin{table}[h!]
	    \scriptsize
		\centering
		\begin{tabular}{|c|c|c|c|c|c|}
			\hline
			\textbf{Method} & \textbf{Bias} & \textbf{Variance} & \textbf{Loop} &\textbf{Metric}& \textbf{\fin{BV-tradeoff}}\\
			\hline
			%[Lin~\emph{et al.} 2020]
			\cite{lin-near-optimal} & $\frac{L}{\sqrt{\mu_x\mu_y}}\ln^3\big(
			%(\frac{L^2}{\mu_x}+\frac{L^2}{\mu_y})
			\frac{1}{\epsilon}\big)$ &\XSolidBrush & 3 & $\cG$ & N/A
			\\
			\hline
			%[Wang and Li~2020]
			\cite{wang-li} & {\footnotesize$\tilde\cO\Big(\big({\frac{L_{xx}}{\mu_x}}+\frac{L \cdot L_{xy}}{\mu_x \mu_y}+{\frac{L_{yy}}{\mu_y}}\big)^{1/2}\Big)
				%\ln ^3\left(\frac{L^2}{\mu_x \mu_y}\right)
				{\ln \left(%\frac{L^2}{\mu_x \mu_y}
					\frac{1}{\epsilon}\right)}$} &\XSolidBrush & 4 & $\cD$ & N/A
			\\
			\hline
			%\xzh{\cite{zhang2021complexity}}
			\cite{yang2020catalyst} &   %\xzh{$\frac{L}{\sqrt{\mu_x\mu_y}}\ln(\frac{1}{\epsilon})$ }
			$\frac{L}{\mu}\ln(\frac{1}{\epsilon})$ &  \XSolidBrush & {2} & {$\cD$} & N/A
			\\
			\hline
			%[Mokhtari~\emph{et al.} 2020]
			\cite{mokhtari2020unified} &    $\frac{L}{\mu}\ln(\frac{1}{\epsilon})$   & \XSolidBrush & 1 & $\cD$ & N/A
			\\
			\hline
			%[Cohen~\emph{et al.} 2021]
			\cite{cohen2020relative,jin2022sharper} &  $(\frac{L_{xx}}{\mu_x}+{\frac{L_{yx}}{\sqrt{\mu_x\mu_y}}}+\frac{L_{yy}}{\mu_y}){\ln(\frac{1}{\epsilon})}$     &  \XSolidBrush & 1 & $\cG$ & N/A
			\\
			\hlinewd{1pt}
			%[Yan \emph{et al.}~2020]
			\cite{yan2020optimal} & $\tilde\cO(\frac{1}{\epsilon}\ln(1/p))$ & $\tilde\cO(\frac{\delta^2}{\mu\epsilon}\ln(1/p))$ & 2 &$\cP_p$ & \XSolidBrush \\
			\hline
			%[Fallah~\emph{et al.} 2020]
			\cite{fallah2020optimal} &   $\frac{L}{\mu}\ln(\frac{1}{\epsilon})$     &  $\frac{\delta^2}{\mu\epsilon}$ & 2 &$\cD$ & \XSolidBrush
			\\
			\hline
			%[Hsieh \emph{et al.}~2020]
			\cite{hsieh2019convergence} & $\frac{1}{\epsilon}$ & $\frac{\delta^2}{\mu\epsilon}$ & 1 & $\cD$ & \XSolidBrush
			\\
			% \hline
			% [Kotsalis~\emph{et al.} 2022] &  $\frac{L}{\mu}\ln(\frac{1}{\epsilon})$     &  $\frac{\delta^2}{\mu\epsilon}$ & x & $\cD$
			%  \\
			\hline
			\textbf{ours} &  $(\frac{L_{xx}}{\mu_x}+{\frac{L_{yx}}{\sqrt{\mu_x\mu_y}}}+\frac{L_{yy}}{\mu_y}){\ln(\frac{1}{\epsilon})}$     &  ${(\frac{\delta_x^2}{\mu_x}+\frac{\delta_y^2}{\mu_y})\frac{1}{\epsilon}}\ln(\frac{1}{\epsilon})$ & 1 & $\cD$ & \CheckmarkBold
			\\
			\hline
		\end{tabular}
            \vspace*{2mm}
		\caption{\rev{Related work: Comparison of methods for solving SCSC saddle point problem in~\eqref{eq:main-problem} with a non-bilinear $\Phi$. \cite{cohen2020relative} requires $\Phi$ to be \textbf{twice} differentiable. Results in $\cD$ metric can be converted into guarantees %duality 
        in the gap metric $\cG$ %by only increasing the logarithmic term by problem parameters
			\sa{while still preserving $\ln(1/\epsilon)+1/\epsilon$ complexity,} see~\cite[Appendix C]{cohen2020relative}. The metric $\cP_p$ denotes the number of oracle calls for $\cG\leq \epsilon$ with probability at least $1-p$. Among single-loop methods, \cite{mokhtari2020unified} employs $\tau=\sigma=\frac{1}{4L}$, \cite{cohen2020relative} employs $\tau=\frac{1}{\mu_x\lambda}$ and $\sigma=\frac{1}{\mu_y\lambda}$, where $\lambda=\frac{L_{xx}}{\mu_x}+\frac{L_{xy}}{\sqrt{\mu_x\mu_y}}+\frac{L_{yy}}{\mu_y}$, \cite{hsieh2019convergence} employs $\tau_k=\sigma_k=\frac{1}{\alpha k+4L}$ for $\alpha\in(0,\mu)$. In the column %``$\cA_{\tau,\sigma}$-\textbf{Set}" we indicated whether a generic set of admissible step sizes with a certifiable rate for bias term is provided for single-loop methods.
			``BV-tradeoff" we indicated whether a systematic analysis is provided for the %rate-robustness 
            \fin{bias-variance} trade-off.}}
			\vspace*{-5mm}
		\label{tab:comparison}
\end{table}}%
\rev{In \cref{tab:comparison}, among the papers we discuss in~\cref{sec:related} we compare %first-order primal-dual 
the deterministic and stochastic methods for solving the SCSC saddle point problem in~\eqref{eq:main-problem} with a \emph{non-bilinear} $\Phi$ --to focus on more relevant papers, we did not include methods for $\Phi$ that is bilinear and/or in the finite-sum form. 
% In the table, 
% %we considered important works closely related to ours, 
% some of them use deterministic oracle for computing $\grad_x\Phi$ and $\grad_y\Phi$ while the others employ stochastic first-order oracles~(SFO) 
% to get noisy estimates $\tilde\grad_x\Phi$ and $\tilde\grad_y\Phi$. 
For deterministic methods, having access to $\grad_x\Phi$ and $\grad_y\Phi$, we only provide the bias term of the oracle complexity --this term represents the work required against the bias introduced due to initialization of the algorithm while computing an $\epsilon$-solution. For methods employing stochastic first-order oracles~(SFO) to get noisy estimates $\tilde\grad_x\Phi$ and $\tilde\grad_y\Phi$, we provide both the bias and variance terms in the oracle complexity result, where variance term denotes the additional oracle calls required due to persistent noise in gradient estimates {compared to the (noiseless) deterministic case}. For all the methods compared, we list how many \emph{nested} loops they employ. Finally, in the last column ``{\fin{BV}-tradeoff}" of Table~\ref{tab:comparison} we indicate whether a systematic analysis is provided for the %rate-robustness
\fin{bias-variance} trade-off for the stochastic methods discussed in the table. While ours %paper 
is achieving near optimal state-of-the art complexities for both bias and variance as a single loop method, it also provides conditions on algorithm parameters describing the dependency between the parameter choice and corresponding certifiable rate --see~\eqref{eq: general SAPD LMI_R1}; hence, our admissibility rule for selecting algorithm parameters allows us to characterize the %rate robustness
\fin{bias-variance} trade-off for SAPD.}
\subsection{Contributions}
\sa{We propose the Stochastic Accelerated Primal-Dual~(SAPD) algorithm which extends %the accelerated primal-dual (APD)
the APD method proposed in~\cite{hamedani2018primal} to the stochastic gradient setting.
%For the stochastic scenario,
%More precisely, 
We assume that the first-order oracles %~(SFO)
$\tilde{\grad}\Phi_x$ and $\tilde{\grad}\Phi_y$ return noisy partial gradients that are \textit{unbiased} and have \textit{finite variance} \fin{bounded by} $\delta_x^2$ and $\delta_y^2$, respectively. \rev{Let $z^*=(x^*,y^*)$ denote the unique saddle point of the SCSC minimax problem in~\eqref{eq:main-problem}.} 
%Let $\Omega_f\triangleq\sup_{x,x'\in\dom f}\norm{x-x'}^2$ and $\Omega_g\triangleq\sup_{y,y'\in\dom g}\norm{y-y'}^2$,
%we show that
For any $\epsilon>0$,
%the iteration complexity on the number iterations $K_\epsilon$
%for computing $\{(x_k,y_k)\}$
%such that
SAPD guarantees \rev{$\mathcal{D}({x}_N,{y}_N)\leq \epsilon$} within \vspace*{-2mm} %can be bounded as
%Suppose stochastic first-order oracles~(SFO) returns partial gradients with additive noise having zero-mean and finite variance,  then the variance term in bound
{\small
\begin{equation*}
    N \leq \mathcal{O}\Big( \Big(
\frac{L_{xx}}{\mu_x} + \frac{L_{yx}}{\sqrt{\mu_x\mu_y}} + \frac{L_{yy}}{\mu_y}
+ \Big( \frac{\delta_x^2}{\rev{\mu_x}} +  \frac{\delta_y^2}{\rev{\mu_y}} \Big)\frac{1}{\epsilon}\Big) \cdot \ln\Big( \frac{
\rev{
%\norm{x_0-x^*}^2+\norm{y_0-y^*}^2}
\cD(x_0,y_0)}}{\epsilon}
\Big)\Big) \vspace*{-2mm}
\end{equation*}}%
iterations. %where $\{(\bar{x}_K,\bar{y}_K)\}$ \sa{is a particular weighted average sequence}.\nsa{Do we need bounded domain for the complexity when we consider distance metric?}
%\xtodo{Do we need to use $\cD'$ here?}
\rev{The oracle complexity bound on the bias term $\cO(\kappa\ln(1/\epsilon))$ is optimal, where $\kappa=L/\mu$}, and the bound on the variance term $\tilde{\cO}((\delta_x^2/\rev{\mu_x}+\delta_y^2/\rev{\mu_y})/\epsilon)$ is optimal up to a log factor}, which can be removed \mg{by} employing \sa{a restarting strategy as in~\cite{fallah2020optimal}} \fin{--see \cref{sec:m-sapd}.}
%\nsa{we might put this to the arxiv version and cite it here?} 
%here, $\delta^2_x$ and $\delta^2_y$ denote  bounds of the variance of estimators of $\grad_x\Phi$ and $\grad_y\Phi$, respectively.
%\sa{Moreover, to the best of our knowledge, our work is the first one showing an iteration complexity result for the gap function $\cG$ on smooth SCSC problems without the bounded domain assumption.}
\sa{Since %we assume that 
the noise is persistent,} linear convergence cannot be achieved -- unlike the finite sum problems where variance reduction-based methods are applicable to obtain linear convergence~\cite{palaniappan2016stochastic}. \fin{However,
%our results show that
for SCSC problems, SAPD with constant step size converges to a neighborhood of the saddle-point at a linear %convergence
rate {$\rho\in(0,1)$}, and the size of the neighborhood, defined as $\limsup_{N\to\infty}\mathbb{E}[\norm{z_N-z^*}^2]$, %will be
%\sa{is} determined by the asymptotic variance of the iterates.
%\nsa{We should emphasize all our contributions in terms of the complexity analysis.}
%It follows from our analysis show that
%\sa{The asymptotic variance of the iterates} 
scales linearly with the gradient noise level; {hence,} we interpret the {ratio} $\limsup_{N\to\infty}\mathbb{E}[\norm{z_N-z^*}^2]/\delta^2$ %asymptotic variance
%the worst asymptotic squared distance of $\{z_N\}$ to $z^*$ in expectation %expected 
%{to the gradient noise variance}
%of the iterates
%(after normalization by the input variance)
as a measure of \emph{robustness},
%to gradient noise,
{which we denote with $\cJ$, where $\delta^2=\max\{\delta^2_x,\delta^2_y\}$.
%In particular,
We evaluate the overall algorithmic performance with two metrics:
%depends on both
SAPD parameters should be tuned to achieve %the convergence
\sa{a} faster rate $\rho$
%as well as %the robustness
with a smaller noise amplification $\cJ$}.}
\mg{Our analysis leads to explicit characterizations of %robustness
$\cJ$ for a particular problem class, and of an upper bound $\cR$ on $\cJ$ for more general problems; both $\cJ$ and $\cR$ are given as functions of SAPD parameters. %\rev{i.e.,} strong convexity and Lipschitz
%parameters (primal stepsize, dual stepsize, strongly convexity constants and Lipschitz
%constants of the problem in~\eqref{eq:main-problem}. 
Based on these characterizations, we develop computationally tractable techniques for optimizing the SAPD parameters %algorithm
to achieve a desired systematic trade-off between
%the convergence rate and the robustness
\sa{$\rho$ and $\cJ$ without assuming the knowledge of 
noise variance bounds, $\delta^2_x$ and $\delta^2_y$.} 
%nor the domain diameters $\Omega_f$ and $\Omega_g$. 
This allows SAPD to enjoy fast convergence
%rate properties
with a robust performance in the presence of stochastic gradient noise. Achieving systematic trade-offs between the 
rate and robustness has been previously studied in \cite{aybat2020robust} in the context of accelerated methods for \sa{smooth} strongly convex %unconstrained
minimization problems. To our knowledge, our work is the first one that can trade-off %the convergence rate with robustness
\sa{$\rho$ with $\cJ$} in a systematic fashion in the context of %saddle-point
\sa{primal-dual algorithms for %solving 
%the SP problem in~
\eqref{eq:main-problem}}}.
%\nsa{Discuss why this extension from \cite{aybat2020robust} is not trivial for primal-dual algorithms.(done in response)}

%\newcomment{
%We also creatively propose a framework of an optimization problem to alyze the trade off between convergence rate and the robustness of the algorithm.
%For deterministic scenarios, the framework provides an wider range of step sizes compared to traditional step size condition, e.g., $\frac{\mu}{L}$, thus leading to an potentially faster convergence rate.
%For stochastic scenarios, users can always ensure a desired convergence rate while benefiting from a high robustness.
%We highlighted that our optimization framework is easy to solve and only depends on the Lipshitz constants and convex modules.
%In other word, it is free of variances $\delta_x,\delta_y$ and problem dimensions while common algorithms require to know the variance of the gradient that is difficult to measure in practice.
%}

%\todo{NSA: This should be refined toinclude stepsize effect}
\rev{For the stochastic MCMC case, SAPD can generate $(\bar{x},\bar{y})$ such that $\cG(\bar x,\bar{y})\leq \epsilon$ within $\cO(L/\epsilon+\delta^2/\epsilon^2)$ oracle calls, which is optimal for this setting in both bias and variance terms. For both SCSC and MCMC scenarios,} the deterministic results\footnote{\rev{In the deterministic scenario, SAPD reduces to APD algorithm~\cite{hamedani2018primal}, which has the optimal rate guarantees for MCMC and SCMC (with $L_{yy=0}$) settings; that said, deterministic SCSC setting was not studied in \cite{hamedani2018primal}.}} can be derived from \mg{our} stochastic results immediately \mg{by setting the noise variances $\delta_x^2=\delta_y^2=0$}. %\mg{We note that} 
\rev{In the deterministic setting,} our algorithm, when applied to \eqref{eq:main-problem} with a bilinear $\Phi$, generates the same iterate sequence with \cite{chambolle2016ergodic} for a specific choice of step size parameters; therefore, SAPD, \rev{being able to handle noisy gradients and non-bilinear couplings,} can be viewed as a general form of the \textit{optimal} method (CP) proposed by Chambolle and Pock \cite{chambolle2016ergodic} for \rev{MCMC and SCSC problems with a bilinear coupling.}
%For this special case,
\rev{\fin{Indeed}, in the deterministic case when $\Phi$ is bilinear, both CP and SAPD hit the lower complexity bounds, $\Omega(L/\epsilon)$ for the MCMC and $\Omega(\frac{L_{yx}}{\sqrt{\mu_x\mu_y}}\ln(1/\epsilon))$ for the SCSC problems, given in~\cite{ouyang2021lower} and \cite{zhang2019lower}, respectively.}
%when $\mathcal{L}$ is strongly convex in $x$ and strongly concave in $y$.
 Moreover, when $\Phi$ is not assumed to be bilinear, \sa{SAPD} guarantees $\cO(L/\epsilon)$ complexity in the MCMC setting and $\cO((\frac{L_{xx}}{\mu_x}+\frac{L_{yy}}{\mu_y}+\frac{L_{yx}}{\sqrt{\mu_x\mu_y}})\cdot \ln(1/\epsilon))$ \rev{complexity in the SCSC setting for the bias term, which are the best bounds shown for \eqref{eq:main-problem}.} \rev{%We also note that 
 \fin{For the SCSC setup, the papers~\cite{cohen2020relative,jin2022sharper} provide} bias guarantees similar to our method; but, they are not applicable to the (noisy) stochastic setting like ours.}
 %-- \sa{we are not aware of any other single-loop methods achieving this complexity bound for \eqref{eq:main-problem}.
 %\nsa{See Niao He's result that she mentioned to Mert.} 
 Furthermore, our framework exploiting block Lipschitz constants $L_{xx}$, $L_{yx}$ and $L_{yy}$, provides %an wider range of
 larger step sizes compared to the traditional step size $\cO(1/L)$.

\rev{Finally, the single-loop design of our algorithm make it suitable for solving large-scale problems efficiently {--usually in methods with nested loops, inner iterations are terminated when a sufficient optimality condition holds and these conditions are usually very conservative, leading to excessive number of inner iterations. Furthermore, solving nonconvex-concave minimax problems using an inexact proximal point method requires solving SCSC subproblems to an increasing accuracy; hence, adopting single-loop algorithms as solvers for SCSC subproblems leads to simple implementations compared to using multi-loop methods as solvers --see~\cite{zhang2022sapd+}.} Indeed, single loop algorithms are preferable compared to multi-loop algorithms in many settings, e.g., see \cite{zhang2020single} for a discussion.}
 \vspace{-1mm}
\sa{
\subsection{Notation} Throughout the paper, $\reals_{++}$ denotes the set of positive real numbers, and $\reals_+=\reals_{++}\cup\{0\}$. We adopted arithmetic using the extended reals with the convention that $\frac{1}{0}\triangleq\infty$, $\frac{0^2}{0}\triangleq 0$, $\frac{0}{0^2}\triangleq \infty$.} We use $\norm{\cdot}$ to denote the Euclidean norm. %; hence, $\norm{\cdot}_{\cX}=\norm{\cdot}$ and $\norm{\cdot}_{\cY}=\norm{\cdot}$.\todo{MG: Is the notation $\|\cdots\|_\cX$ necessary if we use Euclidean norm all the time?}
%\xtodo{ I think not necesssary.}
%\sout{We define
%\begin{equation}\label{BRGM}
%    $\mathbf{D}_{\mathcal{X}}(x,\bar{x}) \triangleq \frac{1}{2}\|x-\bar{x}\|^2$,
    %\quad
%    $\mathbf{D}_{\mathcal{Y}}(y,\bar{y}) \triangleq \frac{1}{2}\|y-\bar{y}\|^2$,
%\end{equation}
%\noindent
%for any $x,\bar{x}\in\cX$ and $y,\bar{y}\in\cY$.} %\nsa{I took Bregman distance equation into the text; unless very important, we should take unnumbered equations into the text to save space.} %Let $\mathbb{R}^+$ denote the set of positive real numbers.
\sa{The proximal operator associated with
a proper, %lower-semicontinuous
closed convex %function $f(x):\mathcal{X}\to \sa{\reals\cup\{+\infty\}}$
$f:\cX\to\reals\cup\{\infty\}$ is given by
$\prox{f}(x)\triangleq\argmin_{v\in\mathcal{X}}f(v)+{\frac {1}{2}}\|x-v\|^{2}$, and $\prox{g}(\cdot)$ is defined similarly. \mg{We let $\mathbb{S}^{d}$ denote the set of symmetric $d\times d$ real matrices.}}

\subsection{\sa{Assumptions and Statement of SAPD Algorithm}}\label{subsec-assump}
%Before starting the analysis, we list part of notations and assumptions as follows.
%\mg{We start with introducing the main notations and assumptions that we will make use of throughout the paper.}
\mg{In the following}, we introduce the assumptions needed throughout this paper.
\begin{assumption}\label{ASPT: lipshiz gradient}
\sa{$f:\mathcal{X}\rightarrow \reals\sa{\cup\{\fin{+}\infty\}}$ and $g:\mathcal{Y}\rightarrow \reals\sa{\cup\{\fin{+}\infty\}}$ are \sa{proper, closed, convex functions} with %modules
moduli $\mu_x,\mu_y\geq 0$.} Moreover,
\sa{$\Phi:\cX\times\cY\to\reals$} is %a continuous function
such that
%\begin{itemize}
    %\item

    (i) for any $y\in \dom g\subset \mathcal{Y}$,
    $\Phi(\cdot,y)$ is convex and differentiable; and $\exists L_{xx} \geq 0$, \sa{$\exists L_{xy} > 0$} such that $\forall x,\bar{x}\in \dom f\subset\cX$ and $\forall y,\bar{y}\in \dom g\subset \mathcal{Y}$,
    \begin{equation}\label{LGX1}
        \| \nabla_x \Phi(x,y) - \nabla_x \Phi(\bar{x},\bar{y})\|
        \leq
        L_{xx}\|x- \bar{x}\|+L_{xy}\|y-\bar{y}\|;
    \end{equation}
    %\item

    (ii) for any $x\in\dom f\subset \mathcal{X}$,
    $\Phi(x, \cdot)$ is \sa{concave} and differentiable; and $\exists$ $L_{yx}{>}0$ and $\mg{\exists} L_{yy} \geq 0$ such that $\forall x,\bar{x}\in \dom f\subset \mathcal{X}$ and $\forall y,\bar{y}\in \dom g\subset \mathcal{Y}$,%\vspace*{-1mm}
    \begin{equation}\label{LGY1}
        \| \nabla_y \Phi(x,y) - \nabla_y \Phi(\bar{x},\bar{y})\|
        \leq
        L_{yx}\|x- \bar{x}\|
        + L_{yy}\|y- \bar{y}\|.
    \end{equation}
%\end{itemize}
\end{assumption}
\begin{remark}
\label{rem:strong_convexity}
\fin{In fact, in terms of strong convexity, we \textit{only} need to assume that $\cL$ defined in \eqref{eq:main-problem} is $\mu_x$-convex in $x$ and $\mu_y$-concave in $y$ ($f,g$ may be \textit{merely convex}, e.g., indicator functions). We argue that %it is important to note that 
\cref{ASPT: lipshiz gradient} holds without loss of generality even for this more general setting.}
%where $\cL(x,y)=f(x)+\Phi(x,y)-g(y)$. 
%Indeed, 
%suppose that $f$ and $g$ are \textit{not} strongly convex (they are merely convex, e.g., indicator functions); but, rather $\Phi$ is $\mu_x$-convex and $\mu_y$-concave. 
\fin{Suppose $\Phi$ is $(L_{xx},L_{xy},L_{yx},L_{yy})$-smooth, i.e., \eqref{LGX1} and \eqref{LGY1} hold, and $\Phi$ is $\mu_x$-strongly convex in $x$ and $\mu_y$-strongly concave in $y$, and $f,g$ are proper, closed, \textit{merely convex} functions. %e.g., indicator functions. 
After properly redefining %the functions
$f,g$ and $\Phi$, %one can argue that 
\cref{ASPT: lipshiz gradient} holds for a different representation of the same problem. Indeed, define $f^0,g^0$ and $\Phi^0$ such that for any $(x,y)\in\dom f\times \dom g$, let
\begin{align*}
   &f^0(x)\triangleq f(x)+\frac{\mu_x}{2}\norm{x}^2,\quad g^0(y)\triangleq g(y)+\frac{\mu_y}{2}\norm{y}^2,\quad  \Phi^0(x,y)\triangleq\Phi(x,y)-\frac{\mu_x}{2}\norm{x}^2+\frac{\mu_y}{2}\norm{y}^2,
\end{align*}
and consider $\min_{x\in\cX}\max_{y\in\cY}f^0(x)+\Phi^0(x,y)-g^0(y)$.
The definition of $\Phi^0$ implies that it is $(L^0_{xx},L^0_{xy},L^0_{yx},L^0_{yy})$-smooth, where $L^0_{xx}\triangleq L_{xx}-\mu_x$, $L^0_{yy}\triangleq L_{yy}-\mu_y$, $L^0_{xy}\triangleq L_{xy}$ and $L^0_{yx}=L_{yx}$. Note that $f^0,g^0$ and $\Phi^0$ satisfy \cref{ASPT: lipshiz gradient}. Furthermore, if $f$ and $g$ are prox-friendly functions, i.e., one can compute $\prox{tf}$ and $\prox{tg}$ efficiently for all $t>0$, then $f^0$ and $g^0$ are also prox-friendly. Indeed, given arbitrary $\bar x\in\cX$,  $\bar y\in\cY$ and $t>0$, one has $\prox{tf^0}(\bar{x})=\prox{\frac{t}{t\mu_x+1}f}\Big(\frac{1}{t\mu_x+1}~\bar{x}\Big)$ and $\prox{tg^0}(\bar{y})=\prox{\frac{t}{t\mu_y+1}g}\Big(\frac{1}{t\mu_y+1}~\bar{y}\Big)$.}
% \begin{align*}
% \prox{tf^0}(\bar{x})
% &=\argmin_{x\in\cX}tf^0(x)+\frac{1}{2}\norm{x-\bar{x}}^2,\\
% &=\argmin_{x\in\cX}tf(x)+\frac{t\mu_x}{2}\norm{x}^2+\frac{1}{2}\norm{x-\bar{x}}^2,\\
% &=\argmin_{x\in\cX}\frac{t}{t\mu_x+1}f(x)+\frac{1}{2}\norm{x-\frac{1}{t\mu_x+1}~\bar{x}}^2\\
% &=\prox{\frac{t}{t\mu_x+1}f}\Big(\frac{1}{t\mu_x+1}~\bar{x}\Big);
% \end{align*}}% 
%therefore, $f^0$ and $g^0$ are also prox-friendly functions.}
\end{remark}
\begin{remark}
\label{rem:mcmc}
\fin{We first analyze the error bounds and oracle complexity of SAPD (\cref{sec:sapd}) and its robustness properties (\cref{sec:robustness}) %, i.e., the worst asymptotic squared distance of the iterate sequence to the unique saddle point $(x^*,y^*)$ 
        under the assumption that $\mu_x\mu_y>0$, i.e., for SCSC minimax problems. Later, in \cref{sec:extensions}, we extend these results to MCMC setting, i.e., $\mu_x=\mu_y=0$.}
\end{remark}
 In many %machine learning
 \sa{ML}~applications,  \fin{as passing over the {whole} %problem
 {dataset} to compute %the 
 \rev{a} full gradient may be {computationally %infeasible or
 impractical},} the {full} gradients are estimated \sa{through sampling from data.} %due to \sa{practical reasons, e.g.,
%\mg{In such settings, stochastic algorithms that are based on estimates the gradients from randomly sampled smaller subsets of data are backbone methods.}
 %computational efficiency reasons. %For instance,
 %This is a typical setting in the saddle-point formulations of
 %when the full gradient is approximated from finitely many samples
 \sa{Within \mg{the context of} SP \mg{problems}, this setting arises in supervised learning tasks, e.g., \cite{bottou2018optimization,palaniappan2016stochastic}}.
 %\mg{We adopt a standard stochastic oracle model from the literature, which says that the stochastic estimates $\tilde\nabla_x{\Phi}$ and $\tilde\nabla_y{\Phi}$ are unbiased and have a finite variance. More specifically, we assume the following:}
 \sa{\mg{In} the rest of the paper, we use $\tilde\nabla_x{\Phi}$ and $\tilde\nabla_y{\Phi}$ to denote \mg{such stochastic estimates of the true gradients $\nabla_x{\Phi}$ and $\nabla_y{\Phi}$}.} \sa{Given %such 
 stochastic oracles
 $\tilde\nabla_x{\Phi}$ and $\tilde\nabla_y{\Phi}$, we propose SAPD algorithm to tackle with \eqref{eq:main-problem}, which}
 \mg{%\mg{In this paper, we consider the}
 %the Stochastic Accelerated Primal-Dual~(SAPD) algorithm which extends the accelerated primal-dual (APD) method proposed in~\cite{hamedani2018primal}.
 %to the computational setting that uses stochastic estimates of the gradients instead of full gradients.
 %SAPD algorithm
 is described in Algorithm~\ref{ALG: SAPD}.} %where the iterates $(x_k, y_k)$ are updated based on the stochastic gradient estimates $\Phi(x_k,y_k;\sa{\omega_k^y})$, $\Phi(x_{k-1},y_{k-1};\sa{\omega_{k-1}^y})$, $\Phi(x_k,y_{k+1};\sa{\omega_k^x})$:
 %Algorithm \ref{ALG: SAPD} below.
\mg{We note that} when $\theta$ is zero, SAPD %algorithm is the same with
\sa{reduces to the well-known} stochastic \mg{(proximal)} gradient descent ascent~\mg{(SGDA)} method.

%{, which roughly says that the gradient noise is unbiased when conditioned on the iterates and has a finite variance.}
%In addition, it is not the Nesterov acceleration as the momentum term involves partial gradients not iterates.}
%\todo{Sep,13th xuan: do we need the independent assumption? Check line 840 and 917.}
%\todo{xuan: We define $\omega_k^y$ and $\omega_k^x$ in assumption 2.1, but before algorithm 2.1, we do not define them.}
%\vspace{-1mm}
% \begin{algorithm}
% \caption{Stochastic\sa{Accelerated Primal-Dual} (SAPD) Algorithm}
% {
% \begin{algorithmic}[1]
% \STATE \textbf{Input:} $\{\tau,\sigma,\theta\}$, $(x_0, y_0) \in \mathcal{X} \times \mathcal{Y}$
% \STATE $(x_{-1},y_{-1})\leftarrow (x_0, y_0)$
% \FOR{$k\geq 0$}
% \STATE $\tilde{q}_k \leftarrow\sa{\tilde\nabla_y} \Phi(x_k,y_k\sa{\omega_k^{y,1}}) -\sa{\tilde\nabla_y}\Phi(x_{k-1},y_{k-1}\sa{\omega_k^{y,2}})$
% \STATE $\tilde{s}_k \leftarrow \sa{\tilde\nabla_y} \Phi(x_k,y_k;\sa{\omega_k^{y,1}}) + \theta \tilde{q}_k$
% \STATE $y_{k+1} \leftarrow \prox{\sigma g}(y_k+\sigma \tilde{s_k})$
% \STATE $x_{k+1} \leftarrow \prox{\tau f}(x_k-\tau\sa{\tilde\nabla_x}\Phi(x_k,y_{k+1};\sa{\omega_k^x}))$
% \ENDFOR
% \end{algorithmic}}
% \end{algorithm}
%\vspace{-1mm}
\mg{We %will also
make the following assumption on the statistical nature of the gradient noise.}
\begin{assumption}\label{ASPT: unbiased noise assumption}
%The gradient information is accessible through a stochastic oracle. %At each iteration
There \mg{exist} $\delta_x,\delta_y\geq 0$ such that for
%\xtodo{Should it be There exists or There exist? MG: Fixed}
all $k\geq 0$, given %the current iterate
the SAPD iterates $(x_k,y_k,y_{k+1})$, %\in \mathcal{X}\times\mathcal{Y}\times\mathcal{Y}$,
the \sa{stochastic} gradients $\sa{\tilde{\nabla}_x}\Phi(x_k,y_{k+1};\omega_k^x)$, $\sa{\tilde{\nabla}_y}\Phi(x_k,y_k;\omega_k^y)$ and \sa{random} sequences
$\{\omega_k^x\}$, $\{\omega_k^y\}$ %\sa{are such that}
satisfy %the following conditions: %\nsa{independence assumption is removed now.}  %\nsa{Filtration discussion was unnecessary here, I moved it to the appendix.}
%\begin{itemize}
    %\item
    % (i) $\{\omega_k^x\}_k$ and $\{\omega_k^y\}_k$ are %sequences of
    % independent sequences %random variables
    % and \sa{also} \mg{are} independent from each other;
    %\todo{MG: Do we need them to be independent? Cond. independence if you fix the past history should be enough?}
    %\xtodo{Don't we say they are uncorrelated? MG: For simplicity we could say independent}
    %\item
    (i) $\mathbb{E}[\sa{\tilde{\nabla}_x}\Phi(x_k,y_{k+1};\omega_k^x)|x_k,y_{k+1} ] = \sa{\nabla_x}\Phi(x_k,y_{k+1})$;
    %\item
    (ii) $\mathbb{E}[\sa{\tilde{\nabla}_y}\Phi(x_k,y_{k};\omega_k^y)|x_k,y_k ] = \sa{\nabla_y}\Phi(x_k,y_k)$;
    %\item
    (iii) $\mathbb{E}[\|\sa{\tilde{\nabla}_x}\Phi(x_k,y_{k+1};\omega_k^x) - \sa{\nabla_x}\Phi(x_k,y_{k+1})\|^2 |x_k,y_{k+1} ]\leq \delta_x^2$;
    %\item
    (iv) $\mathbb{E}[\|\sa{\tilde{\nabla}_y}\Phi(x_k,y_{k};\omega_k^y) - \sa{\nabla_y}\Phi(x_k,y_k)\|^2 |x_k,y_k ]\leq \delta_y^2$.
%\end{itemize}
\end{assumption}
\vspace*{-1mm}
\fin{We should point out that we do not make any independence assumption on the {random} sequences
$\{\omega_k^x\}_k$ and $\{\omega_k^y\}_k$.}
\cref{ASPT: unbiased noise assumption} %can be adaptive
\sa{applies} to most unbiased \rev{estimation} situations.
%\todo{Mert:Say assumptions hold on bounded domains $\mathcal{X},\mathcal{Y}$.(xuan: done)}
For example, when %$\mathcal{X},\mathcal{Y}$
\sa{$\dom f\times\dom g$} is compact, \sa{for $\{\omega_k^x\}\subset\cX^*$ and $\{\omega_k^y\}\subset\cY^*$ having zero-mean and finite-variance, the following additive noise model} is a special case of \cref{ASPT: unbiased noise assumption}: $\Tilde{\nabla}_x\Phi(x_k,y_{k+1};\rev{\omega_k^x})= \nabla_x\Phi(x_k,y_{k+1}) + \omega_k^x$, $\Tilde{\nabla}_y\Phi(x_k,y_k;\rev{\omega_k^y}) = \nabla_y\Phi(x_k,y_k)+\omega_k^y$.
\mg{This type of noise %would
arises %for instance
in the context of privacy-preserving %saddle-point
algorithms.} %where
\sa{Indeed, when %the gradients of
$\grad\Phi$ is associated with} the user data, the user would inject additive noise to the gradients for protecting data privacy, see e.g., \cite[Alg. 1]{privacy-gan}, \cite{kuru2020differentially}.
%This 
\rev{Unbiased noise with a finite variance assumption also %arises
holds for %in the primal-dual
\sa{SP formulation of %empirical risk minimization
ERM} problems if the gradients are estimated from mini-batches on bounded domains,
%i.e., from randomly sampled subset of points
e.g., \cite{NEURIPS2018_08048a9c,zhang2017stochastic}.}
\begin{algorithm}[h]
\caption{Stochastic \sa{Accelerated Primal-Dual} (SAPD) Algorithm}
{\label{ALG: SAPD}
{\small
\begin{algorithmic}[1]
\STATE \textbf{Input:} $\{\tau,\sigma,\theta\}$, $(x_0, y_0) \in \mathcal{X} \times \mathcal{Y}$
\STATE $(x_{-1},y_{-1})\leftarrow (x_0, y_0)$
\FOR{$k\geq 0$}
\STATE $\tilde{q}_k \leftarrow\sa{\tilde\nabla_y} \Phi(x_k,y_k;\sa{\omega_k^y}) -\sa{\tilde\nabla_y}\Phi(x_{k-1},y_{k-1};\sa{\omega_{k-1}^y})$
\STATE $\tilde{s}_k \leftarrow \sa{\tilde\nabla_y} \Phi(x_k,y_k;\sa{\omega_k^y}) + \theta \tilde{q}_k$
\STATE $y_{k+1} \leftarrow \prox{\sigma g}(y_k+\sigma \tilde{s}_k)$
\STATE $x_{k+1} \leftarrow \prox{\tau f}(x_k-\tau\sa{\tilde\nabla_x}\Phi(x_k,y_{k+1};\sa{\omega_k^x}))$
\ENDFOR
\end{algorithmic}}}
\end{algorithm}
 %----------------new section-------------------------------------------------
\section{\sa{Performance Guarantees for SAPD}}
\label{sec:sapd}
Under our noise model (\cref{ASPT: unbiased noise assumption}),
%in the result
we next provide performance guarantees for \mg{the} SAPD algorithm.

\begin{theorem}\label{Thm: main result_R1}
\rev{Suppose $\mu_x,\mu_y>0$ and} \sa{Assumptions~\ref{ASPT: lipshiz gradient}
%,~\ref{ASPT: strongly convex concave},~\ref{ASPT: f and g convex}
and~\ref{ASPT: unbiased noise assumption}}
hold, and $\{ x_k,y_k \}_{k\geq 0}$ are generated by SAPD, stated in \cref{ALG: SAPD}, using %parameters
\sa{$\tau, \sigma>0$ and $\theta\geq 0$} %and $\theta\in[0,\rho]$
that satisfy %the following \sa{conditions}:
{\footnotesize
\begin{equation}
    \label{eq: general SAPD LMI_R1}
    \sa{G\triangleq}
  \begin{pmatrix}
    \frac{1}{\tau}+\mu_x - \frac{1}{\rho\tau} & 0 & 0 & 0 & 0\\
  0 & \frac{1}{\sigma}+\mu_y - \frac{1}{\rho\sigma} & (\frac{\theta}{\rho} - 1)L_{yx} & (\frac{\theta}{\rho} - 1)L_{yy} & 0\\
  0 & (\frac{\theta}{\rho} - 1)L_{yx} & \tfrac{1}{\tau} - L_{xx} & 0 & -  \frac{\theta}{\rho}L_{yx}\\
  0& (\frac{\theta}{\rho} - 1)L_{yy} & 0 & \frac{1}{\sigma} - \alpha & -  \frac{\theta}{\rho}L_{yy}\\
  0 & 0 & - \frac{\theta}{\rho}L_{yx} & -  \frac{\theta}{\rho}L_{yy} & \frac{\alpha}{\rho}
\end{pmatrix}\succeq 0
\end{equation}}%
%where
\sa{for some} \rev{$\alpha \in [0,  \tfrac{1}{\sigma})$} and \rev{$\rho \in (0,1)$}.
%\nsa{I included $\frac{1}{\sigma}$ in the bound.}
%\todo{In the next paragraph we set $\pi_1^y=0$ so replace ${\mathbb{R}_+}$ with $\geq 0?$}
%\nsa{$\reals_+$ include 0, $\reals_{++}$ excludes it. We should give these symbols in the notation above.}%can be arbitrary.
%Moreover, if \cref{ASPT: compact} holds in addition to the existing assumptions,
Then \rev{for any $(x_0,y_0)\in\dom f\times\dom g$ and $N\geq 1$, }
%
%\sa{any compact convex set $X \times  Y\subset \dom f\times \dom g$ such that $x_0\in X$ and $y_0\in Y$,
%\nsa{we do not need compact sets here!}
%and for any %\mg{$\eta_x,\eta_y> 0$}%
%\sa{$\eta_x,\eta_y \geq 0$}},
%\nsa{Mert, $\eta_x$ and $\eta_y$ can be 0, we defined $0^2/0=0$ and you can check the proof that when $\delta^2=0$, then $\eta=0$ choice is valid. MG: Ok but then the upper bounds blow up, and are not well defined. the upper bound becomes like 1/0}
%\sa{the following bound,}
%\nsa{I will revert back to $\eta_x$ and $\eta_y$, as otherwise with this particular choice, we get a loose bound for the deterministic case.}
%\nsa{Do we need to assume $(x^*,y^*)\in X\times Y?$}
%\xtodo{I did not see where we use that }
%\nsa{Below I changed the term inside expectation; original version would cause trouble when $\theta=1$.}
\rev{{\small
\begin{equation}
\label{eq:distance-rate_R1}
    \mathbb{E}[d_N^*]
              \leq \rho^{N}\underbrace{\Big(\tfrac{1}{2\tau}\| x_0 - x^*\|^2
    +  \tfrac{1}{2\sigma}\| y_0 - y^*\|^2\Big)}_{ D_{\tau,\sigma}} + \frac{\rho}{1-\rho}~ %(1-\rho^N)
              \underbrace{\Big(\tfrac{\tau}{1+\tau\mu_x} \Xi^x_{\tau,\sigma,\theta} \delta_x^2
              + \tfrac{\sigma}{1+\sigma\mu_y} \Xi^y_{\tau,\sigma,\theta}\delta_y^2\Big)}_{{\Xi}_{\tau,\sigma,\theta}},\vspace*{-2mm}
\end{equation}}}%
%holds for all $N \geq 1$,
where \rev{$(x^*,y^*)$ is the unique saddle point,}
%{\small
%\begin{align}
    %&
    \rev{$d_N^* \triangleq \tfrac{1}{2\tau}\|x_N-x^*\|^2 +  \tfrac{1}{2\sigma}\left( 1 - \alpha\sigma\right)\|y_N-y^*\|^2$},
    %  $\sa{ D_{\tau,\sigma}}\triangleq \tfrac{1}{2\tau}\| x_0 - x^*\|^2
    % +  \tfrac{1}{2\sigma}\| y_0 - y^*\|^2$,
    %\nonumber
    %\\
    %&
    %\nonumber
    %\\
    %&
    % \rev{${\Xi}_{\tau,\sigma,\theta}
    % \triangleq
    % \tfrac{\tau}{1+\tau\mu_x} \Xi^x_{\tau,\sigma,\theta} \delta_x^2
    %           + \tfrac{\sigma}{1+\sigma\mu_y} \Xi^y_{\tau,\sigma,\theta}\delta_y^2$} 
    %with
     %\nonumber\\
     %&
     $\sa{\Xi^x_{\tau,\sigma,\theta}} \triangleq 1 + \tfrac{\sigma\sa{\theta}(1+\theta)L_{yx}}{2(1+\sigma\sa{\mu_y})}$ and %\nonumber\\
     %&
      $\sa{\Xi^y_{\tau,\sigma,\theta}} \triangleq \tfrac{\tau\theta(1+\theta)L_{yx}}{
     \sa{2(1+\tau\mu_x)}}+\Big(1+2\theta + \tfrac{\theta + \sigma\theta(1+\theta)L_{yy}}{1+\sigma\mu_y}  +  \tfrac{\tau\sigma\theta(1+\theta)L_{yx}L_{xy}}{(1+\tau\mu_x)(1+\sigma\mu_y)} \Big)
    (1+2\theta)$. 
    
    \fin{Moreover, whenever $\delta_x=\delta_y=0$, 
    %\xtodo{maybe $\cG(\bar{x}_N,\bar{y}_N)\leq\frac{1}{K_N(\rho)} D_{\tau,\sigma}$is better?} 
    the gap metric defined in \eqref{eq:gap} satisfies $\cG(\bar{x}_N,\bar{y}_N)\leq %(1-\rho)\rho^{N-1}
     D_{\tau,\sigma}/K_N(\rho)$ for all $N\geq 1$, where $D_{\tau,\sigma}$ is defined in \eqref{eq:distance-rate_R1}, %$\bar{z}_N=(\bar{x}_N,\bar{y}_N)$ such that $\bar{z}_N 
     $(\bar{x}_N,\bar{y}_N)= \tfrac{1}{K_N(\rho)}\sum_{k=1}^{N}\rho^{-k+1}(x_k,y_k)$,}
     %z_k$ \rev{with $z_k=(x_k,y_k)$ for $k\geq 1$} 
     and $K_N(\rho)\triangleq \sum_{k=0}^{N-1}\rho^{-k}=\frac{1-\rho^N}{1-\rho}\rho^{-N+1}$.\looseness=-1
\end{theorem}
\begin{proof}
\rev{See \cref{sec:proof}.}
\end{proof}
\looseness=-1
\begin{remark}
\fin{Consider the stochastic case, i.e., $\delta_x, \delta_y>0$. \cref{eq:distance-rate_R1} implies that when $\alpha\leq \tfrac{c}{\sigma}$ for some $c\in(0,1)$, {the weighted squared distance in expectation, 
$\mathbb{E}[\cD(x_N,y_N)]$},
is bounded by a sum of two terms: {\emph{bias term}} %$\tfrac{1}{K_N(\rho)}\Omega_{\tau,\sigma,\theta}$ 
{$\mathcal{O}(\max\{\tau,\sigma\}\rho^N D_{\tau,\sigma}$)} that goes to zero as $N\to\infty$ and a constant {\emph{variance term}} {$\mathcal{O}(\max\{\tau,\sigma\}\frac{\rho}{1-\rho}\Xi_{\tau,\sigma,\theta})$}, which can be controlled by properly selecting $\tau$, $\sigma$ and $\theta$. Indeed,
{%The variance 
the term $\Xi_{\tau,\sigma,\theta}$ %is determined by two constants $\Xi^x_{\tau,\sigma,\theta}$ and $\Xi^y_{\tau,\sigma,\theta}$ are
depends on algorithm parameters $\{\tau,\sigma,\theta\}$ %and two other free parameters $\eta_x,\eta_y>0$. For this scenario, we set $\eta_x=\frac{1}{\tau}+\mu_x$ and $\eta_y=\frac{1}{\sigma}+\mu_y$. This choice of parameters implies 
in such a way that 
%as the primal and dual step sizes approach 0, %the variance term 
%$\Xi_{\tau,\sigma,\theta}$ vanishes as expected, i.e., 
as $\tau\to 0$ and $\sigma\to 0$, we have $\Xi_{\tau,\sigma,\theta}\to 0$.} Furthermore, there exists $\bar\theta\in(0,1)$ depending only on problem parameters such that for all $\theta\geq \bar\theta$, there is a solution to \eqref{eq: general SAPD LMI_R1} satisfying $\tau=\cO(1-\theta)$, $\sigma=\cO(1-\theta)$ and $\rho=\theta$ --see \cref{sec:parameter_choice}, such that the \emph{variance term}
%$\frac{\rho}{1-\rho}\Xi_{\tau,\sigma,\theta}$ 
can be made arbitrarily small as $\theta\to 1$.} %\xtodo{Should we use $d_N$ instead of $\cD$here?}
\looseness=-1
%Lipschitz constant $\{L_{xy},L_{yx}\}$, and strongly convex modules $\{\mu_x,\mu_y\}$. When $\tau$ and $\sigma$ go to zero, they  vanish as well.
\end{remark}
\begin{remark}
\sa{%Since $\theta\in[0,\rho]$,
If we set $\theta = \rho$} in \cref{eq: general SAPD LMI_R1}, we obtain a simpler \sa{matrix inequality}: %\vspace*{-1mm}
{%\footnotesize
\begin{equation}
\label{Condition: SAPD simple LMI system}
\min\{\tau\mu_x,~\sigma\mu_y\}\geq \frac{1-\theta}{\theta},
\qquad %\sigma\mu_y\geq \frac{1-\theta}{\theta},\quad
% \left(\begin{smallmatrix}
\begin{pmatrix}
  \tfrac{1}{\tau} - L_{xx}  & 0 & -L_{yx}\\
 0 & \tfrac{1}{\sigma}- \alpha  & -L_{yy}\\
 - L_{yx} & - L_{yy} & \frac{\alpha}{ \theta}
\end{pmatrix}
%\end{smallmatrix}\right)
\succeq 0. %\vspace*{-1mm}
\end{equation}}%
% \begin{subequations}
% \label{Condition: SAPD simple LMI system}
% {\begin{equation}
% \label{Condition: SAPD simple LMI 1}
%     \begin{pmatrix}
%   \tfrac{1}{\tau} - L_{xx}  & 0 & - L_{yx}\\
%  0 & \tfrac{1}{\sigma}- \alpha  & - L_{yy}\\
%  - L_{yx} & - L_{yy} & \tfrac{\alpha}{ \theta}
% \end{pmatrix}
% \succeq 0,
% \end{equation}}
% \begin{align}
%     \label{Condition: SAPD simple LMI 2}
%     \sa{\tau\mu_x\geq \frac{1-\theta}{\theta}},\quad \sa{\sigma\mu_y\geq \frac{1-\theta}{\theta}.}
% %\label{Condition: noisy LMI 3  _bounded version}
% \end{align}
% \end{subequations}
%\xtodo{Only complexity bound based on this simple system.}
The %rest of the
SAPD complexity analysis is mainly based on the simpler system in \cref{Condition: SAPD simple LMI system}. \rev{On the other hand, when $\theta = 0$, \sa{SAPD %degenerates to {Stochastic} Gradient Descent Ascent~\sa{(SGDA)} method,
reduces to SGDA,} of which step size conditions can be obtained immediately from \cref{eq: general SAPD LMI_R1} by setting $\theta=\alpha = 0$ --see \cref{SGDA section}.}
\end{remark}
\subsection{Proof of Theorem~\ref{Thm: main result_R1}}
\label{sec:proof}
%\fin{The proofs of lemmas are provided in the appendix.}
We first \fin{provide some key inequalities} %below 
%which derive some key 
%with useful inequalities %below 
for the SAPD iterate sequence $\{x_k,y_k\}_{k\geq 0}$ generated
by
%\todo{We need to make sure APD Alg. is defined}
Algorithm~\ref{ALG: SAPD}, \fin{the omitted proofs are provided in the appendix}. Let %\sa{For $k\geq 0$, define}
%let ${q}_k$ and $s_k$ be the %deterministic
%noiseless versions of $\tilde{q}_k$ and $\tilde{s}_k$, i.e.,
{\small
\begin{align}
\label{eq:qksk}
{q}_k \triangleq \nabla_y \Phi(x_k,y_k) - \nabla_y\Phi(x_{k-1},y_{k-1}),\qquad s_k \triangleq \nabla_y \Phi(x_k, y_k) + \theta q_k,\quad \forall~k\geq 0.
\end{align}}%
\sa{Recall $x_{-1} = x_0,\;y_{-1} = y_0$, thus $q_0=\mathbf{0}$; and for $k\geq 0$, \cref{ASPT: lipshiz gradient} implies that
%\xtodo{This part can be shorten}
{\small
\begin{equation}
           \norm{q_{k+1}} \leq  L_{yx} \|x_{k+1} - x_{k} \|+  L_{yy} \|y_{k+1} - y_{k} \|. \label{INEQ: Cauchy Ineqaulity 1}
           \vspace*{-4mm}
\end{equation}}}%
% \todo{XUAN: do we still use this $\pi_1$ and $\pi_2$?}
% Recall the inequality
% $$
% \|a\|\|b\|\leq \tfrac{c}{2}\|a\|^2+\tfrac{1}{2c}\|b\|^2,\;\forall c>0,\;a,b\in \reals^n.
% $$
% Then for any $\pi_1,\pi_2>0$, \eqref{INEQ: Cauchy Ineqaulity 1} implies that
% \begin{equation}\label{INEQ: Cauchy Ineqaulity 2}
%     \begin{aligned}
%           \theta \langle q_k, y_{k+1} - y_{k} \rangle
%             %\leq &  \theta L_{yx} \|x_{k} - x_{k-1} \|\| y_{k+1} - y_{k} \| + \theta   L_{yy} \|y_{k} - y_{k-1} \| \| y_{k+1} - y_{k} \| \\
%             \leq & \frac{\pi_1\theta L_{yx}}{2}   \|x_{k} - x_{k-1} \|^2 +\frac{\pi_2\theta L_{yy}}{2}   \|y_{k} - y_{k-1} \|^2\\ & + (\frac{ \theta L_{yx}}{2\pi_1} +  \frac{ \theta L_{yy}}{{2}\pi_2}) \|y_{k+1} - y_{k} \|^2.
%     \end{aligned}
% \end{equation}
\vspace*{-2mm}
\begin{lemma}\label{lem: basic lemma for deterministic case}
Let $\{x_k,y_k\}_{k\geq 0}$ be SAPD iterates generated according to Algorithm {\ref{ALG: SAPD}}. %~\ref{APD}.
Then %following inequality holds
for all \sa{$x\in\dom f\subset \cX$, $y\in\dom g\subset \cY$,} and $k\geq 0$,
\vspace*{-2mm}
%\xtodo{It is better if we avoid using inline eqaution}
{\small
\begin{equation}\label{D1}
    \begin{aligned}
       \mathcal{L}( & x_{k+1}, y)  - \mathcal{L}(x, y_{k+1})
        \\
    \leq
    &-\langle q_{k+1}, y_{k+1} - y \rangle + \theta \langle q_k, y_{k} - y \rangle
    + \Lambda_k(x,y) - \Sigma_{k+1}(x,y)+ \Gamma_{k+1}+\sa{\varepsilon^x_k+\varepsilon^y_k},\vspace*{-2mm}
%   \\
%   & -
%       \langle  \tilde{\nabla}_x \Phi(x_k, y_{k+1};\omega_k^x) - \nabla_x \Phi(x_k, y_{k+1}) , x_{k+1} - x \rangle
%         +
%       \langle \tilde{s}_k -s_k,
%       y_{k+1} - y \rangle,
\end{aligned}
\end{equation}}%
where \sa{$\varepsilon^x_k\triangleq
       \langle  \tilde{\nabla}_x \Phi(x_k, y_{k+1};\omega_k^x) - \nabla_x \Phi(x_k, y_{k+1}),~x-x_{k+1} \rangle$ and $\varepsilon^y_k\triangleq\langle \tilde{s}_k -s_k, y_{k+1} - y \rangle$, $q_k$ and $s_k$ are defined as in \eqref{eq:qksk}, and}
%\sa{$q_k = \nabla_y \Phi(x_k, y_k) - \nabla_y \Phi(x_{k-1}, y_{k-1})$, $s_k= \nabla_y \Phi(x_k, y_k) + \theta q_k$, and}
\rev{\small
\begin{equation*}
    \begin{aligned}
        \Lambda_k(x,y) \triangleq  &\frac{1}{2\tau} \|x - x_{k}\|^2 + \frac{1}{2\sigma}\|y-y_{k}\|^2,\quad
        \Sigma_{k+1}(x,y) \triangleq  (\frac{1}{2\tau} + \frac{\mu_x}{2})\|x - x_{k+1}\|^2 + (\frac{1}{2\sigma} + \frac{\mu_y}{2})\|y-y_{k+1}\|^2,\\
     \Gamma_{k+1} \triangleq &
     (\frac{ L_{xx} }{2} - \frac{1}{2\tau}) \| x_{k+1} - x_k \|^2 - \frac{1}{2\sigma}\|y_{k+1}- y_k\|^2  + \theta (L_{yx} \|x_{k} - x_{k-1} \|+ L_{yy} \|y_{k} - y_{k-1} \|) \| y_{k+1} - y_{k} \|.
     \end{aligned}
\end{equation*}}%
\vspace*{-3mm}
\end{lemma}
% \begin{proof}
% \fin{xxxxx}
% \end{proof}
%Next, based on the above \sa{inequality}, 
Next, we %prove
give %\mg{an}
two intermediate results %that assists 
to bound the variance \sa{of the SAPD iterate sequence.}
% \todo[inline]{NSA: It may be better to state it as: $\|\prox{cf}(x) - \prox{cf}(y) \|^2 \leq \frac{1}{1+c\mu}\fprod{\prox{cf}(x) - \prox{cf}(y),~x-y}.$}

\begin{lemma}\label{lemma: intermediate noisy bound}
Let $\{x_k,y_k\}_{k\geq 0}$ be SAPD iterates generated according to Algorithm {\ref{ALG: SAPD}}.
%\sa{Let $q_k \triangleq \nabla_y \Phi(x_k, y_k) - \nabla_y \Phi(x_{k-1}, y_{k-1})$, $s_k\triangleq \nabla_y \Phi(x_k, y_k) + \theta q_k$ for $k\geq 0$.} If
\sa{For $k\geq 0$, let $q_k$ and $s_k$ be defined as in \eqref{eq:qksk}, and
%we also define
let} %auxiliary sequences:}
{\small
\begin{equation*}%\label{eq: intermediate iteration sequence}
    \begin{aligned}
       &{\hat{x}_{k+1}} \triangleq   \prox{\tau f}\left({x_{k}- \tau{\nabla}_x \Phi(x_k, y_{k+1}})\right),\ 
       \hat{y}_{k+1} \triangleq  \prox{\sigma g}\left(y_{k} + \sigma s_k\right),\\
       &{\hat{\hat{x}}_{k+1}  \triangleq   \prox{\tau f}\left({x_{k}- \tau{\nabla}_x \Phi(x_k, \hat{y}_{k+1}})\right)},\  \hat{\hat{y}}_{k+1}\triangleq \prox{\sigma \sa{g}}\left( \hat{y}_{k} + \sigma(1+\theta){\nabla}_y \Phi(\hat{\hat{x}}_{k}, \hat{y}_{k}) - \sigma\theta{\nabla}_y \Phi(x_{k-1}, y_{k-1})  \right),
    \end{aligned}
\end{equation*}}%
% and
% \begin{equation}\label{eq: diff between true and estimate gradient}
%     \begin{aligned}
%               \Delta^{x}_k {\triangleq} \tilde{\nabla}_x \Phi(x_k, y_{k+1};\omega_k^x) - \nabla_x \Phi(x_k, y_{k+1}),\\
%     \Delta^{y}_k {\triangleq} \tilde{\nabla}_y \Phi(x_k, y_{k};\sa{\omega_k^y}) - \nabla_x \Phi(x_k, y_{k}),
%     \end{aligned}
% \end{equation}
then the following inequalities hold for all $k\geq 0$:
{\footnotesize
\begin{subequations}
\begin{align}
    \|x_{k+1} - \hat{x}_{k+1} \| &\leq \frac{\tau}{1+\tau\mu_x}\| \Delta_k^{x}\|,%\label{eq:x_k - x_hat_k}\\
    \qquad \|y_{k+1} - \hat{y}_{k+1}\| \leq \frac{\sigma}{1+\sigma\mu_y}\left((1+\theta)\|\Delta^{y}_k\| + \theta\|\Delta^{y}_{k-1}\|\right),\label{eq:y_k - y_hat_k}\\
    \|{y}_{k+1} - \hat{\hat{y}}_{k+1}\| &\leq  \frac{\sigma}{1+\sigma\mu_y}\left(
    (1+\theta)\|\Delta^{y}_k\| + \theta\|\Delta^{y}_{k-1}\|
    +
    \frac{\tau(1+\theta)L_{yx}}{1+\tau\mu_x}\|\Delta^{x}_{k-1}\|
    \right.\label{eq:y_k - y_hat_hat_k}\\
    & \left.
    + \left(\frac{1 + \sigma(1+\theta)L_{yy}}{1+\sigma\mu_y}  +  \frac{\tau\sigma(1+\theta)L_{yx}L_{xy}}{(1+\tau\mu_x)(1+\sigma\mu_y)} \right)\left((1+\theta)\|\Delta^{y}_{k-1}\| + \theta\|\Delta^{y}_{k-2}\|\right)
    \right), \nonumber
    \end{align}
\end{subequations}}%
where
\rev{\small
%\begin{subequations}
%\begin{align}
    $\Delta^{x}_k {\triangleq} \tilde{\nabla}_x \Phi(x_k, y_{k+1};\omega_k^x) - \nabla_x \Phi(x_k, y_{k+1})$ and $\Delta^{y}_k {\triangleq} {\tilde\nabla_y}{\Phi}(x_k, y_k;\omega_k^y) - \nabla_y \Phi(x_k, y_k)$.
    %\nsa{Consider using $\Delta^{x}_k$ and $\Delta^{y}_k$.}
%\label{eq:Delta_xy}
%\end{align}
%\end{subequations}
}%
\end{lemma}
% \begin{proof}
% \fin{xxxx}
% \end{proof}

%In the next
\fin{The next result, which will be used in the variance analysis for SAPD, follows from \cref{lemma: intermediate noisy bound}.} 
%we provide some inequalities to bound the variance later} in our analysis.
\begin{lemma}\label{Lemma: final noisy bound}
Let $\{x_k,y_k\}_{k\geq 0}$ be SAPD iterates generated according to Algorithm {\ref{ALG: SAPD}}.
% Recall that
% \begin{equation*}
% \begin{aligned}
%       \Delta^{x}_k = \tilde{\nabla}_x \Phi(x_k, y_{k+1};\omega_k^x) - \nabla_x \Phi(x_k, y_{k+1}),\\
%     \Delta^{y}_k = \tilde{\nabla}_y \Phi(x_k, y_{k};\omega_k^x) - \nabla_x \Phi(x_k, y_{k}),
% \end{aligned}
% \end{equation*}
The following inequalities hold for all $k\geq 0$:\vspace*{-2mm}
{\footnotesize
\begin{equation*} %\label{ineq: noisy bound}
   \begin{aligned}
    &\mathbb{E}\left[|\langle \Delta^{x}_k,  \sa{\hat{x}_{k+1}-x_{k+1}} \rangle|\right] \leq \frac{\tau}{1+\tau\mu_x}\delta_x^2,\qquad
     \mathbb{E}\left[|\langle \Delta^{y}_k,  y_{k+1} - \hat{y}_{k+1}\rangle|\right] \leq \frac{\sigma(1+2\theta)}{1+\sigma\mu_y}\delta_y^2,\\
    &\mathbb{E}\left[ |\langle \Delta^{y}_{k-1},  \sa{\hat{\hat{y}}_{k+1}-y_{k+1}} \rangle| \right]
    \leq
    \frac{\sigma}{1+\sigma\mu_y} \left[  \left( \left(1 + \frac{1 + \sigma(1+\theta)L_{yy}}{1+\sigma\mu_y}  +  \frac{\tau\sigma(1+\theta)L_{yx}L_{xy}}{(1+\tau\mu_x)(1+\sigma\mu_y)} \right)\right.\right.
     \\
    &\hspace*{4cm}\cdot\left.\left.
    (1+2\theta)
    + \frac{\tau(1+\theta)L_{yx}}{%(1+\sigma\mu_y)
     \sa{2(1+\tau\mu_x)}}\right)\delta_y^2 + \frac{\tau(1+\theta)L_{yx}}{%(1+\sigma\mu_y)
     \sa{2(1+\tau\mu_x)}}\delta_x^2\right].
    \end{aligned}
    \vspace*{-1mm}
\end{equation*}}%
% Furthermore, if %the parameter $\theta$ for SAPD is equal to zero,
% $\theta=0$, then
% $\langle \Delta^{y}_k,  y_{k+1} - \hat{y}_{k+1}  \rangle \leq \frac{\sigma}{1+\sigma\mu_y}\|\Delta^{y}_k \|^2$.\todo{Why do we state this?}
\end{lemma}
\sa{Before we move on to prove our main result in \cref{Thm: main result_R1}, we give two technical lemmas that help us simplify the SAPD parameter selection rule to \rev{the matrix inequality in \eqref{eq: general SAPD LMI_R1}}.}
\begin{lemma}
\label{lem:equivalent_systems}
\sa{Let $G\in\reals^{5\times 5}$ be the matrix on the left-hand-side of \eqref{eq: general SAPD LMI_R1}, and}
% \begin{equation*}
%   G_1 \triangleq \begin{pmatrix}
%     \frac{1}{\tau}+\mu_x - \frac{1}{\rho\tau} & 0 & 0 & 0 & 0\\
%   0 & \frac{1}{\sigma}+\mu_y - \frac{1}{\rho\sigma} & (\frac{\theta}{\rho} - 1)L_{yx} & (\frac{\theta}{\rho} - 1)L_{yy} & 0\\
%   0 & (\frac{\theta}{\rho} - 1)L_{yx} & \tfrac{1}{\tau} - L_{xx} & 0 & -  \frac{\theta}{\rho}L_{yx}\\
%   0& (\frac{\theta}{\rho} - 1)L_{yy} & 0 & \frac{1}{\sigma} - \alpha & -  \frac{\theta}{\rho}L_{yy}\\
%   0 & 0 & - \frac{\theta}{\rho}L_{yx} & -  \frac{\theta}{\rho}L_{yy} & \frac{\alpha}{\rho}
% \end{pmatrix},
% \end{equation*}
{\scriptsize
\begin{equation*}
 \sa{G'} \triangleq  \begin{pmatrix}
    \frac{1}{\tau}+\mu_x - \frac{1}{\rho\tau} & 0 & 0 & 0 & 0\\
  0 & \frac{1}{\sigma}+\mu_y - \frac{1}{\rho\sigma} & \sa{-|1-\frac{\theta}{\rho}|}~L_{yx} & \sa{-|1-\frac{\theta}{\rho}|}~L_{yy} & 0\\
  0 & \sa{-|1-\frac{\theta}{\rho}|}~L_{yx} & \tfrac{1}{\tau} - L_{xx} & 0 & -  \frac{\theta}{\rho}L_{yx}\\
  0& \sa{-|1-\frac{\theta}{\rho}|}~L_{yy} & 0 & \frac{1}{\sigma} - \alpha & -  \frac{\theta}{\rho}L_{yy}\\
  0 & 0 & - \frac{\theta}{\rho}L_{yx} & -  \frac{\theta}{\rho}L_{yy} & \frac{\alpha}{\rho}
\end{pmatrix},
\end{equation*}}%
then \sa{$G\succeq 0$ %is equivalent to
\rev{if and only if} $G'\succeq 0$}.
\end{lemma}
\begin{proof}
$\forall ~ \mathbf{y} = (y_1,y_2,y_3,y_4,y_5)^\top\in \mathbb{R}^5$, letting $\tilde{\mathbf{y}}= (y_1,-y_2,y_3,y_4,y_5)^\top$, we have
{\footnotesize
\begin{equation*}
    \mathbf{y}^\top G' \mathbf{y} =
    \begin{cases}
    \mathbf{y}^\top \sa{G} \mathbf{y} & \text{if} ~\sa{\theta\leq \rho,}\\
    \tilde{\mathbf{y}}^\top G \tilde{\mathbf{y}} & \text{else;}
    \end{cases}\quad
    \mathbf{y}^\top G \mathbf{y} =
    \begin{cases}
    \mathbf{y}^\top G' \mathbf{y} & \text{if} ~\sa{\theta\leq \rho},\\
    \tilde{\mathbf{y}}^\top G' \tilde{\mathbf{y}} & \text{else.}
    \end{cases}
\end{equation*}}%
% It follows that, if one of $G_1$ and $G_2$ is semi-positive definite, then another one must be also semi-positive definite. Therefore,
\sa{Thus, $G\succeq 0$ is equivalent to $G'\succeq 0$.}
\end{proof}
\begin{lemma}\label{lemma: sub positive matrix}
Suppose %all assumptions in \cref{Thm: main result} and
\sa{the parameters \sa{$\tau, \sigma>0$} and \rev{$\theta\geq 0$}
satisfy \cref{eq: general SAPD LMI_R1} \sa{for some} $\alpha \in [0,  \tfrac{1}{\sigma})$ and $\rho \in (0,1]$, then it follows that}
%\nsa{Add a brief explanation why we need this lemma. xuan:done. stated before lemma A.5.}
%\xtodo{Done before Lemma A.5}
%the following matrix is positive definite:
{\footnotesize
\begin{equation}
    \label{eq: sub matrix of general SAPD LMI}
  \sa{G''} \triangleq \begin{pmatrix}
  \frac{1}{\sigma}(1-\frac{1}{\rho})+\mu_y +\frac{\alpha}{\rho} & (-| 1 - \tfrac{\theta}{\rho} | - \tfrac{\theta}{\rho})L_{yx} & (-| 1 - \tfrac{\theta}{\rho} | - \tfrac{\theta}{\rho})L_{yy} \\
 (-| 1 - \tfrac{\theta}{\rho} | - \tfrac{\theta}{\rho}) L_{yx} & \tfrac{1}{\tau} - L_{xx} & 0 \\
 (-| 1 - \tfrac{\theta}{\rho} | - \tfrac{\theta}{\rho})L_{yy} & 0 & \frac{1}{\sigma} - \alpha \\
\end{pmatrix}\sa{\succeq 0}.
\end{equation}}%
\end{lemma}
\begin{proof}
Since \cref{eq: general SAPD LMI_R1} holds, $\forall \bx=[x_1~x_2~x_3]^\top\in \mathbb{R}^3$, by \cref{lem:equivalent_systems} we have that
%{\small
%\begin{equation}
%\begin{aligned}
% \begin{pmatrix}
%   x_1 \\
%   x_2 \\
%   x_3
% \end{pmatrix}^\top
$\bx^\top
G''
%   \begin{pmatrix}
%   \frac{1}{\sigma}(1-\frac{1}{\rho})+\mu_y +\frac{\alpha}{\rho} & (-| 1 - \tfrac{\theta}{\rho} | - \tfrac{\theta}{\rho})L_{yx} & (-| 1 - \tfrac{\theta}{\rho} | - \tfrac{\theta}{\rho})L_{yy} \\
%  (-| 1 - \tfrac{\theta}{\rho} | - \tfrac{\theta}{\rho}) L_{yx} & \tfrac{1}{\tau} - L_{xx} & 0 \\
%  (-| 1 - \tfrac{\theta}{\rho} | - \tfrac{\theta}{\rho})L_{yy} & 0 & \frac{1}{\sigma} - \alpha \\
% \end{pmatrix}
% \begin{pmatrix}
%   x_1 \\
%   x_2 \\
%   x_3
% \end{pmatrix}
\bx
=
% \begin{pmatrix}
% \sa{0}\\
%   x_1 \\
%   x_2 \\
%   x_3 \\
%   x_1
% \end{pmatrix}^\top
{\bx'}^\top
G'
%   \begin{pmatrix}
%   \frac{1}{\sigma}(1-\frac{1}{\rho})+\mu_y & -| 1 - \tfrac{\theta}{\rho} | L_{yx} & -| 1 - \tfrac{\theta}{\rho} | L_{yy} & 0\\
%   -| 1 - \tfrac{\theta}{\rho} | L_{yx} & \tfrac{1}{\tau} - L_{xx} & 0 & -  \frac{\theta}{\rho}L_{yx}\\
%   -| 1 - \tfrac{\theta}{\rho} | L_{yy} & 0 & \frac{1}{\sigma} - \alpha & -  \frac{\theta}{\rho}L_{yy}\\
%   0 & - \frac{\theta}{\rho}L_{yx} & -  \frac{\theta}{\rho}L_{yy} & \frac{\alpha}{\rho}
% \end{pmatrix}
% \begin{pmatrix}
% \sa{0}\\
%   x_1 \\
%   x_2 \\
%   x_3 \\
%   x_1
% \end{pmatrix}
\bx'\geq 0$,
%\end{aligned}
%\end{equation}}%
where $G'$ is defined in \cref{lem:equivalent_systems} and $\bx'=[0~x_1~x_2~x_3~x_1]^\top$.
%Therefore, matrix \cref{eq: sub matrix of general SAPD LMI} is positive definite.
\end{proof}

% {\begin{equation*}
%     \begin{aligned}
%     & \mathbb{E}\left[\langle \Delta^{x}_k, x_{k+1} - \hat{x}_{k+1} \right]\rangle \leq \frac{\tau}{1+\tau\mu_x}\delta_x^2,\;\;
%      \mathbb{E}\left[\langle \Delta^{y}_k,  y_{k+1} - \hat{y}_{k+1}\rangle\right] \leq \frac{\sigma(1+2\theta)}{1+\sigma\mu_y}\delta_y^2,\\
%     &  \mathbb{E}\left[ \langle \Delta^{y}_{k-1},  y_{k+1} - \hat{\hat{y}}_{k+1}\rangle \right]
%     \leq
%     \frac{\sigma}{1+\sigma\mu_y} \left[  \left( \left(1 + \frac{1 + \sigma(1+\theta)L_{yy}}{1+\sigma\mu_y}  +  \frac{\tau\sigma(1+\theta)L_{yx}L_{xy}}{(1+\tau\mu_x)(1+\sigma\mu_y)} \right)* \right.\right.
%      \\
%     & \left.\left.
%     (1+2\theta)
%     + \frac{\tau(1+\theta)L_{yx}}{%(1+\sigma\mu_y)
%      \sa{2(1+\tau\mu_x)}}\right)\delta_y^2 + \frac{\tau(1+\theta)L_{yx}}{%(1+\sigma\mu_y)
%      \sa{2(1+\tau\mu_x)}}\delta_x^2\right].
%     \end{aligned}
% \end{equation*}}
% which completes the proof.

\sa{Finally, with the following observation, we will be ready to proceed to the proof of Theorem~\ref{Thm: main result_R1}. Let $\{\cF_k^x\}$ and $\{\cF_k^y\}$ be the filtrations such that $\mathcal{F}^x_k \triangleq \mathcal{F}(\{x_i\}_{i=0}^{k},\{y_i\}_{i=0}^{k+1})$ and $\mathcal{F}^y_k \triangleq \mathcal{F}(\{x_i\}_{i=0}^{k},\{y_i\}_{i=0}^{k})$ denote the $\sigma$-algebras generated by the random variables in their arguments.} %It should be noted that,
\mg{A consequence of Assumption~\ref{ASPT: unbiased noise assumption} is that}
%if a random variable
\sa{for $\mathcal{F}^x_k$-measurable random variable $v$, i.e.,} $v\in\mathcal{F}^x_k$, we have that %\nsa{Why do we need compactness? MG: We do not need it, i meant to refer to the other (unbiased noise) assumption}
$\mathbb{E}\left[\langle \tilde{\nabla}\Phi_x(x_k,y_{k+1};\omega_k^x) - \nabla\Phi_x(x_k,y_{k+1}),v \rangle \right]= 0$; similarly,
for $v\in\mathcal{F}^y_k$, it holds that
$\mathbb{E}\left[\langle \tilde{\nabla}\Phi_y(x_k,y_{k};\omega_k^y) - \nabla\Phi_y(x_k,y_k),v \rangle \right]= 0$. We are now ready to give the proof of \cref{Thm: main result_R1}.
\paragraph{Proof of \cref{Thm: main result_R1}}
Fix arbitrary \rev{$(x,y)\in \dom f\times \dom g$}. Since $(x_{k+1},y_{k+1})\in\dom f\times \dom g$, using the \mg{concavity} 
of $\mathcal{L}(x_{k+1}, \cdot)$
and the \mg{convexity} %concavity
of $\mathcal{L}(\cdot, y_{k+1})$, we get
\vspace{-2mm}
{\small
\begin{align}
\label{eq:jensen _bounded version}
 K_{N}(\rho)  \left(\mathcal{L}(\bar{x}_{N}, y)  - \mathcal{L}(x, \bar{y}_{N}) \right)
            \leq
            \sum_{k=0}^{N-1}\rho^{-k} \left( \mathcal{L}(x_{k+1}, y) - \mathcal{L}(x, y_{k+1}) \right),
            \;\sa{\forall}\rho\in \rev{(0,1)}.
\end{align}}%
%\newpage
%\vspace{-4mm}
\noindent Thus, if we multiply $\rho^{-k}$ for both sides of %\eqref{INEQ: pure gap _bounded version}
\eqref{D1} and \sa{sum the resulting inequality from $k=0$ to $N-1$, then using \eqref{eq:jensen _bounded version} we get}
{\small
\begin{equation}\label{INEQ: difference of gap  _bounded version}
    \begin{aligned}
            K_{N} (\rho)  \left(\mathcal{L}(\bar{x}_{N}, y)  - \mathcal{L}(x, \bar{y}_{N}) \right)
            %\leq & \sum_{k=0}^{N-1}\rho^{-k} \left( \mathcal{L}(x_{k+1}, y) - \mathcal{L}(x, y_{k+1}) \right) \\
            \leq &
             \sum_{k=0}^{N-1}\rho^{-k} \Big(
                \underbrace{ -\langle q_{k+1}, y_{k+1} - y \rangle + \theta \langle q_k, y_{k} - y \rangle}_{\text{\bf part 1}} + \Lambda_k \sa{(x,y)} - \Sigma_{k+1} \sa{(x,y)}+ \Gamma_{k+1}
            \\
             & \quad  \underbrace{-
           \langle  \tilde{\nabla}_x \Phi(x_k, y_{k+1};\omega_k^x) - \nabla_x \Phi(x_k, y_{k+1}) , x_{k+1} - x \rangle}_{\text{\bf part 2} }
            +
           \underbrace{\langle \tilde{s}_k -s_k,
          y_{k+1} - y \rangle}_{ \text{\bf part 3} }
          \Big).
    \end{aligned}
\end{equation}}%
\sa{Using Cauchy–Schwarz inequality and \eqref{INEQ: Cauchy Ineqaulity 1}} leads to
{\small
\begin{equation}\label{def: Q_k}
    |\langle q_{k+1}, y_{k+1} - y \rangle|\leq  \sa{S_{k+1}} \triangleq L_{yx}\|x_{k+1}-x_{k}\| \|y_{k+1}-y\|+L_{yy}\|y_{k+1}-y_{k}\| \|y_{k+1}-y\|
\end{equation}}%
for $k \geq -1$.
%\nsa{$Q_k$ was conflicting with a later definition we used: $P_k,Q_k$}
%We will use this definition and the relative bound later.
Recall $x_{-1} = x_0,\;y_{-1} = y_0$, thus \sa{$q_0=\mathbf{0}$; therefore, for \sa{\bf part\;1},}
{\small
\begin{align}
\label{eq:noisy-rate-part1}
         \sum_{k=0}^{N-1}&\rho^{-k}
                 ( \theta \langle q_k, y_{k} - y \rangle-\langle q_{k+1}, y_{k+1} - y \rangle )
                %  = & -\sum_{k=0}^{N-2}\rho^{-k}\langle q_{k+1}, y_{k+1} - y \rangle - \rho^{-N+1}\langle q_{N}, y_{N} - y \rangle + \sum_{k=-1}^{N-2}\rho^{-k-1} \theta \langle q_{k+1}, y_{k+1} - y \rangle
                %  \\
                 = \sum_{k=0}^{N-2}\rho^{-k}\Big(\frac{\theta}{\rho}-1\Big)\langle q_{k+1}, y_{k+1} - y \rangle - \rho^{-N+1}\langle q_{N}, y_{N} - y \rangle
                 \\
                 \leq & \sum_{k=0}^{N-2}\rho^{-k}\sa{|1-\frac{\theta}{\rho}|~S_{k+1}} + \rho^{-N+1}\sa{S_{N}}
                 \rev{\leq  \sum_{k=0}^{N-1}\rho^{-k}{|1-\frac{\theta}{\rho}|~S_{k+1}}
                 + \rho^{-N}{\theta}{S_{N}}},\nonumber
\end{align}}%
where \sa{the first inequality follows from %\sa{$\theta\in[0,\rho]$} and
\cref{def: Q_k}.}
% \sa{Furthermore, for {\bf part 2},} we have that
% \begin{equation}
% \label{eq:noisy-rate-part2}
%     \begin{aligned}
%           & \sum_{k=0}^{N-1}\theta^{-k} (\Lambda_k{(x,y)} - \Sigma_{k+1} \sa{(x,y)} )
%           \\
%             = & \Lambda_ 0 \sa{(x,y)} - \sum_{k=1}^{N-1}\theta^{-k} (\theta \Sigma_{k} \sa{(x,y)} - \Lambda_{k} \sa{(x,y)}  ) - \theta^{-N+1}\Sigma_N \sa{(x,y)}
%             \\
%             \leq & \Lambda_ 0 \sa{(x,y)}  - \theta^{-N+1}\Sigma_N \sa{(x,y)}
%             \leq \frac{1}{2\tau}\Omega_{ X} + \frac{1}{2\sigma}\Omega_{ Y}- \theta^{-N+1}\Sigma_N \sa{(x,y)},
%     \end{aligned}
% \end{equation}
% where the %first and second
% inequalities follow from \eqref{Condition: SAPD simple LMI 2} and
% %the second one is from
% \cref{ASPT: compact}.
% Next, we will analyze  \eqref{INEQ: difference of gap  _bounded version} for $x=x^*$ and $y=y^*$ in the expected sense.
\rev{
Next,
for $\Delta_k^{x}$ and $\hat{x}_{k+1}$ defined as in
\cref{lemma: intermediate noisy bound}, we write \textbf{part 2} as follows:
{\footnotesize
\begin{equation}\label{eq:noisy-rate-part2-dist-metric}
    \begin{aligned}
            \MoveEqLeft \sum_{k=0}^{N-1} -\rho^{-k} \langle \Delta^{x}_k , x_{k+1} -  x \rangle
            =
             \sum_{k=0}^{N-1} \rho^{-k}
            \Big(\langle \Delta^{x}_k , {\hat{x}_{k+1}}-x_{k+1} \rangle  - \langle \Delta^{x}_k , {\hat{x}_{k+1}} -  x \rangle\Big).
    \end{aligned}
\end{equation}}%
Finally, for $\Delta_k^{y}$, $\hat{y}_{k+1}$ and $\hat{\hat{y}}_{k+1}$ defined as in
\cref{lemma: intermediate noisy bound}, we also write \textbf{part 3} as follows:
{\footnotesize
\begin{equation}\label{eq:noisy-rate-part3-dist-metric}
    \begin{aligned}
    \hspace{-4pt}
        \MoveEqLeft \sum_{k=0}^{N-1} \rho^{-k} \langle \tilde{s}_k -s_k,
       y_{k+1} - y \rangle
       \\
      = & \hspace{-4pt}
         \sum_{k=0}^{N-1} \rho^{-k}\Big[(1 +\theta) \langle \Delta^{y}_k,
       y_{k+1} - {\hat{y}_{k+1}}  + {\hat{y}_{k+1}}  - y \rangle
       %\\  &
       - \theta \langle \Delta^{y}_{k-1},
       y_{k+1} - {\hat{\hat{y}}_{k+1}} + {\hat{\hat{y}}_{k+1}}  - y \rangle\Big].
    \end{aligned}
\end{equation}}%
%\xtodo{In the statement of \cref{Thm: main result_R1}, the def of $d_N$ does not have (x,y)}
\rev{Let $d_N(x,y)\triangleq \tfrac{1}{2\tau}\|x-x_N\|^2 +  \tfrac{1}{2\sigma}\left( 1 - \alpha\sigma\right)\|y-y_N\|^2$.} Adding $\rho^{-N}d_N(x,y)$ to both sides of \eqref{INEQ: difference of gap  _bounded version}, then using \eqref{eq:noisy-rate-part1}, \eqref{eq:noisy-rate-part2-dist-metric} and \eqref{eq:noisy-rate-part3-dist-metric},
%for any fixed $(x,y)\in \dom f\times \dom g$,
we get
{\small
\begin{equation}
\label{INEQ: ergordic gap  _bounded version}
\begin{aligned}
        %\MoveEqLeft
        K_{N} ({\rho})  \left(\mathcal{L}(\bar{x}_{N}, y)  - \mathcal{L}(x, \bar{y}_{N}) \right)
        + {\rho}^{-N} d_N(x,y)
        %\\ &
        \leq
      U_N(x,y) + \sum_{k=0}^{N-1}{\rho}^{-k}(P_k(x,y) + Q_k),
\end{aligned}
\end{equation}}%
where $U_N(x,y)$, $P_k(x,y)$ and $Q_k$ for $k=0,\ldots,N-1$ are defined as
{\small
\begin{equation*}
    \begin{aligned}
    {U_N(x,y)}   \triangleq
      &  {\sum_{k=0}^{N-1}{\rho}^{-k}
      \Big(\Gamma_{k+1} + \Lambda_k{(x,y)} - \Sigma_{k+1} {(x,y)} + {{|1-\frac{\theta}{\rho}|~S_{k+1}}} \Big)
      \rev{+ {\rho}^{-N}\Big(d_N(x,y) + {{\theta}{S_{N}}}\Big)},}
      \\
           P_k(x,y)  \triangleq
      & -\langle \Delta^{x}_k , {\hat{x}_{k+1}} -  x \rangle
     +    (1+\theta)\langle \Delta^{y}_k, {\hat{y}_{k+1}} - y \rangle
     -\theta \langle \Delta^{y}_{k-1}, {\hat{\hat{y}}_{k+1}} -y  \rangle,
     \\
     Q_k  \triangleq &
     \langle \Delta^{x}_k , {\hat{x}_{k+1} -  x_{k+1}} \rangle
     +
     (1+\theta)\langle \Delta^{y}_k, {y_{k+1} - \hat{y}_{k+1}}  \rangle
     -\theta \langle \Delta^{y}_{k-1}, {y_{k+1} - \hat{\hat{y}}_{k+1}}  \rangle.
    \end{aligned}
\end{equation*}}}%
\rev{We first uniformly upper bound $\mathbb{E}\left[ Q_k \right]$ for all $k\geq 0$ using \cref{Lemma: final noisy bound}, i.e.,
{\small
\begin{equation}
\begin{aligned}
\label{eq:Q-bound}
     \mathbb{E}[\sum_{k=0}^{N-1}{\rho}^{-k} {Q_k}]  \leq
     \Big[\tfrac{\tau}{1+\tau\mu_x}\Xi_{\tau,\sigma,\theta}^x \delta_x^2  +     \tfrac{\sigma}{1+\sigma\mu_y}\Xi_{\tau,\sigma,\theta}^y \delta_y^2\Big] \sum_{k=0}^{N-1}{\rho}^{-k}.
\end{aligned}
\end{equation}}}%
\rev{Next, for arbitrarily fixed} $(x,y)\in \dom f \times \dom g$, we analyze $U_N(x,y)$. 
After adding and subtracting $\tfrac{\alpha}{2}\|y_{k+1}-y_{k}\|^2$,  and rearranging the terms, we get
{\small
\begin{equation}
\label{eq:U_def _bounded version}
    \begin{aligned}
           U_N(x,y) =  &
           \sa{\frac{1}{2}\sum_{k=0}^{N-1}{\rho}^{-k}\Big(\xi_k^\top A\xi_k-\xi_{k+1}^\top B \xi_{k+1}\Big)}
           \rev{+{\rho}^{-N}(d_N(x,y) + {\theta}{S_{N}})}\\
           =  &
           \frac{1}{2}\xi_0^\top A\xi_0- \frac{1}{2}\sum_{k=1}^{N-1}{\rho}^{-k+1}[\xi_k^\top{( B
           - \tfrac{1}{\xuan{\rho}}A)}\xi_k]
           \rev{-{\rho}^{-N + 1 }\Big( \frac{1}{2}\xi_N^\top   B \xi_N  - \frac{1}{\rho}d_N(x,y) - \frac{\theta}{\rho}{S_{N}}\Big)},
    \end{aligned}
\end{equation}}%
\sa{
%as the second inequality follows from the fact that $\xi_0^\top  B \xi_0=0$\nsa{edit here},
where $A, B \in\reals^{5\times 5}$
%such that $ B $ is a diagonal matrix
and $\xi_k\in\reals^5$ are defined for $k\geq 0$ as follows:}
{\scriptsize
%\begin{align*} &
$A \triangleq
  \begin{pmatrix}
    \frac{1}{\tau} & 0 & 0 & 0 & 0\\
  0 & \frac{1}{\sigma} & 0 & 0 & 0\\
  0 & 0 & 0 & 0 & {\theta L_{yx}}\\
  0& 0 & 0 & 0 & {\theta L_{yy}}\\
  0 & 0 & \theta L_{yx} & \theta L_{yy} & -\alpha
\end{pmatrix}$,
%\quad
$\xi_k \triangleq \left(
    \begin{array}{*{20}{c}}
  \|x_k - x\| \\
  \|y_k - y\| \\
  \|x_k - x_{k - 1}\| \\
  \|y_k - y_{k - 1}\| \\
  \|y_{k+1} - y_k\|
   \end{array}
   \right)$, {\normalsize and}
   %\\&
   $B \triangleq
  \begin{pmatrix}
    \frac{1}{\tau}+\mu_x & 0 & 0 & 0 & 0\\
  0 & \frac{1}{\sigma}+\mu_y & \sa{-|1-\frac{\theta}{\rho}|}~L_{yx} & \sa{-|1-\frac{\theta}{\rho}|}~L_{yy} & 0\\
  0 & \sa{-|1-\frac{\theta}{\rho}|}~L_{yx} & \tfrac{1}{\tau} - L_{xx} & 0 & 0\\
  0& \sa{-|1-\frac{\theta}{\rho}|}~L_{yy} & 0 & \frac{1}{\sigma} - \alpha & 0\\
  0 & 0 & 0 & 0 & 0
\end{pmatrix}$}
%\end{align*}}%
such that $x_{-1} = x_{0}$ and $y_{-1} = y_{0}$.
\sa{In Lemma~\ref{lem:equivalent_systems} we show that \cref{eq: general SAPD LMI_R1}
%implies that
is equivalent to $B -\tfrac{1}{\rho}A\succeq 0$;} therefore, it follows from \eqref{eq:U_def _bounded version} that for any given $(x,y)$, the following inequality holds w.p. 1, 
{\small
\begin{align*}
    U_N(x,y) \leq  \sa{\frac{1}{2}\xi_0^\top A\xi_0}- {\rho}^{-N + 1 } \Big( \frac{1}{2}\xi_N^\top  B \xi_N  -  \frac{1}{\rho}d_N(x,y) -\frac{\theta}{\rho}\sa{S_{N}}\Big).
\end{align*}}%
\rev{Note that $\frac{1}{2}\xi_0^\top A\xi_0=\frac{1}{2\tau}\norm{x-x_0}^2+\frac{1}{2\sigma}\norm{y-y_0}^2$.} Furthermore,
%\nsa{When we allow $\theta>\rho$, we need to edit the following.}
%\todo{Xuan: I edited the following equation in June. 23. If it is correct, please remove my color.}
{\small
\begin{equation*}
%\label{PF3eq:permutation_result_2_change  _bounded version}
\begin{aligned}
%\textstyle
\frac{1}{2}\xi_N^\top    B  \xi_N  -\frac{\theta}{\rho}\sa{S_{N}}=
%\\
%  = &  \frac{1}{2}\xi_N^\top
%       \xuan{ \begin{pmatrix}
%   \frac{1}{\tau}+\mu_x & 0 & 0 & 0 & 0\\
%   0 & \frac{1}{\sigma}+\mu_y & (-| 1 - \tfrac{\theta}{\rho} | - \tfrac{\theta}{\rho})L_{yx} & (-| 1 - \tfrac{\theta}{\rho} | - \tfrac{\theta}{\rho})L_{yy} & 0\\
%   0 & (-| 1 - \tfrac{\theta}{\rho} | - \tfrac{\theta}{\rho})L_{yx} & {\frac{1}{\tau} - L_{xx}} & 0 & 0\\
%   0 & (-| 1 - \tfrac{\theta}{\rho} | - \tfrac{\theta}{\rho})L_{yy} & 0 & \frac{1}{\sigma} - \alpha & 0\\
%   0 & 0 & 0 & 0 & 0
% \end{pmatrix}}
% \xi_N
% \\ = &
&\frac{1}{2\rho\tau} \|x_N-x\|^2 +  \frac{1}{2}{\left(\frac{1}{\rho\sigma} - \frac{\alpha}{\rho} \right)}\|y_N-y\|^2\\
&
+\frac{1}{2}\xi_N^\top
     \begin{pmatrix}
       \frac{1}{\tau}(1-\frac{1}{\rho})+\mu_x & \mathbf{0}_{1\times3} & 0\\
       \mathbf{0}_{3\times 1} & \sa{G''} & \mathbf{0}_{3\times 1} \\
       0 & \mathbf{0}_{1\times3} & 0\\
     \end{pmatrix}
%      \begin{pmatrix}
%   \frac{1}{\tau}(1-\frac{1}{\rho})+\mu_x & 0 & 0 & 0 & 0\\
%   0 &  {\frac{1}{\sigma}+\mu_y -\frac{1}{\rho\sigma} + \frac{\alpha}{\rho}} & (-| 1 - \tfrac{\theta}{\rho} | - \tfrac{\theta}{\rho})L_{yx} & (-| 1 - \tfrac{\theta}{\rho} | - \tfrac{\theta}{\rho})L_{yy} & 0\\
%   0 & (-| 1 - \tfrac{\theta}{\rho} | - \tfrac{\theta}{\rho})L_{yx} & {\frac{1}{\tau} - (-| 1 - \tfrac{\theta}{\rho} | - \tfrac{\theta}{\rho})L_{xx}} & 0 & 0\\
%   0 & (-| 1 - \tfrac{\theta}{\rho} | - \tfrac{\theta}{\rho})L_{yy} & 0 & \frac{1}{\sigma} - \alpha & 0\\
%   0 & 0 & 0 & 0 & 0
% \end{pmatrix}
\xi_N
%\\
  \geq %&
%   \frac{1}{2\rho\tau} \|x_N-x\|^2 +  \frac{1}{2\rho\sigma}(1-\alpha\sigma)\|y_N-y\|^2 =
\rev{\frac{1}{\rho}d_N(x,y)},
\end{aligned}
\end{equation*}}%
%for $(x,y)\in X\times Y$,
which
%\sa{the first inequality follows from Cauchy-Schwarz and the second inequality follows from \cref{ASPT: lipshiz gradient}}; %the third inequality is due to $\pi_3^y>0$; and
{%the inequality
follows from \cref{eq: general SAPD LMI_R1,lemma: sub positive matrix}}, \rev{where} \sa{$G''$} is defined in \cref{eq: sub matrix of general SAPD LMI}.
%Now, we are ready to show \eqref{eq:distance-rate_R1}.
\rev{\sa{Since $(x,y)\in \dom f\times \dom g$ is fixed arbitrarily, all the results we have derived so far hold for any $(x,y)$; thus,}
\begin{equation}\label{eq:U-bound-distance-metric}
    U_N(x,y)\leq \tfrac{1}{2\tau}\norm{x-x_0}^2+\tfrac{1}{2\sigma}\norm{y-y_0}^2,\quad\forall (x,y)\in\dom f\times\dom g\quad \text{w.p. 1.}
\end{equation}}%
Finally, from \cref{ASPT: unbiased noise assumption}, for
$k\geq-1$, we have
%{\small
%\begin{equation*}
%\begin{aligned}&
$\mathbb{E}\left[ \langle \Delta^{x}_k , {\hat{x}_{k+1}} -  x^* \rangle \right] = \mathbb{E}\left[ \langle \Delta^{y}_k, {\hat{y}_{k+1}} - y^* \rangle \right] = \mathbb{E}\left[ \langle \Delta^{y}_{k-1},{\hat{\hat{y}}_{k+1}} - y^* \rangle \right] = 0$.
    %   \mathbb{E}\left[ \|\Delta^{x}_k\|^2 \right] \leq \delta_x^2,
    %   \mathbb{E}\left[ \|\Delta^{y}_k\|^2 \right] \leq \delta_y^2;
%\end{aligned}
%\end{equation*}}%
%\xuan{Indeed, it is easy to see that, $\hat{x}_{k+1}$ and $\tilde{x}_k$ are the deterministic function of $x_k$ and $y_{k+1}$, thus, the first expectation is zero. And, it is similar to the other two expectations.}
\sa{Thus, %it follows that
$\mathbb{E}[P_k(x^*,y^*)]=0$ for any %$k\in\{0,\ldots, N-1\}$
$k\geq 0$.}
%\nsa{Mert, can you also please verify these expectation results?}\mtodo{Looked good to me.}
%\xtodo{It is correct to me}
\sa{Therefore, %combining this result with 
from \eqref{eq:Q-bound}, we get\vspace*{-3mm}
{\small
\begin{equation}
\label{eq:variance-bound}
   \mathbb{E}\Big[ \sum_{k=0}^{N-1}{\rho}^{-k}(P_k(x^*,y^*)+ Q_k)\Big]\leq K_N({\rho})~{\Xi}_{\tau,\sigma,\theta}. \vspace*{-2mm}
\end{equation}}}%
%Hence, $U_N(x^*,y^*) \leq \frac{1}{2\tau}\norm{x^*-x_0}^2+\frac{1}{2\sigma}\norm{y^*-y_0}^2$ w.p. $1$.
\rev{Note $d_N^*=d_N(x^*,y^*)$; hence,} it follows from \eqref{INEQ: ergordic gap  _bounded version}, \eqref{eq:U-bound-distance-metric} and \eqref{eq:variance-bound} that
$$\mathbb{E}[ K_{N} ({\rho})  \left(\mathcal{L}(\bar{x}_{N}, y^*)  - \mathcal{L}(x^*, \bar{y}_{N}) \right)+ {\rho}^{-N}d_N(x^*,y^*)]\leq K_N({\rho})~{\Xi}_{\tau,\sigma,\theta}+ D_{\tau,\sigma}.$$
%\end{aligned}
%\end{equation*}}}%
\rev{Since $\mathcal{L}(\bar{x}_{N}, y^*)  - \mathcal{L}(x^*, \bar{y}_{N})\geq 0$, %for $N\geq 1$,
\eqref{eq:distance-rate_R1} immediately follows from %\eqref{INEQ: distance gap}
above inequality.}
\qed
\vspace{-2mm}
\subsection{Parameter Choices for SAPD}
\label{sec:parameter_choice}
\sa{%The reason
We employ the %complicated
matrix inequality~(MI) in~\eqref{Condition: SAPD simple LMI system} to describe the admissible set of algorithm parameters that guarantee convergence. Our aim is to enjoy  \rev{a} wide range of parameters %as possible %,  which helps to eliminate
%in order
to improve the robustness of SAPD, i.e., to control the noise amplification of the algorithm,} %when dealing with
\sa{in the presence of} noisy gradients. %situation.
Although, it seems difficult to find an \mg{explicit} solution to the
%matrix inequality, %
MI in \cref{Thm: main result_R1},
we can compute a particular solution \sa{to it %due to
 by exploiting} \mg{its} %special
 structure.
%\sa{Indeed}, for some $c_{\tau},c_{\sigma}\in\sa{[0,1)}$, we first let
% \begin{equation}
% \label{eq:pi_selection}
% \pi^x = \frac{c_{\tau}}{\tau},\quad \pi^y_1 = \frac{c_{\sigma}}{3(1+\theta)\sigma},\quad \pi^y_2 = \frac{c_{\sigma}}{3\theta\sigma},\quad \pi^y_3 = \frac{c_{\sigma}}{3\sigma}.
% \end{equation}
%Then we have a special matrix inequality for \cref{Thm: main result} as follows,
% \sa{For this particular choice of $(\pi^x,\pi_1^y,\pi_2^y,\pi_3^y)$, %the matrix inequality in~
% \eqref{Condition: SAPD simple LMI system} reduces to}
% \begin{subequations}
% \label{Condition: SAPD simple LMI system}
% \begin{equation}\label{Condition: SAPD simple LMI 1}
%     \begin{pmatrix}
%   \tfrac{1-c_{\tau}}{\tau} - L_{xx}  & 0 & - L_{yx}\\
%  0 & \tfrac{1-c_{\sigma}}{\sigma}- \alpha  & - L_{yy}\\
%  - L_{yx} & - L_{yy} & \tfrac{\alpha}{ \theta}
% \end{pmatrix}
% \succeq 0,
% \end{equation}
% \begin{align}
%     \label{Condition: SAPD simple LMI 2}
%     \sa{\tau\mu_x}\geq \frac{1-\theta}{\theta},\quad \sa{\sigma\mu_y}\geq \frac{1-\theta}{\theta}.
% \end{align}
% \end{subequations}
Next, in \cref{LEMMA: Noise LMI after young's ineq-R1}, we give \sa{an intermediate condition to help us construct the particular solution provided in \cref{Corollary: explicit solution to noisy LMI-R1} for the SCSC setting.}
%\xtodo{I suggest move  \cref{LEMMA: Noise LMI after young's ineq-R1} to the end of the paper. The referees also have such complains.}
\begin{lemma}\label{LEMMA: Noise LMI after young's ineq-R1}
Let $\tau, \sigma>0$, \rev{$\theta\in(0,1)$} be a solution to the following 
system:
{\small
\begin{equation}
\label{eq:sufficient_cond_noisy_LMI-R1}
\sa{\min\{\tau\mu_x,\sigma\mu_y\}\geq \frac{1-\theta}{\theta}},\quad \frac{1}{\tau} \geq L_{xx} + \pi_1  L_{yx},
\quad \frac{\rev{c}}{\sigma} \geq \frac{ \theta L_{yx}}{\pi_1} + \Big(\frac{\theta}{\pi_2} +\pi_2\Big)  L_{yy},
\end{equation}}
\sa{for some $\pi_1,\pi_2>0$ and \rev{$c\in(0,1]$}.} 
Then $\{\tau, \sigma,\theta,\alpha\}$ is a solution to \eqref{Condition: SAPD simple LMI system} for \sa{$\alpha = \frac{\theta L_{yx}}{\pi_1} + \frac{\theta L_{yy}}{\pi_2}$.}
\end{lemma}
\begin{proof}
\sa{Since the first inequalities in both \eqref{eq:sufficient_cond_noisy_LMI-R1} 
and~\eqref{Condition: SAPD simple LMI system} are the same}, \mg{we} only need to show \sa{the
MI in~\eqref{Condition: SAPD simple LMI system}
 holds}. Substituting $\alpha = \tfrac{\theta L_{yx}}{\pi_1} + \tfrac{\theta L_{yy}}{\pi_2}$ into \eqref{Condition: SAPD simple LMI system}, we get
{\footnotesize
\begin{equation*}
\begin{aligned}
\begin{pmatrix}
  \tfrac{1}{\tau} - L_{xx} & 0 & - L_{yx}\\
 0 &\tfrac{1}{\sigma}- \alpha & - L_{yy}\\
 - L_{yx} & - L_{yy} & \tfrac{\alpha}{ \theta}
\end{pmatrix}
= 
\underbrace{
\begin{pmatrix}
\tfrac{1}{\tau} - L_{xx} & 0 & - L_{yx}\\
 0 & 0 & 0 \\
 - L_{yx} & 0 & \tfrac{L_{yx}}{\pi_1}
\end{pmatrix}}_{M_1}
+
\underbrace{
\begin{pmatrix}
  0 & 0 & 0\\
 0 & \frac{1}{\sigma}- \frac{\theta L_{yx}}{\pi_1} - \frac{\theta L_{yy}}{\pi_2} & - L_{yy}\\
 0 & - L_{yy} & \frac{ L_{yy}}{\pi_2}
\end{pmatrix}}_{M_2}.
% &
%     \begin{pmatrix}
%   \tfrac{1}{\tau} - L_{xx} & 0 & - L_{yx}\\
%  0 & \tfrac{1}{\sigma}-  \tfrac{\theta L_{yx}}{\pi_1} - \tfrac{\theta L_{yy}}{\pi_2} & - L_{yy}\\
%  - L_{yx} & - L_{yy} & \tfrac{L_{yx}}{\pi_1} + \tfrac{ L_{yy}}{\pi_2}
% \end{pmatrix}
% \triangleq &
% M_1 + M_2,
\end{aligned}
\end{equation*}}%
% where
% ${\everymath={\scriptscriptstyle}
% M_1
%  \triangleq
%  \begin{pmatrix}
% \tfrac{1}{\tau} - L_{xx} & 0 & - L_{yx}\\
%  0 & 0 & 0 \\
%  - L_{yx} & 0 & \tfrac{L_{yx}}{\pi_1}
% \end{pmatrix}}$ and
% ${
% \everymath={\scriptscriptstyle}
% M_2 \triangleq
%     \begin{pmatrix}
%   0 & 0 & 0\\
%  0 & \frac{1}{\sigma}-  \frac{\theta L_{yx}}{\pi_1} - \frac{\theta L_{yy}}{\pi_2} & - L_{yy}\\
%  0 & - L_{yy} & \frac{ L_{yy}}{\pi_2}
% \end{pmatrix}}$.
 \sa{Therefore, since $\pi_1,\pi_2>0$, \rev{the second and the third inequalities in~\eqref{eq:sufficient_cond_noisy_LMI-R1}}
 imply $M_1
 \succeq 0$ and $M_2 \succeq 0$, respectively.}
% $$
%         M_1
%  \succeq
% \begin{pmatrix}
% \pi_1 L_{yx} & 0 & - L_{yx}\\
%  0 & 0 & 0 \\
%  - L_{yx} & 0 & \tfrac{L_{yx}}{\pi_1}
% \end{pmatrix}
% \succeq 0,
% \quad
% % $$
% % %  \sa{similarly, \eqref{Condition: noisy LMI sufficient 2} implies that}
% % $$
%         M_2
%  \succeq
% \begin{pmatrix}
%   0 & 0 & 0\\
%  0 &\pi_2 L_{yy} & - L_{yy}\\
%  0 & - L_{yy} & \tfrac{ L_{yy}}{\pi_2}
% \end{pmatrix}
% \succeq 0,\; \text{respectively.}
% $$
Thus, $M_1+M_2\succeq 0$. 
\end{proof} 
\cref{LEMMA: Noise LMI after young's ineq-R1} \mg{shows}
 that every solution to
\eqref{eq:sufficient_cond_noisy_LMI-R1} can be converted to a solution to
\eqref{Condition: SAPD simple LMI system}. Next, based on 
this lemma, we will give an explicit parameter choice for Algorithm~\ref{ALG: SAPD}. 
\begin{corollary}\label{Corollary: explicit solution to noisy LMI-R1}
\sa{Suppose $\mu_x,\mu_y > 0$. If $L_{yy}>0$, for any given $\beta\in(0,1)$ and $\rev{c\in(0,1]}$, let $\tau, \sigma>0$ and $\theta\in (0,1)$ be chosen satisfying
{\small
\begin{equation}
\label{Condition: SP solution to noisy LMI-R1}
%\begin{gather}
\tau = \frac{1-\theta}{\mu_x \theta},\quad \sigma = \frac{1-\theta}{\mu_y\theta},
%\label{Condition: SP solution to noisy LMI 2}\\
\quad
\theta  \geq \overline{\theta}\triangleq\max\{\overline{\theta}_1,~\overline{\theta}_2\},
%\label{Condition: SP solution to noisy LMI 1}
%\end{gather}
\end{equation}}}%
%\todo{MG: in the reply letter, we should say why we change the stepsize choice; we can say we simplified the analysis for the distance function, this caused simplifications in the parameter choice in Corollary 2.11 etc.}
where $\overline{\theta}_1,~\overline{\theta}_2\in\sa{(0,1)}$, depending on the choice of $\beta$ and \rev{$c$}, are defined as
{\footnotesize
\begin{equation}
\label{eq:theta1-R1}
        \sa{\overline{\theta}_1\triangleq 1 -\tfrac{\rev{c}\beta (L_{xx} + \mu_x) \mu_y}{2L_{yx}^2}
          \Big(\sqrt{ 1+ \tfrac{4\mu_xL_{yx}^2}{\rev{c}\beta\mu_y(L_{xx}+\mu_x)^2}}-1\Big),\quad \overline{\theta}_2\triangleq
        1 - \tfrac{\rev{c^2}(1-\beta)^2}{8}\tfrac{\mu_y^2}{L_{yy}^2} \Big( \sqrt{1+\tfrac{16L_{yy}^2}{\rev{c^2}(1-\beta)^2\mu_y^2}}-1\Big).}
\end{equation}}%
%\mg{and depend on the choice of $\beta$.}
\sa{On the other hand, if $L_{yy}=0$, let $\tau, \sigma>0$ and $\theta\in (0,1)$ be chosen as in \eqref{Condition: SP solution to noisy LMI-R1} for $\overline{\theta}_1$ in \eqref{eq:theta1-R1} with $\beta=1$ and $\overline{\theta}_2=0$.\footnote{\sa{Our parameter selection when $L_{yy}=L_{xx}=0$ recovers $(\tau,\sigma,\theta)$ choice in~\cite[Eq.(49)]{chambolle2016ergodic}.}}} \sa{Then $\alpha = \frac{\rev{c}}{\sigma}-\sqrt{\theta}L_{yy}>0$,  and $\{\tau, \sigma,\theta,\alpha\}$ is a solution to \eqref{Condition: SAPD simple LMI system}.} \sa{Moreover, when $L_{yy}>0$, the minimum $\overline{\theta}$ is attained at the unique $\beta^*\in (0,1)$ such that $\overline{\theta}_1=\overline{\theta}_2$.}
\end{corollary}
{%To compute 
\rev{For} this \sa{particular} solution, we set $\rho=\theta$; hence, $\theta$ is not only the momentum parameter, but it also determines the linear rate for \rev{the bias term in \eqref{eq:distance-rate_R1}
%which gives an error bound on $\mathbb{E}[d_N^*]$
.}}
%\nsa{Below discussion is not good. We can say that we can compute $\beta^*$ that minimizes $\overline{\theta}$.}
\vspace{-2mm}
\subsection{Iteration Complexity Bound for SAPD} \label{sec:complexity} %\label{section3.3}
%\cite{palaniappan2016stochastic} proposed a variance reduced version of the stochastic gradient method for solving strongly monotone variational inequalities with linear convergence guarantees.
 \sloppy In this part, we study the %lower
 iteration complexity bound for SAPD,
 %which are given on the number of iterations $N$ for computing
 \rev{to compute $({x}_N,{y}_N)$ such that $\cD({x}_N,{y}_N)\leq \epsilon$,
 %$\mathbb{E}\left[\mathcal{L}\left(\bar{x}_{N}, y\right)-\mathcal{L}\left(x, \bar{y}_{N}\right) + \Delta^{\theta}_N(x,y)\right]\leq \epsilon$,
 where $\epsilon>0$ is a given tolerance and $\cD(\cdot,\cdot)$ denotes the 
 %gap 
 distance function defined in \eqref{eq:gap}.}
 \begin{theorem}\label{Prop: iteration complex for SAPD_R1}
Suppose $\mu_x, \mu_y>0$, and \cref{ASPT: lipshiz gradient,ASPT: unbiased noise assumption} hold. For any $\epsilon>0$, suppose the SAPD parameters $\{\tau,\sigma,\theta\}$ are chosen such that
{\small
\begin{align}
\label{eq:SAPD-parameter-choice-R1}
    \tau = \frac{1-\theta}{\mu_x \theta},\quad \sigma = \frac{1-\theta}{\mu_y\theta},\quad \theta=\max\{\overline{\theta},~\overline{\overline{\theta}}\},
\end{align}}%
where \rev{$\overline{\theta}$ is set as in \eqref{Condition: SP solution to noisy LMI-R1} for some
arbitrary $\beta\in(0,1]$ and $c=\tfrac{1}{2}$}, and
{\small
\begin{align}
\label{eq:theta_bound_2-R1}
    \overline{\overline{\theta}}\triangleq\max\left\{ \overline{\overline{\theta}}_1 ,\; \overline{\overline{\theta}}_2 \right\},\quad
    \sa{\overline{\overline{\theta}}_1}  = \max\Big\{0,1 -
\sa{\frac{1}{12\Xi^x(\beta)}}\frac{\mu_x}{\delta_x^2}\epsilon\Big\},\quad \sa{\overline{\overline{\theta}}_2} = \max\Big\{0,1 -
\frac{1}{12\Xi^y(\beta)}\frac{\mu_y}{\delta_y^2}\epsilon\Big\}\mg{,}
\end{align}}%
\sa{such that ${\Xi}^x(\beta)\triangleq 1+ \Psi(\beta)$ and ${\Xi}^y(\beta) \triangleq 
% 3\Big(5  + \beta \frac{ L_{xy}}{L_{yx}}- \beta\Big) 
\frac{27-3\beta}{2}  + \frac{3\beta L_{xy}}{L_{yx}}
+ \sa{\frac{\mu_y}{\mu_x}\Psi(\beta)}
$
\mg{with} $\Psi(\beta)\triangleq \min\left\{\sqrt{\frac{\beta}{2}\frac{\mu_x}{\mu_y}},\frac{1-\beta}{4}\frac{L_{yx}}{L_{yy}}\right\}$.} 
%and $\mu=\min\{\mu_x,\mu_y\}$.
Then, the iteration complexity of SAPD, \sa{as stated in~\cref{ALG: SAPD},} \sa{to generate \rev{$(x_\epsilon,y_\epsilon)\in\cX\times\cY$ such that $\cD(x_\epsilon,y_\epsilon)=\mathbb{E}[\mu_x\norm{x_\epsilon-x^*}^2+\mu_y\norm{y_\epsilon-y^*}^2]\leq \epsilon$}}
is
\rev{\footnotesize
\begin{equation}
\label{eq:complexity-1-SCSC}
\mathcal{O}\left(\Big[
\frac{L_{xx}}{\mu_x} + \frac{L_{yx}}{\sqrt{\mu_x\mu_y}} + \frac{L_{yy}}{\mu_y}
+ \Big(\Big(1+\sqrt{\frac{\mu_x}{\mu_y}} \Big) \frac{\delta_x^2}{\mu_x} +  \Big(1{+\frac{L_{xy}}{L_{yx}}}+\sqrt{\frac{\mu_y}{\mu_x}}
\Big)\frac{\delta_y^2}{\mu_y} \Big)\frac{1}{\epsilon}\Big]
\cdot\ln\Big(\frac{
%\max\{\frac{\mu_x}{\mu_y},\frac{\mu_y}{\mu_x}\}\|z_0-z^*\|^2
\cD(x_0,y_0)}{\epsilon}\Big)\right).
\end{equation}}%
\rev{Furthermore, choosing $\beta= \min\{\tfrac{1}{2},\frac{\mu_y}{\mu_x},\frac{\mu_x}{\mu_y}\}$ leads to the following iteration complexity:
%\xtodo{change $\cD'$}
{\footnotesize
\begin{equation}\label{eq:complexity-2-SCSC}
\mathcal{O}\left( \Big[
\frac{L_{xx}}{\mu_x} + \frac{L_{yx}}{\min\{\mu_x,\mu_y\}} + \frac{L_{yy}}{\mu_y}
+ \Big(\frac{\delta_x^2}{\mu_x} +  \Big(1+\frac{L_{xy}}{L_{yx}}
\Big)\frac{\delta_y^2}{\mu_y} \Big)\frac{1}{\epsilon}\Big]
\cdot\ln\Big(\frac{
%\max\{\frac{\mu_x}{\mu_y},\frac{\mu_y}{\mu_x}\}\|z_0-z^*\|^2
\cD(x_0,y_0)}{\epsilon}\Big)\right).
\end{equation}}}%
\end{theorem}
%\xtodo{Maybe move this proof to the end or appendix}
\begin{proof}
Given $\beta\in(0,1)$, let $\{\tau, \sigma,\theta,\alpha\}$ be a particular solution to \eqref{Condition: SAPD simple LMI system} constructed according to \cref{Corollary: explicit solution to noisy LMI-R1}. Therefore, using these particular parameter values together with $\rho=\theta$ within \cref{Thm: main result_R1}, we know that  \eqref{eq:distance-rate_R1} holds, \fin{i.e.,} for any  $N\geq0$, it follows  that
{\small
\begin{equation*}
     \mathbb{E}\Big[\frac{1}{2\tau}\|x_N-x^*\|^2 +  \frac{1 - \alpha\sigma}{2\sigma}\|y_N-y^*\|^2\Big]
              \leq \rho^{N}\Big(\frac{1}{2\tau}\| x_0 - x^*\|^2
    +  \frac{1}{2\sigma}\| y_0 - y^*\|^2\Big) + \frac{\rho}{1-\rho}\Xi_{\tau,\sigma,\theta}.
\end{equation*}}%
Using the parameter choice
$
\tau = \frac{1-\theta}{\theta\mu_x},\; \sigma = \frac{1-\theta}{\theta\mu_y},\; \alpha = \rev{\frac{c}{\sigma}}-\sqrt{\theta} L_{yy},\;\rho = \theta,
$
and letting \rev{$c=\frac{1}{2}$}, 
we first obtain that
$
\frac{1-\alpha\sigma}{\sigma} \geq \frac{1}{2\sigma};
$ then this inequality together with our parameter choice leads to
\begin{equation}\label{eq:rate-z-R1}
     \mathbb{E}[\mu_x\|x_N-x^*\|^2+\mu_y\|y_N-y^*\|^2 ]
              \leq 2 %\max\Big\{\frac{\mu_x}{\mu_y},\frac{\mu_y}{\mu_x}\Big\}
              \theta^{N}\Big(\mu_x\| x_0 - x^*\|^2+\mu_y\| y_0 - y^*\|^2\Big)
 + 4\Xi_{\tau,\sigma,\theta}. %{\min\{\mu_x,\mu_y\}}.
\end{equation}
Note \eqref{eq:SAPD-parameter-choice-R1} implies $\Xi_{\tau,\sigma,\theta}=(1-\theta)(\Xi^x_{\tau,\sigma,\theta}\frac{\delta_x^2}{\mu_x} +\Xi^y_{\tau,\sigma,\theta} \frac{\delta_y^2}{\mu_y})$, \rev{where $\Xi^x_{\tau,\sigma,\theta}$ and $\Xi^y_{\tau,\sigma,\theta}$ are defined in the statement of \cref{Thm: main result_R1}.}
% $\frac{\tau}{1+\tau\mu_x}=\frac{1-\theta}{\mu_x}$ and  $\frac{\sigma}{1+\sigma\mu_y}=\frac{1-\theta}{\mu_y}$. 
Thus, %using \rev{$\mu\triangleq\min\{\mu_x,\mu_y\}$}, 
for any $\epsilon>0$,
the right side of \eqref{eq:rate-z-R1}
can be bounded by $\epsilon$ when
{\begin{align}
    \label{eq:n123-R1}
   %\max\Big\{\frac{\mu_x}{\mu_y},\frac{\mu_y}{\mu_x}\Big\}
   \theta^{N}\Big(\mu_x\| x_0 - x^*\|^2+\mu_y\| y_0 - y^*\|^2\Big)
   %\| z_0 - z^*\|^2 
   \leq \frac{\epsilon}{6},\quad (1-\theta)\Xi^x_{\tau,\sigma,\theta}\frac{\delta_x^2}{\mu_x} \leq \frac{\epsilon}{12},\quad
    (1-\theta)\Xi^y_{\tau,\sigma,\theta}\frac{\delta_y^2}{\mu_y} \leq \frac{\epsilon}{12}.
\end{align}}%
\sa{
%; moreover, the other two terms that depend on $\theta$ satisfy $1+\Xi^x_{\tau,\sigma,\theta}>1$ and $1+2\theta+\Xi^y_{\tau,\sigma,\theta}>1$.
Therefore,
%\nsa{I have not check the remaining inequalities carefully.} 
to get a sufficient condition on $\theta$ for the last two inequalities in~\eqref{eq:n123-R1} to hold, we first upper bound $\Xi^x_{\tau,\sigma,\theta}$ and $\Xi^y_{\tau,\sigma,\theta}$
%, and then exploit $\frac{\tau}{1+\tau\mu_x}=\frac{1-\theta}{\mu_x}$ and $\frac{\sigma}{1+\sigma\mu_y}=\frac{1-\theta}{\mu_y}$ to get a lower bound on $\theta$
.}
The parameter choice of $\tau$ and $\sigma$ in \eqref{eq:SAPD-parameter-choice-R1}
implies that
\rev{\footnotesize
$$
     \Xi^x_{\tau,\sigma,\theta} = 1+ \theta \sa{(1-\theta^2)} \frac{L_{yx}}{2\mu_y},\quad \Xi^y_{\tau,\sigma,\theta} = \Big( %1 + 2\theta + \theta^2
     \sa{(1+\theta)^2}+ \sa{\theta(1-\theta^2)}\frac{L_{yy}}{\mu_y} + \theta(1+\theta)(1-\theta)^2\frac{ L_{yx}L_{xy}}{\mu_x\mu_y}\Big)(1+2\theta) + \theta\sa{(1-\theta^2)}\frac{L_{yx}}{2\mu_x}.
$$}%
Since $0<\theta \leq 1$, \sa{we have $1-\theta^2\leq 2(1-\theta)$; thus,}
$\Xi^x_{\tau,\sigma,\theta} \leq \sa{1}+ (1-\theta)\frac{L_{yx}}{\mu_y}$ and
{\small
\begin{equation}
\label{eq:Xi_y_1-R1}
   \Xi^y_{\tau,\sigma,\theta} \leq  6\Big(2+ (1-\theta)\frac{L_{yy}}{\mu_y} + (1-\theta)^2\frac{ L_{yx}L_{xy}}{\mu_x\mu_y}\Big) + (1-\theta)\frac{L_{yx}}{\mu_x}.
\end{equation}}%
On the other hand, since $\theta\geq \overline{\theta}=\max\{\overline{\theta}_1,\overline{\theta}_2\}$ and $\rev{c=\frac{1}{2}}$, the inequality $\sqrt{a+b} \leq \sqrt{a} + \sqrt{b}$ for all $a,b\geq 0$, and \eqref{eq:theta1-R1} 
  together imply that
$1-\theta\leq \sa{\min\Big\{\frac{\sqrt{\sa{\beta}\mu_x\mu_y/2}}{L_{yx}}, (1-\beta)\frac{\mu_y}{4L_{yy}}\Big\}}$.
\sa{Thus,}
{\small
\begin{equation}\label{eq: Xi_x bound-R1}
    \Xi^x_{\tau,\sigma,\theta} \leq 1+(1-\theta)\frac{L_{yx}}{\mu_y}\leq {\Xi}^x(\beta),
\end{equation}}
\sa{and within \eqref{eq:Xi_y_1-R1} bounding $(1-\theta)\frac{L_{yx}}{\mu_x}$ similarly and using $(1-\theta)^2 \leq \rev{\frac{\beta}{2}}\frac{\mu_x\mu_y}{L_{yx}^2}$, we get}
{\small
\begin{equation*}
   \Xi^y_{\tau,\sigma,\theta} \leq  6\Big(2 + (1-\theta)\frac{L_{yy}}{\mu_y} + \frac{\beta L_{xy}}{2L_{yx}}\Big) + \sa{\frac{\mu_y}{\mu_x}\min\Big\{\sqrt{\tfrac{\beta}{2}\tfrac{\mu_x}{\mu_y}},\tfrac{1-\beta}{4}\tfrac{L_{yx}}{L_{yy}}\Big\}}.
\end{equation*}}%
Next, it follows from $1-\theta\leq (1-\beta)\frac{\mu_y}{4L_{yy}}$ that
   $ \Xi^y_{\tau,\sigma,\theta} \leq
    \sa{{\Xi}^x(\beta)}$.
Therefore, this inequality together with \cref{eq: Xi_x bound-R1} and  the definition of $\overline{\overline{\theta}}$
%and our $(\tau,\sigma)$ choice in \eqref{eq:SAPD-parameter-choice} 
imply that
% requiring $\frac{1-\theta}{\mu_x}\frac{\delta_x^2}{2}~\sa{{\Xi}^x(\beta)}\leq \frac{\epsilon}{3}$ and $\frac{1-\theta}{\mu_y}\frac{\delta_y^2}{2}~\sa{{\Xi}^x(\beta)} \leq \frac{\epsilon}{3}$
% is sufficient for
% the last two inequalities in~\cref{eq:n123-R1} to hold. 
% \sa{It follows from the definition of $\overline{\overline{\theta}}$ in~\eqref{eq:theta_bound_2} that the sufficient condition holds for $\theta \geq \overline{\overline{\theta}}$.}
% \sa{Thus, combining this requirement with
% $\theta\geq\overline{\theta}$ in \cref{Coroallary: explicit solution to noisy LMI-DIST}, it follows that
\eqref{eq:SAPD-parameter-choice-R1} provides us with a particular parameter choice for SAPD \sa{such that the last two inequalities in~\eqref{eq:n123-R1} hold}. Indeed, our choice in \eqref{eq:SAPD-parameter-choice-R1} satisfies \eqref{Condition: SAPD simple LMI system} which is a simpler MI obtained by setting $\theta = \rho$ in \cref{eq: general SAPD LMI_R1}; therefore, $\rho=\theta\in(0,1)$ provides us with an upper bound on the actual convergence rate --see \eqref{eq:distance-rate_R1}.
To compute the upper complexity bound for SAPD, we 
%only need to 
\rev{next} analyze how $N$ should grow depending on $\epsilon$ such that the first inequality 
in~\eqref{eq:n123-R1} holds.
%     It follows from \eqref{eq:SAPD-parameter-choice}, \cref{Thm: main result} and
%     \cref{ASPT: compact} that
% {\small
% \begin{equation*}
% \frac{1}{K_N(\theta)}\Omega_{\tau,\sigma,\theta}
%     = \frac{\theta^{N-1}}{1-\theta^N}  \Big( (1+\theta)\mu_x{\Omega_f} + (1+ 3\theta) \mu_y {\Omega_g}\Big)/2
%     \leq \frac{\theta^{N-1}}{1-\theta^N}  \left( \mu_x {\Omega_f} + 2\mu_y{\Omega_g}\right).
% \end{equation*}}
\sa{%It follows that
The first inequality in~\eqref{eq:n123-R1} holds for
%\xtodo{Change $\cD'$}
\rev{$N\geq 1 + \ln(6%\max\{\frac{\mu_x}{\mu_y},\frac{\mu_y}{\mu_x}\}
\cD(x_0,y_0)/\epsilon)/{\ln(\tfrac{1}{\theta})}$}. Thus, SAPD can generate a point $(x_\epsilon,y_\epsilon)\in\cX\times\cY$ such that $\cD(x_\epsilon,y_\epsilon)\leq \epsilon$ within} %$N_\epsilon\in\integers_+$ iterations where
{\small
\begin{equation}\label{INEQ: bound of N, noise case-R1}
 N_\epsilon = \mathcal{O}\Big(  \ln\Big(\frac{
 %\max\{\frac{\mu_x}{\mu_y},\frac{\mu_y}{\mu_x}\}\|z_0-z^*\|^2
 \cD(x_0,y_0)}{\epsilon}\Big)/\ln(\tfrac{1}{\theta}) \Big)
\end{equation}}%
iterations. \sa{In the remaining part of the proof, we will bound the term $\ln(\tfrac{1}{\theta})^{-1}$ in terms of given $\epsilon>0$. According to \eqref{eq:SAPD-parameter-choice-R1}, $\theta=\max\{\overline{\theta},~\overline{\overline{\theta}}\}$; hence, it follows from \eqref{Condition: SP solution to noisy LMI-R1} and \eqref{eq:theta_bound_2-R1} that $\theta\in\{\overline{\theta}_1,\overline{\theta}_2,\overline{\overline{\theta}}_1,\overline{\overline{\theta}}_2\}\subset (0,1)$. Since $\ln(1/\theta)$ is convex for $\theta\in\reals_{++}$, we immediately get $\frac{1}{\ln(\frac{1}{\theta})}  \leq \frac{1}{1-\theta}$ for $\theta\in(0,1)$. Therefore, we trivially get the bound
\begin{equation}
\label{eq:simple_theta_bound-R1}
\frac{1}{\ln(\frac{1}{\theta})}\leq \max\{(1-\overline{\theta}_1)^{-1}, (1-\overline{\theta}_2)^{-1}, (1-\overline{\overline{\theta}}_1)^{-1}, (1-\overline{\overline{\theta}}_2)^{-1}\}.
\end{equation}
First, we equivalently rewrite $(1-\overline{\theta}_1)^{-1}$
and $(1-\overline{\theta}_2)^{-1}$
as follows:}
{\footnotesize
\begin{equation*}
         \sa{(1-\overline{\theta}_1)^{-1}}
         =  \frac{1}{2}\left(\frac{L_{xx}}{\mu_x}+1\right) + \sqrt{\frac{1}{4}\left(\frac{L_{xx}}{\mu_x} + 1\right)^2 + \frac{\rev{2}L_{yx}^2}{\beta \mu_x\mu_y}},\quad
         \sa{(1-\overline{\theta}_2)^{-1}}
     = \frac{1}{2} + \sqrt{\frac{1}{4} +\frac{\rev{16}L_{yy}^2}{\left(1-\beta\right)^2\mu_y^2}};
\end{equation*}}%
\sa{finally, $(1-\overline{\overline{\theta}}_1)^{-1}=\tfrac{1}{12}{\Xi}^x(\beta)\frac{\delta_x^2}{\mu_x}\frac{1}{\epsilon}$ and $(1-\overline{\overline{\theta}}_2)^{-1}=\tfrac{1}{12}{\Xi}^x(\beta)\frac{\delta_y^2}{\mu_y}\frac{1}{\epsilon}$}.
\sa{Thus, using four identities we derived above within \eqref{eq:simple_theta_bound-R1} and combining it with \eqref{INEQ: bound of N, noise case-R1}, we achieve the desired %upper complexity
bound for SAPD.}
\end{proof}
\begin{remark}
\label{rem:SCMC}
\rev{Whenever $\mu_x\gg\mu_y$ or $\mu_x\gg\mu_y$, the variance bound in~\eqref{eq:complexity-2-SCSC} is better than~\eqref{eq:complexity-1-SCSC}. There is a bias-variance trade-off for this improvement, i.e., $\frac{L_{yx}}{\sqrt{\mu_x\mu_y}}$ term in bias
%{\color{red} for the }bias {\color{red} term in \eqref{eq:complexity-1-SCSC} }
degrades to $\frac{L_{yx}}{\min\{\mu_x,\mu_y\}}$. However, in certain scenarios, the improvement in variance justifies this degradation in bias. For instance, suppose $\cL$ is $\mu_x$-strongly convex in $x$ for $\mu_x=\cO(1)$, there exists $\cD_y$ such that $\norm{y}\leq \cD_y$ for $y\in\dom g$, and $\Phi$ affine in $y$; hence, $L_{yy}=0$ --see DRO problem in~\cref{sec:real-data}. Let $h(x)\triangleq\max_y\cL(x,y)$ denote the primal function. Using Nesterov's smoothing technique in~\cite{nesterov2005smooth}, one can smooth $h$, which leads to an SCSC problem: $\min_x\{ h_{\mu_y}(x)\triangleq\max_y\cL(x,y)-\tfrac{\mu_y}{2}\norm{y}^2$\}, for which choosing the smoothing parameter $\mu_y=\frac{\epsilon}{2\cD_y^2}$ implies $|h(\cdot)-h_{\mu_y}(\cdot)|\leq \epsilon$. To compute an $\epsilon$-solution for the regularized problem with $\mu_y=\Theta(\epsilon)$, \eqref{eq:complexity-1-SCSC} implies 
$\tilde\cO(\frac{\delta_x^2}{\epsilon^{3/2}}+\frac{\delta_y^2}{\epsilon^{2}})$ while \eqref{eq:complexity-2-SCSC}
%\todo{MG: This should be \eqref{eq:complexity-2-SCSC} instead?} 
gives us $\tilde\cO(\frac{\delta_x^2}{\epsilon}+\frac{\delta_y^2}{\epsilon^{2}})$.}
\end{remark}
\begin{remark}
\sa{Given $\epsilon>0$, for sufficiently small $\delta_x^2>0$, \eqref{eq:theta_bound_2-R1} implies that $\overline{\overline{\theta}}_1 = 0$; similarly, $\overline{\overline{\theta}}_2 = 0$ for sufficiently small $\delta_y^2>0$. \rev{Therefore, $\delta_x^2=\delta_y^2=0$ implies $\overline{\overline{\theta}}_1 = \overline{\overline{\theta}}_2 = 0$.}}
% when there is no noise for $\grad_x\Phi$, i.e., $\delta_x^2= 0$, we get $\overline{\overline{\theta}}_1 = 0$; similarly, $\overline{\overline{\theta}}_2 = 0$ if $\delta_y^2 = 0$.}
%\todo{MG:this is not clear; is it the convention we are using?}
\end{remark}
\rev{{Our} bound's variance term {(the term that depends on the noise levels $\delta_x^2$ and $\delta_y^2$)} in \cref{Prop: iteration complex for SAPD_R1} is optimal {with respect to its dependency on $\epsilon$} up to a log factor, which can further be
eliminated through employing a restarting strategy in the lines of our previous work~\cite{aybat2019universally}} \fin{(see \cref{sec:m-sapd}).}

\section{Robustness and Convergence Rate \sa{Trade-off}}\label{sec:robustness}
In this section, \sa{assuming $\mu_x,\mu_y>0$}, we study \sa{the trade-off\mg{s} between robustness-to-gradient noise and the convergence rate depending on the choice of SAPD parameters, \fin{i.e., bias-variance trade-off for SAPD.}}
Given the saddle point $z^*\triangleq (x^*,y^*)\in \mathcal{X}\times \mathcal{Y}$ of  \eqref{eq:main-problem}, we first define the \emph{robustness} as follows:
{\small
\begin{equation}\label{def-robustness-measure}
 {\cJ\triangleq\limsup_{N\to\infty}\cJ_N},\quad \hbox{where}\quad \mathcal{J}_N \triangleq {\mathbb{E}\Big[ 
 %\frac{1}{\delta_x^2}\|x_N -x^* \|^2+ \frac{1}{\delta_y^2}\| y_N - y^*\|^2
 \|z_N-z^*\|^2\Big]/\delta^2,}
\end{equation}}%
{where $z_N\triangleq (x_N,y_N)$ and $\delta\triangleq\max\{\delta_x,\delta_y\}$.} \mg{The quantity $\mathcal{J}_N$ %reflects the variance 
\fin{is the expected squared distance of \sa{$z_N$} to $z^*$}, %in each iteration
normalized by the level of gradient noise {$\delta^2$.} %$\delta_x^2$ and $\delta_y^2$. %each step.
%It can therefore be viewed as a measure of
\sa{Thus, using $\cJ$ we measure how much SAPD %iterates
amplifies the gradient noise asymptotically.} Due to the persistent stochastic noise, %the iterates
%$(x_N, y_N)$
\sa{$\{z_N\}$ does} not typically converge to $z^*$ but oscillate around it with a positive variance. The limit $\mathcal{J}$ provides a bound on the \mg{expected} size of \sa{neighborhood %that the iterates
$\{z_N\}$ \mg{accumulates in},} %converge to,
\sa{i.e., from Jensen's lemma, we get}
%it follows that in the sense that we have
{\small
$${\limsup_{N\to\infty}\mathbb{E}[\norm{z_N-z^*}] \leq \limsup_{N\to\infty} \sqrt{\mathbb{E}\left[\|x_N -x^*\|^2+ \| y_N - y^*\|^2\right]}  
%\leq \max\{\delta_x, \delta_y\} \sqrt{\mathcal{J}}
=\delta\sqrt{\mathcal{J}}.}$$}}%
Therefore, smaller values of $\mathcal{J}$ will lead to better robustness to noise and will give a better asymptotic performance. % %Next, we will study the robustness of SAPD algorithms as a function of the choice of the parameters.
\sa{Below we derive an explicit characterization of $\mathcal{J}$ for
%some particular quadratic
\sa{a particular class of} SCSC problems; %on the other hand,
and we will obtain an upper bound on $\mathcal{J}$ for more general SCSC problems in~\cref{sec:robustness-bound}.}
%Although it seems abstract, we analyze it from two aspects, i.e., an upper bound and an explicit computation.
\subsection{Explicit Estimates for Robustness to Noise}\label{sec:explicit-robustness}
We consider the special case {of \eqref{eq:main-problem}} %\mg{when $\Phi$ is a quadratic of the form}
{when $\Phi$ is bilinear and  $f,g$ have simple quadratic forms, i.e.,}\vspace*{-2mm}
{\small
\begin{equation}
{\Phi(x,y) = \langle Kx, y \rangle,\ f(x)=\frac{\mu_x}{2} \|x\|^2,\ g(y)=\frac{\mu_y}{2}\|y\|^2,}
\label{eq:special-Phi}
\vspace*{-2mm}
\end{equation}}%
where $K\in \mathbb{R}^{d\times d}$ is a symmetric matrix, {and the noise is additive, i.e.,}
{\small
\begin{equation}
\label{assump-additive-noise}
\begin{aligned}
         \Tilde{\nabla}_x\Phi(x_k,y_{k+1}{;\omega^x_k}) = \nabla_x\Phi(x_k,y_{k+1}) + \omega_k^x,\quad
         \Tilde{\nabla}_y\Phi(x_k,y_k{;\omega^y_k}) = \nabla_y\Phi(x_k,y_k)+\omega_k^y,
\end{aligned}
\end{equation}}%
%which %is a special case of
{satisfying} \cref{ASPT: unbiased noise assumption}. We also assume that {there exists $\delta>0$ such that} %for each $k\geq 0$, 
\rev{$\{w_k^x\}$ and $\{w_k^y\}$ are i.i.d Gaussian with zero mean and %an
 isotopic covariance,} i.e.,
% \todo{we use $\delta_x=\delta_y=\delta$, max is also true here}
%that is a multiple of the identity, i.e.,
{\small
\begin{equation}\label{assump-gaussian-noise} 
\rev{\mathbb{E}\left[w_k^x\right]=0_d,\quad \mathbb{E}\left[w_k^y\right]=0_d,}\quad
\mathbb{E}\left[w_k^x (w_k^x)^\top\right] =  \frac{{\delta^2}}{d} I_d, \quad  \mathbb{E}\left[w_k^y (w_k^y)^\top\right] = \frac{{\delta^2}}{d} I_d.
\end{equation}}%
Clearly, the unique saddle point to \eqref{eq:special-Phi} is the origin, i.e., $(x^*,y^*)={(0_d, 0_d)}$.
% \todo{Sep 7th, xuan: I add this notation for simplicity}
% \newcomment{\begin{definition}\label{def: quadratic function parameter}
%  $F(\mu_x,\mu_y,\|K\|_2,d,\delta_x,\delta_y)(x,y)\triangleq  \cL(x,y)$ s.t. \cref{eq:special-Phi,assump-additive-noise,assump-gaussian-noise}.
% \end{definition}}
% \xuan{For simplicity, we use the notation $F(\mu_x,\mu_y,\|K\|_2,d,\delta_x,\delta_y)$ to claim the properties of the quadratic function \cref{eq:special-Phi}.}
{We will show that the robustness measure $\mathcal{J}$ defined in \eqref{def-robustness-measure} is finite,
%$$ \mathcal{J}= \limsup_{k\to\infty} \mathbb{E}\norm{x_k}^2/\delta^2+ \mathbb{E}\norm{y_k}^2/\delta^2,
%$$
%which reflects the asymptotic variance of the iterates. Furthermore,
and {that it admits a closed form solution.}}
{We first note that} according to \cref{ALG: SAPD}, {for $k\geq 0$,} \vspace*{-2mm}%we have the following iterations
{\small
\begin{equation}
\begin{aligned}
            \rev{x_{k}} &= \rev{\frac{1}{1+\tau\mu_x}\big(x_{k-1} - \tau K^\top y_{k} - \tau \omega_{k-1}^x\big)},\\
             y_{k+1} &= \frac{1}{1+\sigma\mu_y}\left(y_k + \sigma(1+\theta)K x_{k} - \sigma\theta K x_{k-1} + \sigma(1+\theta)\omega^y_k - \sigma\theta\omega^y_{k-1}\right).
\end{aligned}
\end{equation}}%
Next, {for $k\geq 0$, we %consider the vector
define}
%\nsa{$z_k$ is not a good notation as we reserve $z=(x,y)$}
\fin{${\tilde z_k} \triangleq [x_{k-1}^\top~ y_{k}^\top]^\top\in\reals^{2d}$} and
% which is a {vertical} concatenation {of iterates from %the last
% two consecutive iterations,
% %at step $k$. We also introduce the vector
% and}
\fin{$w_k \triangleq
%\begin{bmatrix}
[(w^x_{k-1})^\top (w^y_{k-1})^\top  (w^y_k)^\top]^\top
%\end{bmatrix}
\in\reals^{3d}$},
which is the {vertical} concatenation of the %realization of the
noise realization at {step} $k-1$ and $k$. %After a straightforward computation, we see that
{The %column
vector} {$\tilde{z}_k$} satisfies
%the following \mg{affine} recursion:
\begin{equation} \rev{{\tilde z_{k+1}} = A {\tilde z_k} + B w_k,} %\quad \hbox{where}
\label{lin-dyn-sys}
\end{equation}
\rev{\footnotesize
$$A\triangleq\begin{bmatrix}
\frac{1}{1+\tau\mu_x} I_d & \frac{-\tau}{(1+\tau\mu_x)} K^\top  \\
\frac{1}{1+\sigma\mu_y}\Big(\frac{\sigma(1+\theta)}{1+\tau\mu_x} - \sigma\theta \Big)K
&
\frac{1}{1+\sigma\mu_y}\Big( I_d -\frac{\tau\sigma(1+\theta)}{1+\tau\mu_x}KK^\top\Big)
\end{bmatrix},\ 
B\triangleq\begin{bmatrix}
\frac{-\tau}{1+\tau\mu_x}I_d & 0_d & 0_d \\
\frac{-\tau\sigma(1+\theta)}{(1+\tau\mu_x)(1+\sigma\mu_y)}K & \frac{-\sigma\theta}{1+\sigma\mu_y}I_d  & \frac{\sigma(1+\theta)}{1+\sigma\mu_y}  I_d
\end{bmatrix}.$$}%
% \newcomment{Therefore, $z_{k}=A^{k}z_0 + \sum_{i=0 }^{k-1}A^{i}B\omega_{k-1-i}$, thus the convergence rate is }
{From \eqref{lin-dyn-sys}, \rev{using the noise model}, it is easy to see that %the covariance matrix
$\Sigma_k \triangleq \mathbb{E} {[\tilde z_k \tilde z_k^\top]}$}
%at time $k$
satisfies %the recursion
\rev{\footnotesize
\begin{equation*}
    \begin{aligned}
        \Sigma_{k+1}  &= A \Sigma_k A^\top + \frac{\delta^2}{d} BB^\top  +  \mathbb{E}[B\omega_k\tilde{z}_k^\top A^\top + A\tilde{z}_k\omega_k^\top B^\top ]\\
        & = A \Sigma_k A^\top + \frac{\delta^2}{d} BB^\top  +  \mathbb{E}[B\omega_k(A \tilde z_{k-1}+ B w_{k-1})^\top A^\top + A(A \tilde z_{k-1} + B w_{k-1})\omega_k^\top B^\top]
        % \\
        % &
        % =A \Sigma_k A^\top + \frac{\delta^2}{d} BB^\top  +  \mathbb{E}[B\omega_k \tilde z_{k-1}^\top A^\top A^\top+ B\omega_k w_{k-1}^\top B^\top A^\top + A A \tilde z_{k-1} \omega_k^\top B^\top
        % +
        %  A B w_{k-1}\omega_k^\top B^\top]
        %  \\
         =  A \Sigma_k A^\top + \frac{\delta^2}{d}R
    \end{aligned}
\end{equation*}}%
for $k\geq 0$, where 
% $R \triangleq R_1 + R_2+ R_2^\top$ such that
% {\footnotesize
% $$
% W \triangleq \begin{bmatrix}
%     0_d & 0_d & 0_d \\
%     0_d & 0_d & I_d \\
%     0_d & 0_d & 0_d \\
% \end{bmatrix}, R \triangleq R_1 + R_2+ R_2^\top
% $$}%
% and
% {\footnotesize
% $$
% \begin{aligned}
%     & R_1 \triangleq BB^\top =
% \begin{bmatrix}
%     \frac{\tau^2}{(1+\tau\mu_x)^2}I_d & \frac{\tau^2\sigma(1+\theta)}{(1+\tau\mu_x)^2(1+\sigma\mu_y)}K^\top \\
%     \frac{\tau^2\sigma(1+\theta)}{(1+\tau\mu_x)^2(1+\sigma\mu_y)}K &  \frac{\tau^2\sigma^2(1+\theta)^2}{(1+\tau\mu_x)^2(1+\sigma\mu_y)^2}KK^\top + \Big(\frac{\sigma^2\theta^2}{(1+\sigma\mu_y)^2} + \frac{\sigma^2(1+\theta)^2}{(1+\sigma\mu_y)^2}\Big)I_d
% \end{bmatrix},\\
% & R_2 \triangleq BWB^\top A^\top =
% \begin{bmatrix}
%     0_d & 0_d \\
%     0_d & -\frac{\sigma^2\theta(1+\theta)}{(1+\sigma\mu_y)^2}I_d
% \end{bmatrix}A^\top = \begin{bmatrix}
%     0_d & 0_d \\
% \frac{\sigma^2\tau\theta(1+\theta)}{(1+\tau\mu_x)(1+\sigma\mu_y)^2}K
% & -\frac{\sigma^2\theta(1+\theta)}{(1+\sigma\mu_y)^3}(I_d - \frac{\tau\sigma(1+\theta)}{1+\tau\mu_x}KK^{\top})
% \end{bmatrix}
% \end{aligned}
% $$}% 
%Therefore, we can further compute
\rev{\footnotesize
$R \triangleq \begin{bmatrix}
    c_1 I_d & c_2 K^\top \\
    c_2 K & c_3KK^\top + c_4 I_d
\end{bmatrix}
$} \rev{and $\{c_i\}_{i=1}^4$ are some constants.}\footnote{\label{fnt:constants}
%\begin{equation}\label{eq:R-c}
%\begin{aligned}
\rev{These constants can be computed explicitly as follows: $c_1 \triangleq \frac{\tau^2}{(1+\tau\mu_x)^2}$, $c_2 \triangleq c_1\frac{\sigma(1+\theta)}{1+\sigma\mu_y} +  \sqrt{c_1}\theta(1+\theta)\frac{\sigma^2}{(1+\sigma\mu_y)^2}$, $c_3 \triangleq %\frac{\tau^2\sigma^2(1+\theta)^2}{(1+\tau\mu_x)^2(1+\sigma\mu_y)^2} 
    (1+\theta)^2\frac{\sigma^2}{(1+\sigma\mu_y)^2}(c_1+ \sqrt{c_1}\frac{2\sigma\theta}{1+\sigma\mu_y})$, and $c_4 \triangleq 
     %\frac{\sigma^2\theta^2}{(1+\sigma\mu_y)^2} + \frac{\sigma^2(1+\theta)^2}{(1+\sigma\mu_y)^2} - \frac{2\sigma^2\theta(1+\theta)}{(1+\sigma\mu_y)^3}
     \frac{\sigma^2}{(1+\sigma\mu_y)^2}(1+2\theta(1+\theta)\frac{\sigma\mu_y}{1+\sigma\mu_y})$.}}
%\end{aligned}
%\end{equation}}%
Linear dynamical systems subject to Gaussian noise {such as \eqref{lin-dyn-sys}} have been well studied in the robust control literature. {In fact, %it is known that 
the limit
$ \Sigma_\infty\triangleq \lim_{k\to\infty}\Sigma_{k}
$
exists if the spectral radius %\footnote{{The spectral radius of a square matrix $A$ is defined as the largest absolute value of its eigenvalues.}} 
of $A$, {denoted by $\rho(A)$}, is less than one}, {and it} satisfies the Lyapunov equation:
%if we take the limit in \eqref{eq-cov-matrix-dynamics} we obtain
%the Lyapunov equation
{\small
%\begin{equation}
$     \Sigma_\infty \triangleq {A \Sigma_\infty A^\top} + \frac{\delta^2}{d} \rev{R},
$
 %    \label{eq-lyap}
%\end{equation}
}%
%whose solution is unique as long as the spectral radius $\rho(A)<1$.
%where the solution can be expressed as the series sum}
%the series sum}
whose solution is given {in the form of an infinite series} %by the series sum
{\small
%\begin{equation}\label{eq-variance-formula}
%\Sigma_\infty = \sum_{j=0}^\infty (A^T)^j B B^T A^j
$\rev{\Sigma_\infty = \frac{\delta^2}{d}\sum_{k=0}^\infty A^k R (A^k)^\top}$}
%\end{equation}}%
(see~\cite{zhou1996robust}). %It is also easy to see
Note $\mathcal{J} ={\limsup_{k\to\infty}\frac{1}{\delta^2}\mathbb{E}[\norm{z_k-z^*}^2]=\lim_{k\to\infty}\frac{1}{\delta^2}\Tr(\Sigma_k)}=\frac{1}{\delta^2} %\mbox{trace}
\Tr (\Sigma_\infty)=\sa{\frac{1}{d}\sum_{k=0}^\infty \Tr(A^kR(A^k)^{\top})}$.
{We also observe that $\mathcal{J}$ is {\it invariant} under orthogonal transformations, i.e., for any orthogonal matrix $Z$, $\tilde{A}\triangleq Z^\top A Z$ and \rev{$\tilde{R}\triangleq Z^\top R {Z}$}
%\todo{$\tilde{B}\triangleq Z^\top B Z$?(corrected already)}
satisfy\vspace*{-2mm}
{\small
\begin{equation}
\mathcal{J} =\frac{1}{\sa{\delta^2}}
\Tr(\tilde{\Sigma}_\infty), \quad \mbox{{\normalsize where}} \quad  \tilde{\Sigma}_\infty \triangleq \frac{\delta^2}{d}\sum_{k=0}^\infty \rev{\tilde{A}^k \tilde{R} (A^k )^\top} = Z^\top \Sigma_\infty Z
\label{eq-transformed-J}\vspace{-2mm}
\end{equation}}%
solves the transformed Lyapunov equation
{%\small
%\begin{equation}
$\tilde{\Sigma}_\infty = \tilde{A} \tilde{\Sigma}_\infty \tilde{A}^\top + \frac{\delta^2}{d}  \rev{\tilde{R}}$. %
%\label{eq-transformed-Lyap}
%\end{equation}
}}%
In order to compute $\mathcal{J}$ explicitly, 
%in the following 
we will choose \sa{a particular} orthogonal matrix $Z$ %in a special way
so that solving {the transformed Lyapunov equation} %\eqref{eq-transformed-Lyap}
explicitly will be simple. 
%than solving the Lyapunov equation \mg{without any transformation}. %\eqref{eq-lyap}.
First, we consider the eigenvalue decomposition of $K\in\mathbb{S}^d$, i.e., $K = U \Lambda U^\top$, 
%of the $K$ matrix
where $\Lambda$ is a diagonal matrix such that $\Lambda_{ii} = \lambda_i$, \sa{and \mg{$\{\lambda_i\}_{i=1}^d$} are the} eigenvalues %of $K$ 
in increasing order:
$ \lambda_1 \leq \lambda_2 \leq \dots \leq \lambda_d$. Then, 
%we can write 
$A = V A_\Lambda V^\top$, where
\rev{\footnotesize
$V \triangleq \begin{bmatrix}
U & 0_d \\
0_d &  U
\end{bmatrix}$} and
\rev{\footnotesize
$
% A_\Lambda \triangleq \begin{bmatrix}
% \frac{1}{1+\tau\mu_x} I_d & \frac{-\tau}{(1+\tau\mu_x)} \Lambda  \\
% \frac{1}{1+\sigma\mu_y}\Big(\frac{\sigma(1+\theta)}{1+\tau\mu_x}  - \sigma\theta \Big)\Lambda
% &
% \frac{1}{1+\sigma\mu_y}\Big( I_d -\frac{\tau\sigma(1+\theta)}{1+\tau\mu_x}\Lambda ^2\Big)
% \end{bmatrix},
A_\Lambda \triangleq \begin{bmatrix}
a_1 I_d & a_2 \Lambda  \\
a_3 \Lambda
&
a_4\Lambda ^2+a_5 I_d
\end{bmatrix}$}%
 \rev{for constants
{\footnotesize
$
%\begin{aligned}
a_1 = \frac{1}{1+\tau\mu_x}$, $a_2 = \frac{-\tau}{1+\tau\mu_x}$, $a_3 = \frac{\sigma}{1+\sigma\mu_y}\Big(\frac{1+\theta}{1+\tau\mu_x}  - \theta \Big)$, $a_4 = \frac{-\tau\sigma(1+\theta)}{(1+\tau\mu_x)(1+\sigma\mu_y)}$, $a_5 = \frac{1}{1+\sigma\mu_y}$.}
%\end{aligned}
}%
%and \rev{$V\in\mathbb{S}^{2d}$ is a block-diagonal matrix with each block on the diagonal equal to $U$, i.e., $V=\diag(\{U\}_{i=1}^2)$.}
% {
% $$ V = \begin{bmatrix}
% U & 0_d & 0_d & 0_d \\
% 0_d &  U & 0_d & 0_d\\
% 0_d & 0_d & U & 0_d \\
% 0_d & 0_d & 0_d & U
% \end{bmatrix}.$$}
\mg{Furthermore, %it is straightforward to see that
we can permute the entries of $A_\Lambda$ so that it becomes a %$4\times 4$
block diagonal matrix}, %In other words,
\sa{i.e.,} there exists a permutation matrix $P$ such that %\todo{MG: It may be good to write down this perm. matrix explicitly?}
$PA_\Lambda P^\top =  \diag(\{\sa{\tilde{A}_i}\}_{i=1}^d)\triangleq \mg{\tilde{A}}$, where for each $i\in\{1,\ldots, d\}$, \rev{$\mg{\tilde{A}_i}\in\reals^{2\times 2}$ is
%a $4\times 4$ matrix
defined by
{\footnotesize
$
\tilde{A}_i\triangleq
% \begin{bmatrix}
% \frac{1}{1+\tau\mu_x}  & -\frac{\tau}{(1+\tau\mu_x)} \lambda_i \\
% \frac{1}{1+\sigma\mu_y}\Big(\frac{\sigma(1+\theta)}{1+\tau\mu_x}\lambda_i  - \sigma\theta \lambda_i \Big)
% &
% \frac{1}{1+\sigma\mu_y}\Big( 1 -\frac{\tau\sigma(1+\theta)}{1+\tau\mu_x}\lambda_i ^2\Big)
% \end{bmatrix} = 
\begin{bmatrix}
    a_1 & a_2\lambda_i \\
    a_3\lambda_i & a_4\lambda_i^2 +a_5
\end{bmatrix}.
$
}}%
%\xtodo{Do we need those $a_1,a_2,a_3,a_4$ in article?}
%where
\sa{Thus, for $Z = VP^\top$, we have $\tilde{A} = Z^\top A Z$, and \rev{$\tilde{R} =Z^\top R Z=\mbox{diag} \{ \tilde{R} _i\}_{i=1}^d$}} \sa{such that}
\rev{\footnotesize
$\tilde{R}_i\triangleq \begin{bmatrix}
c_1 & c_2\lambda_i \\
c_2\lambda_i & c_3\lambda_i^2 + c_4
\end{bmatrix}
$,
}%
where \rev{$c_1,c_2,c_3$ and $c_4$ are explicitly given in~\cref{fnt:constants}.} %related to the definition of $R$.}
%Note that 
\sa{Both $\tilde{A}$ and $\tilde{R}$} have a block diagonal structure; therefore, %we see
$\tilde{\Sigma}_\infty = \diag( \{\tilde{S}_i\}_{i=1}^d)$,
where for each $i\in\{1,\ldots, d\}$, $\tilde{S}_i$ %solves
is the unique solution to
\begin{align}
\label{eq:Lyapunov-eq}
    \tilde{S}_i = \tilde{A}_i \tilde{S}_i \tilde{A}_i^\top + \frac{\sa{\delta^2}}{d}  \rev{\tilde{R}_i}.
\end{align}
\rev{This Lyapunov equation} %~\eqref{closed form H2 equation}
is a $2\times 2$ system, which %we can solve
\sa{can be solved} for $\tilde{S}_i$ \rev{explicitly by inverting a $3\times 3$ symbolic matrix --since $\tilde{S}_i$ is symmetric, one needs to %compute
\fin{solve for 1 off-diagonal and 2 diagonal elements.}} Using \eqref{eq-transformed-J} and %\eqref{eq-soln-Lyapunov-eqn}
$\tilde{\Sigma}_\infty = \diag( \{\tilde{S}_i\}_{i=1}^d)$ will yield us an explicit formula for $\mathcal{J}$. 
% \mg{%We solve \eqref{closed form H2 equation}
% \sa{$\tilde{S}_i$ can be computed using} the ``Matlab Symbolic Toolbox",
% %this allows
% \sa{allowing} us to compute the robustness measure $\mathcal{J}$ explicitly}. 

\sa{Next, for $A$ in \cref{lin-dyn-sys},
we define $\rho_{\rm{true}} \triangleq (\rho(A))^2$,
which determines the exact \mg{(asymptotic)} convergence rate of \rev{$\mathbb{E}[\norm{\tilde z_N-z^*}^2]$; hence, {$\mathbb{E}[d_N^*]$} in Theorem~\ref{Thm: main result_R1} also converges with this asymptotic rate.}
Furthermore, \rev{it can also be shown for this quadratic model that $\mathbb{E}[\sup_{(x,y)\in\mathcal{X}\times \mathcal{Y}}\{\mathcal{L}(x_k,y)-\mathcal{L}(x,y_k)\}]$ converges with the same rate {(see \cref{sec:gap-rate} for more details)}}}.
%\nsa{Check if this is still true.}
%\xtodo{It is still True, but the appendix has been removed. Maybe we can say because of the good property of bilinear quadratic}
%It should be noted that
%This way \mg{of computing} 
Robustness measure $\cJ$ and convergence rate $\rho_{\rm true}$ computed in this section are independent of our theoretical analysis of the SAPD algorithm; they reflect the exact \rev{asymptotic} behavior of the algorithm for a quadratic function \sa{in~\eqref{eq:special-Phi}}, and establishing the interaction between $\cJ$ and $\rho_{\rm true}$ helps us understand some fundamental relations for SAPD.\looseness=-1

\begin{figure}[h]
\vspace{-6mm}
\label{fig:tradeoff}
\centering
\subfigure[]{
\label{fig:fundamental trade-off curve}
\includegraphics[width = 0.45\textwidth]{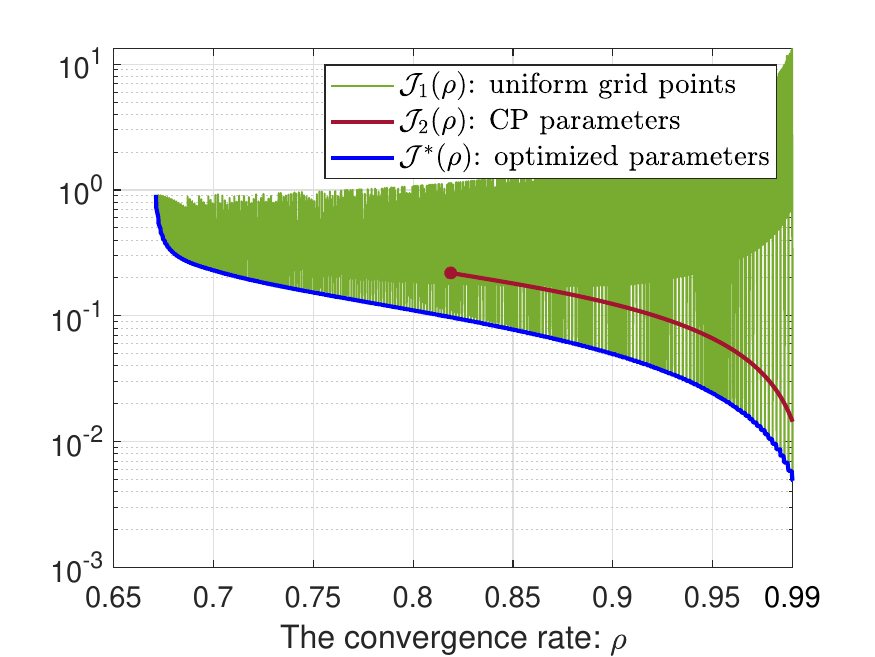}
}
\centering
\subfigure[]{
\label{fig:trade-off bewtween R and rho_star}
\includegraphics[width = 0.45\textwidth]{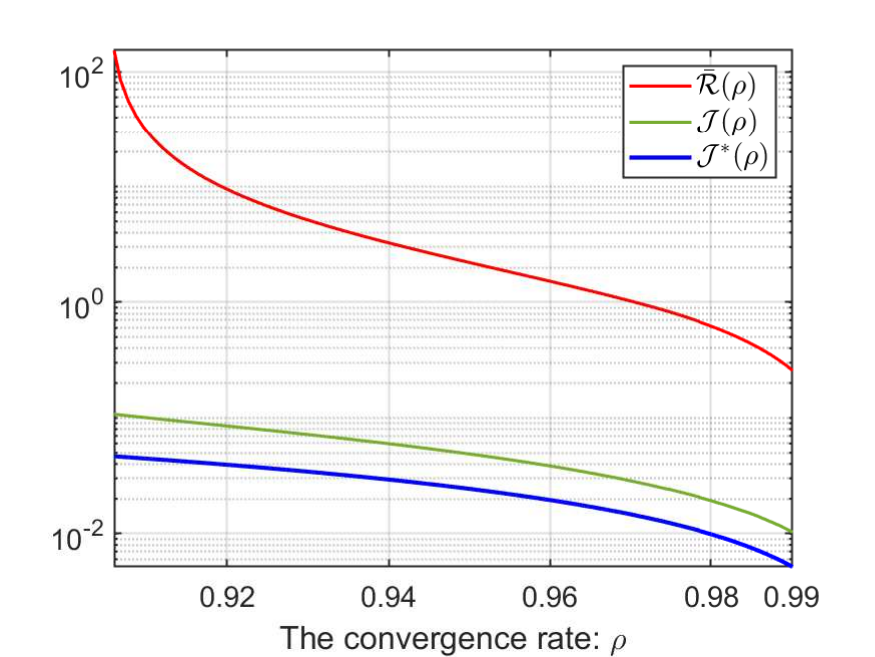}}
\vspace*{-2mm}
\caption{The rate-robustness trade-off \sa{for \eqref{eq:special-Phi} when $\mu_x=\mu_y=1$, $\|K\|_2=10$, $d=30$ and $\delta_x=\delta_y=\delta=10$.} The best \sa{achievable} rate is $0.67$. \sa{The point indicated with a red $``*"$} in \cref{fig:fundamental trade-off curve} is the particular \mg{choice} of CP parameters given in~\cite[Eq.(49)]{chambolle2016ergodic}. \sa{\cref{fig:trade-off bewtween R and rho_star} illustrates that employing the SAPD parameters obtained through minimizing $\bar{\cR}$, an upper bound on $\cJ$ defined in~\cref{sec:robustness-bound}, one can closely track the efficient frontier $\cJ^*$. The best certifiable rate is $0.9049$.}} \vspace*{-5mm}
%The minimum convergence rate is 0.6716 in the left picture.}
\end{figure}

{We numerically illustrate the fundamental rate-robustness trade-off \sa{in \cref{fig:tradeoff} for \eqref{eq:special-Phi} through plotting 3 curves: $\cJ_1$, $\cJ_2$ and $\cJ^*$}.
For $\cJ_1$, we %equally
\sa{uniformly} grid the %step size
\sa{parameter} space $(\tau,\sigma,\theta)\in[0,0.5]\times[0,0.5]\times[0,2]$ \sa{using} {$500\times 500 \times 200$} points;
then, for each \sa{grid} point, we compute the corresponding $(\rho_{\text{true}},\cJ)$ values and plot it. For $\cJ_2$, we employ the step sizes \mg{suggested} in \cite[Algorithm 5]{chambolle2016ergodic} for the CP method\footnote{{Although our method SAPD generalizes the CP method beyond the bilinear problem, SAPD coincides with CP on this particular problem as it has a \sa{bilinear coupling function $\Phi$}.}} and plot $(\rho_{\text{true}},\cJ)$ \mg{(see the \cref{section: CP parameters explanation} for more details)}.
For $\cJ^*$, defining $\cJ^*(\rho)\triangleq\min_{\tau,\sigma,\theta\geq0}\{\cJ:\rho_{\text{true}} = \rho\}$, we plot $(\rho,\cJ^*(\rho))$ \mg{which illustrates the best robustness that can be achieved for a given rate}.}
\sa{In \mg{the} $\cJ_1$ plot, there are vertical lines as there exist many points in the grid} sharing the same %convergence
rate while they have %totally
\sa{very} different robustness values. As seen in~\cref{fig:fundamental trade-off curve}, for great majority of parameter choices \sa{from} the uniform grid, the corresponding robustness is very poor, i.e., very high $\cJ$ value. As a consequence, we infer that it \mg{is} necessary to control the robustness through properly tuning the algorithm parameters.
\sa{The $\cJ_2$ plot demonstrates that} for %the same 
\rev{a fixed rate} \sa{CP parameter choice} ensures relatively lower $\cJ$ values compared to the majority of points in the uniform grid; \sa{but, $\cJ_2$} is still far away from the efficient frontier $\cJ^*$.
\sa{%From the
As indicated in $\cJ_2$ plot}, %we observe that
the best convergence rate CP parameters can achieve is only around {$0.83$}, while the best rate \sa{achieved among the uniform grid} is %less than
\sa{$0.67$}. \looseness=-1

\sa{%We emphasize that
While the parameter optimization problem to compute $\cJ^*(\rho)$ for a given convergence rate $\rho\in(0,1)$ can be done for the special case of \eqref{eq:main-problem} %defined by
%\fin{corresponding to} 
given in \eqref{eq:special-Phi}, this is not a trivial task for a more general coupling function $\Phi$; therefore,} we provide an alternative %formulation
\fin{model} to achieve a similar trade-off result between an \textit{upper bound} on $\cJ$ and a \textit{bound on the convergence rate} in \cref{proposition robustness bound nonconvex}.\vspace*{-5mm}

%\todo{Xuan: I add a comment for the special point in \cref{fig:fundamental trade-off curve} in the caption.}
%\nsa{In \cref{fig:trade-off bewtween R and rho_star}, it says $\cR$; but, should not it be $\bar{\cR}$?(xuan: fixed)}
\begin{figure}[htbp]
\centering
\subfigure[]{
\label{fig:rho_level}
\includegraphics[width = 0.4\textwidth]{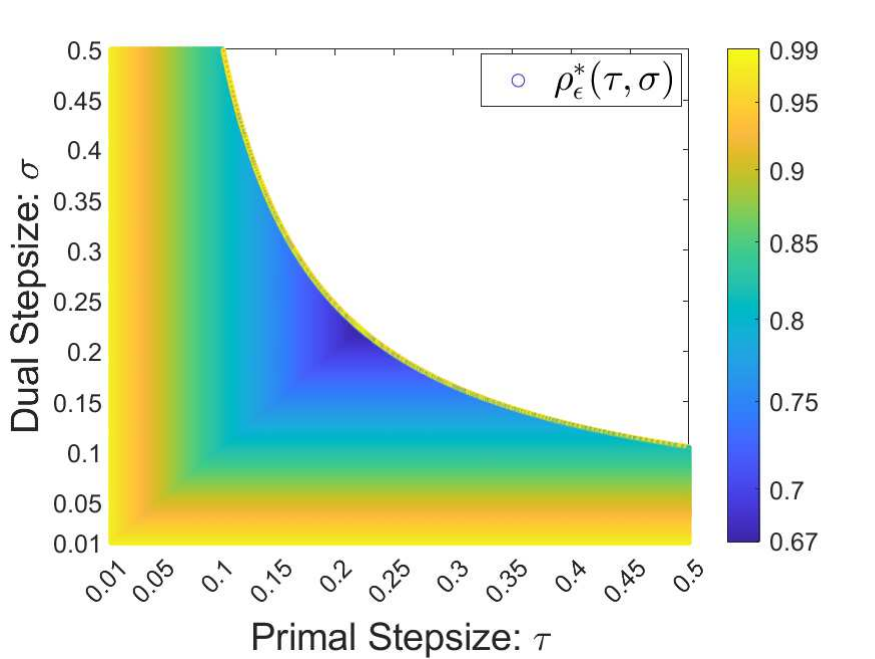}
}
\centering
\subfigure[]{
\label{fig:robustness_level}
\includegraphics[width = 0.4\textwidth]{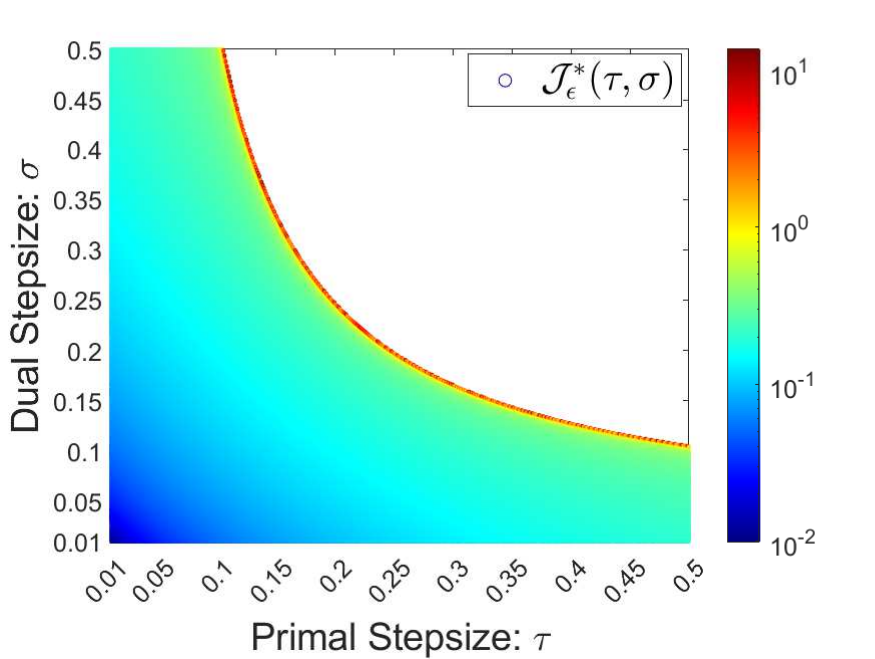}
}
\vspace*{-3mm}
\caption{\sa{The effect of the step sizes on rate and robustness of SAPD running on \eqref{eq:special-Phi} when $\mu_x=\mu_y=1$, $\|K\|_2=10$, $d=30$ and $\delta_x=\delta_y=\delta=10$}. The best \sa{achievable} rate is $0.67$.
%The stepsizes $\tau,\sigma$ satisfy $\exists \theta, s.t.~ \rho_{\text{true}} < 1-\epsilon$, where $\epsilon=0.01$. The minimum convergence rate is 0.6716 in the left picture.
}\vspace*{-3mm}
\end{figure}
\rev{Next, we analyze how primal-dual step sizes, $\tau$ and $\sigma$, affect the convergence rate and the robustness level.
{For any given $\epsilon\in(0,1)$,} we define
$$\rho_\epsilon^*(\tau,\sigma)\triangleq\min_{{\theta}\geq 0}\{ \rho_{\text{true}} :\rho_{\text{true}}{\leq} 1 -\epsilon\},\quad \cJ^*_{\epsilon}(\tau,\sigma)\triangleq \min_{\theta\geq 0}\{\cJ:\rho_{\text{true}} {\leq} 1 -\epsilon\}.$$}%
\sa{We consider
the same experiment described in the caption of Figure~\ref{fig:tradeoff},} setting x-axis as
$\tau$, and y-axis as %the dual stepsize
$\sigma$, we plot %the convergence rate
$\rho_{\epsilon}^*(\tau,\sigma)$ in \cref{fig:rho_level} and $\cJ^*_{\epsilon}$ in \cref{fig:robustness_level}, for $\epsilon = 0.01$.
\sa{We observe that except for the boundary points, simultaneously increasing $\tau$ and $\sigma$
leads to a faster convergence rate
at the expense of a decrease in robustness level --  as one approaches the boundary, there is a significant increase in both convergence rate coefficient $\rho^*_\epsilon(\tau,\sigma)$ and $\cJ$ values.}
These results illustrate the fundamental trade-offs between the convergence rate %versus
\fin{and} robustness for SAPD. 
\vspace{-2mm}
\subsection{An Upper Bound for the Robustness \sa{Measure $\cJ$}}
\label{sec:robustness-bound}
%As we mentioned earlier,
$\cJ$ is hard to compute \sa{in general}; to alleviate this issue, we can alternatively minimize an upper bound on $\cJ$ to control \mg{the robustness level.} %$\cJ$.
\mg{We start with a proposition that provides an upper bound} %\sa{on $\cJ$}.
%of the robustness $\mathcal{J}_N$.
\sa{on the robustness measure $\mathcal{J}$.}
\begin{theorem}\label{proposition robustness bound nonconvex}
Suppose \sa{Assumptions~\ref{ASPT: lipshiz gradient}, %\ref{ASPT: strongly convex concave}, \ref{ASPT: f and g convex}
and \ref{ASPT: unbiased noise assumption}} hold,  and $\{ x_k,y_k \}_{k\geq 0}$ are generated by SAPD stated in Algorithm \ref{ALG: SAPD}, the parameters $\{ \tau, \sigma,\theta \}$ satisfy the conditions in \cref{Thm: main result_R1} %\xuan{some $\rho\in(0,1)$}
\sa{for some $\alpha \in [0,  \tfrac{1}{\sigma})$ and $\rho \in (0,1)$}. Then, for 
{$\delta \triangleq \max\{\delta_x,\delta_y\}$}, $\cJ$ can be bounded as follows:\vspace*{-2mm}
\rev{{\small
\begin{equation}
\label{eq:robustness-bound}
            {\cJ=\limsup_{N\rightarrow\infty} \mathcal{J}_N
            \leq
             \frac{2\rho}{1-{\rho}}\cdot\max\Big\{\tau,\frac{\sigma}{1-\alpha\sigma}\Big\}\cdot B_{\tau,\sigma,\theta},}
             %\Xi_{\tau,\sigma,\theta},
\end{equation}}
where {$B_{\tau,\sigma,\theta}\triangleq \tfrac{\tau}{1+\tau\mu_x} \Xi^x_{\tau,\sigma,\theta}
              + \tfrac{\sigma}{1+\sigma\mu_y} \Xi^y_{\tau,\sigma,\theta}$} and $\Xi^x_{\tau,\sigma,\theta}$, $\Xi^y_{\tau,\sigma,\theta}$ are defined in \cref{Thm: main result_R1}.}%
\end{theorem}
%\mtodo{Can we support $\alpha=1/\sigma$ here in this proposition or not?}
%\xtodo{No, we can't because we need the bound for $y_N-y^*$.}
\begin{proof}
\rev{Let $C_{\tau,\sigma} \triangleq \min\{{\frac{1}{2\tau},\frac{1}{2\sigma}(1-\alpha\sigma)}\}$}. Then, from \cref{Thm: main result_R1}, we have that
\rev{\small
\begin{equation*}
{C_{\tau,\sigma} \delta^2 \mathcal{J}_N
%\leq
=
 C_{\tau,\sigma}\mathbb{E}\left[ \|x_N -x^* \|^2+\| y_N - y^*\|^2 \right]
 \leq 
 \rho^{N} D_{\tau,\sigma} + \frac{\rho}{1-\rho}~  {\Xi}_{\tau,\sigma,\theta},} \vspace*{-3mm}
\end{equation*}}%
\sa{which implies \eqref{eq:robustness-bound} since \rev{$\rho^{N} D_{\tau,\sigma} \to 0$} %goes to zero
as $N\to\infty$, and we have {$\frac{\delta_x^2}{\delta^2}\leq 1$ and $\frac{\delta_y^2}{\delta^2}\leq 1$}.} %goes to infinity.
%\mg{This completes the proof}.
% Letting $N$ go to infinity and dividing the above inequality by $\tfrac{1-\xuan{\rho}}{1-\xuan{\rho}^N}*C_{\tau,\sigma}*\delta^2$, it follows that
% $$
% \mg{\limsup_{N\rightarrow\infty}} \,
% \mathcal{J}_N \leq \frac{1}{1-\xuan{\rho}} *\frac{1}{C_{\tau,\sigma}} *\frac{1}{\delta^2}*\Xi_{\tau,\sigma,\theta}
% $$
%\mg{where we used the fact that }
\end{proof}
\sa{%Note that the
This upper bound %is theoretically correct only 
holds for the SAPD parameters satisfying our step size conditions in \cref{eq: general SAPD LMI_R1}, %Recall that our stepsize condition
\rev{which are only sufficient for ensuring a linear rate; but, %these conditions are 
they may not be necessary.}}

\rev{Next we %can study
investigate the trade-off between the convergence rate bound implied by the matrix inequality in~\cref{eq: general SAPD LMI_R1} and %the upper bound for the robustness.
the robustness upper bound provided in \cref{proposition robustness bound nonconvex}.
%the above result. %proposition. %Furthermore,
}
%\xtodo{Do we still need lemma 3.2 if $\theta$ can be greater than $\rho$?}
 \sa{
 \begin{lemma}
 \label{lem:SDP}
 Given $\rho\in(0,1)$, let $(t_\rho,s_\rho,\theta_\rho,\alpha_\rho)$ be an element of $\cP_\rho$, where
{\footnotesize
\begin{align*}
%\min_{\substack{\theta,t,s\geq 0\\ 0\leq\alpha\leq s}}
&\cP_\rho\triangleq\{(t,s,\theta,\alpha):\ t,s,\theta,\alpha\geq 0,~ \alpha\leq s,~G_\rho(t,s,\theta,\alpha) \succeq 0\}\ \hbox{\normalsize and}\\
&G_\rho(t,s,\theta,\alpha) \triangleq
  \begin{pmatrix}
 (1- \frac{1}{\rho}) t+\mu_x & 0 & 0 & 0 & 0\\
  0 & (1- \frac{1}{\rho}) s+\mu_y & (\frac{\theta}{\rho} - 1)L_{yx} & (\frac{\theta}{\rho} - 1)L_{yy} & 0\\
  0 & (\frac{\theta}{\rho} - 1)L_{yx} & t - L_{xx} & 0 & -  \frac{\theta}{\rho}L_{yx}\\
  0& (\frac{\theta}{\rho} - 1)L_{yy} & 0 & s - \alpha & -  \frac{\theta}{\rho}L_{yy}\\
  0 & 0 & - \frac{\theta}{\rho}L_{yx} & -  \frac{\theta}{\rho}L_{yy} & \frac{\alpha}{\rho}
\end{pmatrix}.
\end{align*}}%
%\todo{Xuan(8.24): Shall we change the order into $t,s,\theta,\alpha\geq 0$ in line 543}
%If $\theta_\rho\leq \rho$,
Then \rev{the bias term, i.e., $\rho^N D_{\tau,\sigma}$ defined in~\eqref{eq:distance-rate_R1}, converges to $0$ with rate $\rho$} for SAPD employing $\tau=1/t_\rho$, $\sigma=1/s_\rho$ and $\theta=\theta_\rho$ . If %$\theta_p>\rho$ or
%the SDP is infeasible,
$\cP_\rho=\emptyset$, then \eqref{eq: general SAPD LMI_R1} does not have a solution for the given $\rho$ value.
\end{lemma}
\begin{proof}
This result immediately follows from \cref{Thm: main result_R1}.
\end{proof}
With the help of Lemma~\ref{lem:SDP}, one can do a binary search on $(0,1)$ interval to compute the best rate $\rho^*$ we can justify using the \mg{matrix inequality} %LMI system
in \eqref{eq: general SAPD LMI_R1}, i.e., $\rho^*\triangleq\min_{\rho\geq 0}\{\rho: \cP_\rho\neq\emptyset\}$. For any $\rho\in(0,1)$, checking whether $\cP_\rho$ is nonempty or not requires solving \rev{a} 4-dimensional SDP.}

\sa{\mg{Next}, we numerically illustrate that} \mg{the explicit upper bound we derived in \cref{proposition robustness bound nonconvex} provides a reasonable approximation to the actual robustness measure $\cJ$}. \mg{For this purpose, we consider the same %quadratic
example from Section~\ref{sec:explicit-robustness} where $\cJ$ can be explicitly computed, and compare $\cJ$ to its upper bound \rev{given in~\eqref{eq:robustness-bound}}.}
%for quadratics considered in Section~\ref{section: Explicit Estimates for Robustness to Noise
%mimics the same behavior of robustness measure $\cJ$ as a function of convergence rate.
%For a numerical illustration of \sa{the trade-off between $\cJ$ and the convergence rate},
%we consider the SP problem \eqref{eq:main-problem}
In %the quadratic problem
\eqref{eq:special-Phi}, we take $\mu_x=\mu_y=1$,
where \mg{we generate the symmetric matrix $K\in\mathbb{R}^{d\times d}$ randomly with $\|K\|_2 = 10$ and $d=30$}. We assume the noise model given in \eqref{assump-additive-noise} and \eqref{assump-gaussian-noise} with $\delta_x = \delta_y = 10$. \mg{Consequently}, we have
 $ L_{xx} = L_{yy} = 0$, $L_{yx} = 10$.
 %for Assumptions \ref{ASPT: strongly convex concave}-\ref{ASPT: f and g convex}.
 Employing the particular parameter choice in \cref{Corollary: explicit solution to noisy LMI-R1}, \sa{i.e., setting $\beta=1$, $\overline{\theta}_2=0$ and \rev{$c=1$}, we can certify that SAPD converges with rate $\rho \approx 0.9049$ using $\tau$ and $\sigma$ as in \eqref{Condition: SP solution to noisy LMI-R1} and $\theta=\rho$. We %realized
 have found out that $\rho^*$ obtained using the binary search for this example was also equal to $0.9049$, i.e., our special solution in Corollary~\ref{Corollary: explicit solution to noisy LMI-R1} leads to the optimal rate \rev{bound} $\rho^*$.}%\looseness=-1
%\xuan{In addition, we provide a method to find an optimal convergence rate in \cref{search theta}. We note that the optimal convergence rate obtained by \cref{search theta} is also approximately equal to 0.9049. In other word, the convergence rate $\rho_*$ obtained by our special solution is optimal for this example.}
%compute an explicit value of the convergence rate $\rho$,
%$\rho_{p}\triangleq \rho=\theta\approx 0.9382$.
%\nsa{How do you compute this? If you set $\beta=1$ and set $\overline{\theta}_2=0$, we get $0.9049$} \mtodo{Instead of $\rho_p$ should we use $\rho_*$ if it is the best rate in our grid?}

% \todo{MG: $\{r_k\}$ can be confused with $\theta^*_R$. Maybe call $r_k$ or something else instead?(XUAN:DONE)}
\sa{For any $\rho\in[\rho^*,1)$, we can optimize SAPD parameters minimizing the %upper 
bound for robustness %provided 
in \eqref{eq:robustness-bound} while ensuring that the bias term converges linearly with %the given
rate \rev{not worse than} $\rho$,
%through solving
i.e.,}\vspace*{-2mm}
%the best robustness level while asking for the convergence with $r_k$ is obtained by solving the following optimization problem,
%\todo{Important10.15.20201 :Mention that this optimization dont need to know $\delta_x$ and $\delta_y$}
%\xtodo{Where should we say that?} \mtodo{Added a line.}
{\small
\begin{equation}\label{robustness problem}
\begin{aligned}
     \sa{\mathcal{R}(\rho)} \triangleq & \min_{\tau,\sigma,\theta,\alpha\geq 0}\left\{ \frac{2\rho}{1-\rho}\cdot\max\Big\{\tau,\frac{\sigma}{1-\alpha\sigma}\Big\}\cdot{B_{\tau,\sigma,\theta}}
     %\Xi_{\tau,\sigma,\theta}\sa{/\delta^2}
     ~:~\eqref{eq: general SAPD LMI_R1} \mg{\text{ holds for } (\tau,\sigma,\theta,\rho,\alpha)}\right\}.
\end{aligned}
\end{equation}}%
%\mtodo{The problem above is minimized over both $\tau$ and $\theta$ but also over $\alpha$?}
%Note that for solving this optimization problem, 
\rev{
To be able to solve~\eqref{robustness problem}, we do not need to know $\delta_x$ or $\delta_y$.
%It is a
\eqref{robustness problem} is
non-convex; however, it has some structure. 
%namely \xuan{the objective function is an strictly increasing function of $\tau$.}  \xuan{Furthermore, a smaller $\tau$ will yield a wider ranger of other variables. Therefore, the optimal solution to \cref{robustness problem} must require $\tau$ attains its minimal value. }
% Then, if we add the following extra \mg{constraint} to \cref{robustness problem}
% \todo{Xuan: constrain w.r.t $\tau$}
% % \todo{MG: Here, we should emphasize that $c_\tau,c_\sigma$ is fixed, but not variables on $(0,1)$; otherwise if they are variables then the feasibility set would be open and the minimum may not be attained.(XUAN:DONE)}
% \begin{equation}\label{EQ: extra constrains}
% \begin{aligned}
%           & \tau = \frac{1-\xuan{\rho}}{\mu_x\xuan{\rho}},
% \end{aligned}
% \end{equation}
%  \xuan{we would not change the optimal value of the problem \cref{robustness problem}.}
% %Moreover, adding the extra constrains \cref{EQ: extra constrains} to the problem \cref{robustness problem} leads to
% % $$
% % \frac{1}{C_{\tau,\sigma}} = \max\{\frac{2(1-\xuan{\rho})}{\mu_x},\frac{2(1-\theta)}{\mu_y - (\frac{1}{\theta} -1)\alpha}\},
% % $$
% % and we define $C_{\mu} = \min\{\frac{2}{\mu_x},\frac{2}{\mu_y - (\tfrac{1}{\theta}-1)\alpha}\}$.
%As a \mg{consequence},
%\sa{Therefore, we can alternatively solve a simpler optimization problem.}
In the next lemma, we provide a simpler optimization problem exploiting this structure.}\looseness=-1
%as stated in the following lemma
\begin{lemma}\label{lemma: equivalence between robustness problem}
\sa{Given $\rho\in(0,1)$ and $\tau>0$, %we define a set
let $S_\rho(\tau) \triangleq \{ (\sigma,\theta,\alpha): \cref{eq: general SAPD LMI_R1} \text{ holds for } (\tau,\sigma,\theta,\rho,\alpha)\}$. Suppose $\cup_{\tau>0}S_\rho(\tau)\neq\emptyset$. Then, $\frac{1-\rho}{\mu_x\rho}=\min\{\tau:\ S_\rho(\tau)\neq\emptyset\}$. Moreover, for any $\tau_1 \geq \tau _2 \geq \frac{1-\rho}{\mu_x\rho}$, %such that $S_\rho(\tau_2)\neq\emptyset$,
we have $S_{\rho}(\tau_1) \subset S_\rho(\tau_2)$.} Finally, %the optimal value of
\sa{$\cR(\rho)$ defined in \cref{robustness problem} %is equal to that of the problem \cref{relaxed robustness problem 1}  as follows
can also be computed as}\vspace*{-2mm}
{\small
\begin{equation}\label{relaxed robustness problem 1}
        \sa{\cR(\rho)=}\min_{\sigma,\theta,\alpha\geq 0}\Big\{ \frac{2\rho}{1-\rho}\cdot\max\Big\{\tau,\frac{\sigma}{1-\alpha\sigma}\Big\}\cdot
        %\Xi_{\tau,\sigma,\theta}\sa{/\delta^2}
        {B_{\tau,\sigma,\theta}}~:~
        %\\
        % & \text{s.t}.\;\eqref{Condition: SAPD simple LMI 1},~\eqref{Condition: SAPD simple LMI 2},~\eqref{EQ: extra constrains}, ~\theta = r_k.
        %&\ \ \ \ \text{s.t.}\quad
            \tau = \frac{1-\rho}{\mu_x\rho},\quad \cref{eq: general SAPD LMI_R1} \mg{\text{ holds}}\Big\}.
\end{equation}}%
\end{lemma}
%\todo{Xuan(8.24): I notice that some equation use tfrac while other use frac. Shall we make it consistent?}
%\nsa{Sometimes, I use tfrac for constant fractions or to save some space.}
\begin{proof}
\sa{Since $\cup_{\tau>0}S_\rho(\tau)\neq\emptyset$, there exists $(\tau,\sigma,\theta,\rho,\alpha)$ satisfying \cref{eq: general SAPD LMI_R1}; thus, $\tfrac{1}{\tau}(1-\tfrac{1}{\rho})+\mu_x\geq 0$, i.e., $\tau\geq \bar{\tau}\triangleq\frac{1-\rho}{\mu_x\rho}$. Say $\tau\geq \bar{\tau}$, then we have $\tfrac{1}{\bar{\tau}} -L_{xx}\geq \tfrac{1}{\tau} - L_{xx}$, which implies that $(\sigma,\theta,\alpha)\in S_\rho(\bar{\tau})$. Therefore, we can conclude that $\frac{1-\rho}{\mu_x\rho}=\min\{\tau:\ S_\rho(\tau)\neq\emptyset\}$ because \cref{eq: general SAPD LMI_R1} requires that $\tau\geq \bar{\tau}$.}

\sa{Suppose $(\sigma,\theta,\alpha)\in S_{\rho}(\tau_1)$ for some $\tau_1>0$.
% First, since $S_\rho(\tau_2)\neq\emptyset$, we have $\tfrac{1}{\tau_2}(1-\tfrac{1}{\rho})+\mu_x\geq 0$. Second $\tau_1 \geq \tau_2$ implies that $\tfrac{1}{\tau_2} -L_{xx}\geq \tfrac{1}{\tau_1} - L_{xx}$.
% \begin{equation}\label{eq: monotionic tau}
%      \tfrac{1}{\tau_2}+\mu_x -\tfrac{1}{\rho}\geq\tfrac{1}{\tau_1}+\mu_x -\tfrac{1}{\rho},\qquad \tfrac{1}{\tau_2} -L_{xx}\geq \tfrac{1}{\tau_1} - L_{xx}.
% \end{equation}
% If we let $M_{\tau_1}$ and $M_{\tau_2}$ be the matrix in equation $\cref{eq: general SAPD LMI_R1}$ with parameters $(\tau_1, \sigma,\theta,\rho,\alpha)$ and $(\tau_2, \sigma,\theta,\rho,\alpha)$, respectively; then it follows from \cref{eq: monotionic tau} that $M_{\tau_2} \succeq M_{\tau_1}\succeq 0$.
The same arguments %above 
also show that for any  $\tau_2\in[\frac{1-\rho}{\mu_x\rho}, \tau_1]$, \mg{we have}
%Therefore,
$(\sigma,\theta,\alpha)\in S_\rho(\tau_2)$; hence, $S_{\rho}(\tau_1)\subset S_\rho(\tau_2)$. Furthermore, the objective %function 
in \cref{robustness problem} is strictly increasing in $\tau$;
% and given $\rho\in(0,1]$, the smallest value of $\tau$ for which \cref{eq: general SAPD LMI_R1} holds is $\frac{1-\rho}{\mu_x\rho}$, i.e., $\frac{1-\rho}{\mu_x\rho}=\min\{\tau:\ S_\rho(\tau)\neq\emptyset\}$.
thus,}
%\nsa{This is not trivial: how do we know that \cref{eq: general SAPD LMI_R1} has a solution for $\tau=\frac{1-\rho}{\mu_x\rho}$}
%we can conclude that 
the optimal values of \cref{robustness problem} 
%is equal to that of 
and \cref{relaxed robustness problem 1} are equal.
\end{proof}

%Although  \cref{relaxed robustness problem 1} is non-convex, we further provide a method to approximately solve it as follows.

\rev{For any %given
$\rho\in(0,1]$, %we have an lower bound of $\tau$, which is
we consider two necessary conditions for \cref{eq: general SAPD LMI_R1}: i) $\tau \geq \frac{1-\rho}{\mu_x\rho}$, ii)
%\sa{Furthermore, another necessary condition for \cref{eq: general SAPD LMI_R1} is %that
{%\everymath={\scriptstyle}
$
\begin{psmallmatrix}
  \mu_y & \big(\tfrac{\theta}{\rho} - 1\big)L_{yx} \\
  \big(\tfrac{\theta}{\rho} - 1\big)L_{yx} & \tfrac{1}{\tau} - L_{xx}
\end{psmallmatrix}\succeq 0
$,} {which further implies %that
}
%We can further compute that
{\small $(\tfrac{\theta}{\rho} - 1)^2L_{yx}^2\leq \mu_y(\frac{1}{\tau} - L_{xx})$}. Thus, for fixed $\rho$, any solution to \cref{eq: general SAPD LMI_R1} satisfies $\theta\leq\bar{\theta}_\rho$, which is defined in~\eqref{eq:theta_bound}. Indeed,  
%When $\theta\leq \rho$, the upper bound of $\theta$ is trivial.
either $\theta\in[0,\rho]$, or when $\theta \geq \rho$, the necessary conditions %above
imply that}\vspace*{-2mm}
{\small
\begin{equation}
\label{eq:theta_bound}
\theta
\leq
\rho\Bigg( 1 + \frac{\sqrt{\mu_y(\frac{1}{\tau} - L_{xx})}}{L_{yx}}\Bigg)
\leq
\rho\Bigg( 1 + \frac{\sqrt{\mu_y(\frac{\mu_x\rho}{1-\rho} - L_{xx})}}{L_{yx}}\Bigg)\sa{\triangleq \overline{\theta}_\rho}.\vspace*{-1mm}
\end{equation}}%
%\xtodo{how to state that the definition in the following lemma are well defined?}
% %%%%%% OLD DEFINITION %%%%%%%%%%%%
% \begin{definition}
% \sa{Given $c \in (0,1)$, define $L_c\triangleq\{(t,s,\theta,\alpha)\in\reals^4_+:\ \alpha=c s\}$ and $\rho^*_c\triangleq\min_{\rho\geq 0}\{\rho:\ \cP_\rho\cap L_c\neq \emptyset\}$, which is well defined as $\cP_\rho\cap L_c\neq\emptyset$ for all $c\in(0,1)$ when $\rho=1$. Clearly, $0<\rho^*\leq \rho^*_c\leq 1$, and $\cP_\rho\cap L_c\neq \emptyset$ for all $\rho \in [\rho_c^*,1]$.}
% \end{definition}
% %%%%%% OLD DEFINITION %%%%%%%%%%%%

\mg{%In the following discussion,
Next, we %will
discuss how an upper bound %to robustness
\sa{on $\cJ$} can be computed efficiently 
%by a bisection search 
through bisection over the rate parameter $\rho$ and a grid search on $\theta$.
%We start with a definition and provide a lemma for obtaining a bound on robustness.
}
% \xtodo{we already use $c$ in our particular solution. Is it ok if we still $c$ here. For me, I think it is acceptable because we don't use particular solution here and both $c$ are ratios.}
\begin{definition}\label{def-C-rho}
\sa{For $\rho \in (0,1)$, %define
let $\cC_\rho\triangleq\{c\in(0,1):\ \cP_\rho\cap L_c\neq \emptyset\}$, where \mg{$\cP_\rho$ is as in \cref{lem:SDP}} and $L_c\triangleq\{(t,s,\theta,\alpha)\in\reals^4_+:\alpha=c s\}$. \fin{The definition implies that} $\cC_\rho\neq\emptyset$ for all
$\rho\in[\rho^*,1)$.}%\mtodo{Recall the def of $\rho_*$ here? Is it the fastest rate that can be certified?}
%\xtodo{I am not sure if I can edit this version,  but one can copy the following to remark 3.5: recall that $\rho^*$ is the smallest $\rho$ admissible to \cref{eq: general SAPD LMI_R1}.}
\end{definition}
%\xtodo{I add the reference to the  footnote.}
\begin{remark}
\sa{For any $\rho\in[\rho^*,1)$, $\cC_\rho\neq\emptyset$ is a convex set, see \cref{lemma: c rho connected}.
 hence, $\cC_\rho \subset [0,1]$ is an interval. Thus, $\bar{c}_\rho\triangleq\sup \cC_\rho$ and $\barbelow{c}_\rho\triangleq\inf \cC_\rho$ can be computed %using 
via bisection. Each bisection iteration is a $3$-dimensional SDP checking the feasibility of $\{(t,s,\theta)\in\reals^3_+: G_\rho(t,s,\theta,c s)\succeq\mathbf{0}\}$ for a given $c\in(0,1)$.}
\end{remark}
% \todo{xuan: Do we need to say, for all $\rho \in [\rho_c^*,1]$, $\cP_\rho\cap L_c\neq \emptyset$?}
\begin{lemma}\label{lemma: upper bound of robustness problem}
\sa{Given $\rho\in [\rho^*,1)$ and %any
$c\in\cC_\rho$, let  $\barbelow{\theta}_{c,\rho}\triangleq\inf \Theta_{c,\rho}$ and $\overline{\theta}_{c,\rho}\triangleq\sup \Theta_{c,\rho}$, where $\Theta_{c,\rho}\triangleq\{\theta:\ \exists (s,\theta)\in S_{c,\rho}\}$ and $S_{c,\rho}\triangleq\{(s,\theta):\ \exists (t,s,\theta,\alpha)\in\cP_\rho \mbox{ s.t. } t=\frac{\mu_x\rho}{1-\rho},~\alpha=c s\}$.
%$L_{c,\rho}\triangleq\{(t,s,\theta,\alpha)\in\reals^4_+:\ t=\frac{\mu_x\rho}{1-\rho},~\alpha=c s\}$,
% $\barbelow{\theta}_{c,\rho}\triangleq\min_{s,\theta}\{\theta:~\cP_\rho\cap L_{c,\rho}\neq\emptyset\}$ and $\overline{\theta}_{c,\rho}\triangleq\max_{s,\theta}\{\theta:~\cP_\rho\cap L_{c,\rho}\neq\emptyset\}$.
%\mtodo{the fact that we are taking max over $s$ is a bit confusing here because it seems that the term inside the paranthesis does not depend on $s$.}
%\xtodo{Both $\cP_{\rho}$ and $L_{c,\rho}$ depend on s. \textbf{M}: I get that, but are we first maximzing over $s$ and then maximize over $\theta$ for the maximized $s$? why do we need max wrt $s$?\\ X: when we find the range of $\theta$, the only two free variables are $\theta$ and $s$. After we know the range of $\theta$, then we find the max of $s$ for each $\theta$. $\min_{s,\theta}$ means the free variable is $\theta$ and $s$. But, it is indeed a little bit confusing. M: Ok then I suggest we write it $\max_\theta \max_s ..$ rather than $\max_{\theta,s}$ because the order of taking max usually matters.}
For fixed $K_\theta\in \mathbb{Z}_+$, let $\{\theta_k\}_{k=1}^{K_\theta}\subset[\barbelow{\theta}_{c,\rho}, \overline{\theta}_{c,\rho}]\subset[0,\bar\theta_\rho]$ be an arbitrary set of grid points such that $\theta_1=\barbelow{\theta}_{c,\rho}$ and $\theta_{K_\theta}=\overline{\theta}_{c,\rho}$.
%\todo{Sep 7th xuan: I change $\bar{\cR}_{c}$ into $\bar{\cR}_{c,m}$ to show the grid size of theta}
Define {$\bar{\cR}_{c}(\rho)\triangleq \min_{k=1,\ldots K_\theta}\bar{\cR}_c(\rho,\theta_k)$}, where \vspace*{-1mm}
%\nsa{I rewrote this lemma, to clarify the ambiguity on previously defined sets.}
{\small
\begin{align}
    \label{eq:theta_problem}
    \bar{\cR}_c(\rho,\theta)\triangleq\min_{\sigma\geq 0}\left\{\max\Big\{\frac{2}{\mu_x},~\frac{2\rho\sigma}{(1-c)(1-\rho)}\Big\}\cdot
    %\Xi_{\tau,\sigma,\theta}/\delta^2
    {B_{\tau,\sigma,\theta}}:~\tau = \frac{1-\rho}{\mu_x\rho},~\alpha=\frac{c}{\sigma},~\cref{eq: general SAPD LMI_R1} \mg{\text{ holds}}\right\}.
\end{align}}}%
\sa{Then, $\bar{\cR}_{c}(\rho)\geq \cR(\rho)$. Furthermore, for any fixed $\rho\in[\rho^*,1)$, $c\in\cC_\rho$ and $\theta\in[\barbelow{\theta}_{c,\rho}, \overline{\theta}_{c,\rho}]$, $\sigma_c(\rho,\theta)\triangleq 1/\max\{s: (s,\theta)\in S_{c,\rho}\}$ is the unique optimal solution to \eqref{eq:theta_problem}.
%$\sigma_c(\rho,\theta)\triangleq 1/\max\{s:~s\in S_{c,\rho,\theta}\}$ is an optimal solution to \eqref{eq:theta_problem}, where %$S_{c,\rho,\theta}\triangleq\{s:~\exists (t,s,\theta,\alpha)\in\cP_\rho\cap L_{c,\rho}\}\subset\reals_+$.
%$S_{c,\rho,\theta}\triangleq\{s:~\exists (s,\theta)\in S_{c,\rho}\}\subset\reals_+$.
}
\end{lemma}
\begin{proof}
\sa{Given $\rho\in[\rho^*,1)$ and $c \in\cC_\rho$, since we fix $\theta$ and $\alpha=c/\sigma$ while
deriving \cref{eq:theta_problem}, we immediately get $\bar{\cR}_{c}(\rho)\geq \cR(\rho)$ due to \cref{lemma: equivalence between robustness problem}. Lastly, after fixing $\rho\in[\rho^*,1)$, $c\in\cC_\rho$ and $\theta\in[\barbelow{\theta}_{c,\rho}, \overline{\theta}_{c,\rho}]$, the objective 
%function 
in \cref{eq:theta_problem} is increasing in $\sigma>0$, and $\sigma\mapsto 1/\sigma=s$ is a bijection between the feasible region of \eqref{eq:theta_problem} and %$S_{c,\rho,\theta}$
$\{s: (s,\theta)\in S_{c,\rho}\}$. Therefore, the unique %optimal
solution $\sigma_c(\rho,\theta)$
can be computed by solving a \mg{one-dimensional} SDP, i.e., $\max\{s: (s,\theta)\in S_{c,\rho}\}$ for fixed $\rho$, $c$ and $\theta$.}
\end{proof}
\sa{Given $K_c,K_\rho\in\integers_+$, let  $P\triangleq\{\rho_k\}_{k=1}^{K_\rho}\subset[\rho^*,1]$ and $C_\rho\triangleq\{c_k\}_{k=1}^{K_c}\subset\cC_\rho$ be the grid points.
%i.e., for any given $c\in\cC$, one only needs to consider the interval of certifiable rates, i.e., $P\cap [\rho_c^*,1]$.
%$m$ large enough, we grid  the interval of all possible certifiable rates, i.e., let $\{\rho_k\}_{k=1}^K$ be the grid points such that $\rho_1=\rho_c^*$ and $\rho_K=1$, where $K\in\integers_+$ is the grid size.
%Indeed, compared to \eqref{relaxed robustness problem 1}, we have an additional constraint $\alpha = \frac{c}{\sigma}$ in \eqref{eq:theta_problem}, that is why we compute $\bar{\cR}_{c}(\rho)$ only for $\rho\in P$ such that $\rho\geq \rho^*_c\geq\rho^*$.
%In \cref{fig:influence of c_alpha}, setting x-axis as $c\in(0,1)$, we plot the  convergence rate $\rho_c^*$ and the robustness upper bound $\bar{\cR}_{c}(\rho_{c}^*)$.
Finally, for $\rho\in P$, we define $\bar{\cR}(\rho)\triangleq \min_{c\in C_\rho} \bar{\cR}_c(\rho)$, where %for $c\in\cC$ such that $\rho\geq \rho^*_c$,
$\bar{\cR}_{c}(\rho)$ can be computed based on Lemma~\ref{lemma: upper bound of robustness problem} for any $c\in C_\rho$.
%; otherwise, for $\rho<\rho^*_c$, we set $\bar{\cR}_{c}(\rho)=\infty$.
Therefore, for any $\rho\in P$, computing $\bar{\cR}(\rho)$ using Lemma~\ref{lemma: upper bound of robustness problem} will yield $(\tau_\rho, \sigma_\rho, \theta_\rho)$ achieving $\bar{\cR}_{c}(\rho)$ for some $c\in C_\rho$ such that $\bar{\cR}(\rho)=\bar{\cR}_c(\rho)$.
%\xtodo{Shall we mention that $\bar{\cR}(\rho)=\cR(\rho)$ as the grid size goes to infinity?}
Thus, for the quadratic model assumed in \cref{sec:explicit-robustness},
we can compute the robustness measure, defined in \eqref{def-robustness-measure}, corresponding to $(\tau_\rho, \sigma_\rho, \theta_\rho)$, which we call $\cJ(\rho)$.
% \sa{Given $c\in(0,1)$ and $m$, for each $\rho = \rho_k\in[\rho_c^*,1]$, computing $\bar{\cR}_{c}(\rho_k)$ based on Lemma~\ref{lemma: upper bound of robustness problem} will yield $(\tau_k, \sigma_k, \theta_k)$. For the quadratic model assumed in \cref{sec:explicit-robustness},
% we can can compute the robustness measure, defined in \eqref{def-robustness-measure}, corresponding to $(\tau_k, \sigma_k, \theta_k)$, which we call $\cJ_{k}$.
Recall that in Section~\ref{sec:explicit-robustness}, we defined $\cJ^*(\rho)\triangleq\min_{\tau,\sigma,\theta\geq0}\{\cJ:\rho_{\text{true}} = \rho\}$.
To numerically illustrate the rate vs robustness trade-off and also to demonstrate that we can control robustness through optimizing $\bar{\cR}$, in~\cref{fig:trade-off bewtween R and rho_star}, we plot robustness measure $\mathcal{J}(\rho)$, corresponding to %parameters found
$(\tau_\rho, \sigma_\rho, \theta_\rho)$ computed by minimizing its upper bound $\bar{\cR}(\rho)$, against the convergence rate values $\rho\in P$ in the x-axis, and compare $\cJ(\rho)$ with $\cJ^*(\rho)$ and $\bar{\cR}(\rho)$, where we set $K_\rho=K_\theta=100$ and $K_c=50$.}

{In~\cref{fig:trade-off bewtween R and rho_star}, we observe that $\cJ(\rho)$ computed for SAPD parameters optimizing $\bar{\cR}(\rho)$ closely tracks $\cJ^*(\rho)$.
%have the same pattern, and Moreover, the true robustness is very closed to the optimized robustness.
Therefore, we infer that %optimizing
\rev{minimizing} the upper bound 
%is practical to 
helps us optimize the robustness for the problem class used in these experiments.}
\section{Extensions} \label{sec:extensions}
We now
show that SAPD admits
the optimal oracle complexity bound for the stochastic %merely convex 
\fin{MCMC case}, i.e., {when} $\mu_x = \mu_y=0$. This result can be viewed as a 
nontrivial extension of the deterministic complexity result in \cite{hamedani2018primal} to the stochastic gradient setting.} 
\begin{remark}
Suppose $\mu_x=\mu_y=0$, and the parameters $\tau, \sigma>0$ and $\theta\in(0,1]$ satisfy \eqref{Condition: SAPD simple LMI system}. The \sa{first} condition \eqref{Condition: SAPD simple LMI system} implies that $\theta=1$.
\end{remark}
\vspace*{-2mm}
\begin{theorem}\label{Prop: Iteration bound for mcmc noise}
Suppose $\mu_x=\mu_y=0$, \cref{ASPT: lipshiz gradient,ASPT: unbiased noise assumption}
%~\ref{ASPT: strongly convex concave},~\ref{ASPT: f and g convex},
%~\ref{ASPT: unbiased noise assumption} and~\ref{ASPT: compact} 
hold. Assume that $\Omega_x\triangleq\sup_{x_1,x_2 \in \dom f } \|x_1 - x_2 \|<\infty$ and $\Omega_y\triangleq\sup_{y_1,y_2 \in \dom g } \|y_1 - y_2 \|<\infty$. For any $\epsilon>0$, suppose %the parameters
$\{\tau,\sigma,\theta\}$ are chosen such that
{\small
\begin{align}
\label{eq:SAPD-parameter-choice-MC}
\tau =
        \min
        \Big\{\frac{1}{L_{yx}+L_{xx}},~
        %\frac{c_{\tau}\epsilon}{3\delta_x^2}
        \sa{ {\frac{\sa{2}}{15}}\cdot\frac{\epsilon}{\delta_x^2}}\Big\},\quad
        \sigma =
        \min\Big\{\frac{1}{\sa{L_{yx}+2L_{yy}}},~\sa{\frac{1}{L_{xy}}},~
        %\frac{2c_{\sigma}\epsilon}{(54+9c_{\sigma})\delta_y^2}
        %\xuan{\frac{2}{57 + 9 \frac{L_{xy}}{L_{yx}}}}
        \sa{\frac{1}{\sa{72}}}\cdot\frac{\epsilon}{\delta_y^2}\Big\},\quad \theta=1.
\end{align}}%
%{
%---------------------------------------------------------------------
%xuan's comment
%---------------------------------------------------------------------
% \xtodo{xuan(8.23): \cref{eq:SAPD-parameter-choice-MC} should be
% $$
% \tau =  \min
%         \bigg\{\frac{1 }{L_{yx}+L_{xx}},~
%         \frac{1}{6}\cdot\frac{\epsilon}{\delta_x^2}\bigg\},\quad
%         \sigma =
%         \min\bigg\{\frac{1}{2(L_{yx}+L_{yy})},~\frac{1}{2L_{xy}},~
%         {\frac{2}{123}}\cdot\frac{\epsilon}{\delta_y^2}\bigg\}
% $$}
% %---------------------------------------------------------------------
%}
\sa{Then \fin{for the gap metric $\cG(\cdot,\cdot)$, defined in~\eqref{eq:gap},} %iteration complexity of SAPD method {to generate a point $(\bar{x},\bar{y})\in\cX\times\cY$ such that
$\cG(\bar{x}_N,\bar{y}_N)\leq \epsilon$ for all $N\geq N_\epsilon$ such that}
{\small
\begin{equation*}%\label{Iteration bound for noisy MCMC}
    \sa{N_\epsilon} = \mathcal{O}\Big( \frac{(L_{yx}+L_{xx}){\Omega_x^2} + \sa{\max\{L_{yx}+L_{yy},~L_{xy}\}}{\Omega_y^2}}{\epsilon}  + \frac{\delta_x^2{\Omega_x^2} + %\xuan{(1+\frac{L_{xy}}{L_{yx}})}
    \delta_y^2{\Omega_y^2}}{\epsilon^2} \Big).
\end{equation*}}%
\end{theorem}
\begin{proof}
Since $\mu_x=\mu_y=0$, the \sa{first} condition in~\eqref{Condition: SAPD simple LMI system} trivially holds for $(\tau,\sigma,\theta)$ as in \eqref{eq:SAPD-parameter-choice-MC}. Furthermore,  \eqref{eq:SAPD-parameter-choice-MC} implies that
 $\frac{1}{\tau} \geq L_{yx}+L_{xx}$ and \sa{$\frac{1}{\sigma}\geq L_{yx}+ 2L_{yy}$;} therefore, $(\tau,\sigma,\theta)$ in \eqref{eq:SAPD-parameter-choice-MC} with $\pi_1=\pi_2=1$ satisfy the conditions in Lemma~\ref{LEMMA: Noise LMI after young's ineq-R1}. \sa{Thus, $(\tau,\sigma,\theta)$ %together 
 with $\alpha = L_{yx}+L_{yy}$} 
 %is a solution to 
 solves \eqref{Condition: SAPD simple LMI system}.
 
 %\todo{``\eqref{eq:SAPD-parameter-choice-MC} is a solution to \eqref{Condition: SAPD simple LMI 1} and \eqref{Condition: SAPD simple LMI 2} when $\mu_x=\mg{\mu_y} = 0,\;\theta=1$."\\ MG:This sounds right, but not entirely clear to me or to the readers I believe}
% \begin{equation}\label{particular solution for noisy LMI when MCMC}
% \begin{aligned}
%         & \tau =
%         \min
%         \bigg\{\frac{1 - c_{\tau}}{L_{yx}+L_{xx}},~
%         \frac{c_{\tau}\epsilon}{3\delta_x^2}\bigg\}
%         ,\quad  \sigma =
%         \min\bigg\{\frac{(1-c_{\sigma})/2}{L_{yx}+L_{yy}},~
%         \frac{2c_{\sigma}\epsilon}{(54+9c_{\sigma})\delta_y^2}\bigg\},
%         % &\alpha = L_{yx}+L_{yy},\; \pi^x = \frac{c_{\tau}}{\tau},\; \pi^y_1 = \frac{c_{\sigma}}{3(1+\theta)\sigma},\;\pi^y_2 = \frac{c_{\sigma}}{3\theta\sigma},\;\pi^y_3 = \frac{c_{\sigma}}{3\sigma},
%         % \\
%         % & c_{\tau} = \frac{1}{1+ \frac{(L_{yx}+L_{yy})\epsilon}{3\delta_x^2}},\; c_{\sigma} = \frac{-(5+\frac{4(L_{yx}+L_{yy})\epsilon}{9\delta_y^2})+\sqrt{(5+\frac{4(L_{yx}+L_{yy})\epsilon}{9\delta_y^2})^2+24}}{2},
% \end{aligned}
% \end{equation}
% %\mg{with the convention that} $\frac{0}{0}\triangleq 1$.
% From \eqref{particular solution for noisy LMI when MCMC},
% by a trivial computation, we know that \eqref{particular solution for noisy LMI when MCMC} is a solution to \eqref{Condition: SAPD simple LMI 1} and \eqref{Condition: SAPD simple LMI 2} when $\mu_x=\mg{\mu_y} = 0,\;\theta=1$.
\rev{The analysis in the proof of~\cref{Thm: main result_R1} until the end of \eqref{eq:U-bound-distance-metric} is valid for our choice of parameters in~\eqref{eq:SAPD-parameter-choice-MC}. To get a bound for the expected gap, we next analyze $\bar{P}\triangleq\sup\{\sum_{k=0}^{N-1}P_k(x,y): (x,y)\in\dom f\times\dom g\}$.} \sa{For some arbitrary $\eta_x> 0$}, define \sa{$\{\tilde{x}_k\}$} sequence as follows:
%\begin{equation}
%\label{eq:xtilde}
        $\tilde{x}_0 \triangleq x_0$, %\quad
        and
        $\tilde{x}_{k+1} \triangleq \argmin_{x'\in \dom f}  - \langle \Delta^{x}_k , x' \rangle + \frac{\sa{\eta_x}}{2}\| x' - \tilde{x}_k\|^2$, for $k\geq 0$, where $\Delta_k^{x}$ is defined as in 
\cref{lemma: intermediate noisy bound}.
%\end{equation}}%
Then, from \cite[Lemma 2.1]{nemirovski2009robust}, for all $x\in\dom f$
%\nsa{These are the same $\eta_x$ and $\eta_y$ we have used before, right?.}
%\todo{Xuan: Yes, those are same ones.(June 24)}
% \begin{equation*}
%     \frac{1}{2}\| x - \tilde{x}_{k+1}\|^2 \leq \frac{1}{2}\| x - \tilde{x}_{k}\|^2 - \tau\langle \Delta^{x}_k , x -  \tilde{x}_{k} \rangle +\frac{\tau^2}{2}\| \Delta^{x}_k\|^2,\;\forall x\in X.
% \end{equation*}
% Thus,
%\sa{\small
%\begin{equation*}
    %   $\langle \Delta^{x}_k , x -  \tilde{x}_{k} \rangle   \leq \frac{\eta_x}{2}\| x - \tilde{x}_{k}\|^2 -\frac{\eta_x}{2}\| x - \tilde{x}_{k+1}\|^2  +\frac{1}{2\eta_x}\| \Delta^{x}_k\|^2$.
%\end{equation*}}%
%Furthermore, we have that
%Thus, using $\tilde{x}_0 = x_{0}$ 
we get
{\footnotesize
\begin{equation}\label{INEQ: bound of x - tilde x }
    \begin{aligned}
          \sum_{k=0}^{N-1} \langle \Delta^{x}_k , x -  \tilde{x}_{k} \rangle   \leq
          \sum_{k=0}^{N-1} \frac{\eta_x}{2}\| x - \tilde{x}_{k}\|^2 -\frac{\eta_x}{2}\| x - \tilde{x}_{k+1}\|^2  +\frac{1}{2\eta_x}\| \Delta^{x}_k\|^2
            %= & \frac{1}{2\tau}\| x - {x}_{0}\|^2 -  \sum_{k=0}^{N-2}\frac{\rho^{-k}}{2\tau}(\| x - \tilde{x}_{k+1}\|^2 - \rho\| x - \tilde{x}_{k+1}\|^2) -\frac{\rho^{-N+1}}{2\tau}\| x - \tilde{x}_{N}\|^2\\
            % = & \sa{\frac{\eta_x}{2}(\| x - {{x}_{0}}\|^2 -\rho^{-N+1}\| x - \tilde{x}_{N}\|^2)+ \sum_{k=0}^{N-1}\frac{\rho^{-k}}{2\eta_x} \| \Delta^{x}_k\|^2 - \sum_{k=1}^{N-1}\frac{\eta_x}{2}\rho^{-k}(1-\rho)\| x - \tilde{x}_{k}\|^2}
            % \\
            %\leq &  \frac{1}{2\tau}\| x - x_{0}\|^2  -\frac{\rho^{-N+1}}{2\tau}\| x - \tilde{x}_{N}\|^2+ \frac{\tau}{2} \sum_{k=0}^{N-1}\rho^{-k} \| \Delta^{x}_k\|^2\\
            \leq {\frac{\eta_x}{2}} \sa{\norm{x-x_0}^2} %\Omega_X
            + {\frac{1}{2\eta_x}}\sum_{k=0}^{N-1}\| \Delta^{x}_k\|^2;
    \end{aligned}
\end{equation}}%
hence, using $\hat{x}_{k+1}$ defined in \cref{lemma: intermediate noisy bound}, we get
{\footnotesize
\begin{equation}\label{eq:noisy-rate-part2}
    \begin{aligned}
            %\MoveEqLeft
            \mathbb{E}\Big[\sup_{x\in\dom f}\Big\{ \sum_{k=0}^{N-1} -\langle \Delta^{x}_k , \hat x_{k+1} -  x \rangle\Big\}\Big]
            % \\
            % =&    \sum_{k=0}^{N-1} - \langle \Delta^{x}_k , {\hat{x}_{k+1}} -  \tilde{x}_k \rangle
            % + \langle \Delta^{x}_k , x -\tilde{x}_k   \rangle
            % \\
            \leq \sum_{k=0}^{N-1} \mathbb{E}\Big[
            \langle \Delta^{x}_k , \tilde{x}_k-{\hat{x}_{k+1}} \rangle
          +\frac{1}{2\sa{\eta_x}}\| \Delta^{x}_k\|^2\Big] + \frac{\sa{\eta_x}}{2} \Omega_x^2.
          %\sa{\norm{x-x_0}^2}, %\Omega_X
    \end{aligned}
\end{equation}}%
%Next, we consider \textbf{part 3}, let $\Delta_k^{y}$ be defined as in \cref{lemma: intermediate noisy bound}.
\rev{Similarly, for arbitrary $\eta_y> 0$, we construct two auxiliary sequences:
%as follows:
let $\tilde{y}_0^+  = \tilde{y}_0^- = y_0$, and we define}
{\small
%\begin{equation*}
          $\tilde{y}^+_{k+1} \triangleq \argmin_{y'\in \dom g}   \langle \Delta^{y}_k , y' \rangle + \frac{\sa{\eta_y}}{2}\| y' - \tilde{y}^+_k\|^2$, and %\quad
          $\tilde{y}^-_{k+1} \triangleq \argmin_{y'\in \dom g}   -\langle \Delta^{y}_k , y' \rangle + \frac{\sa{\eta_y}}{2}\| y' - \tilde{y}^-_k\|^2$,}
%\end{equation*}}%
\sa{for $k\geq 0$. Thus, as in as in~\eqref{INEQ: bound of x - tilde x }, it follows from \cite[Lemma 2.1]{nemirovski2009robust} that for all $y\in \dom g$,}
%   $\langle \Delta^{y}_k ,   \tilde{y}_{k}^+ - y  \rangle   \leq \frac{\eta_y}{2}\| y - \tilde{y}_{k}^+\|^2 -\frac{\eta_y}{2}\| y - \tilde{y}_{k+1}^+\|^2  +\frac{1}{2\eta_y}\| \Delta^{y}_k\|^2$, %\\
%   and
%       $\langle \Delta^{y}_k ,   y-\tilde{y}_{k}^-  \rangle   \leq \frac{\eta_y}{2}\| y - \tilde{y}_{k}^-\|^2 -\frac{\eta_y}{2}\| y - \tilde{y}_{k+1}^-\|^2  +\frac{1}{2\eta_y}\| \Delta^{y}_k\|^2$.
% Therefore, %similar to
% as in~\eqref{INEQ: bound of x - tilde x }, %using \eqref{eq:y-inner-bound}
we get\footnote{\fin{$\delta_x=0$ implies $\Delta^{x}_k=\mathbf{0}$; %\textbf{part 2} is equal to $0$ and
%we can set 
hence, for $\eta_x=0$, %for which
\eqref{eq:noisy-rate-part2} %trivially
becomes $0\leq 0$.} Similarly, when $\delta_y=0$, we can set $\eta_y=0$.}
{\footnotesize
\begin{equation*}
%\label{INEQ:bound of y-tilde y}
\begin{aligned}
            \MoveEqLeft \sum_{k=0}^{N-1}  2\langle \Delta^{y}_k , \tilde{y}_{k}^+ -  y \rangle- \langle \Delta^{y}_{k-1} ,\tilde{y}_{k-1}^- - y  \rangle
            \leq
            \frac{3\eta_y}{2} \sa{\norm{y-y_0}^2} %\Omega_Y
            + \frac{1}{2\eta_y}\sum_{k=0}^{N-1} 2\| \Delta^{y}_k\|^2 + \| \Delta^{y}_{k-1}\|^2;
\end{aligned}
\end{equation*}}%
hence, using $\hat{y}_{k+1}$ and $\hat{\hat{y}}_{k+1}$ defined in \cref{lemma: intermediate noisy bound}, we get
{\footnotesize
\begin{equation*}%\label{eq:noisy-rate-part3}
    \begin{aligned}
    \hspace{-4pt}
        \MoveEqLeft 
        \mathbb{E}\Big[\sup_{y\in\dom g}\Big\{\sum_{k=0}^{N-1} 2 \langle \Delta^{y}_k,
       {\hat{y}_{k+1}} 
       %- \tilde{y}_k^+ + \tilde{y}_k^+ 
       - y \rangle
      - \langle \Delta^{y}_{k-1},
       {\hat{\hat{y}}_{k+1}}  
       %- \tilde{y}_{k-1}^- + \tilde{y}_{k-1}^- 
       - y \rangle\Big\}\Big]
      \\
      \leq &
        \sum_{k=0}^{N-1} \mathbb{E}\Big[2 \langle \Delta^{y}_k,
       {\hat{y}_{k+1}} - \tilde{y}_k^+  \rangle
      %\\ &
      - \langle \Delta^{y}_{k-1},
       {\hat{\hat{y}}_{k+1}} - \tilde{y}_{k-1}^-  \rangle 
      + \frac{1}{2{\eta_y}} \Big(2\| \Delta^{y}_k\|^2 +  \| \Delta^{y}_{k-1}\|^2 \Big)\Big] +\frac{\sa{3\eta_y}}{2} %\sa{\norm{y-y_0}^2}. 
      \Omega_y^2.
    %   \\
    %   \leq
    %   & \sa{\sum_{k=0}^{N-1} \theta^{-k}\Big[(1 +\theta)\Big( \frac{1}{2\pi^y_1}\| \Delta^{y}_k\|^2 +
    %   \frac{\pi^y_1}{2}\|y_{k+1} - y_k\|^2 +    \langle \Delta^{y}_k, y_k - \tilde{y}_k^+  \rangle \Big)}
    %   \\
    %   & \sa{+ \frac{\theta}{2} \Big(\big(\frac{1}{\pi^y_2}+\frac{1}{\pi^y_3}\big)\| \Delta^{y}_{k-1}\|^2
    %   +\pi^y_2\|y_{k+1} - y_k\|^2
    %   %+\frac{1}{2\pi^y_3}\| \Delta^{y}_{k-1}\|^2
    %   + \pi^y_3\| y_k - y_{k-1}\|^2-2\langle \Delta^{y}_{k-1}, y_{k-1} - \tilde{y}_{k-1}^-  \rangle\Big)\Big]}
    %   \\
    % %   & -\theta \langle \Delta^{y}_{k-1}, y_{k-1} - \tilde{y}_{k-1}^-  \rangle \Big]
    % %   \\
    %   & \sa{+ \frac{1}{2\eta_y}\sum_{k=0}^{N-1} \theta^{-k}\Big((1+\theta)\| \Delta^{y}_k\|^2 +  \theta\| \Delta^{y}_{k-1}\|^2\Big)  +\frac{\eta_y}{2}(1+2\theta)\Omega_Y.}
    \end{aligned}
\end{equation*}}%
Thus, combining this bound with \eqref{eq:noisy-rate-part2} we get $\mathbb{E}[\bar{P}]\leq N(\frac{1}{2}\frac{\delta_x^2}{\eta_x}+\frac{3}{2}\frac{\delta_y^2}{\eta_y})+\frac{\eta_x}{2}\Omega_x^2+\frac{3\eta_y}{2}\Omega_y^2$, where we used $\mathbb{E}[\langle \Delta^{x}_k , \tilde{x}_k-{\hat{x}_{k+1}} \rangle]=\mathbb{E}[\langle \Delta^{y}_k,
       {\hat{y}_{k+1}} - \tilde{y}_k^+  \rangle]=\mathbb{E}[\langle \Delta^{y}_{k-1},
       {\hat{\hat{y}}_{k+1}} - \tilde{y}_{k-1}^-  \rangle]=0$ for $k\geq 0$. Therefore, setting $\eta_x=1/\tau$ and $\eta_y=1/\sigma$ and using the fact that $\theta=1$ implies $K_N(\theta)=N$, it follows from  \eqref{INEQ: ergordic gap  _bounded version}, \eqref{eq:Q-bound} and \eqref{eq:U-bound-distance-metric} that 
{\small
\begin{equation} \label{bound2}
\begin{aligned}
\mathbb{E}[ &\sup_{(x,y)\in \cX\times \cY}\{ \mathcal{L}(\bar{x}_{N}, y)  - \mathcal{L}(x, \bar{y}_{N})\}]
            %+ \frac{1}{N}\Delta_N(x,y)
               \\
 &\leq \frac{1}{N}\Big(\frac{1}{\tau}{\Omega_x^2} + \frac{2}{\sigma}{\Omega_y^2}\Big)  +  {\tau(1+\sa{2}\Xi^x_{\tau,\sigma,\theta}) \frac{\delta_x^2}{2}
     +    \sigma(3+\sa{2}\Xi^y_{\tau,\sigma,\theta}) \frac{\delta_y^2}{2}}
    %   &= \frac{1}{N}\Big(\frac{1}{\tau}{\Omega_f} + \frac{2}{\sigma}{\Omega_g}\Big)  +  {\tau(3 + \sa{2}\sigma L_{yx}) \frac{\delta_x^2}{2}
    %  +    \sigma(27 + 12\sigma L_{yy} + 12\tau\sigma L_{yx}L_{xy}+2\tau L_{yx}) \frac{\delta_y^2}{2}}
    %  \\
      \leq \frac{1}{N}\Big(\frac{1}{\tau}{\Omega_x^2} + \frac{2}{\sigma}{\Omega_y^2}\Big)  +  {\frac{5\tau\delta_x^2}{2}}
     +   \sa{24 \sigma\delta_y^2}, %(19 + 3\frac{L_{xy}}{L_{yx}}) \frac{\sigma\delta_y^2}{2}.
\end{aligned}
\end{equation}}%
%{
%---------------------------------------------------------------------
%xuan's comment
%---------------------------------------------------------------------
% \xtodo{xuan(June. 30): \cref{{bound2}} should be
% \small{$$
% \begin{aligned}
% &\mathbb{E}[ \sup_{(x,y)\in X\times Y}\{ \mathcal{L}(\bar{x}_{N}, y)  - \mathcal{L}(x, \bar{y}_{N})
% + \frac{1}{N}\Delta_N(x,y)\} ]
% \leq
% {\frac{1}{N}}\Omega_{\tau,\sigma,\theta}  + \Xi_{\tau,\sigma,\theta} \\
%  =&  \frac{1}{N}\Big(\frac{1}{\tau}{\Omega_f} + \frac{2}{\sigma}{\Omega_g}\Big)  +  {\tau(1+2\Xi^x_{\tau,\sigma,\theta}) \frac{\delta_x^2}{2}
%      +    \sigma(3+2\Xi^y_{\tau,\sigma,\theta}) \frac{\delta_y^2}{2}}
% \\
% =& \frac{1}{N}\Big(\frac{1}{\tau}{\Omega_f} + \frac{2}{\sigma}{\Omega_g}\Big)  +  {\tau(3 + 2\sigma L_{yx}) \frac{\delta_x^2}{2}
% +    \sigma(27 + 12\sigma L_{yy} + 12\tau\sigma L_{yx}L_{xy}+2\tau L_{yx}) \frac{\delta_y^2}{2}}
% \\
% \leq& \frac{1}{N}\Big(\frac{1}{\tau}{\Omega_f} + \frac{2}{\sigma}{\Omega_g}\Big)  +
% 2\tau\delta_x^2
% +   {\frac{41}{2} \sigma\delta_y^2}
% \end{aligned}
% $$}}
%---------------------------------------------------------------------
%}
% \todo{Xuan(8.24): Does $24$ come from $\tfrac{47}{2}\leq 24$ in equation 2.26?}
where for the last inequality we first substitute $\Xi^x_{\tau,\sigma,\theta}$ and $\Xi^y_{\tau,\sigma,\theta}$ defined in \cref{Thm: main result_R1}, and then use $\tau L_{yx} \leq 1$, $\sigma \max\{L_{yx}, L_{xy}\} \leq \sa{1}$ and $\sigma L_{yy} \leq \frac{1}{2}$ due to  \cref{eq:SAPD-parameter-choice-MC}. For any $\epsilon>0$, %if we desire
requiring
%{
%---------------------------------------------------------------------
%xuan's comment
%---------------------------------------------------------------------
% \todo{xuan(8.23): We need to add $\sigma L_{xy}\leq \frac{1}{2}$.}
%---------------------------------------------------------------------
%}
% \mg{it suffices to}
%\xuan{a sufficient condition is that}%have %one way is to separate them into three parts:
% \sa{impose \eqref{n23} together with}
{\small
\begin{equation}
\label{n1 2}
    \frac{1}{N}\Big(\frac{1}{\tau}{\Omega_x^2} + \frac{2}{\sigma}{\Omega_y^2}\Big) {\leq} %<
    \frac{\epsilon}{3},\quad
   \frac{5\tau\delta_x^2}{\sa{2}} \leq \frac{\epsilon}{3},\quad
   %(19 + 3\frac{L_{xy}}{L_{yx}}) \frac{\sigma\delta_y^2}{2}
   \sa{24\sigma\delta_y^2}\leq \frac{\epsilon}{3},
\end{equation}}%
implies that \eqref{bound2} can be bounded 
%from above 
by $\epsilon$.
%{
%---------------------------------------------------------------------
%xuan's comment
%---------------------------------------------------------------------
% \todo{xuan(June. 30): \cref{n1 2} should be $$
%  2\tau\delta_x^2 \leq \frac{\epsilon}{3},\quad
%   %(19 + 3\frac{L_{xy}}{L_{yx}}) \frac{\sigma\delta_y^2}{2}
%   \frac{41}{2}\sigma\delta_y^2\leq \frac{\epsilon}{3}.
% $$I dont know why the front is red. This is not edited by Dr. Aybat.}
%---------------------------------------------------------------------
%}
\sa{Our parameter choice in \eqref{eq:SAPD-parameter-choice-MC} implies that the second and the third inequalities in \eqref{n1 2} trivially hold.}
It suffices to choose $N$ large enough depending on given $\epsilon$ so that \eqref{n1 2} holds, i.e., $N\geq \frac{3}{\epsilon}(\frac{1}{\tau}{\Omega_x^2} + \frac{2}{\sigma}{\Omega_y^2})$.
%Note that
From \eqref{eq:SAPD-parameter-choice-MC} we have
$\frac{1}{\tau} \leq L_{yx}+L_{xx} + {\frac{15}{\sa{2}}\frac{\delta_x^2}{\epsilon}}$ and
$\frac{1}{\sigma}\leq \sa{\max\{L_{yx}+2L_{yy},~L_{xy}\}} + %\frac{57+9\frac{L_{xy}}{L_{yx}}}{2}
\sa{72}\frac{\delta_y^2}{\epsilon}$;
% \todo{Xuan(8.24):Should it be $
% \frac{1}{\sigma}\leq {\max\{L_{yx}+2L_{yy},~L_{xy}\}} + %\frac{57+9\frac{L_{xy}}{L_{yx}}}{2}
% {72}\frac{\delta_y^2}{\epsilon}$
% }
\sa{thus, $\cG(\bar{x}_N,\bar{y}_N)\leq \epsilon$ holds for all $N\geq N_\epsilon$.}
\end{proof}
\subsection{Robustness Measure for the MCMC Setting}\label{subsec:robustness-MCMC}
\rev{In MCMC setting, based on the %expected 
gap result in \cref{Prop: Iteration bound for mcmc noise},  one can adopt $\tilde{J}\triangleq \limsup_{N\to\infty}\mathbb{E}[\sup\{\cL(\bar{x}_N,y)-\cL(x,\bar{y}_N):(x,y)\in\cX\times\cY\}]/\delta^2$ as the corresponding robustness metric --this definition %of robustness 
	would be %in the same spirit with
	parallel to the definition in~\cite{aybat2020robust}, where the authors consider first-order stochastic algorithms for smooth strongly convex minimization $f^*=\min_x f(x)$ and defined the robustness as $\limsup_{N\to\infty}\mathbb{E}[f(x_N)-f^*]/\delta^2$.}

\rev{Alternatively, one can extend the ideas of SCSC setting to MCMC setting in the following way based on {Tikhonov} regularization. Assume $\mu_x = \mu_y = 0$ and consider the MCMC {saddle point problem} $\min_x \max_y \mathcal{L}(x,y)$. Given a regularization parameter $\mu>0$, let $\mathcal{J}_\mu$ be the robustness of the {SAPD iterate sequence generated when SAPD is implemented on the following} regularized problem: %i.e. the robustness of the regularized problem 
		\begin{align}
			\label{eq:Tikhonov}
			\min_{x\in\cX} \max_{y\in\cY} \mathcal{L}_\mu(x,y)\triangleq \mathcal{L}(x,y)+\frac{\mu}{2}\|x\|^2 - \frac{\mu}{2}\|y\|^2,    
		\end{align}
		where \fin{$\cL_\mu$ is $\mu$-convex in $x$ and $\mu$-concave in $y$}. Using the results in~\cite{ferris1991finite}, under some technical conditions on $\cL$, e.g., $\Phi$ is smooth convex-concave and $f,g$ are indicator functions of some polyhedra, one can show that there exists $\bar\mu>0$ and $({x}^*,{y}^*)\in\cX\times\cY$ such that $({x}^*,{y}^*)$ is the unique saddle point of $\cL_\mu(\cdot,\cdot)$ for all $\mu\in(0,\bar\mu]$; moreover, $({x}^*,{y}^*)=\argmin\{\norm{x}^2+\norm{y}^2:\ (x,y)\in \cZ^*\}$ where $\cZ^*\subset\cX\times\cY$ denotes the set of saddle points of the original MCMC problem $\min_x\max_y\cL(x,y)$. Therefore, rather than directly solving the MCMC problem with SAPD %as explained in \cref{sec:extensions} 
		using the parameters as stated in~\cref{Prop: Iteration bound for mcmc noise} and use the alternative robustness measure $\tilde{J}$ based on the expected gap defined above, one can instead solve the regularized SCSC problem in~\eqref{eq:Tikhonov} for $\mu>0$ sufficiently small, which would generate a least-norm solution of the original MCMC problem, and directly use the originally defined robustness metric $\mathcal{J}_\mu={\limsup_{N\to\infty}\mathbb{E}[\norm{z_N-z^*}^2]/\delta^2}$ corresponding to the SAPD iterate sequence generated while solving the SCSC problem in~\eqref{eq:Tikhonov}, {where $z_N=(x_N,y_N)$, $z^*=(x^*,y^*)$ and $\delta=\max\{\delta_x,\delta_y\}$.}}
\section{Numerical Experiments}
\label{sec:numerics}
\sa{In this section, we compared SAPD against %other methods
S-OGDA~\cite{fallah2020optimal}, SMD~\cite{nemirovski2009robust} and SMP~\cite{juditsky2011solving} for solving \eqref{eq:main-problem} with synthetic and real-data.\footnote{The code is made available at \url{https://github.com/XuanZhangg/SAPD}.}}

\vspace{-1mm}
\subsection{\sa{Regularized Bilinear SP Problem with Synthetic Data}}
\label{sec:random-data}
\begin{figure}[t!]
  \centering
  \includegraphics[width=.4\linewidth]{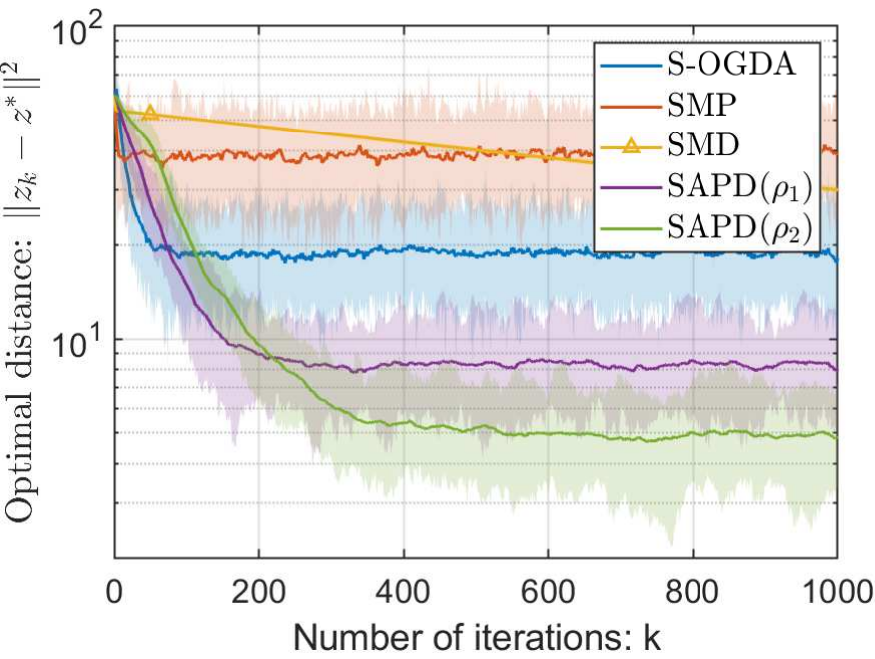}
\includegraphics[width=.4\linewidth]{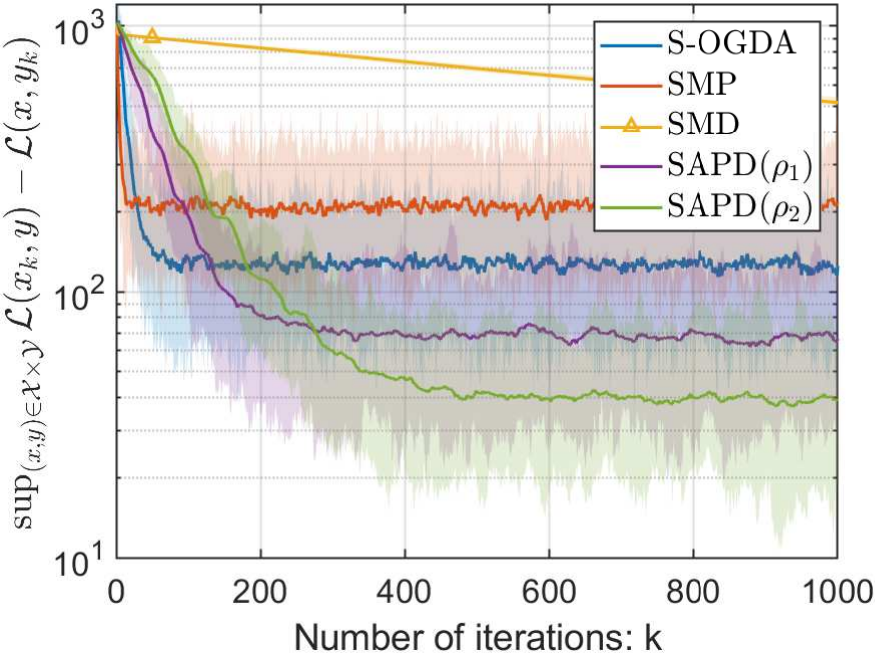}
  \label{fig:trade-off in practice}
\caption{%The rate-robustness trade-off of different algorithms in practice for a quadratic function.
\sa{Comparison of SAPD against SOGDA, SMP and SMD on a synthetic toy problem using the distance (left) and gap (right) metrics.} The rates of SAPD are $\rho_1 = 0.99$ and $\rho_2=0.995$.
%\textbf{shall we use a lower $\delta_x$? Todo cannot work for this figure so I note it by this way.?}
}\vspace*{-5mm}
\end{figure}

{We first tested \sa{SAPD, S-OGDA and SMP} %with other primal-dual algorithms for a quadratic function $F(2,2,10,30,5,5)$
on the regularized bilinear SP problem defined in~\eqref{eq:special-Phi}. In this experiment, we set $\mu_x=\mu_y=1$, $\|K\|_2=10$, $d=30$ and $\delta_x=\delta_y=\delta=5$. \sa{Since SMD step size condition requires a bound on the stochastic gradients, SMD is implemented on \eqref{eq:special-Phi} with additional $(x,y)\in X\times Y$ constraint where $X=\{x\in\cX:\norm{x}\leq \sqrt{d}\}$ and $Y=\{y\in\cY:\norm{y}\leq \sqrt{d}\}$.} Letting \mg{$x$-axis} as the iteration counter, we plot the $50$ sample paths %of different
\sa{for each} algorithm in \cref{fig:trade-off in practice}. \sa{The step sizes for S-OGDA, SMD and SMP are selected as in~\cite{fallah2020optimal}, \cite{nemirovski2009robust} and \cite{juditsky2011solving}, respectively.
Specifically, except for SAPD, all algorithms use primal and dual step sizes that are set equal, %depending on
and their value is a function of $L = \max\{\mu_x, \mu_y,L_{xy},L_{yx}\}$; indeed, S-OGDA uses $\frac{1}{8L}$, SMP uses $\frac{1}{\sqrt{3}L}$, and SMD uses $\frac{2}{\sqrt{5G N}}$, where $N$ denotes the total iteration budget for SMD, and $G>0$ is 
%a fixed constant 
such that
%$\mathbb{E}\left[2\cD_x\|\tilde{\nabla}_x\cL_{\mu_y}\|^2+2\cD_y\|\tilde{\nabla}_y\cL_{\mu_y} \|^2\right]$ where $\cD_x$ and $\cD_y$ are the diameters for the domains of $x$ and $y$.
$\mathbb{E}[2\|\tilde{\nabla}\cL(x,y;\omega^x,\omega^y)\|^2]\leq G$ uniformly for all $(x,y)\in X\times Y$.}
%\nsa{What happened to diameters here?}
%\xtodo{The diameter of $\cX$ is its dimension $d=30$. The diameter of $\cY$ is also its dimension $d=30$. I only use this for SMD.}
%letting $\tilde{\nabla}\cL$ be the estimator of gradients of \eqref{eq:special-Phi} wit $(x,y)$.
The step sizes for SAPD are determined by
minimizing $\bar{\cR}(\rho)$ for $\rho\in\{\rho_1,\rho_2\}$, where $\rho_1=0.99$ and $\rho_2=0.995$. \sa{This process leads to
% $(\tau,\sigma,\theta) = (0.0101, 0.0115,0.6447)$ for $\text{SAPD}(\rho_1)$, and to $(\tau,\sigma,\theta) = (0.0050,0.0076, 0.1739)$ for $\text{SAPD}(\rho_2)$
$(\tau,\sigma,\theta) = (0.010, 0.012,0.645)$ for $\text{SAPD}(\rho_1)$, and to $(\tau,\sigma,\theta) = (0.005,0.008, 0.174)$ for $\text{SAPD}(\rho_2)$.}}
%\xtodo{How should we write the step size of SMD in a nice way? In fact, paper \cite{nemirovski2009robust} use $G_{SMD} = \mathbb{E}\left[2\cD_x\|\tilde{\nabla}_x\cL_{\mu_y}\|^2+2\cD_y\|\tilde{\nabla}_y\cL_{\mu_y} \|^2\right]$, where $\cD_x$ and $\cD_y$ are the diameters of the regions of $x$ and $y$. }
\mg{In Figure \ref{fig:trade-off in practice}, %we observe that
SAPD outperforms the others %algorithms 
in both metrics, i.e., $\cD$ and $\cG$. 
%(distance to the saddle point and the gap function). 
Since $\rho_1 < \rho_2$, %we also observe that 
SAPD with $\rho=\rho_1$ leads to a faster decay of the bias term than that with $\rho=\rho_2$. However, due to \sa{rate and robustness trade-off}, the choice of $\rho=\rho_2$ is more robust to noise, leading to a smaller asymptotic variance of %the iterates
\sa{$\{z_k\}$} as expected.}%\vspace*{-4mm}
\vspace{-2mm}
\subsection{\sa{Distributionally Robust Optimization with Real Data}}
\label{sec:real-data}
\mg{Next, we consider \sa{$\ell_2$-regularized variant of the
%following problem from \cite{namkoong2016stochastic} which arises in
distributionally robust %stochastic
optimization problem from~\cite{namkoong2016stochastic}, i.e.,}
%\nsa{$\rho$ is not a good symbol, we use it for rate. MG: Right.}
%\begin{equation}
(DRO): $\min_{\sa{x\in\cS}} \max_{\sa{y} \in \mathcal{P}_{\sa{r}}}\tfrac{\sa{\mu_x}}{2}\norm{x}^2+\sum_{i=1}^n \sa{y_i} \sa{\phi_i(x)}$ where %f(x)+\Phi(x,y),\quad\hbox{where}\quad \Phi(x,y) =\sum_{i=1}^n \sa{y_i} \sa{\phi_i(x)},
%\label{pbm-dist-robust}
%\end{equation}
%we assume the functions
$\sa{\phi_i}: \mathbb{R}^d \to \mathbb{R}$ is a strongly convex %and twice continuously differentiable, representing the cost
\sa{smooth loss} function corresponding to the $i$-th data point, \sa{$\mu_x>0$ is a regularization parameter}, \sa{$\cS\triangleq \{x\in\reals^d:\ \norm{x}^2\leq \cD_x\}$ for some given model diameter $\cD_x>0$} and $\mathcal{P}_{\sa{r}} \triangleq \{y\in \mathbb{R}^n_+:~\mathbf{1}^\top y = 1,~\|y - \mathbf{1}/n\|^2 \leq \frac{r}{n^2}\}$ -- here, $\mathbf{1}$ denotes the vector with all entries equal to one, and $\mathcal{P}_{r}$ is the uncertainty set around the uniform distribution $\mathbf{1}/n$ whose radius is determined by the parameter $r$. In the special case when $r=0$, the problem recovers the %empirical risk minimization
\sa{ERM} problems arising in supervised learning from labeled data which assigns uniform weights $y_i = 1/n$ to all data points. When $r>0$, the problem is to minimize a worst-case objective %that is
\sa{to be} robust %with respect to
\sa{against} uncertainty in the underlying data distribution. %This min-max formulation
\sa{(DRO)} has several advantages to construct confidence intervals for the parameters of predictive models in supervised learning, see \cite{namkoong2016stochastic}. \sa{This SP problem} %\eqref{pbm-dist-robust}
is affine in the dual variable $y$; therefore, it is not strongly convex with respect to $y$. However, we can approximate it, %this problem,
in a similar spirit to Nesterov's smoothing technique in~\cite{nesterov2005smooth}, with the following SCSC problem:
  \begin{equation}
%   \sa{\min_{\sa{x\in\reals^d}} \max_{y\in \reals^n} \cL_{\mu_y}(x,y)\triangleq \mathbbm{1}_{\cS}(x)+\tfrac{\sa{\mu_x}}{2}\norm{x}^2+\Phi(x,y)-\tfrac{{\mu_y}}{2} \|y\|^2-\mathbbm{1}_{\cP_r}(y)},
  \sa{\min_{x\in\cS} \max_{y\in \cP_r} \cL_{\mu_y}(x,y)\triangleq \tfrac{\sa{\mu_x}}{2}\norm{x}^2+\Phi(x,y)-\tfrac{{\mu_y}}{2} \|y\|^2,}
  \label{pbm-dist-robust-scsc}
  \end{equation}
  where %\sa{$\mathbbm{1}_{\cS}(\cdot)$, $\mathbbm{1}_{\cP_r}(\cdot)$ denote the indicator functions, and}
  $\Phi(x,p) = \sum_{i=1}^n y_i \phi_i(x)$
  %and $g(y) = \frac{{\mu_y}}{2} \|y\|^2$
  %by choosing $\newcomment{\mu_y}$ \newcomment{properly}.
  \rev{for some properly chosen smoothing parameter $\mu_y>0$ --see~\cref{rem:SCMC}}.}
  %This is similar
  %from the literature
  %and we also provide the iteration complexity bound in \cref{Prop: iteration complex for SAPD for SCMC}. %\newcomment{It should be noted that although \cref{pbm-dist-robust-scsc} does not have the strongly convex term $f(x)$ as that in \cref{eq:main-problem}, one can construct it by adding an $L_2$ regularizer.}
%\xtodo{I add the definition of $\cD_y$}
In our tests, we consider the binary logistic loss with an $l_2$ regularizer, i.e., $\phi_i(x)=\ln(1+\exp(-b_i \sa{a_i}^\top x))$ and \sa{set $r=2\sqrt{n}$}. We can then apply SAPD to the SCSC problem in \eqref{pbm-dist-robust-scsc} which admits the Lipschitz constants $L_{xy} = L_{yx} = \|A\|_2$, $L_{xx} = \max_{\sa{i=1,\ldots,n}}\{\tfrac{1}{4}\|a_i\|_2^2$\}, and $L_{yy}=0$, \sa{where $A\in\reals^{n\times d}$ is the data matrix with rows $\{a_i\}_{i=1}^n$ and columns $\{A_j\}_{j=1}^d$.} %Recall the gap function $\cG$ defined in \eqref{eq:gap}.
Since
$\cD_y\triangleq\sup_{y\in\cP_r}\norm{y}=1$,
%\xtodo{need to change statement here since we dont have such result for $\cD$}, 
{for any given $\epsilon>0$, 
we set $\mu_y = \frac{\epsilon}{2\cD_y^2}$ according to~\cref{rem:SCMC}. 
}
%up to a log factor.
%\end{aligned}
%\end{equation}}%
\begin{figure}[t]
\label{figs: LG Regression}
\centering
\subfigure[Dry Bean: $\rho_1 =0.9986$, $\rho_2=0.9997$.]{
\label{fig: Dry Bean}
\includegraphics[width = 0.31\textwidth]{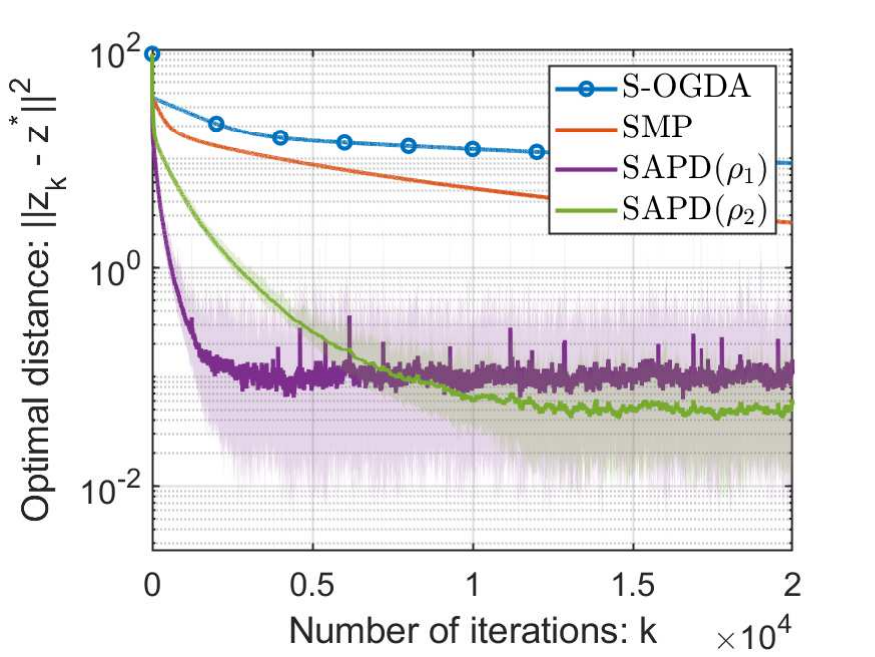}
\includegraphics[width = 0.31\textwidth]{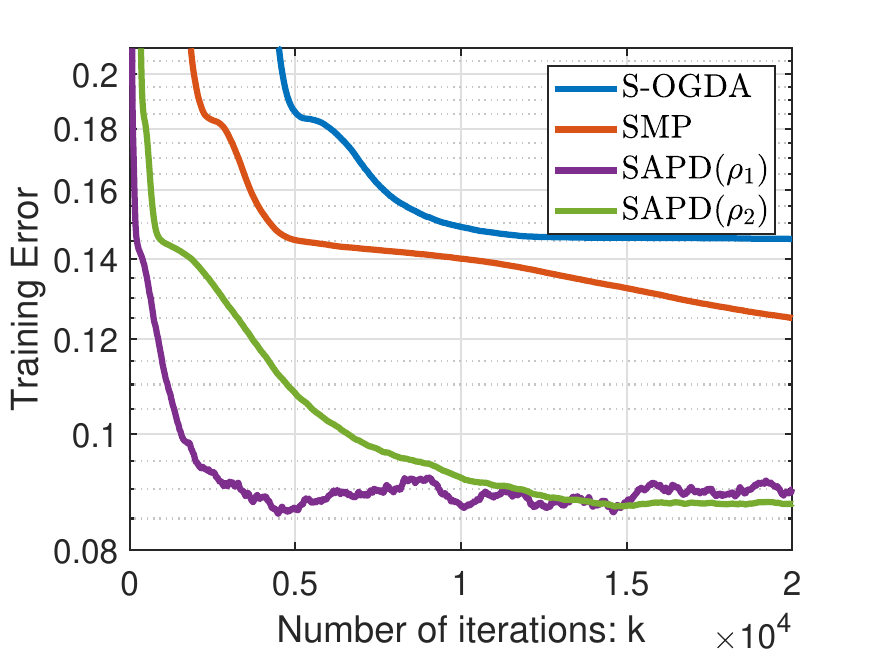}
\includegraphics[width = 0.31\textwidth]{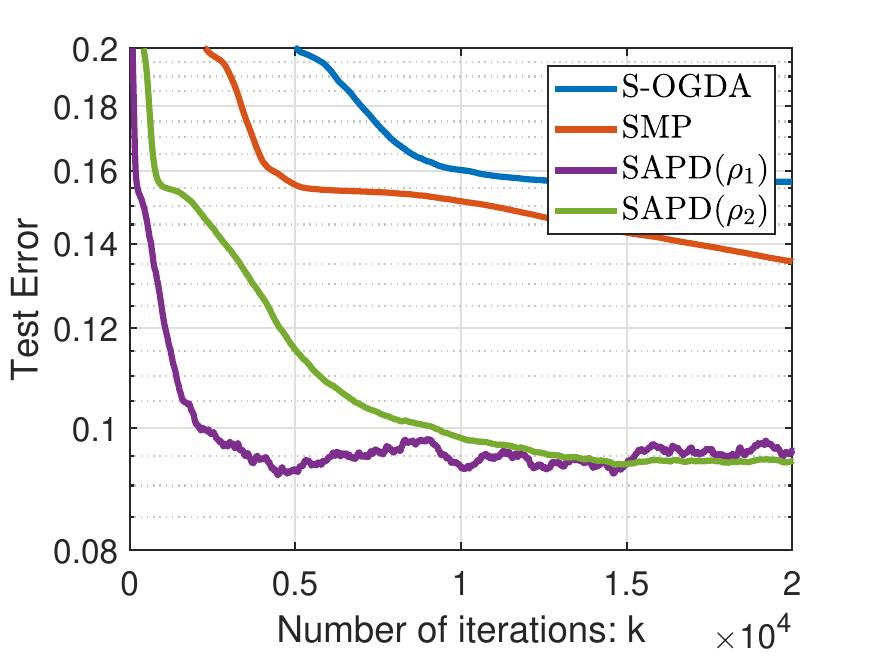}
}
% \vspace{-4mm}
\centering
% \vspace{-4mm}
\subfigure[Arcene: $\rho_1 =0.989$, $\rho_2=0.992$.]{
\label{figs: Arcene}
\vspace{-3mm}
\includegraphics[width = 0.31\textwidth]{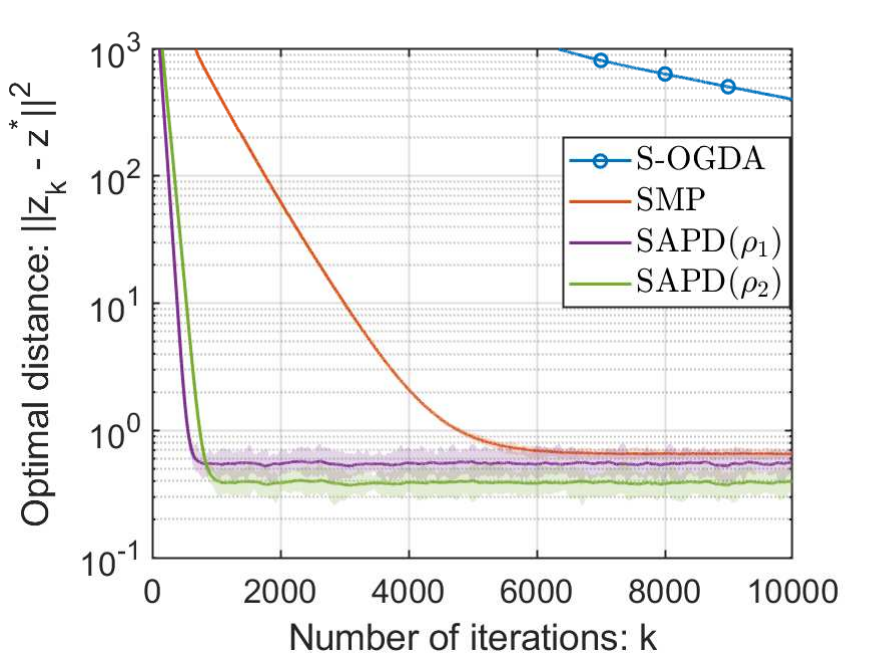}
\includegraphics[width = 0.31\textwidth]{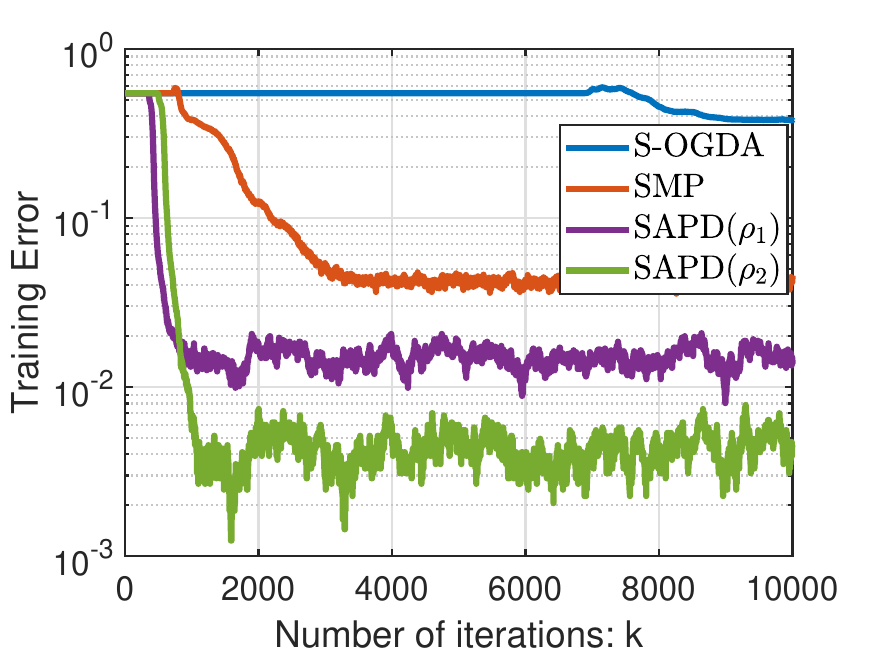}
\includegraphics[width = 0.31\textwidth]{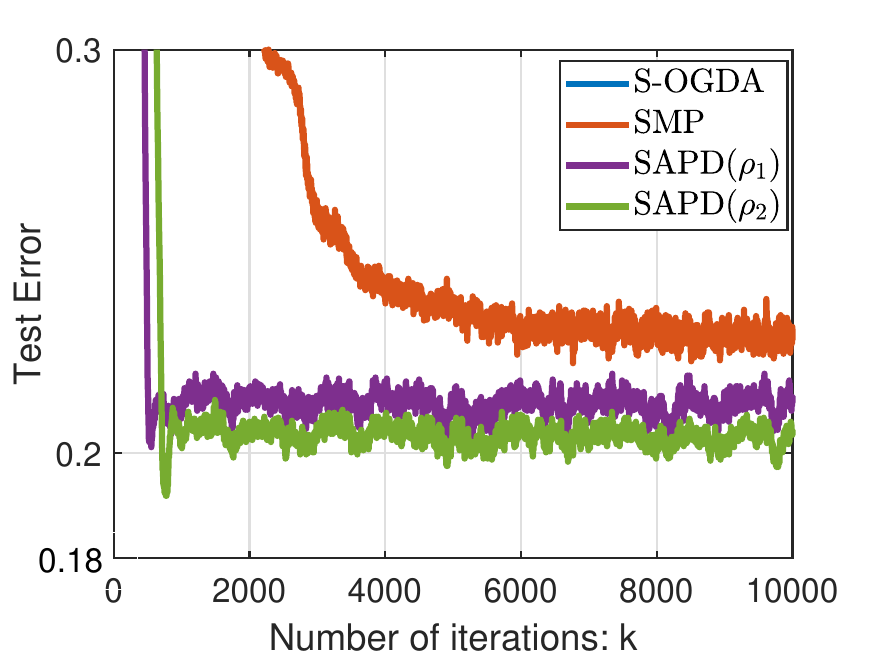}
}
\vspace*{-4mm}
\caption{Comparison of SAPD against S\newcomment{-}OGDA~\cite{fallah2020optimal} and SMP~\cite{juditsky2011solving} on real-data sets.}
\vspace*{-5mm}
\end{figure}

\sa{We perform experiments on two data sets: 1) Dry Bean data set~\cite{koklu2020multiclass} with $d=16$, $n=9528$ with a test data set of $4083$ \mg{points}; 2) Arcene data set~\cite{guyon2004result} with $d=10,000$, $n=97$ and test size of
$96$ \mg{points}. In these experiments we set the regularization parameter $\mu_x$ through cross-validation. The source of noise in gradient computations is mini-batch sampling of data.
% \xtodo{Important note: For the dry bean data set. The influence of noise is slight for SMP and S-OGDA because their step sizes is too conservative. Intuitively, one should expect curves of SMP and S-OGDA go below SAPD. However, it follows from the training error and test error that even after SMP and S-OGDA show a convergence behavior, they almost stop and can not improve more in acceptable future. In other word, our method to optimize robustness is efficient and guarantee a fast convergence rate. fOR other data set, When running SMP, OGDA for the other data set, not only the influence of noise is huge but also the convergence is slow.}
For the Dry Bean data set, we set $\mu_x = 0.01$, $\mu_y = 10$ and use batch size $1$, and normalize each feature column $A_j\in\reals^n$ using $A_j \gets \frac{A_j - \min(A_j)}{\max(A_j) - \min(A_j)}$, where both min and max are taken over the elements of $A_j$.
%\nsa{Xuan, check this one. This normalization does not look right.}\xtodo{I find this normalization method from some machine learning blogs}\mtodo{We can add a reference from the blogs and say that the denominator is non-zero in the normalization}
%\xtodo{It is the first method in wiki pedia ''https://en.wikipedia.org/wiki/Feature_scaling'', shall we refer to it?}\mtodo{We can check with Prof. Aybat, but to me because this is unusual scaling it is better to give a reference to Wikipedia (even though that is not a textbook/perfect reference).}\nsa{Xuan, sometimes you do not read carefully. Please check the webpage you suggest and see if what you wrote here is the same with what the first method says.}\xtodo{It should be $A_j = \frac{A_j - \min(A_j)}{\max(A_j) - \min(A_j)}$}
For the Arcene data set, we set $\mu_x = 0.02$ and $\mu_y = 10$, and use batch size 10, and normalize the data matrix such that $A\gets A/\min\{\sqrt{d},\sqrt{n}\}$.}
%Then employing those parameters, we solve \cref{robustness problem}
\sa{As described in Section~\ref{sec:robustness-bound}, given a desired rate $\rho\geq \rho^*$, we compute $(\tau,\sigma,\theta)$ for SAPD that achieves $\bar{\cR}(\rho)$}. \sa{We plot SAPD statistics for two different rates to illustrate that our framework can trade-off rate and robustness in an effective manner. The other methods we tested set the primal and dual step sizes equal. Indeed, for S-OGDA \cite{hsieh2019convergence,fallah2020optimal} the step size is $\frac{1}{8L}$, and
for SMP~\cite{juditsky2011solving}, it is $\frac{1}{\sqrt{3}L}$, where $L = \max\{L_{xx} + \mu_x, \mu_y,L_{xy},L_{yx}\}$. In~\cref{figs: LG Regression}, we plotted the optimization error using the distance metric $\cD$, training and the test errors. We reported the results for 50 sample paths. \mg{Our results show that for both \sa{test and training errors,} reported in terms of distances to the solution, SAPD achieves a good performance \sa{for both rate and robustness}.
%in terms of convergence rate as well as robustness to noise.
}}%
% $\min\{\frac{1}{\sqrt{3}L},\sqrt{\frac{\cD_x+\cD_y}{7N\max\{(\delta_x^2,\delta_y^2\}}}\}$,
% where $N$ is the number of total iterations;
% for SMD \cite{nemirovski2009robust}, letting $\tilde{\nabla}_x\cL_{\mu_y},\tilde{\nabla}_y\cL_{\mu_y}$ be the estimator of gradients of $\cL_{\mu_x}$ w.r.t. $x,y$, respectively, and $\mathbb{E}\left[2\cD_x\|\tilde{\nabla}_x\cL_{\mu_y}\|^2+2\cD_y\|\tilde{\nabla}_y\cL_{\mu_y} \|^2\right]$ bounded by $G_{SMD}$,the step size is $\frac{2}{\sqrt{5G_{SMD}N}}$.
\rev{
\section{Future Work}
In a follow-up paper~\cite{zhang2022sapd+}, we have considered SARAH variance reduction~\cite{nguyen2017sarah} on weakly convex-strongly concave~(WCSC) problems, \fin{and proposed 
%In~\cite{zhang2022sapd+}, we study 
an inexact proximal point method based on SAPD, which serves as a subroutine for inexactly solving SCSC sub-problems.} We have implemented a variance reduction framework within SAPD, which not only improved the oracle complexity from $\cO(\epsilon^{-4})$ to
$\cO(\epsilon^{-3})$; but, we have also improved the best condition number dependency from $\cO(L\kappa^3/\epsilon^3)$ to $\cO(L\kappa^2/\epsilon^3)$, where $\kappa\triangleq L/\mu_y$ with $L$ being the Lipschitz constant of $\grad \Phi$ and $\mu_y$ \fin{being} the strong concavity constant of $\cL(x,\cdot)$ uniformly for all $x\in\dom f$. \fin{While incorporating SARAH within SAPD helps for WCSC problems,} using the same variance reduction analysis for SAPD on SCSC problems does not help in improving the complexity results we established in \cref{Prop: iteration complex for SAPD_R1}. {That being said, for the SCSC case, in a recent relevant paper, we have applied a bias reduction strategy called Richardson-Romberg extrapolation to SAPD \cite{bugra-richardson} and in our experiments we have observed that this technique has not only created an improved bias performance but also exhibits an improved dependency to gradient noise variance.} As a future work on the SCSC setting with noisy gradients, it would be interesting to design \fin{efficient variance/bias reduction techniques} for SAPD. 
%to improve the $\cO(1/\epsilon)$ worst-case complexity due to variance-term while still achieving the linear rate for the bias term. 
The method in~\cite{zhang2022sapd+} has two nested loops, another important future research direction involves establishing convergence guarantees for SAPD as a single-loop method when implemented for solving \fin{WCSC and weakly convex-merely-concave problems.}}

\bibliographystyle{siamplain}
{%\scriptsize
 \bibliography{references} %{arXiv_final}

\begin{thebibliography}{10}

\bibitem{aybat2019universally}
{\sc N.~S. Aybat, A.~Fallah, M.~Gurbuzbalaban, and A.~Ozdaglar}, {\em A
  universally optimal multistage accelerated stochastic gradient method}, in
  Advances in Neural Information Processing Systems, 2019, pp.~8525--8536.

\bibitem{aybat2020robust}
{\sc N.~S. Aybat, A.~Fallah, M.~Gurbuzbalaban, and A.~Ozdaglar}, {\em Robust
  accelerated gradient methods for smooth strongly convex functions}, SIAM
  Journal on Optimization, 30 (2020), pp.~717--751.

\bibitem{beck2017first}
{\sc A.~Beck}, {\em First-order methods in optimization}, Society for
  Industrial and Applied Mathematics, Philadelphia, PA, 2017.

\bibitem{ben2009robust}
{\sc A.~Ben-Tal, L.~El~Ghaoui, and A.~Nemirovski}, {\em Robust optimization},
  vol.~28, Princeton University Press, 2009.

\bibitem{bottou2018optimization}
{\sc L.~Bottou, F.~E. Curtis, and J.~Nocedal}, {\em Optimization methods for
  large-scale machine learning}, Siam Review, 60 (2018), pp.~223--311.

\bibitem{bugra-richardson}
{\sc B.~Can, M.~Gurbuzbalaban, and N.~S. Aybat}, {\em A variance-reduced
  stochastic accelerated primal dual algorithm}, 2022,
  \url{https://doi.org/10.48550/ARXIV.2202.09688},
  \url{https://arxiv.org/abs/2202.09688}.

\bibitem{chambolle2011first}
{\sc A.~Chambolle and T.~Pock}, {\em A first-order primal-dual algorithm for
  convex problems with applications to imaging}, Journal of {M}athematical
  {I}maging and {V}ision, 40 (2011), pp.~120--145.

\bibitem{chambolle2016ergodic}
{\sc A.~Chambolle and T.~Pock}, {\em On the ergodic convergence rates of a
  first-order primal--dual algorithm}, Mathematical Programming, 159 (2016),
  pp.~253--287.

\bibitem{chen2017accelerated}
{\sc Y.~Chen, G.~Lan, and Y.~Ouyang}, {\em Accelerated schemes for a class of
  variational inequalities}, Mathematical Programming, 165 (2017),
  pp.~113--149.

\bibitem{cohen2020relative}
{\sc M.~B. Cohen, A.~Sidford, and K.~Tian}, {\em Relative lipschitzness in
  extragradient methods and a direct recipe for acceleration}, in 12th
  Innovations in Theoretical Computer Science Conference (ITCS 2021), Schloss
  Dagstuhl-Leibniz-Zentrum f{\"u}r Informatik, 2021.

\bibitem{condat2016fast}
{\sc L.~Condat}, {\em Fast projection onto the simplex and the $l_1$ ball},
  Mathematical Programming, 158 (2016), pp.~575--585.

\bibitem{cui2016analysis}
{\sc S.~Cui and U.~V. Shanbhag}, {\em On the analysis of reflected gradient and
  splitting methods for monotone stochastic variational inequality problems},
  in 2016 IEEE 55th Conference on Decision and Control (CDC), IEEE, 2016,
  pp.~4510--4515.

\bibitem{fallah2020optimal}
{\sc A.~Fallah, A.~Ozdaglar, and S.~Pattathil}, {\em An optimal multistage
  stochastic gradient method for minimax problems}, in 2020 59th IEEE
  Conference on Decision and Control (CDC), IEEE, 2020, pp.~3573--3579.

\bibitem{ferris1991finite}
{\sc M.~C. Ferris and O.~L. Mangasarian}, {\em Finite perturbation of convex
  programs}, Applied Mathematics and Optimization, 23 (1991), pp.~263--273.

\bibitem{gidel2018variational}
{\sc G.~Gidel, H.~Berard, G.~Vignoud, P.~Vincent, and S.~Lacoste-Julien}, {\em
  A variational inequality perspective on generative adversarial networks}, in
  International Conference on Learning Representations, 2019,
  \url{https://openreview.net/forum?id=r1laEnA5Ym}.

\bibitem{golub1996matrix}
{\sc G.~H. Golub and C.~F. Van~Loan}, {\em Matrix computations,. johns},
  Hopkins Studies in Mathematical Sciences, 3rd edition edition,  (1996).

\bibitem{guyon2004result}
{\sc I.~Guyon, S.~Gunn, A.~Ben-Hur, and G.~Dror}, {\em Result analysis of the
  {N}ips 2003 feature selection challenge}, Advances in {N}eural {I}nformation
  {P}rocessing {S}ystems, 17 (2004).

\bibitem{hamedani2018primal}
{\sc E.~Y. Hamedani and N.~S. Aybat}, {\em A primal-dual algorithm with line
  search for general convex-concave saddle point problems}, SIAM Journal on
  Optimization, 31 (2021), pp.~1299--1329.

\bibitem{hsieh2019convergence}
{\sc Y.-G. Hsieh, F.~Iutzeler, J.~Malick, and P.~Mertikopoulos}, {\em On the
  convergence of single-call stochastic extra-gradient methods}, in Advances in
  Neural Information Processing Systems, 2019, pp.~6938--6948.

\bibitem{jin2022sharper}
{\sc Y.~Jin, A.~Sidford, and K.~Tian}, {\em Sharper rates for separable minimax
  and finite sum optimization via primal-dual extragradient methods}, in
  Conference on Learning Theory, PMLR, 2022, pp.~4362--4415.

\bibitem{juditsky2011first}
{\sc A.~Juditsky and A.~Nemirovski}, {\em First order methods for nonsmooth
  convex large-scale optimization, i: general purpose methods}, Opt. for
  Machine Learning, 30 (2011), pp.~121--148.

\bibitem{juditsky2011solving}
{\sc A.~Juditsky, A.~Nemirovski, and C.~Tauvel}, {\em Solving variational
  inequalities with stochastic mirror-prox algorithm}, Stochastic Systems, 1
  (2011), pp.~17--58.

\bibitem{koklu2020multiclass}
{\sc M.~Koklu and I.~A. Ozkan}, {\em Multiclass classification of dry beans
  using computer vision and machine learning techniques}, Computers and
  Electronics in Agriculture, 174 (2020).

\bibitem{kotsalis2020simple}
{\sc G.~Kotsalis, G.~Lan, and T.~Li}, {\em Simple and optimal methods for
  stochastic variational inequalities, i: Operator extrapolation}, SIAM Journal
  on Optimization, 32 (2022), pp.~2041--2073.

\bibitem{kuru2020differentially}
{\sc N.~Kuru, S.~Ilker~Birbil, M.~Gurbuzbalaban, and S.~Yildirim}, {\em
  Differentially private accelerated optimization algorithms}, SIAM Journal on
  Optimization, 32 (2022), pp.~795--821.

\bibitem{lin-near-optimal}
{\sc T.~Lin, C.~Jin, and M.~I. Jordan}, {\em Near-optimal algorithms for
  minimax optimization}, in Conference on Learning Theory, PMLR, 2020,
  pp.~2738--2779.

\bibitem{privacy-gan}
{\sc Y.~{Liu}, J.~{Peng}, J.~J.~Q. {Yu}, and Y.~{Wu}}, {\em Ppgan:
  Privacy-preserving generative adversarial network}, in IEEE Conference on
  Parallel and Distributed Systems, 2019, pp.~985--989.

\bibitem{mokhtari2020unified}
{\sc A.~Mokhtari, A.~Ozdaglar, and S.~Pattathil}, {\em A unified analysis of
  extra-gradient and optimistic gradient methods for saddle point problems:
  Proximal point approach}, in International Conference on Artificial
  Intelligence and Statistics, 2020, pp.~1497--1507.

\bibitem{namkoong2016stochastic}
{\sc H.~Namkoong and J.~C. Duchi}, {\em Stochastic gradient methods for
  distributionally robust optimization with f-divergences}, in Advances in
  Neural Information Processing Systems, 2016, pp.~2208--2216.

\bibitem{nedic2009subgradient}
{\sc A.~Nedi{\'c} and A.~Ozdaglar}, {\em Subgradient methods for saddle-point
  problems}, Journal of Optimization Theory and Applications, 142 (2009),
  pp.~205--228.

\bibitem{nemirovski2004prox}
{\sc A.~Nemirovski}, {\em Prox-method with rate of convergence o (1/t) for
  variational inequalities with lipschitz continuous monotone operators and
  smooth convex-concave saddle point problems}, SIAM Journal on Optimization,
  15 (2004), pp.~229--251.

\bibitem{nemirovski2009robust}
{\sc A.~Nemirovski, A.~Juditsky, G.~Lan, and A.~Shapiro}, {\em Robust
  stochastic approximation approach to stoc. programming}, SIAM {J}ournal on
  {O}ptimization, 19 (2009), pp.~1574--1609.

\bibitem{nesterov2005excessive}
{\sc Y.~Nesterov}, {\em Excessive gap technique in nonsmooth convex
  minimization}, SIAM Journal on Optimization, 16 (2005), pp.~235--249.

\bibitem{nesterov2005smooth}
{\sc Y.~Nesterov}, {\em Smooth minimization of non-smooth functions},
  Mathematical {P}rogramming, 103 (2005), pp.~127--152.

\bibitem{nesterov2009primal}
{\sc Y.~Nesterov}, {\em Primal-dual subgradient methods for convex problems},
  Mathematical {P}rogramming, 120 (2009), pp.~221--259.

\bibitem{nguyen2017sarah}
{\sc L.~M. Nguyen, J.~Liu, K.~Scheinberg, and M.~Tak{\'a}{\v{c}}}, {\em Sarah:
  A novel method for machine learning problems using stochastic recursive
  gradient}, in International Conference on Machine Learning, PMLR, 2017,
  pp.~2613--2621.

\bibitem{ouyang2021lower}
{\sc Y.~Ouyang and Y.~Xu}, {\em Lower complexity bounds of first-order methods
  for convex-concave bilinear saddle-point problems}, Mathematical Programming,
  185 (2021), pp.~1--35.

\bibitem{palaniappan2016stochastic}
{\sc B.~Palaniappan and F.~Bach}, {\em Stochastic variance reduction methods
  for saddle-point problems}, in Advances in Neural Information Processing
  Systems, 2016, pp.~1416--1424.

\bibitem{strang2005linear}
{\sc G.~Strang}, {\em Linear algebra and its applications.: Thomson {B}rooks},
  Cole, Belmont, CA, USA,  (2005).

\bibitem{NEURIPS2018_08048a9c}
{\sc C.~Tan, T.~Zhang, S.~Ma, and J.~Liu}, {\em Stochastic primal-dual method
  for empirical risk minimization with $\mathcal{O}(1)$ per-iteration
  complexity}, in Advances in Neural Information Processing Systems, S.~Bengio,
  H.~Wallach, H.~Larochelle, K.~Grauman, N.~Cesa-Bianchi, and R.~Garnett, eds.,
  vol.~31, Curran Associates, Inc., 2018.

\bibitem{thekumparampil2022lifted}
{\sc K.~K. Thekumparampil, N.~He, and S.~Oh}, {\em Lifted primal-dual method
  for bilinearly coupled smooth minimax optimization}, in International
  Conference on Artificial Intelligence and Statistics, PMLR, 2022,
  pp.~4281--4308.

\bibitem{wang2017exploiting}
{\sc J.~Wang and L.~Xiao}, {\em Exploiting strong convexity from data with
  primal-dual first-order algorithms}, in International Conference on Machine
  Learning, PMLR, 2017, pp.~3694--3702.

\bibitem{wang-li}
{\sc Y.~Wang and J.~Li}, {\em Improved algorithms for convex-concave minimax
  optimization}, Advances in Neural Information Processing Systems, 33 (2020),
  pp.~4800--4810.

\bibitem{pmlr-v32-wen14}
{\sc J.~Wen, C.-N. Yu, and R.~Greiner}, {\em Robust learning under uncertain
  test distributions: Relating covariate shift to model misspecification}, in
  Proceedings of the 31st International Conference on Machine Learning, E.~P.
  Xing and T.~Jebara, eds., vol.~32 of Proceedings of Machine Learning
  Research, Bejing, China, 22--24 Jun 2014, PMLR, pp.~631--639.

\bibitem{xie2018differentially}
{\sc L.~Xie, K.~Lin, S.~Wang, F.~Wang, and J.~Zhou}, {\em Differentially
  private generative adversarial network}, arXiv preprint arXiv:1802.06739,
  (2018).

\bibitem{xu2005maximum}
{\sc L.~Xu, J.~Neufeld, B.~Larson, and D.~Schuurmans}, {\em Maximum margin
  clustering}, in Advances in {N}eural {I}nformation {P}rocessing systems,
  2005, pp.~1537--1544.

\bibitem{yan2020optimal}
{\sc Y.~Yan, Y.~Xu, Q.~Lin, W.~Liu, and T.~Yang}, {\em Optimal epoch stochastic
  gradient descent ascent methods for min-max optimization}, Advances in Neural
  Information Processing Systems, 33 (2020), pp.~5789--5800.

\bibitem{yang2020catalyst}
{\sc J.~Yang, S.~Zhang, N.~Kiyavash, and N.~He}, {\em A catalyst framework for
  minimax optimization}, Advances in Neural Information Processing Systems, 33
  (2020), pp.~5667--5678.

\bibitem{zhang2019lower}
{\sc J.~Zhang, M.~Hong, and S.~Zhang}, {\em On lower iteration complexity
  bounds for the convex concave saddle point problems}, Mathematical
  Programming, 194 (2022), pp.~901--935.

\bibitem{zhang2020single}
{\sc J.~Zhang, P.~Xiao, R.~Sun, and Z.~Luo}, {\em A single-loop smoothed
  gradient descent-ascent algorithm for nonconvex-concave min-max problems},
  Advances in Neural Information Processing Systems, 33 (2020), pp.~7377--7389.

\bibitem{zhang2022sapd+}
{\sc X.~Zhang, N.~S. Aybat, and M.~Gurbuzbalaban}, {\em {SAPD+}: An accelerated
  stochastic method for nonconvex-concave minimax problems}, in Advances in
  Neural Information Processing Systems, S.~Koyejo, S.~Mohamed, A.~Agarwal,
  D.~Belgrave, K.~Cho, and A.~Oh, eds., vol.~35, Curran Associates, Inc., 2022,
  pp.~21668--21681.

\bibitem{zhang2017stochastic}
{\sc Y.~Zhang and L.~Xiao}, {\em Stochastic primal-dual coordinate method for
  regularized empirical risk minimization}, The Journal of Machine Learning
  Research, 18 (2017), pp.~2939--2980.

\bibitem{zhao2019optimal}
{\sc R.~Zhao}, {\em Accelerated stochastic algorithms for convex-concave
  saddle-point problems}, Mathematics of Operations Research, 47 (2022),
  pp.~1443--1473.

\bibitem{zhong2020improving}
{\sc J.~Zhong, X.~Liu, and C.-J. Hsieh}, {\em Improving the speed and quality
  of {GAN} by adversarial training}, arXiv preprint arXiv:2008.03364,  (2020).

\bibitem{zhou1996robust}
{\sc K.~Zhou, J.~C. Doyle, K.~Glover, et~al.}, {\em Robust and optimal
  control}, vol.~40, Prentice hall New Jersey, 1996.

\end{thebibliography}
}%
%\inputencoding{utf8}
\newpage
\appendix
\section{Proofs of Lemmas}
\subsection{Auxiliary Results}
\begin{lemma}\label{lemma: proximal ineq}
	Let $f$ be proper, closed and strongly convex with %convexity
	\sa{modulus $\mu>0$}. Then for any $x,x'\in \dom f$, and %constant
	$c>0$,
	{\small
		$
		\|\prox{cf}(x) - \prox{cf}(x') \| \leq \frac{1}{1+c\mu}\|  x- x' \|.
		$}
\end{lemma}
Because \cref{lemma: proximal ineq} %lemma
is a simple extension of  \cite[Theorem 6.42]{beck2017first} to the strongly convex scenario, we 
%do not %show the
%\sa{provide its} 
omit its proof. 
\subsection{Proof of \cref{lem: basic lemma for deterministic case}}\sa{Fix $x\in\dom f$ and $y\in\dom g$.} %arbitrary.
\rev{Invoking \cite[Lemma 7.1]{hamedani2018primal} for the $y-$ and $x-$subproblems in Algorithm~\ref{ALG: SAPD}, and using the definitions of $\varepsilon^x_k$ and $\varepsilon^y_k$, we get} %two inequalities:
%that hold for any  $x \in \mathcal{X}$ and $y\in \mathcal{Y}$:
% \rev{\small
% \begin{equation*}
%     \begin{aligned}
%     f(x_{k+1}) + \langle & \tilde{\nabla}_x \Phi(x_k, y_{k+1};\omega^x_k) , x_{k+1} - x \rangle\\
%      &\leq 
%       f(x) + \frac{1}{2\tau} \left[  \|x- x_k\|^2 -\|x- x_{k+1}\|^2 - \|x_{k+1}- x_k\|^2
%       \right]  - \frac{\mu_x}{2}\|x - x_{k+1}\|^2,\\
%     %\end{aligned}
% %\end{equation*}
% %\begin{equation*}
%     %\begin{aligned}&
%     g(y_{k+1}) - \langle &\tilde{s}_k, %- \theta \omega_{k-1}^y,
%       y_{k+1} - y \rangle \\
%       & \leq 
%       g(y) + \frac{1}{2\sigma} \left[  \|y-y_k\|^2- \|y -  y_{k+1}\|^2 -
%       \|y_{k+1}- y_k\|^2
%       \right] - \frac{\mu_y}{2}\|y - y_{k+1}\|^2.
%     \end{aligned}
%     \vspace*{-2mm}
% \end{equation*}}%
% {Using} the %following
% two equalities:
% {\footnotesize
% \begin{eqnarray*}
% \lefteqn{\langle \tilde{\nabla}_x \Phi(x_k, y_{k+1};\omega_k^x), x_{k+1} - x \rangle}\\
% &&\hbox{} =\langle \tilde{\nabla}_x \Phi(x_k, y_{k+1};\omega_k^x) - \nabla_x \Phi(x_k, y_{k+1}) , x_{k+1} - x \rangle + \langle \nabla_x \Phi(x_k, y_{k+1}), x_{k+1} - x \rangle,
% \end{eqnarray*}}%
% and
% $\langle \tilde{s}_k,
%       y_{k+1} - y \rangle
%         =
%       \langle s_k,
%       y_{k+1} - y \rangle  +   \langle \tilde{s}_k -s_k,
%       y_{k+1} - y \rangle$,
%\xtodo{The referees complain the inner line equations. Should we move it out?}
%Thus, we get %two inequalities that hold for any  $x \in \cX$ and $y\in \cY$:
\rev{\small
	\begin{subequations}
		\begin{equation} \label{IX}
			\begin{aligned}
				f(x_{k+1}) &+ \langle {{\nabla}_x \Phi(x_k, y_{k+1}%;\omega_k^x
					)},~x_{k+1} - x \rangle\\
				\leq &  f(x) + \tfrac{1}{2\tau} (  \|x- x_k\|^2 - \|x- x_{k+1}\|^2 - \|x_{k+1}- x_k\|^2
				){- \tfrac{\mu_x}{2}\|x - x_{k+1}\|^2}+{\varepsilon^x_k},
			\end{aligned}
		\end{equation}
		\vspace{-1mm}
		\begin{equation}\label{IY2}
			\begin{aligned}
				-g(y)&+ g(y_{k+1})
				\\
				\leq &\langle s_k,  y_{k+1} - y \rangle + \tfrac{1}{2\sigma} \left[  \|y-y_k\|^2- \|y -  y_{k+1}\|^2 -
				\|y_{k+1}- y_k\|^2
				\right] - \tfrac{\mu_y}{2}\|y - y_{k+1}\|^2{+\varepsilon^y_k}.
			\end{aligned}
		\end{equation}%
		% \begin{equation}\label{IY1}
		%     \begin{aligned}
		%     g(y_{k+1})  &- \langle {s}_k,
		%       y_{k+1} - y \rangle
		%       \\
		%       \leq & g(y) + \tfrac{1}{2\sigma} (  \|y- y_k\|^2 - \|y- y_{k+1}\|^2 -
		%       \|y_{k+1}- y_k\|^2
		%       )\sa{ - \tfrac{\mu_y}{2}\|y - y_{k+1}\|^2}\sa{+\varepsilon^y_k}.
		%     \end{aligned}
		% \end{equation}
\end{subequations}}%
% For notation convenience, we temporarily ignore the noisy parts of \cref{IX,IY1}, i.e.
% \begin{equation*}
%       - \langle  \tilde{\nabla}_x \Phi(x_k, y_{k+1};\omega_k^x) - \nabla_x \Phi(x_k, y_{k+1}),~x_{k+1} - x \rangle
%       \text{ and }
%       \langle \tilde{s}_k -s_k,
%       y_{k+1} - y \rangle.
% \end{equation*}
%\xtodo{ suggestion for shorten paper: can we say by lipshitchz cts and strongly convex, we can get A.7 frin A.5 directly.}
%Rearranging \sa{the terms in \eqref{IY1}, we get}\vspace*{-2mm}
% {\footnotesize
% \begin{equation}\label{IY2}
% \begin{aligned}
%     -g(y)&+ g(y_{k+1})
%       \\
%       \leq &\langle s_k,  y_{k+1} - y \rangle + \tfrac{1}{2\sigma} \left[  \|y-y_k\|^2- \|y -  y_{k+1}\|^2 -
%       \|y_{k+1}- y_k\|^2
%       \right] - \tfrac{\mu_y}{2}\|y - y_{k+1}\|^2\sa{+\varepsilon^y_k}.
% \end{aligned}
% \end{equation}}%
\sa{Since $y_{k+1}\in \dom g$}, the inner product in (\ref{IX}) can be lower bounded using convexity of \sa{$\Phi(\cdot, y_{k+1})$} in \cref{ASPT: lipshiz gradient} as follows:
{\small
	\begin{align*}
		\langle \nabla_x  \Phi(x_k, y_{k+1}), x_{k+1} - x \rangle
		= &
		\langle \nabla_x \Phi(x_k, y_{k+1}), x_k - x \rangle
		+
		\langle \nabla_x \Phi(x_k, y_{k+1}), x_{k+1} - x_k \rangle
		\\
		\geq &
		\Phi(x_k, y_{k+1}) - \Phi(x, y_{k+1})  %_{\mathcal{X}}
		+ \langle \nabla_x \Phi(x_k, y_{k+1}),~x_{k+1} - x_k \rangle.
\end{align*}}%
Using this inequality after adding $\Phi(x_{k+1}, y_{k+1})$ to both sides of \eqref{IX}, we get
{\small
	\begin{equation}
		\begin{aligned}
			\Phi( x_{k+1}, &y_{k+1}) + f(x_{k+1}) \\
			\leq &  \Phi(x, y_{k+1}) + f(x)
			+ \Phi(x_{k+1}, y_{k+1}) - \Phi(x_k, y_{k+1}) - \langle \nabla_x \Phi(x_k, y_{k+1}), x_{k+1} - x_k \rangle
			\\
			&
			+  \tfrac{1}{2\tau} \left[ \| x - x_k \| ^ 2 - \| x - x_{ k + 1 }  \| ^ 2 -
			\|x_{k+1}- x_k\|^2\right] - \tfrac{\mu_x}{2}\|x - x_{k+1}\|^2+\sa{\varepsilon^x_k}
			\\
			\leq & \Phi(x, y_{k+1}) +  f(x) + \frac{L_{xx}}{2} \| x_{k+1} - x_k \|^2
			\\
			&    +  \tfrac{1}{2\tau} \left[ \| x - x_k \| ^ 2 - \| x - x_{ k + 1 }  \| ^ 2 -
			\|x_{k+1}- x_k\|^2\right]- \tfrac{\mu_x}{2}\|x - x_{k+1}\|^2+\sa{\varepsilon^x_k},
		\end{aligned}
\end{equation}}%
where \rev{the last step follows from \cref{ASPT: lipshiz gradient}, i.e.,
	%\eqref{LGX2}.
	$\grad_x\Phi(\cdot,y_{k+1})$ is Lipschitz with constant $L_{xx}$.} Rearranging \sa{the} terms gives us
{\small
	\begin{equation}\label{IX2}
		\begin{aligned}
			f(&x_{k+1}) - f(x) - \Phi(x,y_{k+1})
			\leq -\Phi(x_{k+1},y_{k+1}) + \frac{L_{xx}}{2} \| x_{k+1} - x_k \|^2
			\\
			&    +  \frac{1}{\sa{2}\tau} \left[ \| x - x_k \| ^ 2 - \| x - x_{ k + 1 }  \| ^ 2 -
			\|x_{k+1}- x_k\|^2\right]- \frac{\mu_x}{2}\|x - x_{k+1}\|^2+\sa{\varepsilon^x_k}.
		\end{aligned}
\end{equation}}%
Then, for $k \geq 0$, by summing \eqref{IY2} and \eqref{IX2}, we obtain
{\small
	\begin{equation}
		\label{eq:one-step-aux1}
		\begin{aligned}
			\mathcal{L}( x_{k+1}, y)  - \mathcal{L}(x, y_{k+1})
			&=
			f(x_{k+1}) + \Phi(x_{k+1},y) - g(y) - f(x) - \Phi(x,y_{k+1}) + g(y_{k+1})
			\\
			&\leq \Phi(x_{k+1}, y) - \Phi(x_{k+1}, y_{k+1}) + \langle s_k, y_{k+1} - y \rangle + \frac{L_{xx}}{2} \| x_{k+1} - x_k \|^2
			\\
			&\ +  \frac{1}{\sa{2}\sigma} \left[  \|y-y_k\|^2- \|y -  y_{k+1}\|^2 -
			\|y_{k+1}- y_k\|^2 \right]- \frac{\mu_y}{2}\|y - y_{k+1}\|^2+ {\varepsilon_k^y}
			%   \right]
			\\
			& \ +  \frac{1}{\sa{2}\tau} \left[ \| x - x_k \| ^ 2 - \| x - x_{ k + 1 }  \| ^ 2 -
			\|x_{k+1}- x_k\|^2\right]
			- \frac{\mu_x}{2}\|x - x_{k+1}\|^2 + {\varepsilon_k^x}. %_{\mathcal{X}}.
		\end{aligned}
\end{equation}}%
%For $ k \geq 0$, recall that we set $ q_k = \nabla_y \Phi(x_k, y_k) - \nabla_y \Phi(x_{k-1}, y_{k-1})$; thus, $ s_k= \nabla_y \Phi(x_k, y_k) + \theta q_k$.
%\xtodo{suggestion for shorten paper: we dont need those many steps to get this relation}
\sa{From \cref{ASPT: lipshiz gradient},} %by the
the concavity of $\Phi(x, \cdot)$ for \sa{fixed $x \in \dom f\subset \mathcal{X}$} \sa{implies}
{\small
	\begin{equation*}
		\begin{aligned}
			\Phi( x_{k+1}, y) -  &\Phi(x_{k+1}, y_{k+1}) + \langle s_k, y_{k+1} - y \rangle
			\\
			\leq & \langle \nabla_y\Phi(x_{k+1},y_{k+1}), y - y_{k+1} \rangle  + \langle \nabla_y \Phi(x_k, y_k) + \theta q_k, y_{k+1} - y \rangle
			\\
			% = & \langle \nabla_y\Phi(x_{k+1},y_{k+1}) - \nabla_y\Phi(x_{k},y_{k}) , y - y_{k+1} \rangle  + \theta\langle  q_k, y_{k+1} - y \rangle
			%\\
			= & {-\langle q_{k+1}, y_{k+1} - y \rangle + \theta \langle q_k, y_{k} - y \rangle +\theta \langle q_k, y_{k+1}-y_k\rangle.}
		\end{aligned}
\end{equation*}}%
%Recall that we temporarily ignore the noisy parts of \cref{IX,IY1}, we add them back and
\sa{Thus,} using the above inequality within~\eqref{eq:one-step-aux1}, we get
{\small
	\begin{equation*}
		\begin{aligned}
			\mathcal{L}( & x_{k+1}, y)  - \mathcal{L}(x, y_{k+1}) \leq
			-\langle q_{k+1}, y_{k+1} - y \rangle + \theta \langle q_k, y_{k} - y+y_{k+1}-y_k \rangle
			%+\theta \langle q_k, y_{k+1}-y_k\rangle
			+ \frac{L_{xx}}{2} \| x_{k+1} - x_k \|^2
			\\
			& +  \frac{1}{\sa{2}\sigma} \left[  \|y-y_k\|^2- \|y -  y_{k+1}\|^2 -
			\|y_{k+1}- y_k\|^2
			\right]- \frac{\mu_y}{2}\|y - y_{k+1}\|^2
			\\
			& +  \frac{1}{\sa{2}\tau} \left[ \| x - x_k \| ^ 2 - \| x - x_{ k + 1 }  \| ^ 2 -
			\|x_{k+1}- x_k\|^2\right]
			- \frac{\mu_x}{2}\|x - x_{k+1}\|^2+\sa{\varepsilon^x_k}+\sa{\varepsilon^x_y}.
			%       \\
			%   &  - \langle  \tilde{\nabla}_x \Phi(x_k, y_{k+1};\omega_k^x) - \nabla_x \Phi(x_k, y_{k+1}),~x_{k+1} - x \rangle  + \langle \tilde{s}_k -s_k,  y_{k+1} - y \rangle.
		\end{aligned}
\end{equation*}}%
\sa{Finally, \eqref{D1} follows from using Cauchy-Schwarz %inequality
	for $\fprod{q_k,y_{k+1}-y_k}$ and \eqref{INEQ: Cauchy Ineqaulity 1}.}
\subsection{Proof of \cref{lemma: intermediate noisy bound}}
The first inequality in~\cref{eq:y_k - y_hat_k} is %a direct result of
from \cref{lemma: proximal ineq}; for the second, we have
{\footnotesize
	\begin{equation*}
		\begin{aligned}
			\|y_{k+1} - \hat{y}_{k+1}\| \leq \frac{\sigma}{1+\sigma\mu_y}\|\tilde{s}_k - s_k \| \leq \frac{\sigma}{1+\sigma\mu_y}\left((1+\theta)\|\Delta^{y}_k\| + \theta\|\Delta^{y}_{k-1}\|\right),
		\end{aligned}
\end{equation*}}%
%where the first inequality is
\sa{which follows from \cref{lemma: proximal ineq} and the %second inequality is from
	triangle inequality. %Before showing
	To show \cref{eq:y_k - y_hat_hat_k}, %first prove the following relations.
	we bound $\| y_{k+1} - \hat{y}_{k+1}\|$ and $\| \hat{y}_{k+1}- \hat{\hat{y}}_{k+1}\|$ separately.} It follows from \cref{lemma: proximal ineq} that
{\small
	%\begin{equation*}
	%\begin{aligned}
	$$\|x_{k+1} - \hat{\hat{x}}_{k+1} \| \leq \frac{\tau}{1+\tau\mu_x}\|\tilde{\nabla}_x\Phi(x_k,y_{k+1};\omega_k^x) - \nabla_x\Phi(x_k,\hat{y}_{k+1})\|.$$}%
%\end{aligned}
%\end{equation*}}%
\sa{After adding and subtracting $\nabla_x\Phi(x_k,{y}_{k+1})$, \cref{ASPT: lipshiz gradient} implies that}
% \begin{equation*}
%     \begin{aligned}
%     \|x_{k+1} - \hat{\hat{x}}_{k+1} \| \leq \frac{\tau}{1+\tau\mu_x}\|\Delta_k^{x} + \nabla_x\Phi(x_k,{y}_{k+1}) -  \nabla_x\Phi(x_k,\hat{y}_{k+1})\|.
%     \end{aligned}
% \end{equation*}
% Furthermore, if we use the triangle inequality and \cref{ASPT: lipshiz gradient}, it follows that:
{\small
	\begin{equation}\label{ineq: x_k - x_hat_hat_k}
		\begin{aligned}
			\|x_{k+1} - \hat{\hat{x}}_{k+1} \| \leq \frac{\tau}{1+\tau\mu_x}\left(\|\Delta_k^{x}\| + L_{xy}\|y_{k+1} - \hat{y}_{k+1}\|\right).
		\end{aligned}
\end{equation}}%
% Then, by \cref{eq:y_k - y_hat_k}, we have that
% \begin{equation}
%     \begin{aligned}
%     \|x_{k+1} - \hat{\hat{x}}_{k+1} \| \leq \frac{\tau}{1+\tau\mu_x}\left(\|\Delta_k^{x}\| +\frac{\sigma L_{xy}}{1+\sigma\mu_y}\left((1+\theta)\|\Delta^{y}_k\| + \theta\|\Delta^{y}_{k-1}\|\right)\right).
%     \end{aligned}
% \end{equation}
%We will use above relation later.
\sa{We will use this relation to bound $\|\hat{y}_{k+1} - \hat{\hat{y}}_{k+1}\|$.
	%Next, we separate \cref{eq:y_k - y_hat_hat_k} into two parts. First,
	Indeed,} using \cref{lemma: proximal ineq}, we have %that
{\footnotesize
	\begin{align*}
		%\label{eq: y_hat_k - y_hat_hat_k}
		\|\hat{y}_{k+1} &- \hat{\hat{y}}_{k+1}\| \leq  \frac{1}{1+\sa{\sigma}\mu_y}\left\|y_k - \hat{y}_k + \sigma(1+\theta) \left( \nabla_y\Phi(x_k,y_k) - \nabla_y\Phi(\hat{\hat{x}}_k,\hat{y}_k) \right)\right\|\\
		%\end{align*}
		%Furthermore, if we employ \cref{ASPT: lipshiz gradient} and the triangle inequality, it follows that
		%\begin{align*}
		%\|\hat{y}_{k+1} - \hat{\hat{y}}_{k+1}\|
		\leq &\frac{1}{1+\sigma\mu_y}\left((1 + \sigma(1+\theta)L_{yy}) \|y_k - \hat{y}_k\| + \sigma(1+\theta)L_{yx}\|x_k - \hat{\hat{x}}_k \|\right)\nonumber\\
		%\end{align*}
		%If we use \cref{ineq: x_k - x_hat_hat_k}, it follows that
		%\begin{align*}
		%\|\hat{y}_{k+1} - \hat{\hat{y}}_{k+1}\|
		\leq & \frac{1}{1+\sigma\mu_y}\left(\Big(1 + \sigma(1+\theta)L_{yy} + \frac{\tau\sigma(1+\theta)L_{yx}L_{xy}}{1+\tau\mu_x}\Big) \|y_k - \hat{y}_k\| %\right.\nonumber\\ &\left.
		+ \frac{\tau\sigma(1+\theta)L_{yx}}{1+\tau\mu_x}\| \Delta_{k-1}^{x}\|\right)\nonumber\\
		%\end{align*}
		%Moreover, it follows from \cref{eq:y_k - y_hat_k}
		%that
		%\begin{equation}
		%\begin{aligned}
		%\|\hat{y}_{k+1} - \hat{\hat{y}}_{k+1}\|
		\leq & \frac{\sigma}{1+\sigma\mu_y}\left(\Big(1 + \sigma(1+\theta)L_{yy}  +  \frac{\tau\sigma(1+\theta)L_{yx}L_{xy}}{1+\tau\mu_x} \Big) \cdot 
		%\right. \cdot \nonumber\\
		%& \left.
		\frac{(1+\theta)\|\Delta^{y}_{k-1}\| + \theta\|\Delta^{y}_{k-2}\|}{1+\sigma\mu_y}
		+ \frac{\tau(1+\theta)L_{yx}}{1+\tau\mu_x}\|\Delta^{x}_{k-1}\|\right), \nonumber
\end{align*}}%
%\end{equation}
\sa{where the second, third and fourth inequalities follow from \cref{ASPT: lipshiz gradient}, \cref{ineq: x_k - x_hat_hat_k} and the second inequality in~\cref{eq:y_k - y_hat_k}, respectively.  Combining this %inequality
	with
	%\mg{The triangle inequality}
	%Next,  using triangle inequality, we have that
	%\begin{align*}
	$\| y_{k+1} - \hat{\hat{y}}_{k+1}\|
	%= \| y_{k+1} - \hat{y}_{k+1} +\hat{y}_{k+1}- \hat{\hat{y}}_{k+1} \|
	\leq \| y_{k+1} - \hat{y}_{k+1}\| + \| \hat{y}_{k+1}- \hat{\hat{y}}_{k+1}\|$,
	%\end{align*}
	%Finally, applying
	and the second one in~\cref{eq:y_k - y_hat_k} %and \cref{eq: y_hat_k - y_hat_hat_k}
	give us the desired bound.} %into the above inequality completes the proof.
\subsection{Proof of \cref{Lemma: final noisy bound}} %First,
With the convention that $y_{-2} = y_{-1} = y_0$, and $x_{-2} = x_{-1} = x_0$, \cref{lemma: intermediate noisy bound} and Cauchy-Schwarz inequality
%, it is easy to obtain that,
imply for all $k\geq 0$ that
%the following inequalities holds:
{\footnotesize
	\begin{equation*}
		\begin{aligned}
			& \langle \Delta^{x}_k, x_{k+1} - \hat{x}_{k+1} \rangle \leq \frac{\tau}{1+\tau\mu_x}\|\Delta^{x}_k \|^2, \\
			& \langle \Delta^{y}_k,  y_{k+1} - \hat{y}_{k+1}\rangle \leq \frac{\sigma}{1+\sigma\mu_y}\left((1+\theta)\|\Delta^{y}_k\|^2 + \theta\|\Delta^{y}_{k-1}\|\|\Delta^{y}_k\|\right),\\
			&  \langle \Delta^{y}_{k-1},  y_{k+1} - \hat{\hat{y}}_{k+1}\rangle
			\leq
			\frac{\sigma}{1+\sigma\mu_y}\Bigg(
			(1+\theta)\|\Delta^{y}_k\|\|\Delta^{y}_{k-1}\|
			+
			\theta\|\Delta^{y}_{k-1}\|^2+
			\frac{\tau(1+\theta)L_{yx}}{1+\tau\mu_x}\|\Delta^{x}_{k-1}\|
			\|\Delta^{y}_{k-1}\|\\
			&
			+
			\Big(\frac{1 + \sigma(1+\theta)L_{yy}}{1+\sigma\mu_y}  +  \frac{\tau\sigma(1+\theta)L_{yx}L_{xy}}{(1+\tau\mu_x)(1+\sigma\mu_y)} \Big)\cdot
			\Big((1+\theta)\|\Delta^{y}_{k-1}\|^2
			+
			\theta\|\Delta^{y}_{k-2}\|\|\Delta^{y}_{k-1}\|\Big)
			\Bigg).
		\end{aligned}
\end{equation*}}
Next, using \cref{ASPT: unbiased noise assumption} and %the inequality
$\| a\|\|b\|\leq  \frac{1}{2}\| b\|^2
+ \frac{1}{2}\| b\|^2$,
\sa{which holds for $a,b \in \mathbb{R}^n$}, and taking the expectation leads to the desired result.
\subsection{Proof of \cref{Corollary: explicit solution to noisy LMI-R1}}
\sa{Consider arbitrary $\tau,\sigma,\pi_1,\pi_2>0$ and $\theta\in(0,1)$. By a \mg{straightforward} %trivial
	calculation,  %\eqref{Condition: noisy LMI sufficient 1}-\eqref{Condition: noisy LMI sufficient 4}
	$\{\tau, \sigma,\theta,\pi_1,\pi_2\}$ is a solution to \eqref{eq:sufficient_cond_noisy_LMI-R1}
	%in \cref{LEMMA: Noise LMI after young's ineq-R1}
	if and only if}
{\small
	\begin{subequations}
		\label{Condition: SP solution sufficient-R1}
		\begin{gather}
			\tau\geq \frac{1-\theta}{\theta\mu_x},\quad \sigma\geq \frac{1-\theta}{\theta\mu_y},\quad
			%\label{Condition: SP solution sufficient 1}\\
			%,\label{Condition: SP solution sufficient 2}\\
			\pi_1 \geq  \frac{\sigma\theta L_{yx}/\rev{c}}{1- \sigma(\pi_2+\frac{\theta}{\pi_2})L_{yy}/\rev{c}},\label{Condition: SP solution sufficient 123-R1}\\
			\sigma(\pi_2+\frac{\theta}{\pi_2})L_{yy}/\rev{c} < 1,\quad
			\frac{1}{\tau}-L_{xx} \geq\pi_1 L_{yx}. \label{Condition: SP solution sufficient 4-R1}
		\end{gather}
\end{subequations}}%
%Next, we employ \cref{Condition: SP solution to noisy LMI 3-R1} as follows
% $$
% \pi_1 = \frac{\sigma}{1-c_{\sigma}}\theta L_{yx} \frac{1}{1-\frac{\sigma}{1-c_{\sigma}}\left(\pi_2+\frac{\theta}{\pi_2}\right)L_{yy}},\quad \pi_2 = \sqrt{\theta}.
% $$
%As a result, \eqref{Condition: SP solution sufficient 3} is satisfied.
\sa{In the remainder of the proof, we fix $(\pi_1,\pi_2)$ as follows:}
{\small
	\begin{align}
		\pi_1 =   \frac{\sigma\theta L_{yx}/\rev{c}}{1-\sigma\left(\pi_2+\frac{\theta}{\pi_2}\right)L_{yy}/\rev{c}}=\rev{\frac{\sigma\theta L_{yx}/\rev{c}}{1-2\sigma\sqrt{\theta}L_{yy}/\rev{c}}},\quad \pi_2 = \sqrt{\theta}.
		\label{Condition: SP solution to noisy LMI 3-R1}
\end{align}}%

%{This choice of $\pi_1$ corresponds to the smallest value of $\pi_1$ that satisfies \eqref{Condition: SP solution sufficient 3}, whereas the choice of $\pi_2$ minimizes the left-hand side of the first inequality in \eqref{Condition: SP solution to noisy LMI 3-R1}.}
\sa{Note the definition of $\overline{\theta}$ implies that $\overline{\theta}\in (0,1)$. Next, we show that $\theta\in [\overline{\theta},1)$ implies $\pi_1,\pi_2>0$; furthermore, we also show that $\tau,\sigma>0$ defined as in \eqref{Condition: SP solution to noisy LMI-R1} for $\theta\in[\overline{\theta}, 1)$ together with $(\pi_1, \pi_2)$ as in \eqref{Condition: SP solution to noisy LMI 3-R1} is a solution to \eqref{Condition: SP solution sufficient-R1}.}

\sa{First, setting $\tau,\sigma$ as in~\eqref{Condition: SP solution to noisy LMI-R1} and $\pi_1,\pi_2$ as in~\cref{Condition: SP solution to noisy LMI 3-R1} imply that \eqref{Condition: SP solution sufficient 123-R1}
	is trivially satisfied}. Next, by substituting $\{\tau,\sigma,\pi_1,\pi_2\}$, \sa{chosen as in \eqref{Condition: SP solution to noisy LMI-R1} and \cref{Condition: SP solution to noisy LMI 3-R1},} into \eqref{Condition: SP solution sufficient 4-R1}, we conclude that $\{\tau,\sigma,\theta,\pi_1,\pi_2\}$ satisfies \eqref{Condition: SP solution sufficient-R1} for any $\theta\in (0,1)$ \rev{satisfying}
{\small
	\begin{gather}
		\frac{2L_{yy}}{\rev{c}\mu_y}\cdot\frac{1-\theta}{\sqrt{\theta}}  \leq 1-\beta,
		\label{Condition: SP solution from sigma 2-R1}\\
		\frac{\theta\mu_x}{1-\theta} - L_{xx} \geq (1-\theta)\frac{L_{yx}^2}{\rev{c}\mu_y}\cdot\Big(1- \frac{2L_{yy}}{\rev{c}\mu_y}\cdot\frac{1-\theta}{\sqrt{\theta}}\Big)^{-1}, \label{Condition: SP solution from tau-R1}
\end{gather}}%
for some \rev{$\beta\in(0,1)$}.
Clearly, a sufficient condition for \eqref{Condition: SP solution from tau-R1} is
{\small
	\begin{equation}\label{Condition: SP solution from tau 2-R1}
		\frac{\theta\mu_x}{1-\theta} - L_{xx} \geq (1-\theta)\frac{L_{yx}^2}{\mu_y}\cdot\frac{1}{\rev{c}\beta}.
\end{equation}}%
\sa{Note that \eqref{Condition: SP solution from sigma 2-R1} implies that $\pi_1>0$. We also have $\pi_2=\sqrt{\theta}>0$ trivially.}

\sa{When $L_{yy}>0$, given any $\beta\in(0,1)$, solving \cref{Condition: SP solution from sigma 2-R1,Condition: SP solution from tau 2-R1} for $\theta\in(0,1)$, we get the \sa{third} condition in \eqref{Condition: SP solution to noisy LMI-R1}. Indeed, it can be checked that $\theta\in[\overline{\theta}_2,1)$ satisfies \eqref{Condition: SP solution from sigma 2-R1} and $\theta\in[\overline{\theta}_1,1)$ satisfies \eqref{Condition: SP solution from tau 2-R1}; thus, $\theta\in[\overline{\theta},1)$ satisfies \eqref{Condition: SP solution from sigma 2-R1} and  \eqref{Condition: SP solution from tau 2-R1} simultaneously. Moreover, when $L_{yy}=0$, one does not need to solve \cref{Condition: SP solution from sigma 2-R1} as the first inequality in \eqref{Condition: SP solution sufficient 4-R1} holds trivially; thus, the only condition on $\theta$ comes from \eqref{Condition: SP solution from tau-R1} which is equivalent to \eqref{Condition: SP solution from tau 2-R1} with $\beta=1$.}  
\sa{The rest follows from  Lemma~\ref{LEMMA: Noise LMI after young's ineq-R1} by setting $\alpha = \frac{\theta L_{yx}}{\pi_1} + \frac{\theta L_{yy}}{\pi_2}$. Indeed, the particular choice of $(\pi_1,\pi_2)$ in \eqref{Condition: SP solution to noisy LMI 3-R1} gives us $\alpha = \frac{\rev{c}}{\sigma}-\sqrt{\theta}L_{yy}$.}
Finally, {it can be verified that $\overline{\theta}_1:[0,1]\to\reals$ and $\overline{\theta}_2:[0,1]\to\reals$ are monotonically decreasing and monotonically increasing functions of $\beta$, respectively. Since $\overline{\theta}_1(0)=1>\overline{\theta}_2(0)$ and $\overline{\theta}_2(1)=1>\overline{\theta}_1(1)$, %we can conclude that
	$\overline{\theta}$ obtains its minimum at \sa{the unique $\beta^*\in(0,1)$ such that $\overline{\theta}_1(\beta^*) = \overline{\theta}_2(\beta^*)$}.} %\looseness=-1

\section{Extensions and Special Cases}
In this section, we \sa{discuss the deterministic case, i.e., $\delta_x=\delta_y=0$, and we also go over a special case of SAPD when $\theta=0$, i.e., SGDA.}

%\xtodo{Do we still need this deterministic result? Are they redundent?}
\subsection{A Deterministic Primal-Dual Method~(APD)}
\label{sec:deterministic}
\sa{When $\delta_x=\delta_y=0$, i.e., $\tilde\nabla_x\Phi=\nabla_x\Phi$ and $\tilde\nabla_y\Phi=\nabla_y\Phi$, we call this deterministic variant of SAPD as APD.} APD, when applied to \eqref{eq:main-problem} with a bilinear $\Phi$, 
%\xtodo{For the consistence, shall we use iterates instead of iteration sequence? MG: I think it is ok the way it is now.} 
generates the same iterate sequence with \cite{chambolle2016ergodic} for a specific choice of step sizes; therefore, APD can be viewed as a general form of the method proposed by Chambolle and Pock \cite{chambolle2016ergodic} for bilinear SP problems. \sa{For %this special case, 
bilinear problems as in~\cite{chambolle2016ergodic},} APD hits the lower complexity bound when $\mathcal{L}$ is strongly convex in $x$ and strongly concave in $y$. Moreover, when $\Phi$ is not assumed to be bilinear, APD has the best iteration complexity bound shown for \sa{single-loop primal-dual first-order algorithm applied to} \eqref{eq:main-problem}. \sa{The convergence guarantees for the deterministic scenario follows directly from \xzh{the proof in \cref{sec:proof} by setting $\delta_x=\delta_y= 0$}.}
\begin{corollary}
\label{Corollary: determinitic SCSC}
Suppose \cref{ASPT: lipshiz gradient} hold, $\delta_x=\delta_y=0$, and $\{ x_k,y_k \}_{k\geq 0}$ be the iterates generated by \sa{APD, which is} the deterministic version of \cref{ALG: SAPD}. The parameter  $ \tau, \sigma>0$ and %$\theta\in(0,\rho]$
\sa{$\theta\geq 0$} satisfy \cref{eq: general SAPD LMI_R1} for some $\alpha \in \sa{[0, \tfrac{1}{\sigma}]}$ and $\rho\in(0,1]$.
%\nsa{I included $\alpha=\frac{1}{\sigma}$ in the bound, which may be helpful for the deterministic case.}
%\xtodo{Numerically, we can. But if we do that, $\Delta_N$ will lose the information of $y$.}
Then, for any $(x, y)$ and $(x_0,y_0) \in \mathcal{X} \times \mathcal{Y}$,  
%the following bound holds for $N\geq 1$:
{\small
\begin{equation*}
\begin{aligned}
    & \mathcal{L}(\sa{\bar{x}_{N}, y)-\mathcal{L}(x, \bar{y}_{N}})
               +
             {\frac{\rho^{\xzh{-N}}}{K_N(\rho)}d_N(x,y)}
           \leq 
           \frac{1}{K_N(\sa{\rho})}(\frac{1}{2\tau}\norm{x-x_0}^2+\frac{1}{2\sigma}\norm{y-y_0}^2), \\
\end{aligned}
\end{equation*}}%
for all $N \geq 1$, where 
% $D_{\tau,\sigma}(x,y) = \tfrac{1}{2\tau} \norm{x-x_0}^2+\tfrac{1}{2\sigma}\norm{y-y_0}^2$, 
\xzh{$d_N(x,y)=\tfrac{1}{2\tau}\|x_N-x\|^2 +  \tfrac{1}{2\sigma}\left( 1 - \alpha\sigma\right)\|y_N-y\|^2$}, $(\bar{x}_N,\bar{y}_N)$ and $K_N({\rho})$ are defined in \cref{Thm: main result_R1}.
\end{corollary}
% \sa{\cref{Corollary: determinitic SCSC} immediately follows from the proof for \cref{Thm: main result_R1}.}
\begin{remark}
\sa{This result extends the Accelerated Primal-Dual (APD) method proposed in~\cite{hamedani2018primal} for MCMC and SCMC SP problems to cover the SCSC scenario as well. Indeed, the result for the MCMC case in~\cite{hamedani2018primal} can be recovered from \cref{Corollary: determinitic SCSC} immediately by setting} $\theta=\rho = 1$ and $\mu_x = \mu_y = 0$.
%\nsa{Why did you define $\tfrac{0}{0}=1$.}
\sa{Furthermore, the step sizes suggested in \cite[Remark~2.3]{hamedani2018primal} satisfy \cref{eq: general SAPD LMI_R1} for a particular choice of $\alpha>0$. Finally,}
\sa{since $K_N(1)=N$, APD achieves the sublinear rate of $\cO(1/N)$ for the MCMC scenario.} 
\end{remark}

\sa{In the rest, we consider SP problem in~\cref{eq:main-problem} under SCSC scenario. Let $(x^*,y^*)$ denote the unique saddle point of \cref{eq:main-problem}. Next, we account for the individual effects of $L_{xx}, L_{yx}, L_{yy}$ as well as $\mu_x, \mu_y$ on the iteration complexity of APD. When $\Phi$ is bilinear, APD requires ${\mathcal{O}\left( \sqrt{1 + \frac{L_{yx}^2}{\mu_x\mu_y}}\cdot \ln(1/\epsilon)\right)}$ iteration to compute $(\bar{x},\bar{y})$ such that $\cG(\bar{x},\bar{y})\leq \epsilon$; this complexity is shown to be optimal in \cite{zhang2019lower}. Moreover, for the general case, i.e., $\Phi$ may not be bilinear, the iteration complexity of APD is ${\mathcal{O}\left((  \frac{L_{xx}}{\mu_x} + \frac{L_{yx}}{\sqrt{\mu_x\mu_y}} + \frac{L_{yy}}{\mu_y}\right)\cdot \ln(1/\epsilon))}$.} 
%We first display  the result   for Non-bilinear coupling function $\Phi(x,y)$.
\begin{proposition}\label{Prop:iteration bound for APD scsc}
{Suppose $\mu_x, \mu_y >0$, and \cref{ASPT: lipshiz gradient} hold.} \sa{Let $(x^*,y^*)\in \mathcal{X}\times \mathcal{Y}$ denote the unique SP of \eqref{eq:main-problem}.}
% Consider the APD method for solving SCSC SP problems with a general coupling function $\Phi(x,y)$ where  \cref{ASPT: lipshiz gradient,ASPT: strongly convex concave,ASPT: f and g convex} are satisfied with $\mu_x,\mu_y>0$. 
For any $\epsilon>0$, and for any given $\beta \in (0,1)$, suppose the \sa{APD} parameters $\{ \tau, \sigma, \theta\}$ are chosen such that
{\small
\begin{equation}
    \tau = \frac{1-\theta}{\mu_x \theta},\quad \sigma = \frac{1-\theta}{\mu_y \theta},\quad \theta = \bar{\theta}
\end{equation}}%
where $\bar{\theta}$ is defined in \eqref{Condition: SP solution to noisy LMI-R1}. Then, the iteration complexity of APD to generate a point $(\bar{x},\bar{y})\in\mathcal{X}\times\mathcal{Y}$ \sa{such that $\cD(\bar{x},\bar{y})\leq \epsilon$} is
{\small
\begin{align}
\label{eq:APD_complexity}
    \mathcal{O}\left(\left( 1 +  \frac{L_{xx}}{\mu_x} + \frac{L_{yx}}{\sqrt{\mu_x\mu_y}} + \frac{L_{yy}}{\mu_y}\right)\cdot\ln(1/\epsilon)\right).
\end{align}}%
Moreover, when $\Phi(x,y)$ is a bilinear function, the iteration  complexity of APD reduces to
$ {\mathcal{O}\left( \left(1 + \frac{L_{yx}}{\sqrt{\mu_x\mu_y}}\right)\cdot \ln(1/\epsilon)\right)}$.
\sa{Furthermore, assuming $\dom f\times \dom g$ is compact, APD can compute $(\bar{x},\bar{y})\in\cX\times\cY$ such that %$\sup\{\cL(\bar{x},y)-\cL(x,\bar{y}):\ (x,y)\in\dom f\times \dom g\}\leq \epsilon$ 
$\cG(\bar{x},\bar{y})\leq \epsilon$ with the same iteration complexity stated above for both bilinear and general cases of $\Phi$.}
\end{proposition}
\begin{proof}
Using the particular parameters \sa{given} in \cref{Corollary: explicit solution to noisy LMI-R1}  within \cref{Corollary: determinitic SCSC}, and following \sa{the similar arguments as in} the proof of \cref{Prop: iteration complex for SAPD_R1}, we immediately get the result. When $\Phi(x,y)$ is bilinear, we only need to set $L_{xx} = L_{yy} = 0$ \sa{in the general result to get the complexity for the bilinear case.}
\end{proof}

According to \cite{zhang2019lower}, \sa{the complexity of APD for the} bilinear case is optimal \sa{in terms of $\mu_x,\mu_y,L_{yx}$ and $\epsilon$}. \sa{Furthermore, the complexity in~\eqref{eq:APD_complexity} for the general case \xzh{obtains} the best we know for a single-loop first-order primal-dual algorithm.} 
%In the next section, we move to the special case SGDA.
\subsection{Stochastic Gradient Descent Ascent Method~\sa{(SGDA)}}\label{SGDA section}
\sa{The SGDA algorithm can be analyzed as a special case of SAPD with $\theta=0$. Our analysis leads to a wider range of admissible step sizes and establishes the iteration complexity bound for SGDA and shows its dependence on $L_{xx},L_{yx},L_{yy}$ and $\mu_x,\mu_y$ explicitly.} 
%\todo{MG: Give a reference to existing conditions on stepsize for GDA and mention if this generalizes them or not. When primal stepsize=dual stepsize, are we recovering known results for the stepsize choice of GDA?}
%\xtodo{Next section will discuss that}
\begin{corollary}\label{THM: gda theta= 0}
Suppose Assumptions~\ref{ASPT: lipshiz gradient}, \ref{ASPT: unbiased noise assumption} hold, and $\{ x_k,y_k \}_{k\geq 0}$ are generated by SAPD, stated in \cref{ALG: SAPD}, using parameters $\theta = 0$ and $\tau,\sigma>0$ satisfying
{\small
\begin{equation}\label{Condition: SGDA LMI}
        \begin{pmatrix}
        \frac{1}{\sigma} + \mu_y -\frac{1}{\rho\sigma} &  -L_{yx} & -L_{yy} \\
         -L_{yx} & \frac{1}{\tau}- L_{xx} & 0\\
      -L_{yy} & 0 &  \frac{1}{\sigma} 
    \end{pmatrix}\succeq 0,\qquad \sa{\tau \mu_x\geq \frac{1-\rho}{\rho}},
\end{equation}}%
for some  
$\rho\in\sa{(0,1)}$. %is related to convergence rate.  
%Then SAPD degenerates to SGDA. Moreover, 
Then, for any compact set $X \times  Y\subset \dom f\times \dom g$ such that $x_0\in X$ and $y_0\in Y$, and for any $\eta_x,\eta_y\geq 0$, the following bound holds for $N \geq 1$:
\xzh{
{\small
\begin{equation*}
            \mathbb{E}\Big[ \frac{1}{2\tau}\|x_N-x^*\|^2 + \frac{1}{2\sigma}\|y_N-y^*\|^2 \Big] 
           \leq 
           \rho^N (\frac{1}{2\tau}\|x_0-x^*\|^2 + \frac{1}{2\sigma}\|y_0-y^*\|^2 ) + \frac{\rho}{1-\rho}\Xi'_{\tau,\sigma},
\end{equation*}}}%
where 
%$\Delta'_{N}(x,y) \triangleq \tfrac{1}{2\rho\tau}\|x_N-x\|^2 +  \tfrac{1}{2\rho\sigma} \|y_N-y\|^2$, $\Omega'_{\tau,\sigma} \triangleq \tfrac{1}{2}(\tfrac{1}{\tau}+\eta_x)\Omega_X + \tfrac{1}{2}(\tfrac{1}{\sigma}+\eta_y)\Omega_Y$, and
\xzh{$\Xi'_{\tau,\sigma} \triangleq \tfrac{\tau}{1+\tau\mu_x}\delta_x^2
 +   \tfrac{\sigma}{1+\sigma\mu_y} \delta_y^2$.}%
\end{corollary}
\begin{proof}
Setting \sa{$\theta$ and $\alpha$ to $0$} in \cref{eq: general SAPD LMI_R1} immediately leads to the above result.
\end{proof}
\begin{remark}
\sa{When $\mu_x=\mu_y=0$, unlike SAPD, SGDA does not have an admissible $(\tau,\sigma)$ pair with convergence guarantees. Indeed, from \cref{Condition: SGDA LMI}, $\mu_x = 0$ implies that $\rho=1$ %has to be $1$ 
so that the second inequality is satisfied; furthermore, $\rho=1$ and $\mu_y=0$ imply that first diogonal element in the matrix inequality~(MI) becomes $0$; thus, there is no $(\tau,\sigma)$ such that the MI holds. It is worth emphasizing that \cref{Condition: SGDA LMI} not having a solution when $\mu_x=\mu_y=0$ is not because our analysis is not tight enough; indeed, there are examples for which SGDA iterate sequence does not converge to a saddle point when $\mu_x=\mu_y=0$.}
\end{remark}

\subsubsection{Parameter Choices for SGDA}
We provide a %special 
\sa{particular} solution to the %complicated 
matrix inequality \cref{Condition: SGDA LMI} following \sa{a similar technique we used for deriving a particular parameter choice for SAPD. } 
%Then we have a special matrix inequality for \cref{Thm: main result} as follows,
%\xuan{Note that the symbol $\pi_x$ has been used and it shares the same expression here.} With the assistant of the 
% \sa{For this} particular choice of $(\pi^x,\pi^y)$, 
% \eqref{Condition: SGDA LMI} reduces to
% %\begin{subequations}
% %\label{Condition: SGDA LMI}
% \begin{equation}
% \label{Condition: SGDA LMI}
% %\label{Condition: SAPD noisy LMI 1 GDA}
%     \begin{pmatrix}
%         \frac{1-c_\tau}{\tau}- L_{xx} & 0  & -L_{yy}  \\
%       0 & \frac{1-c_\sigma}{\sigma}  & -L_{yx}  \\
%         -L_{yy} & -L_{yx} & \sa{\frac{1}{\sigma} (1-\frac{1}{\rho}) + \mu_y}  
% \end{pmatrix}\succeq 0,\qquad 
% % \end{equation}
% % \begin{align}
% %     \label{Condition: SAPD noisy LMI 2 GDA}
%     \tau\mu_x\geq \frac{1-\rho}{\rho}.
% %\end{align}
% \end{equation}
% %\end{subequations}
Next, in \cref{LEMMA: Noise LMI after young's ineq GDA}, we 
%give an intermediate condition 
\sa{provide an auxiliary system, simpler than \cref{Condition: SGDA LMI},} to construct the particular solution given in \cref{Coroallary: explicit solution to noisy LMI GDA}.
\begin{lemma}\label{LEMMA: Noise LMI after young's ineq GDA}
%Suppose  \cref{ASPT: lipshiz gradient,ASPT: strongly convex concave,ASPT: f and g convex,ASPT: unbiased noise assumption,ASPT: compact} hold. 
Let $\tau,\sigma>0$, $\rho\in(0,1)$, and $\pi_1,\pi_2>0$ \sa{satisfy}
%be a solution to the following \sa{inequality system}:
{\small
\begin{subequations}
\label{eq:sufficient_cond_noisy_LMI GDA}
\begin{gather}
\frac{1}{\tau} - L_{xx} - \pi_1  L_{yx}\geq 0,\label{Condition: noisy LMI sufficient 1 GDA}\\  
 \frac{1}{\sigma}   -\pi_2  L_{yy}\geq 0,\label{Condition: noisy LMI sufficient 2 GDA}\\
 \sa{\frac{1}{\sigma} (1-\frac{1}{\rho})} + \mu_y \geq \frac{L_{yx}}{\pi_1} + \frac{L_{yy}}{\pi_2}, \label{Condition: noisy LMI sufficient 3 GDA}\\
 \tau\mu_x\geq \frac{1-\rho}{\rho}. \label{Condition: noisy LMI sufficient 4 GDA}
\end{gather}
\end{subequations}}%
Then $\{\tau, \sigma,\rho\}$ is a solution to %\eqref{Condition: SAPD simple LMI 1}  and \eqref{Condition: SAPD simple LMI 2}, 
\eqref{Condition: SGDA LMI}.
%; thus, a solution to \eqref{Condition: SAPD simple LMI system}.
%\eqref{Condition: SAPD simple LMI 1} and \eqref{Condition: SAPD simple LMI 2}.
\end{lemma}
\begin{proof}
We only need \sa{to verify that} the matrix inequality in \cref{Condition: SGDA LMI} %is correct
holds. \sa{Permuting the rows and columns in} \eqref{Condition: noisy LMI sufficient 3 GDA}, it follows that
{\small
\begin{equation*}
\begin{aligned}
\begin{pmatrix}
      \frac{1}{\sigma} \sa{(1-\frac{1}{\rho})}+ \mu_y  & -L_{yx}  & -L_{yy}  \\
       -L_{yx} & \frac{1}{\tau}- L_{xx}  & 0  \\
        -L_{yy} & 0 & \frac{1}{\sigma}  
\end{pmatrix}
\succeq & 
    \begin{pmatrix}
       \frac{L_{yx}}{\pi_1} + \frac{L_{yy}}{\pi_2}  & -L_{yx}  & -L_{yy}  \\
       -L_{yx} & \frac{1}{\tau}- L_{xx}  & 0  \\
        -L_{yy} & 0 & \frac{1}{\sigma}  
\end{pmatrix}
\sa{\triangleq \tilde{M}}.
\end{aligned}
\end{equation*}}%
\sa{Note $\tilde{M} = \tilde{M}_1 + \tilde{M}_2$ for}
{\small
\begin{equation*}
    \tilde{M}_1 
 \triangleq 
 \begin{pmatrix}
       \frac{L_{yx}}{\pi_1} & -L_{yx}  & 0 \\
       -L_{yx} & \frac{1}{\tau}- L_{xx}  & 0  \\
        0 & 0 & 0 
\end{pmatrix},
\qquad
\tilde{M}_2 \triangleq 
    \begin{pmatrix}
       \frac{L_{yy}}{\pi_2}  & 0  & -L_{yy}  \\
      0 & 0 & 0  \\
        -L_{yy} & 0 & \frac{1}{\sigma}  
\end{pmatrix}.
\end{equation*}}%
The condition in \eqref{Condition: noisy LMI sufficient 1 GDA} and  \eqref{Condition: noisy LMI sufficient 2 GDA} imply that
{\small
$$
        \tilde{M}_1 
 \succeq 
\begin{pmatrix}
       \frac{L_{yx}}{\pi_1} & -L_{yx}  & 0 \\
       -L_{yx} & \pi_1 L_{yx}  & 0  \\
        0 & 0 & 0 
\end{pmatrix}
\succeq 0,
\qquad
\tilde{M}_2 
 \succeq 
\begin{pmatrix}
       \frac{L_{yy}}{\pi_2}  & 0  & -L_{yy}  \\
      0 & 0 & 0  \\
        -L_{yy} & 0 &\pi_2 L_{yy}  
\end{pmatrix}
\succeq 0, %\;\text{respectively.}
$$}%
respectively. Thus $\tilde{M}_1+\tilde{M}_2\succeq 0$, which completes the proof.
\end{proof}

\sa{Lemma~\ref{LEMMA: Noise LMI after young's ineq GDA} helps us describe a subset of solutions to the matrix inequality system in \cref{Condition: SGDA LMI} using the solutions of an inequality system in~\cref{eq:sufficient_cond_noisy_LMI GDA} that is easier to deal with. Next, based on based on \cref{LEMMA: Noise LMI after young's ineq GDA}, we will construct a family of admissible parameters for SGDA, i.e., SAPD with $\theta = 0$, such that the iterate sequence will exhibit the desired convergence behavior.}
\begin{corollary}\label{Coroallary: explicit solution to noisy LMI GDA}
\sa{Suppose  $\mu_x, \mu_y > 0$. For any $\beta_1,\beta_2\in(0,1)$ such that $\beta_1+\beta_2< 1$, $\{\tau, \sigma,\rho\}$ chosen satisfying}
{\small
\begin{equation}\label{parachocie for GDA}
%\max\Big\{ \frac{1}{1 + \frac{\beta_1\mu_x\mu_y}{ \beta_1\mu_y L_{xx} + L_{yx}^2}},\quad  \frac{1}{1+ \frac{(1-c_\sigma)\beta_2\beta_3\mu_y^2}{L_{yy}^2}}\Big\} \label{parachocie for GDA}\\
        \tau = \frac{1-\rho}{\rho\mu_x },\quad \sigma =\sa{\tfrac{1}{1-\beta_1-\beta_2}}\cdot\frac{1-\rho}{\rho\mu_y} ,\quad \sa{\rho \geq \bar{\rho} \triangleq \Big(1+\frac{1}{L(\beta_1,\beta_2)}\Big)^{-1}}
        % \\
        % \pi_1 = \frac{L_{yx}}{\beta_1\mu_y},\quad \pi_2 = \frac{L_{yy}}{\beta_2\mu_y}\label{Condition: SP solution to noisy LMI 2 GDA},
        % \\
        % \pi^x = \frac{c_{\tau}}{\tau},\; \pi^y_1 =  \frac{c_{\sigma}}{\sigma}\label{Condition: SP solution to noisy LMI 4 GDA}.
\end{equation}}%
is a solution to \eqref{Condition: SGDA LMI}, where \sa{$L(\beta_1,\beta_2)\triangleq\max\Big\{\tfrac{L_{xx}}{\mu_x}+\tfrac{1}{\beta_1}\cdot\tfrac{L_{yx}^2}{\mu_x\mu_y},\ \tfrac{1}{\beta_2(1-\beta_1-\beta_2)}\cdot\tfrac{L_{yy}^2}{\mu_y^2}\Big\}$}.
\end{corollary}
\begin{proof}
 \sa{The proof is based on the result in Lemma~\ref{LEMMA: Noise LMI after young's ineq GDA}.} 
 %\cref{LEMMA: Noise LMI after young's ineq}. 
\sa{Let $\beta_1,\beta_2,\beta_3\in(0,1)$ such that $\beta_1+\beta_2<1$. Given any $\rho\in(0,1)$, let $\tau=\frac{1-\rho}{\rho\mu_x}$, $\sigma=\frac{1}{\beta_3}\cdot\frac{1-\rho}{\rho\mu_y}$ and 
%$\pi_1,\pi_2$ be as in \eqref{Condition: SP solution to noisy LMI 2 GDA}. 
\sa{let $\pi_1 = \frac{L_{yx}}{\beta_1\mu_y}$, $\pi_2 = \frac{L_{yy}}{\beta_2\mu_y}$.} If we substitute $(\tau,\sigma,\rho,\pi_1,\pi_2)$} into \eqref{Condition: noisy LMI sufficient 1 GDA}-\eqref{Condition: noisy LMI sufficient 3 GDA}, we get
 %it follows that
{\small
\begin{subequations}
\begin{gather}
\mu_x\frac{\rho}{1-\rho} - L_{xx} - \frac{L_{yx}^2}{\beta_1\mu_y}  \geq 0,\quad
 \beta_3\mu_y\frac{\rho}{1-\rho}   -\frac{L_{yy}^2}{\beta_2\mu_y}  \geq 0,\label{eq:rho_cond}\\
 \beta_3\mu_y\frac{\rho}{1-\rho} (1-\frac{1}{\rho})+ \mu_y %-\frac{1}{\rho\frac{1}{\beta_3}\cdot\frac{1-\rho}{\rho\mu_y }} 
 \geq \sa{(\beta_1+\beta_2)\mu_y}.
 %\frac{L_{yx}}{\frac{L_{yx}}{\beta_1\mu_y}} + \frac{L_{yy}}{\frac{L_{yy}}{\beta_2\mu_y}}.
 \label{eq:beta_cond}
 %\frac{1-\rho}{\rho\mu_x }\geq \frac{1-\rho}{\rho\mu_x}, 
\end{gather}
\end{subequations}}%
\sa{Next, we solve this inequality system in terms of $\rho\in(0,1)$. Note \eqref{eq:rho_cond} holds for
{\small
\begin{gather*}
    \rho \geq\max\left\{ \left(1 + \frac{1}{ \tfrac{L_{xx}}{\mu_x} + \tfrac{L_{yx}^2}{\beta_1\mu_x\mu_y}}\right)^{-1},\ \left(1+ \frac{\beta_2\beta_3}{L_{yy}^2/\mu_y^2}\right)^{-1}\right\},
\end{gather*}}%
and \eqref{eq:beta_cond} holds whenever $\beta_1 + \beta_2 + \beta_3 \leq 1$. To minimize the lower bound on $\rho$, the optimal choice for $\beta_3$ is $\beta_3=1-\beta_1-\beta_2>0$. Thus, $\{\tau,\sigma,\rho,\pi_1,\pi_2\}$ satisfying \eqref{parachocie for GDA} is a solution to \cref{eq:sufficient_cond_noisy_LMI GDA}, which 
%Finally, according to \cref{LEMMA: Noise LMI after young's ineq GDA},  
implies that} $\{\tau,\sigma,\rho\}$ is a solution to \cref{Condition: SGDA LMI}.
\end{proof}
\sa{To determine the best certifiable convergence rate, i.e., the smallest $\rho$, one can optimize $\beta_1$  and $\beta_2$.} Finally, using the above parameter choice, we %conclude 
\sa{establish} the iteration complexity bound for SGDA in the next subsection.
\subsubsection{Iteration Complexity Bound for SGDA}
\sloppy In this part, we study the %lower 
iteration complexity bound for SGDA
%, i.e., $\theta=0$ for SAPD, which are given on the number of iterations $N$ 
to generate a point $(\bar{x},\bar{y})\in\cX\times\cY$ such that \xzh{$\cD(x_\epsilon,y_\epsilon)\leq \epsilon$}. The proof technique is very similar to that for SAPD.
 
\begin{proposition}\label{Propostion: iteration complex for GDA} Suppose $\mu_x,\mu_y>0$, and  \cref{ASPT: lipshiz gradient,ASPT: unbiased noise assumption} hold. 
%\sa{Furthermore, we assume that $\dom f$ and $\dom g$ are compact sets with $\Omega_f\triangleq \sup\{\norm{x-x'}^2:~x,x'\in\dom f\}$ and $\Omega_g\triangleq \sup\{\norm{y-y'}^2:~y,y'\in\dom g\}$.}
For any $\epsilon>0$, and for any given $\beta_1,\beta_2\in(0,1)$ \sa{satisfying} $\beta_1+\beta_2 < 1$, suppose the parameters $\{\tau,\sigma\}$ are chosen such that
{\small
\begin{equation}\label{eq: SP parameter for SGDA}
    \tau = \frac{1-\rho}{\mu_x\rho},\quad \sigma = \frac{1}{1-\beta_1-\beta_2}\cdot\frac{1-{\rho}}{\mu_y\rho}, \quad \rho = \max\{ \overline{\rho}, \overline{\overline{\rho}}\},
\end{equation}}%
where $\overline{\rho}$ is defined in \cref{parachocie for GDA} and
    $\overline{\overline{\rho}} \triangleq \max \{ \overline{\overline{\rho}}_1, \overline{\overline{\rho}}_2\}$
such that
\xzh{{\small
\begin{equation}
    \label{eq:barbar-rho-gda}
    \sa{\overline{\overline{\rho}}_1=\max\Big\{0,~1 - \frac{(1-\beta_1-\beta_2)\mu_x}{6\delta_x^2}~\epsilon\Big\}, \quad
\overline{\overline{\rho}}_2=\max\Big\{0,~\frac{(1-\beta_1-\beta_2)^2\mu_y}{6\delta_y^2}~\epsilon\Big\}}
\end{equation}}
}%
{with the convention that $ \overline{\overline{\rho}}_1 = 0$ if $\delta_x^2 = 0$ and $ \overline{\overline{\rho}}_2 = 0$ if $\delta_y^2 = 0$}. Then the iteration complexity of SGDA method, i.e., SAPD with $\theta = 0$, {as stated in~\cref{ALG: SAPD},} to generate a point $(x_\epsilon,y_\epsilon)\in\cX\times\cY$ such that \xzh{$\cD(x_\epsilon,y_\epsilon)\leq \epsilon$} is
{\small
\begin{equation}\label{Iteration bound for noisy SCSC GDA}
    \mathcal{O}\left( \left( 
\frac{L_{xx}}{\mu_x} + \frac{L_{yx}^2}{\mu_x\mu_y} + \frac{L_{yy}^2}{\mu_y^2}
+ \left( \frac{\delta_x^2}{\mu_x} +  \frac{\delta_y^2}{\mu_y} \right)\frac{1}{\epsilon}\right) \cdot \ln\left( %\frac{\mu_x{\Omega_f} + \mu_y{\Omega_g}}{\epsilon}+ 1
\xzh{\frac{\cD(x_0,y_0)}{\epsilon}}
\right)\right).
\end{equation}}%
\end{proposition}
\begin{proof}
Given $\beta_1,\beta_2\in(0,1)$ \sa{such that}
$
\beta_1+\beta_2< 1$, \xzh{letting %$\{\tau,\sigma,\rho,\pi^x,\pi^y\}$
\sa{$\{\tau,\sigma,\rho\}$ be chosen according to \cref{eq: SP parameter for SGDA}}, we know that \cref{Condition: SGDA LMI} is satisfied by \cref{Coroallary: explicit solution to noisy LMI GDA}.
%and we set $\eta_x = \frac{1}{\tau} + \mu_x$, 
%$\eta_y = \frac{1}{\sigma} + \mu_y$. 
Therefore}, using these particular parameter values,
%\sa{and setting $X=\dom f$ and $Y=\dom g$ within \cref{THM: gda theta= 0}}, 
it follows  from \cref{THM: gda theta= 0} that
\xzh{{\small
\begin{equation*}
\begin{aligned}
          &  \mathbb{E}\Big[\mu_x\|x_N-x^*\|^2 + (1-\beta_1-\beta_2)\mu_y\|y_N-y^*\|^2\Big] \\
        \leq 
        & \rho^N \Big( \mu_x\|x_0-x^*\|^2 + (1-\beta_1-\beta_2)\mu_y\|y_0-y^*\|^2 \Big) +  \frac{2(1-\rho)}{\mu_x}\delta_x^2
     +  \frac{2(1-\rho)}{(1-\beta_1-\beta_2)\mu_y}\delta_y^2.
\end{aligned}
\end{equation*}}}%
Because $1-\beta_1-\beta_2 \in(0,1)$, we further know that
\xzh{{\small
\begin{equation}\label{bound gda}
\begin{aligned}
          \mathbb{E}\Big[\cD(x_N,y_N)\Big] 
        \leq 
        \frac{1}{ (1-\beta_1-\beta_2)}\rho^N\cD(x_0,y_0)
        +  \frac{2(1-\rho)}{ (1-\beta_1-\beta_2)\mu_x}\delta_x^2
     +  \frac{2(1-\rho)}{(1-\beta_1-\beta_2)^2\mu_y}\delta_y^2.
\end{aligned}
\end{equation}}}%
For any $\epsilon>0$, the right side of (\ref{bound gda}) can be bounded by $\epsilon>0$ when
\xzh{{\small
\begin{equation}
    \frac{1}{ (1-\beta_1-\beta_2)}\rho^N \cD(x_0,y_0) \leq \frac{\epsilon}{3},
    \quad  \frac{2(1-\rho)}{ (1-\beta_1-\beta_2)\mu_x}\delta_x^2\leq \frac{\epsilon}{3},\quad
  \frac{2(1-\rho)}{(1-\beta_1-\beta_2)^2\mu_y}\delta_y^2 \leq \frac{\epsilon}{3}. \label{n1 gda} %\label{n2 gda}
\end{equation}}}%
\xzh{\sa{Substituting $\overline{\overline{\rho}}$ values given in \cref{parachocie for GDA} into the second and the third conditions in \eqref{n1 gda}, we have that these two conditions will hold when $\rho\geq\overline{\overline{\rho}}$. Moreover,
the first inequality in~\eqref{n1 gda} holds for
{$N\geq 1 + \ln(\frac{3}{1-\beta_1-\beta_2}
\cD(x_0,y_0)/\epsilon)/{\ln(\tfrac{1}{\rho})}$}. Thus, SGDA can generate a point $(x_\epsilon,y_\epsilon)\in\cX\times\cY$ such that $\cD(x_\epsilon,y_\epsilon)\leq \epsilon$ within}
{\small
\begin{equation*}
 N_\epsilon = \mathcal{O}\Big(  \ln\Big(\frac{
 \cD(x_0,y_0)}{\epsilon}\Big)/\ln(\tfrac{1}{\rho}) \Big)
\end{equation*}}
iterations.} Then the rest of the proof is repeating the proof of \cref{Prop: iteration complex for SAPD_R1}.
%, and one can easily show that the number of iterations for \cref{n1 gda} to be satisfied is bounded by \eqref{Iteration bound for noisy SCSC GDA}.
\end{proof}

\sa{Since we \xzh{adopt Gauss-Seidel type update rather than a Jacobi-type}, the effect of Lipschitz constants in the complexity bound are different, i.e., compare $\frac{L_{xx}}{\mu_x}$ with $\frac{L_{yy}^2}{\mu_y^2}$. Furthermore, we also observe that adopting a momentum term as in SAPD, i.e., $\theta>0$, the $\cO(1)$ constant improves from    $\frac{L_{xx}}{\mu_x}+\frac{L_{yx}^2}{\mu_x\mu_y}+\frac{L_{yy}^2}{\mu_y^2}$ for SGDA to $\frac{L_{xx}}{\mu_x}+\frac{L_{yx}}{\sqrt{\mu_x\mu_y}}+\frac{L_{yy}}{\mu_y}$ for SAPD.}

\section{{Supporting Results for the Robustness Analysis}}
\sa{In this section, we provide some details about our robustness analysis.}\label{appendix-C}
\subsection{CP parameters}
\label{section: CP parameters explanation}
%For the bilinear objective function \cref{eq: bilinear example},
\sa{Consider \eqref{eq:main-problem} with $\Phi$, $f$ and $g$ defined as in~\eqref{eq:special-Phi}.} Using the notations in our paper, the step size condition in \cite[Algorithm 5]{chambolle2016ergodic} can be summarized as
{\small
\begin{equation}
    \label{eq: CP stepsize}
    1+\mu_x\tau = 1 + \mu_y\sigma = \frac{1}{\theta},\qquad \frac{1}{\tau} \geq \theta L_{yx}^2\sigma.
\end{equation}}%
In fact, the above condition is a quadratic inequality of $\theta$, which is
$
\frac{L_{yx}^2}{\mu_y\mu_y}(1-\theta)^2 - \theta \leq 0;
$
thus, $\theta\in [1+\frac{\mu_x\mu_y}{2L_{yx}^2} - \sqrt{(1+\frac{\mu_x\mu_y}{2L_{yx}^2})^2 - 1},1]$. In \cref{fig:fundamental trade-off curve}, we compute and plot the $(\rho_{\rm true},\cJ)$ for all possible $(\tau,\sigma,\theta)$ satisfying \cref{eq: CP stepsize}. Moreover, condition \cref{eq: CP stepsize} holds with equality \sa{at the point indicated with $``*"$ in red color.}
\subsection{\sa{Convergence of the Gap Function Bias Term for \eqref{eq:special-Phi}}}\label{sec:gap-rate}
\sa{Consider \eqref{eq:main-problem} for $\Phi$, $f$ and $g$ as defined in~\eqref{eq:special-Phi}. We will show $\cG(x_k,y_k)$ and \xzh{$\mathbb{E}[\|z_N-z^*\|^2]$} converge with the same rate, thus $\cG(x_k,y_k)$ has the same rate with $d_N^*$, where $\cG$ is defined in~\eqref{eq:gap} and $d_N^*$ is defined in \cref{Thm: main result_R1}.} 
%The notations refer to \cref{sec:explicit-robustness}.
\sa{First, we can compute $\cG(x_{k},y_{k})$ explicitly, i.e.,}
{\small
\begin{equation}
\label{eq:Gk-bilinear}
    \cG(x_{k},y_{k}) = \mathbb{E}\left[\frac{\mu_x}{2}\|x_{k}\|^2 + \frac{1}{2\mu_y}\|K x_{k}\|^2 + \frac{\mu_y}{2}\|y_{k}\|^2 + \frac{1}{2\mu_x} \|\sa{K^\top} y_{k}\|^2\right].
\end{equation}}%
\sa{Recall the \mg{augmented} vector \xzh{$\tilde{z}_{k}=[x_{k-1}; y_k]$} obtained by vertical concatenation for all $k\geq 0$ such that $x_{-1}=x_0$ and $y_{-1}=y_0$, and $(x_0,y_0)\in\cX\times\cY$ is a given initial point. Let $P_x$ and $P_y$ be matrices with appropriate dimensions such that \xzh{$x_k=P_x \tilde{z}_{k+1}$} and $y_k=P_y \tilde{z}_k$. Note %that 
\eqref{lin-dyn-sys} implies $\tilde{z}_k=A^k\tilde{z}_0+\sum_{i=1}^{k}A^{i-1}B\omega_{k-i}$; thus,
% \begin{equation*}
%     \|x_{k+1}\|^2 = \|\sa{A_1 \tilde{z}_k} + \newcomment{B_1\omega_k} \|^2 = \|A_1A^k\tilde{z}_0 + \sum_{i=1}^{k}\sa{A_1A^{i-1}}B\omega_{k-i} +\newcomment{B_1 \omega_k} \|^2.
% \end{equation*}
% %where we use the relation $z_{k} = A^kz_0 + \sum_{i=0}^{k-1}A^iB\omega_{k-1-i}$, and $\{a_i^T\}$ are block rows of $A$.
% \newcomment{If we define $\mathbf{1}_{d,x} \triangleq diag(\mathbf{I}_{d\times d}, \mathbf{0}_{3d\times 3d}),\; \mathbf{1}_{d,y} \triangleq diag(\mathbf{0}_{d\times d},\mathbf{I}_{d\times d}, \mathbf{0}_{2d\times 2d}) \in \mathbb{R}^{d\times4d}$. Then}
%{\small
%\begin{equation*}
   \xzh{$\|x_{k}\|^2 = \|P_x\tilde{z}_{k+1} \|^2 = \|P_x A^{k+1}\tilde{z}_0 + \sum_{i=1}^{k+1}P_x A^{i-1}B\omega_{k+1-i} \|^2$.}
%\end{equation*}}%
%For the sake of simplicity, we assume $\delta_x = \delta_y = \delta$ for some $\delta>0$. 
The noise model assumed in \cref{assump-additive-noise} and \eqref{assump-gaussian-noise} implies that
% \xtodo{I think the expectation of inner product of $\omega_0$ and $\tilde{z}_0$ is 0. Or we can assume it without loss of generality.}
{\footnotesize
\begin{align*}
    \mathbb{E}\left[  \|x_{k}\|^2  \right] = &  \|P_x A^{\xzh{k+1}}\tilde{z}_0\|^2 + \delta^2\sum_{i=1}^{\xzh{k+1}}\Tr((P_x A^{i-1}B)^\top P_x A^{i-1}B) =\|P_x A^{\xzh{k+1}}\tilde{z}_0\|^2 + \delta^2\sum_{i=1}^{\xzh{k+1}}\|P_xA^{i-1}B\|_F^2.
\end{align*}}%
We can also write the other terms in~\eqref{eq:Gk-bilinear} using the same argument as above:
{\footnotesize
\begin{align*}
    \mathbb{E}\left[  \|Kx_{k}\|^2  \right] &=  \|K P_x A^{\xzh{k+1}}\tilde{z}_0\|^2 + \delta^2\sum_{i=1}^{\xzh{k+1}}\|KP_x A^{i-1}B\|_F^2,\\
    \mathbb{E}\left[  \|y_{k}\|^2  \right] &=  \|P_y A^k\tilde{z}_0 \|^2 + \delta^2\sum_{i=1}^{k}\|P_yA^{i-1}B\|_F^2,\\
    \mathbb{E}\left[  \|K^\top y_{k}\|^2  \right] &=   \|K^\top P_y A^k\tilde{z}_0 \|^2 + \delta^2\sum_{i=1}^{k}\|K^\top P_y A^{i-1}B\|_F^2.
\end{align*}}%
Therefore,
{\footnotesize
\begin{equation*}
\begin{aligned}
     & \cG(x_{k},y_{k}) \\
     = & \frac{\mu_x}{2}\|P_x A^{\xzh{k+1}}\tilde{z}_0 \|^2 +\frac{\mu_y}{2}\|P_y A^{\xzh{k+1}}\tilde{z}_0 \|^2 + \frac{1}{2\mu_y} \|KP_x A^k\tilde{z}_0 \|^2+  \frac{1}{2\mu_x}\|K^\top P_y A^k\tilde{z}_0 \|^2 
     +
     \xzh{\delta^2\Big( \frac{\mu_x}{2}\|P_xA^{\xzh{k}}B\|_F^2 + \frac{\mu_y}{2}\|P_y A^{\xzh{k}}B\|_F^2}
     \\
    &  +\sum_{i=1}^{k} \big(\frac{\mu_x}{2}\|P_xA^{i-1}B\|_F^2 + \frac{\mu_y}{2}\|P_y A^{i-1}B\|_F^2 +\frac{1}{2\mu_y} \|K P_x A^{i-1}B\|_F^2 +\frac{1}{2\mu_x}\|K^\top P_y A^{i-1}B\|_F^2\big)\Big). \label{bias-gap-decay}
\end{aligned}
\end{equation*}}%
\mg{The matrix $A$ is non-symmetric in general. By considering the Jordan decomposition of the \xzh{$2d\times 2d$} matrix $A$, it is known that there exists a positive constant $c_1$ and a non-negative integer \xzh{$0\leq m_1< 2d$} such that}
%By Gelfand's formula (see e.g. \cite{kozyakin2009accuracy}), there exists a sequence $\varepsilon_j\geq 0$ such that for every $j\geq 1$,
\mg{$\| A^{k}\| \leq c_1 k^{m_1} \rho(A)^k, %\quad \mbox{for every} \quad k\geq 1, %(\rho(A) + \varepsilon_j)^j \quad \mbox{and} \quad \lim_{j\to\infty} \varepsilon_j = 0
$ for all $k\geq 1$
(see e.g. \cite{golub1996matrix,strang2005linear}).
%Applying this inequality to \eqref{bias-gap-decay}, we observe that 
%$\cG(x_{k},y_{k}) \leq c_2  $
%This shows that 
\sa{Therefore,} the bias term of $\cG(x_{k},y_{k})$ is bounded by $c_2 \xzh{(k+1)}^{2m_1}\rho(A)^{2k}$ for some positive constant $c_2$. \sa{Thus, we} conclude that the bias diminish\mg{es} exponentially with rate $\rho(A)^2$.}
}
% \begin{figure}[h]
%   \centering
%   \includegraphics[width=.5\linewidth]{picture/robustness_analysis_with_different_alpha/influence_of_calpha-eps-converted-to.pdf}
%   \label{fig:influence of c_alpha}
% \caption{The trade-off between the robustness level \sa{$\bar{\cR}_{c,m}(\rho^*_c)$ and the convergence rate $\rho^*_{c}$} for a quadratic function.}
% \end{figure}
% \begin{figure}[h]
%   \centering
%   \includegraphics[width=.32\linewidth]{picture/robustness_analysis_with_different_alpha/alpha0point10-eps-converted-to.pdf}
%   \includegraphics[width=.32\linewidth]{picture/robustness_analysis_with_different_alpha/alpha0point50-eps-converted-to.pdf}
%   \includegraphics[width=.32\linewidth]{picture/robustness_analysis_with_different_alpha/alpha0point90-eps-converted-to.pdf}
% %   \includegraphics[width=.7\linewidth]{picture/Robustness_analysis/robustness_compare.png}
% %   \includegraphics[width=.7\linewidth]{picture/Robustness_analysis/robustness_compare_sp.png}
%   \label{fig:quadratic robustness level figure}
% \caption{The rate-robustness trade-off based on different $c$ values for a quadratic function}
% \end{figure}

\subsection{$\cC_{\rho}$ is a connected set}
\label{sec:connected}
{Given $\rho$, we \mg{next show that the set $\cC_{\rho}$ is connected.} This result allows us to use \cref{lemma: upper bound of robustness problem} for optimizing the robustness.}
\begin{lemma}
\label{lemma: c rho connected}
For any $\rho\in[\rho^*,1)$, $\cC_\rho \subset (0,1)$ is a \mg{non-empty} convex set \mg{where $\cC_\rho$ is defined by \eqref{def-C-rho}}. Hence, it is connected.
\end{lemma}
\begin{proof}
Since $\cC_\rho\neq\emptyset$, let $c_1,c_2\in\cC_\rho$. Without loss of generality, \mg{suppose} $c_2\geq c_1$. We aim to show that for any $\beta\in [0,1]$, we have $c^*\triangleq\beta c_1 + (1-\beta )c_2\in \cC_{\rho}$. Since $c_1,c_2\in\cC_{\rho}$, there exist $t_i,s_i,\theta_i>0$ such that $G_{\rho}(t_i,s_i,\theta_i,c_i s_i)\succeq \mathbf{0}$ for $i=1,2$. For a given %some free
$\lambda\in[0,1]$, let $(t^*,s^*,\theta^*) \triangleq \lambda (t_1,s_1,\theta_1) + (1-\lambda )(t_2,s_2,\theta_2)$. \mg{It suffices to} %In the following, we will 
construct $\lambda\in[0,1]$ such that $G_{\rho}(t^*,s^*,\theta^*,c^*s^*)\succeq \mathbf{0}$.  
%\mg{where $(t^*,s^*,\theta^*) \triangleq \lambda (t_1,s_1,\theta_1) + (1-\lambda )(t_2,s_2,\theta_2)$}. 
This \mg{will show} that $c^*\in\cC_\rho$. 

It follows from the definition of $G_\rho$ that
{\footnotesize
\allowdisplaybreaks
\begin{equation*}
\begin{aligned}
&G_{\rho}(t^*,s^*,\theta^*,c^*s^*) =
\\
&
      \lambda\begin{pmatrix}
 (1- \frac{1}{\rho}) t_1+\mu_x & 0 & 0 & 0 & 0\\ 
  0 & (1- \frac{1}{\rho}) s_1+\mu_y & (\frac{\theta_1}{\rho} - 1)L_{yx} & (\frac{\theta_1}{\rho} - 1)L_{yy} & 0\\ 
  0 & (\frac{\theta_1}{\rho} - 1)L_{yx} & t_1 - L_{xx} & 0 & -  \frac{\theta_1}{\rho}L_{yx}\\ 
  0& (\frac{\theta_1}{\rho} - 1)L_{yy} & 0 & (1-c_1)s_1 & -  \frac{\theta_1}{\rho}L_{yy}\\
  0 & 0 & - \frac{\theta_1}{\rho}L_{yx} & -  \frac{\theta_1}{\rho}L_{yy} & \frac{c_1s_1}{\rho}
\end{pmatrix}
\\
  & + 
  \lambda \begin{pmatrix}
  0 & 0 & 0 & 0 & 0 \\
  0 & 0 & 0 & 0 & 0 \\
  0 & 0 & 0 & 0 & 0 \\
  0 & 0 & 0 & (c_1-c^*)s_1 & 0\\
  0 & 0 & 0 & 0 & \frac{(c^* - c_1)s_1 }{\rho}
  \end{pmatrix}\\
     &  + (1-\lambda)\begin{pmatrix}
 (1- \frac{1}{\rho}) t_2+\mu_x & 0 & 0 & 0 & 0\\ 
  0 & (1- \frac{1}{\rho}) s_2+\mu_y & (\frac{\theta_2}{\rho} - 1)L_{yx} & (\frac{\theta_2}{\rho} - 1)L_{yy} & 0\\ 
  0 & (\frac{\theta_2}{\rho} - 1)L_{yx} & t_2 - L_{xx} & 0 & -  \frac{\theta_2}{\rho}L_{yx}\\ 
  0& (\frac{\theta_2}{\rho} - 1)L_{yy} & 0 & (1-c_2)s_2 & -  \frac{\theta_2}{\rho}L_{yy}\\
  0 & 0 & - \frac{\theta_2}{\rho}L_{yx} & -  \frac{\theta_2}{\rho}L_{yy} & \frac{c_2s_2}{\rho}
  \end{pmatrix}
  \\
  & + 
  (1-\lambda) \begin{pmatrix}
  0 & 0 & 0 & 0 & 0 \\
  0 & 0 & 0 & 0 & 0 \\
  0 & 0 & 0 & 0 & 0 \\
  0 & 0 & 0 & (c_2-c^*)s_2 & 0\\
  0 & 0 & 0 & 0 & \frac{(c^* - c_2)s_2 }{\rho}
    \end{pmatrix}\\
% \end{aligned}
% \end{equation*}
% \begin{equation*}
%     \begin{aligned}
    =& \lambda G_{\rho}( t_1 , s_1 ,\theta_1, c_1 s_1) + (1-\lambda) G_{\rho}( t_2 , s_2 ,\theta_2, c_2 s_2)
    \\
    & + 
  \begin{pmatrix}
  0 & 0 & 0 & 0 & 0 \\
  0 & 0 & 0 & 0 & 0 \\
  0 & 0 & 0 & 0 & 0 \\
  0 & 0 & 0 & \lambda(c_1-c^*)s_1 + (1-\lambda) (c_2-c^*)s_2 & 0\\
  0 & 0 & 0 & 0 & \frac{\lambda(c^* - c_1)s_1 + (1 - \lambda)(c^* - c_2)s_2 }{\rho}
  \end{pmatrix} 
\end{aligned}
\end{equation*}}%
Since $c_i\in\cC_\rho$ implies $G_{\rho}( t_i , s_i ,\theta_i, c_i s_i) \succeq \mathbf{0}$ for $i=1,2$, we have $G_{\rho}(t^*,s^*,\theta^*,c^*s^*)\succeq 0$ if
{\small
\begin{equation*}
    \begin{aligned}
        &\lambda(c_1 -\beta c_1 - (1 - \beta) c_2)s_1 + (1-\lambda)(c_2 -\beta c_1 - (1 - \beta) c_2)s_2 \geq 0 \\
        & \lambda(\beta c_1 +(1 - \beta) c_2 - c_1)s_1 + (1-\lambda)(\beta c_1 +(1 - \beta) c_2 - c_2)s_2 \geq 0.
    \end{aligned}
\end{equation*}}%
Therefore, to show the desired result, it is sufficient to find $\lambda\in[0,1]$ such that
{\small
\begin{equation*}
    \begin{aligned}
        &\lambda(1 - \beta) (c_1- c_2)s_1 + (1-\lambda)\beta(c_2 - c_1)s_2 \geq 0 \\
        & \lambda(1 - \beta) (c_2 - c_1)s_1 + (1-\lambda)\beta(c_1 - c_2)s_2 \geq 0.
    \end{aligned}
\end{equation*}}%
This system is equivalent to $\lambda(1 - \beta) s_1 - (1 - \lambda) \beta s_2 = 0$, which yields $\lambda = \frac{\beta s_2}{(1-\beta)s_1 +\beta s_2}$ and this completes the proof.
\end{proof}

\section{Multi-stage SAPD (M-SAPD)}
\label{sec:m-sapd}
%\todo{Multi-stage SAPD should we given as an appendix/supp. material to the reply letter otherwise it will be distractive. say we can do it, and see the appendix of this reply letter.}
\sa{Consider running \texttt{SAPD} in stages as shown in~\cref{ALG:M-SAPD}. \mg{The main idea is to run each stage $t$ for $n_t$ iterations, where within each stage constant primal and dual stepsize $\tau_t$, $\sigma_t$ and momentum parameter $\theta_t$ that depends on the stage is used. By choosing these constants $n_t$,  $\tau_t$, $\sigma_t$ and $\theta_t$ in a particular fashion, we will show that we can improve the complexity of SAPD by a logarithmic factor.}}
%\xtodo{I use subscript for $\tau_t,\sigma_t,\theta_t$ because we have $2^t$.}
%\xtodo{We should add subscript here. MG: Ok i did that.}
\begin{algorithm}[h]
\caption{Multi-stage Stochastic Accelerated Primal-Dual (M-SAPD) Algorithm}
{\label{ALG:M-SAPD}
{\small
\begin{algorithmic}[1]
\STATE Initial point $(x^0_0,y^0_0)$, \sa{parameter sequence $\{\tau_t,\sigma_t,\theta_t\}$, the stage-length sequence $\{n_t\}$}.
Set $n_0 = 0$.
\FOR{$t\geq 0$}
\STATE $(x^{t+1}_0,y^{t+1}_0)\leftarrow \texttt{SAPD}(x^t_0,y^t_0,\tau_t,\sigma_t,\theta_t,n_t)$
\ENDFOR
\end{algorithmic}}}
\end{algorithm}

%To running M-SAPD, we provide a particular solution to our LMI. It can be verified that this particular solution satisfies MI by the same way in the main article.
\mg{The following result is a simple consequence of our \cref{Corollary: explicit solution to noisy LMI-R1}, which builds on a particular choice of stepsize and momentum in our framework.}
\begin{corollary}\label{Coroallary: explicit solution to noisy LMI-RS}
Suppose $\mu_x,\mu_y > 0$. If $L_{yy}>0$, for any given $\beta\in(0,1)$, let $\tau, \sigma>0$ and $\theta\in (0,1)$ be chosen satisfying
{\small
\begin{equation}
\label{Condition: SP solution to noisy LMI-RS}
\tau = \frac{1-\theta}{\mu_x \theta},\quad \sigma = \frac{1-\theta}{\mu_y\theta},
\quad
\theta  \geq \bar{\theta}\triangleq\max\{\bar{\theta}_1,~\bar{\theta}_2\},
\end{equation}}
where $\bar{\theta}_1,~\bar{\theta}_2\in(0,1)$, depending on the choice of $\beta$, are defined as
{\footnotesize
\begin{equation}
\label{eq:theta1-RS}
        \bar{\theta}_1\triangleq 1 -\tfrac{\beta (L_{xx} + \mu_x) \mu_y}{4L_{yx}^2}
          \Big(\sqrt{ 1+ \tfrac{8\mu_xL_{yx}^2}{\beta\mu_y(L_{xx}+\mu_x)^2}}-1\Big),\quad \bar{\theta}_2\triangleq
        1 - \tfrac{(1-\beta)^2}{32}\tfrac{\mu_y^2}{L_{yy}^2} \Big( \sqrt{1+\tfrac{64L_{yy}^2}{(1-\beta)^2\mu_y^2}}-1\Big).
\end{equation}}%
If $L_{yy}=0$, let $\tau, \sigma>0$ and $\theta\in (0,1)$ be chosen as in \eqref{Condition: SP solution to noisy LMI-RS} for $\bar{\theta}_1$ in \eqref{eq:theta1-RS} with $\beta=1$ and $\bar{\theta}_2=0$. Then $\alpha = \frac{1}{2\sigma}-\sqrt{\theta}L_{yy}>0$,  and $\{\tau, \sigma,\theta,\alpha\}$ is a solution to MI \cref{eq: general SAPD LMI_R1}.%\todo{MG:Here, lets cite the exact MI equation(Sure, I did that)} 
Moreover, when $L_{yy}>0$, the minimum $\bar{\theta}$ is attained at unique $\beta^*\in (0,1)$ such that $\bar{\theta}_1=\bar{\theta}_2$.
\end{corollary}
\begin{proof}
It directly follows from \cref{Corollary: explicit solution to noisy LMI-R1} by letting $c=\frac{1}{2}$.
\end{proof}
\mg{We recall that in \cref{Thm: main result_R1}, we obtained the performance bound
\begin{equation}
\label{ineq-perf-bound}
    \mathbb{E}[d_N^*]
              \leq \rho^{N}\underbrace{\Big(\tfrac{1}{2\tau}\| x_0 - x^*\|^2
    +  \tfrac{1}{2\sigma}\| y_0 - y^*\|^2\Big)}_{ D_{\tau,\sigma}} + \frac{\rho}{1-\rho}~ %(1-\rho^N)
              \underbrace{\Big(\tfrac{\tau}{1+\tau\mu_x} \Xi^x_{\tau,\sigma,\theta} \delta_x^2
              + \tfrac{\sigma}{1+\sigma\mu_y} \Xi^y_{\tau,\sigma,\theta}\delta_y^2\Big)}_{{\Xi}_{\tau,\sigma,\theta}},\vspace*{-2mm}
\end{equation}
where $\mathbb{E}[d_N^*]$ denotes the weighted expected distance squared to the saddle point at the $N$-th step},
$${\Xi^x_{\tau,\sigma,\theta}} \triangleq 1 + \tfrac{\sigma{\theta}(1+\theta)L_{yx}}{2(1+\sigma{\mu_y})},
\quad
{\Xi^y_{\tau,\sigma,\theta}} \triangleq \tfrac{\tau\theta(1+\theta)L_{yx}}{
     {2(1+\tau\mu_x)}}+\Big(1+2\theta + \tfrac{\theta + \sigma\theta(1+\theta)L_{yy}}{1+\sigma\mu_y}  +  \tfrac{\tau\sigma\theta(1+\theta)L_{yx}L_{xy}}{(1+\tau\mu_x)(1+\sigma\mu_y)} \Big)
    (1+2\theta).$$
\mg{With the choice of parameters given in \cref{Coroallary: explicit solution to noisy LMI-RS}}, we can also provide the following \mg{explicit} bound for the ``variance term" ${\Xi}_{\tau,\sigma,\theta}$ \mg{on the right hand-side of \eqref{ineq-perf-bound}}. %convergence
%result
%\cref{Thm: main result_R1}.
%In the following lemma, we provide a bound for the variance term ${\Xi}_{\tau,\sigma,\theta}$ in the convergence result \cref{Thm: main result_R1} given certain paramater choice.
\begin{lemma}\label{lemma:Xi-bound-Response}
Suppose $\{\tau,\sigma,\theta\}$ are choose according to \cref{Condition: SP solution to noisy LMI-RS}. In addition, let  $$\theta\geq 1-\min\Big\{ \frac{\mu_y}{L_{yx}}, \frac{\mu_y}{L_{yy}},\sqrt{\frac{\mu_x\mu_y}{L_{yx}L_{xy}}}, \frac{\mu_x}{L_{yx}}\Big\}.$$ Then we have ${\Xi}_{\tau,\sigma,\theta}\leq 25(1-\theta)\Big(\frac{\delta_x^2}{\mu_x} +\frac{\delta_y^2}{\mu_y}\Big).$
\end{lemma}%\todo{MG: In the new version, we seem to have ${\Xi}_{\tau,\sigma,\theta}$ instead of ${\Xi}_{\tau,\sigma,\theta}$.(fixed it)}
\begin{proof}
\mg{Substituting} $\tau = \frac{1-\theta}{\theta\mu_x}$ and $\sigma=\frac{1-\theta}{\theta\mu_y}$ into ${\Xi}_{\tau,\sigma,\theta}$, \mg{after straightforward computations,}%it follows that
\begin{equation}\label{eq:Xi-simple}
{\Xi}_{\tau,\sigma,\theta} = (1-\theta)\Big(\Xi^x_{\tau,\sigma,\theta}\frac{\delta_x^2}{\mu_x} + \Xi^y_{\tau,\sigma,\theta}\frac{\delta_y^2}{\mu_y}\Big),
\end{equation}
with $\Xi^x_{\tau,\sigma,\theta}= 1+ \frac{\theta(1+\theta)(1-\theta)}{2}\frac{L_{yx}}{\mu_y}$, and
{\small
$$
\Xi^y_{\tau,\sigma,\theta} = \left(1+2\theta +\theta^2+ \theta(1+\theta)(1-\theta)\frac{L_{yy}}{\mu_y} + \theta(1+\theta)(1-\theta)^2\frac{L_{yx}}{\mu_x}\frac{L_{xy}}{\mu_y}\right)(1+2\theta)+\frac{\theta(1+\theta)(1-\theta)}{2}\frac{L_{yx}}{\mu_x}.
$$}%
Moreover, using the fact that $\theta\leq 1$, we obtain that
\begin{equation}\label{eq:Xi-bound-M}
    \Xi^x_{\tau,\sigma,\theta} \leq 1 + (1-\theta)\frac{L_{yx}}{\mu_y}, \;\Xi^y_{\tau,\sigma,\theta} \leq 12 + 6(1-\theta)\frac{L_{yy}}{\mu_y} +6(1-\theta)^2\frac{L_{yx}}{\mu_x}\frac{L_{xy}}{\mu_y} +  (1-\theta)\frac{L_{yx}}{\mu_x}.
\end{equation}
On the other hand, since $1-\theta\leq \min\{\frac{\mu_y}{L_{yx}},\frac{\mu_y}{L_{yy}},\sqrt{\frac{\mu_x\mu_y}{L_{yx}L_{xy}}}, \frac{\mu_x}{L_{yx}}\}$, using \cref{eq:Xi-bound-M} %\todo{MG: mention which inequalities more specifically? you mean the last one?(I add the equation number)} 
within \cref{eq:Xi-simple} completes the proof.
\end{proof}
The following corollary states the convergence result of SAPD by using
\mg{our} particular parameter choice. It will help us to establish convergence \mg{bounds} for M-SAPD in each stage.
\begin{corollary}\label{Corollary:one-step-relation-for-M-SAPD}
Suppose Assumptions~\ref{ASPT: lipshiz gradient}, \ref{ASPT: unbiased noise assumption} hold, and $\{z_k\}_{k\geq0}=\{ (x_k,y_k) \}_{k\geq 0}$ are generated by SAPD stated in \cref{ALG: SAPD}. Let $z^*=(x^*,y^*)$ be the unique saddle point of $\cL(x,y)$.
\mg{Suppose that the} parameters $\{\tau,\sigma,\theta\}$ are chosen according to \cref{Condition: SP solution to noisy LMI-RS}. In addition, let  $$\theta\geq 1-\min\{\frac{\mu_y}{L_{yx}},\frac{\mu_y}{L_{yy}},\sqrt{\frac{\mu_x\mu_y}{L_{yx}L_{xy}}}, \frac{\mu_x}{L_{yx}}\}.$$ Then, for any $N\geq 0$, it follows that
\begin{equation}\label{INEQ: M-SAPD relation 4}
\begin{aligned}
              \cD(x_N,y_N)  \leq 2\theta^N  \cD(x_0,y_0) + (1-\theta)\delta_\mu,
\end{aligned}
\end{equation}
where $\delta_\mu \triangleq100(\frac{\delta_x^2}{\mu_x} +\frac{\delta_y^2}{\mu_y}) $, and $\cD(x,y)=\mathbb{E}\Big[\mu_x\|x-x^*\|^2 + \mu_y\|y-y^*\|^2\Big]$.
\end{corollary}
\begin{proof}
For any  $N\geq0$, it follows from \cref{Thm: main result_R1} that
\begin{equation*}
     \mathbb{E}[\frac{1}{2\tau}\|x_N-x^*\|^2 +  \frac{1 - \alpha\sigma}{2\sigma}\|y_N-y^*\|^2]
              \leq \rho^{N}\Big(\frac{1}{2\tau}\| x_0 - x^*\|^2
    +  \frac{1}{2\sigma}\| y_0 - y^*\|^2\Big) + \frac{\rho}{1-\rho}{\Xi}_{\tau,\sigma,\theta}.
\end{equation*}
Using the parameter choice
$$
\tau = \frac{1-\theta}{\theta\mu_x},\; \sigma = \frac{1-\theta}{\theta\mu_y},\; \alpha = \frac{1}{2\sigma}-\sqrt{\theta} L_{yy},\;\rho = \theta,
$$
we first obtain that
$
\frac{1-\alpha\sigma}{\sigma} \geq \frac{1}{2\sigma};
$ then this inequality together with our parameter choice leads to
\begin{equation*}
     \cD(x_N,y_N)
              \leq 2 \theta^{N}\cD(x_0,y_0)
     + 4{\Xi}_{\tau,\sigma,\theta}.
\end{equation*}
Then the desired result follows directly from \cref{lemma:Xi-bound-Response}.
\end{proof}

Next, in the following \mg{result}, we choose the number of steps $n_t$ and parameters $\tau_t, \sigma_t$ for each stage $t$ of M-SAPD in a particular fashion, and obtain performance bounds for each stage.
%give the stepsize choice for M-SAPD in each stage and shows the convergence result at that stage.
\begin{theorem}\label{THM: general noise - err relation for M-SAPD}
Suppose Assumptions~\ref{ASPT: lipshiz gradient}, \ref{ASPT: unbiased noise assumption} hold. Let  $\{\{z^t_k=(x^t_k, y^t_k)\}^{n_t}_{k=0}\}_{t\geq0}$ be the iterates generated by M-SAPD stated in \cref{ALG:M-SAPD} with the following parameters
\begin{equation*}\label{Condtion: M-SAPD paramter 1}
\begin{aligned}
         &\theta_t =  1 - \frac{1-\theta}{2^t},\;\theta=\max\Big\{ 1-\min\{\frac{\mu_y}{L_{yx}},\frac{\mu_y}{L_{yy}},\sqrt{\frac{\mu_x\mu_y}{L_{yx}L_{xy}}}, \frac{\mu_x}{L_{yx}}\}, \bar{\theta}\Big\}
        \\
        & \tau_t = \frac{1-\theta_t}{\mu_x\theta_t },
        \sigma_t = \frac{1-\theta_t}{\mu_y\theta_t}, n_t =
        \begin{cases}
        n_0, & t=0 \\
        \lceil \frac{p 2^t log(2)}{1-\theta} \rceil,& t\geq 1
        \end{cases},
\end{aligned}
\end{equation*}
where $p\geq 3$ is an arbitrary real number and $\bar{\theta}$ is defined in \cref{Condition: SP solution to noisy LMI-RS}.
Then for each $t\geq 0$,
\begin{equation}\label{INEQ: M-SAPD final relation}
            \cD(x_0^{t+1},y_0^{t+1})
            \leq \frac{\exp(-n_0(1-\theta))}{2^{t(p-1)-1}}\cD(x_0^{0},y_0^{0}) + \frac{1}{2^{t-1}}(1-\theta) \delta_\mu.
\end{equation}
\end{theorem}
\begin{proof}
For each $t\geq0$, it is easy to \mg{see} %obtain
that $\theta_t\geq \theta$.
% If we let
% \begin{equation}\label{Condtion: M-SAPD paramter 2}
%     \begin{aligned}
%         & \pi_{1,t} = 2\sigma_t\theta_t L_{yx} \frac{1}{1-2\sigma_t\left(\pi_{2,t}+\frac{\theta_t}{\pi_{2,t}}\right)L_{yy}},
%         \;
%         \pi_{2,t} = \sqrt{\theta_t},
%         \;
%         \alpha_t = \frac{\theta_t L_{yx}}{\pi_
%          {1,k}} + \frac{\theta_t L_{yy}}{\pi_{2,t}},
%     \end{aligned}
% \end{equation}
Then it follows from \cref{Coroallary: explicit solution to noisy LMI-RS} that $\{\theta_t,\tau_t,\sigma_t\}$
% together with \eqref{Condtion: M-SAPD paramter 2}
is a solution to MI \cref{eq: general SAPD LMI_R1}.%\todo{Mention MI eqn. number directly or provide the MI. Otherwise proof cannot be easily followed. (I add it.)}
Recall that $z^1_0=z^0_{n_0}$; therefore, it follows from \cref{Corollary:one-step-relation-for-M-SAPD}  that
\begin{equation}\label{INEQ: induction k = 0}
    \begin{aligned}
           \cD(x_0^{1},y_0^{1}) \leq & 2\theta_0^{n_0} \cD(x_0^{0},y_0^{0})  +  (1-\theta)\delta_\mu
            \\
            = & 2(1-\frac{1-\theta}{2^0})^{n_0} \cD(x_0^{0},y_0^{0})  +  (1-\theta)\delta_\mu
            \\
            \leq & 2\exp( - n_0(1-\theta)) \cD(x_0^{0},y_0^{0})  +   (1-\theta)\delta_\mu
            \\
            = &  \frac{\exp( - n_0(1-\theta)) \cD(x_0^{0},y_0^{0}) }{2^{0*(p-1)-1}} + \frac{1}{2^0}   (1-\theta) \delta_\mu
            \\
            \leq &
           \frac{\exp( - n_0(1-\theta)) \cD(x_0^{0},y_0^{0}) }{2^{0*(p-1)-1}}+ \frac{1}{2^{-1}} (1-\theta) \delta_\mu,
    \end{aligned}
\end{equation}
where the first inequality is from \cref{Corollary:one-step-relation-for-M-SAPD}; the second inequality is from the fact that $(1-x)^n\leq\exp(-nx)$; the last inequality uses the fact $\nu\geq 2$. Thus, \eqref{INEQ: M-SAPD final relation} is true for $t = 0$. Then, we suppose that \eqref{INEQ: M-SAPD final relation} is true for $t = i$. When  $t=i+1$, it also follows from \cref{Corollary:one-step-relation-for-M-SAPD} that
\begin{equation}\label{INEQ: induction k = i+1}
    \begin{aligned}
            \cD(x_0^{i+2},y_0^{i+2})
            \leq &
            2\theta_{i+1}^{n_{i+1}}\cD(x_0^{i+1},y_0^{i+1}) +  (1-\theta_{i+1}) \delta_\mu
            \\
            = & 2(1-\frac{1-\theta}{2^{i+1}})^{n_{i+1}}\cD(x_0^{i+1},y_0^{i+1})  +  \frac{1}{2^{i+1}}(1-\theta) \delta_\mu
            \\
            \leq & 2\exp(-\frac{n_{i+1}(1-\theta)}{2^{i+1}})\cD(x_0^{i+1},y_0^{i+1}) +  \frac{1}{2^{i+1}}(1-\theta) \delta_\mu
            \\
            \leq &  \frac{1}{2^{p-1}}\cD(x_0^{i+1},y_0^{i+1})  + \frac{1}{2^{i+1}}(1-\theta) \delta_\mu,
    \end{aligned}
\end{equation}
where the last inequality additionally uses the fact that
$n_{i+1} = \lceil \frac{p2^{i+1}log(2)}{1-\theta} \rceil$. If we substitute \eqref{INEQ: M-SAPD final relation} for $t =i$ into \eqref{INEQ: induction k = i+1}, it follows that
\begin{equation*}
    \begin{aligned}
           \cD(x_0^{i+2},y_0^{i+2})
            \leq & \frac{1}{2^{p-1}}\left[\frac{\exp(-n_0(1-\theta))}{2^{i(p-1)-1}} \cD(x_0^{0},y_0^{0})  + \frac{1}{2^{i-1}}(1-\theta) \delta_\mu\right]
            \\
            & + \frac{1}{2^{i+1}}(1-\theta) \delta_\mu
            \\
            = & \frac{\exp(-n_0(1-\theta))}{2^{(i+1)(p-1)-1}} \cD(x_0^{0},y_0^{0})  +(\frac{1}{2^{i+p-2}}+ \frac{1}{2^{i+1}}) (1-\theta) \delta_\mu
            \\
            \leq &\frac{\exp(-n_0(1-\theta))}{2^{(i+1)(p-1)-1}} \cD(x_0^{0},y_0^{0})  +\frac{1}{2^{i}} (1-\theta) \delta_\mu,
    \end{aligned}
\end{equation*}
where the last inequality is \mg{due to the fact that} $p\geq3$. Then, by an induction argument, we conclude.
\end{proof}
Finally, in the following corollary, we combine \mg{our previous results} to obtain an iteration complexity result for M-SAPD \mg{given in \cref{ALG:M-SAPD}}. \mg{This corollary shows that it is possible to remove the logarithmic factor in the iteration complexity bounds we provided for SAPD, by using the multi-stage variant M-SAPD with parameters given in \cref{THM: general noise - err relation for M-SAPD}.}
\begin{corollary} Suppose $\mu_x,\mu_y>0$, and Assumptions~\ref{ASPT: lipshiz gradient}, \ref{ASPT: unbiased noise assumption} hold. For any $\epsilon>0$, suppose the parameters $\{\tau_t,\sigma_t,\theta_t,n_t\}_{t\geq0}$ and $p$ are chosen according to \cref{THM: general noise - err relation for M-SAPD} and let
$
n_0 = \mathcal{O}\Big(\frac{1}{1-\theta}\ln (\frac{2}{\epsilon})\Big),
$
Then, the  complexity of M-SAPD, as stated in Algorithm \cref{ALG:M-SAPD}, to generate $z_\epsilon=(x_\epsilon,y_\epsilon)\in\cX\times \cY$ such that $\cD(x_\epsilon,y_\epsilon)\leq\epsilon$ is
\begin{equation*}
\mathcal{O}\Big(  (
\frac{\max\{L_{xx},L_{yx}\}}{\mu_x} + \sqrt{\frac{L_{yx}L_{xy}}{\mu_x\mu_y}} + \frac{\max\{L_{yy},L_{yx}\}}{\mu_y}) \ln (\frac{1}{\epsilon})+ p( \frac{\delta_x^2}{\mu_x} + \frac{\delta_y^2}{\mu_y})\frac{1}{\epsilon} \Big).
\end{equation*}
\begin{proof}
First, we define $N(t) \triangleq \sum_{i=0}^{t} n_i$. Note that, for $t\geq 1$, it follows from the fact $\lceil x \rceil < 2x $ that
$$
N(t)-n_0 = \sum_{i=1}^{t} n_i =\sum_{i=1}^{t} \lceil \frac{p2^{i}\ln(2)}{1-\theta} \rceil\leq \frac{2 p\ln(2)}{1-\theta}\sum_{i=1}^{t} 2^i=\frac{4p (2^t-1)\ln(2)}{1-\theta}.
$$
Furthermore, given an arbitrary positive integer $n$, there exists an unique $T$ such that $N(T)<n\leq N(T+1)$. For such pair of $(n,T)$, it follow that
$$
n - n_0 \leq N(T+1) - n_0\leq \frac{4p (2^{T+1}-1)\ln(2)}{1-\theta}.
$$
Then, we can obtain that
\begin{equation}\label{INEQ: 2^k bound}
    2^{T} \geq\frac{(1-\theta)(n-n_0)}{8p\ln(2)} +\frac{1}{2}\geq \frac{(1-\theta)(n-n_0)}{8p\ln(2)}.
\end{equation}
Moreover, letting $\hat{z}_n = z^{T+1}_{n-N(T)}$,
according to stage $T+1$ of M-SAPD, it follows from \cref{Corollary:one-step-relation-for-M-SAPD} that
\begin{equation*}
    \cD(\hat{x}_n,\hat{y}_n) \leq\nu\theta_{T+1}^{n-N(T)} \cD(x_0^{T+1},y_0^{T+1}) + \frac{1}{\nu^{T+1}}(1-\theta) \delta_\mu.
\end{equation*}
If we use \eqref{INEQ: M-SAPD final relation} within the above equation, it follows that
\begin{equation*}
\begin{aligned}
         \cD(\hat{x}_n,\hat{y}_n) & \leq2\theta_{T+1}^{n-N(T)} \left[\frac{\exp(-n_0(1-\theta))}{2^{T(p-1)-1}} \cD(x_0^{0},y_0^{0})
         +\frac{1}{2^{T-1}}(1-\theta) \delta_\mu\right]
          +  \frac{1}{2^{T+1}}(1-\theta) \delta_\mu.
         \\
         & \leq\frac{\exp(-n_0(1-\theta))}{2^{T(p-1)-2}} \cD(x_0^{0},y_0^{0})  +\frac{1}{2^{T-2}} (1-\theta) \delta_\mu+ \frac{1}{2^{T+1}}(1-\theta) \delta_\mu
         \\
         & \leq\exp(-n_0(1-\theta)) \cD(x_0^{0},y_0^{0})  + \frac{1}{2^{T-3}}(1-\theta) \delta_\mu,
\end{aligned}
\end{equation*}
where we \mg{used} the fact that $\theta\leq 1$ in the second inequality and $p\geq3$ in the third inequality.
Furthermore, if we use \eqref{INEQ: 2^k bound} within above inequality, it follows that
\begin{equation*}
\begin{aligned}
    \cD(\hat{x}_n,\hat{y}_n)  &\leq\exp(-n_0(1-\theta)) \cD(x_0^{0},y_0^{0})  +\frac{64p\ln(2)}{n-n_0} \delta_\mu
\end{aligned}
\end{equation*}
For $\epsilon>0$, a sufficient condition for $\cD(\hat{x}_n,\hat{y}_n) \leq \epsilon$ is
$$
\exp(-n_0(1-\theta)) \cD(x_0^{0},y_0^{0})  \leq \frac{\epsilon}{2},\;\frac{64p\ln(2)}{n-n_0} \delta_\mu\leq \frac{\epsilon}{2}.
$$
Since we let $n_0 = \mathcal{O}\Big( \frac{1}{1-\theta}\ln (\frac{1}{\epsilon}) \Big)$, the \mg{first} inequality on \mg{the left hand-side} is trivially satified. \mg{This} means that after at most $n_\epsilon$ iterations of M-SAPD, it will generate $\hat{z}_{n_\epsilon}$ s.t. $\cD(\hat{x}_{n_\epsilon},\hat{y}_{n_\epsilon}) \leq \epsilon$, where
$$
n_0 = \mathcal{O}\Big( \frac{1}{1-\theta}\ln (\frac{1}{\epsilon}) \Big),\; n_\epsilon = \mathcal{O}\Big(n_0 + p( \frac{\delta_x^2}{\mu_x} + \frac{\delta_y^2}{\mu_y})\frac{1}{\epsilon}\Big).
$$
Then, \mg{using} the choice of $\theta$, we conclude that
$$
 n_\epsilon= \mathcal{O}\Big(  (
\frac{\max\{L_{xx},L_{yx}\}}{\mu_x} + \sqrt{\frac{L_{yx}L_{xy}}{\mu_x\mu_y}} + \frac{\max\{L_{yy},L_{yx}\}}{\mu_y}) \ln (\frac{1}{\epsilon})+ p( \frac{\delta_x^2}{\mu_x} + \frac{\delta_y^2}{\mu_y})\frac{1}{\epsilon} \Big).
$$
\end{proof}
\end{corollary}

\section{\sa{Euclidean projection} onto the Intersection of the Simplex and \mg{the} $f$-divergence Ball}
In this section, we show an efficient method to solve \sa{the proximal problems $\min_{y\in\cP_r}\frac{\mu_y}{2}\norm{y}^2+\tfrac{1}{2\sigma}\norm{y-(y_k+\sigma\tilde{s}_k)}^2$ \mg{arising} when SAPD \mg{is} applied to \eqref{pbm-dist-robust-scsc}. In the rest, we consider a generic form of this problem. Indeed, given some $\bar{p}\in \mathbb{R}^n$ and $R>0$, we aim to solve
{\small
\begin{equation}
\label{eq: projection problem}
        p^* \triangleq \argmin_{p\in\cP}  \|p - \bar{p}\|^2,\quad \hbox{where}\quad \mathcal{P} \triangleq \{p\in \mathbb{R}^n_+:~\mathbf{1}^\top p = 1,~\|p - \mathbf{1}/n\|^2 \leq R^2\}. 
\end{equation}}}%
Next, we construct an equivalent problem to \eqref{eq: projection problem}, mainly because computing a dual optimal solution for the new formulation would be easier. Let $\cS\triangleq\{p\in\reals^n_+:\ \mathbf{1}^\top p=1 \}$. For $p\in\cP\subset\cS$, we have $\norm{p-\mathbf{1}/n}^2=\norm{p}^2-1/n$ since $\mathbf{1}^\top p=1$. Therefore, \eqref{eq: projection problem} is equivalent to
{\small
\begin{align}
\label{eq:projection-reformulation}
p^*=\argmin_{p\in\cS}\{\frac{1}{2}\norm{p-\bar{p}}^2:~\norm{p}^2\leq R^2+\frac{1}{n}\triangleq \bar{R}^2\}
\end{align}}%
In the \mg{literature}, many \sa{efficient} methods are provided to compute \sa{the Euclidean projection of a given point onto a unit simplex}, e.g., see~\cite{condat2016fast}. \sa{Therefore, we assume that $\| {p^s}\|>\bar{R}$, where $p^s \triangleq \argmin\{\|p-\bar{p}\|^2:~p\in\cS\}$; otherwise, i.e., $\norm{p^s}\leq \bar{R}$, we trivially have $p^*=p^s$; thus, $p^*$ can be efficiently computed with one of these simplex projection methods \mg{from} the literature.} 
%\xtodo{Do we need to say more about $p^*$ must satisfy $\| p^*- \frac{\mathbf{1}}{n}\|=r$?}
\sa{Since we assume that $\| {p^s}\|>\bar{R}$, $p^*$ must satisfy $\| p^*\|=\bar{R}$.}
%\mtodo{Is it $\leq R$? this is not clear to me. for instance if $\bar{p}=\mathbf{1}/n$, then $p*=\bar{p}$?}
%\nsa{No, it is with equality. In your example, $\bar{p}=p^s=p^*$. The reason it would be equality comes from a simple argument in convex optimization stating that the optimal solution should be on the boundary when you add the new set.}
\sa{The Lagrangian function for the problem in~\cref{eq:projection-reformulation} can be written as
{\small
\begin{equation}
     \cL(p,\lambda) \triangleq \mathbbm{1}_{\cS}(p)+\frac{1}{2}\|p-\bar{p}\|^2 + \frac{\lambda}{2}\Big(\|p\|^2 - \bar{R}^2\Big),
     %\; \text{where}\;p^T\mathbf{1}=1, p\geq 0,
\end{equation}}%
where $\mathbbm{1}_{\cS}(\cdot)$ denotes the indicator function of $\cS$.\footnote{\mg{The indicator function $\mathbbm{1}_{\cS}(\cdot)$ is defined as $\mathbbm{1}_{\cS}(x)=0$ if $x\in \cS$, $\mathbbm{1}_{\cS}(x)=+\infty$ otherwise.}} %Define
% The dual problem can be written as
% \begin{equation}
% \label{eq: Lagrange dual}
%     \max_{\lambda \geq 0}  \cL(p^*(\lambda),\lambda),
% \end{equation}
% where
{\small
\begin{equation}
\label{eq:p-lambda}
    p^*(\lambda) \triangleq \argmin_{p\in\reals^n} \cL(p,\lambda) = \argmin_{p\in\cS} \|p-\frac{\bar{p}}{1+\lambda}\|^2.
\end{equation}}% 
The aim is to compute $\lambda^*\geq 0$ such that $\|p^*(\lambda^*)\| = \bar{R}$, considering that $(p^*(\lambda^*),\lambda^*)$ is a KKT point; thus, $p^*=p^*(\lambda^*)$. It is essential to observe three critical points: \emph{i)} $p^*(0)=p^s$, which implies  $\norm{p^*(0)}>\bar{R}$; \emph{ii)} $\norm{p^*(\lambda_1)}\leq\norm{p^*(\lambda_2)}$ for all $\lambda_2\geq \lambda_1\geq 0$; \emph{iii)} $p^*(\lambda)\to \mathbf{1}/n$ as $\lambda\nearrow\infty$, which also implies that $\norm{p^*(\lambda)}<\bar{R}$ for sufficiently large $\lambda>0$ since $\mathbf{1}/n\in\cP$. These observations show that we can start from $\lambda=0$ and keep gradually increasing it until the first time $\norm{p^*(\lambda)}=\bar{R}$.}
%\nsa{I will continue to edit the rest.}

\sa{For numerical stability, i.e., %not to let 
\mg{for avoiding} $\lambda\to \infty$, instead of \eqref{eq:p-lambda}, we will consider an equivalent problem:} for $\gamma\in(0,1]$,
{\small
\begin{equation}
    \sa{p^*_\gamma \triangleq \argmin_{p\in\cS}} \frac{1}{2}\|p - \gamma\bar{p}\|^2.
\end{equation}}%
\sa{Our aim is to compute $\bar{\gamma}\in(0,1)$ such that $\|p^*_{\bar{\gamma}}\| = \bar{R}$.}%
\sa{Let $u\in\reals^n$ consist of elements of $\bar{p}$ in the descending order, i.e., $u_1\geq u_2\geq...\geq u_n$. Define
$K_\gamma \triangleq  \max_{k\in[n]}\{k:(\sum_{i=1}^k \gamma u_i-1)/k<\gamma u_k\}$ and $q_\gamma \triangleq (\sum_{i=1}^{K_\gamma}\gamma u_i-1)/K_\gamma$, where $[n]\triangleq\{1,\ldots,n\}$. Note that $K_\gamma\geq 1$ is well-defined since $\gamma u_1-1<\gamma u_1$. From \cite[Algorithm 1]{condat2016fast}, we know that}
{\small
\begin{equation} 
\label{eq: projection solution}
    p^*_\gamma = \begin{cases}
    \gamma u_i - q_\gamma & i = 1,2,...,K_\gamma, \\
    0 & \text{otherwise.}
    \end{cases}
\end{equation}}%
\sa{It follows from \cref{eq: projection solution} that the equation $\|p^*_{\bar{\gamma}}\| = \bar{R}$ has a unique positive solution,}
{\small
\begin{equation}
\label{eq: gamma solution}
    \bar{\gamma} =\sqrt{\frac{\bar{R}^2 -1/K_{\bar{\gamma}} }{\sum_{i=1}^{K_{\bar{\gamma}}}u_i^2 + (\sum_{i=1}^{K_{\bar{\gamma}}} u_i)^2/K_{\bar{\gamma}}}}.
\end{equation}}%
\sa{Since $K_{\bar{\gamma}}$ depends on $\bar{\gamma}\in(0,1)$, we cannot solve \eqref{eq: gamma solution} for $\bar{\gamma}$ immediately.}

\sa{At this point, it is essential to observe that for any $\gamma\in(0,1)$, the definitions of $K_\gamma$ and $q_\gamma$ imply that $K_\gamma = k\in[n]$ if and only if}
{\small
\begin{equation}
    \label{eq: projection lemma condition 1}
    \gamma u_{k + 1}\leq q_\gamma < \gamma u_{k},
\end{equation}}%
\sa{where we define $u_{n+1} \triangleq 0$. Since $K_{\bar{\gamma}}\in[n]$, we can set $K_{\bar{\gamma}}=k$ for $k=1,2,...,n$ and check whether $K_{\bar{\gamma}}$ satisfies the condition in~\cref{eq: projection lemma condition 1} for $\bar{\gamma}$ computed by 
\cref{eq: gamma solution}. Then substituting such $K_{\bar{\gamma}}$ into \cref{eq: projection solution} yields the solution $p^*$.}

\end{document}